\documentclass[12pt,twoside]{article}
\usepackage[T1]{fontenc}\usepackage[latin1]{inputenc}
\usepackage{french,multicol,amsmath,amsfonts}

\font\eightrm=ecrm0800 scaled\magstep1 
\font\bfone=ecbx1000 scaled\magstep2   

\makeatletter
\def\ps@headings{
   \let\@mkboth\@gobbletwo
   \def\@oddhead{\ifnum\value{page}=1\else\hfill\thepage\fi}
   \def\@oddfoot{}
   \def\@evenhead{\thepage\hfill}
   \def\@evenfoot{}
   }
\renewenvironment{abstract}{\narrower\footnotesize\bf
   Abstract.\quad\footnotesize\rm}{\par\bigskip}

\def\section{\@startsection{section}{1}{\z@}{0.8truecm}{\medskipamount}
   {\centering\bf}}
\def\trivlist{\bigbreak\vskip-\parskip
   \@trivlist \labelwidth\z@ \leftmargin\z@
   \itemindent\z@ \def\makelabel##1{##1}}
\def\@thmcounter#1{\noexpand\arabic{#1}}
\def\@thmcountersep{}
\def\@begintheorem#1#2{\it \trivlist \item[\hskip 
   \labelsep{\bf #1\ #2.\quad}]}
\def\@opargbegintheorem#1#2#3{\it \trivlist
   \item[\hskip \labelsep{\bf #1\ #2.\quad{\rm #3}}]}
\makeatother

\newtheorem{Theoreme}{Th\'eor\`eme}[section]
\long\def\theoreme#1{\begin{Theoreme} #1 \end{Theoreme}}
\newtheorem{Propositionjyd}[Theoreme]{Proposition}
\long\def\proposition#1{\begin{Propositionjyd} #1 \end{Propositionjyd}}
\newtheorem{Definitionjyd}[Theoreme]{D\'efinition}
\long\def\definition#1{\begin{Definitionjyd} #1 \end{Definitionjyd}}
\newtheorem{Lemme}[Theoreme]{Lemme}
\long\def\lemme#1{\begin{Lemme} #1 \end{Lemme}}
\newtheorem{Remarque}[Theoreme]{Remarque}
\long\def\remarque#1{\begin{Remarque} \rm~\nopagebreak#1 \end{Remarque}}

\newenvironment{Dem}[1]{\begin{trivlist}\item[\hskip\labelsep{\sc #1}]}
{\end{trivlist}}
\long\def\dem#1#2{\begin{Dem}{#1} ~\nopagebreak\smallskip #2 \end{Dem}}

\newtheorem{IntroDefinition}{D\'efinition}[part]
\long\def\Introdefinition#1{\begin{IntroDefinition} #1 \end{IntroDefinition}}
\newtheorem{Hypothese}{Hypoth\`ese\!\!} 
\long\def\hypothese#1{\begin{Hypothese} #1 \end{Hypothese}}
\newtheorem{Resultats}[Hypothese]{R\'esultats\!\!}
\long\def\resultats#1{\begin{Resultats} #1 \end{Resultats}}

\def\thebibl#1#2#3{
\section*{{\normalsize\bf #2}}
\list
  {[\arabic{enumi}]}{\settowidth\labelwidth{[#1]}\leftmargin\labelwidth
    \advance\leftmargin\labelsep
    \setlength{\itemsep}{0cm}
    \usecounter{enumi}\setcounter{enumi}{#3}}
    \def\newblock{\hskip .11em plus .33em minus -.07em}
    \sloppy
    \sfcode`\.=1000\relax}

\newenvironment{index des notations}[1]
{\begin{multicols}{4}[\section*{Index des notations}\smallskip]\begin{trivlist}
   \addcontentsline{toc}{section}{Index des notations}}
{\end{trivlist}\end{multicols}}
\let\labind\label 

\newskip\cqfdpostskipamount\cqfdpostskipamount=8pt plus 3pt minus 1.5pt
\def\cqfdpostskip{\vskip\cqfdpostskipamount}

\def\cqfd{\nolinebreak\quad\nolinebreak\hfill 
         \vbox{\hrule
         \hbox to 6pt{\vrule height 5,2pt \hfil \vrule}
         \hrule}
         \cqfdpostskip}
\def\cqfdpartiel{\nolinebreak\quad\nolinebreak\hfill
          \hbox{\ldots}
          \cqfdpostskip}
\def\cqfr{\nolinebreak\quad\nolinebreak\hfill
          \rule{4pt}{4pt} \cqfdpostskip}

\def\negsmallskip{\vskip-\smallskipamount}
\def\negmedskip{\vskip-\medskipamount}

\newcommand{\cala}{{\cal A}} \newcommand{\calb}{{\cal B}}
\newcommand{\calc}{{\cal C}} \newcommand{\cald}{{\cal D}}
\newcommand{\cale}{{\cal E}} \newcommand{\calf}{{\cal F}}
\newcommand{\calh}{{\cal H}} \newcommand{\call}{{\cal L}}
\newcommand{\calo}{{\cal O}} \newcommand{\calr}{{\cal R}}
\newcommand{\cals}{{\cal S}} \newcommand{\calv}{{\cal V}}
\newcommand{\calw}{{\cal W}}
\def\scalo{{\scriptstyle \calo}}

\let\Bbb\mathbb
\newcommand{\C}{{\Bbb C}} \newcommand{\G}{{\Bbb G}}
\newcommand{\M}{{\Bbb M}} \newcommand{\N}{{\Bbb N}}
\newcommand{\R}{{\Bbb R}} \newcommand{\Z}{{\Bbb Z}}

\let\goth\mathfrak
\newcommand{\gotha}{{\goth a}} \newcommand{\gothb}{{\goth b}}
\newcommand{\gothc}{{\goth c}} \newcommand{\gothg}{{\goth g}}
\newcommand{\gothh}{{\goth h}} \newcommand{\gothj}{{\goth j}}
\newcommand{\gothk}{{\goth k}} \newcommand{\gothm}{{\goth m}}
\newcommand{\gothn}{{\goth n}} \newcommand{\gothq}{{\goth q}}
\newcommand{\gotht}{{\goth t}} \newcommand{\gothu}{{\goth u}}
\newcommand{\gothv}{{\goth v}}

\def\Sum{\displaystyle\sum}
\def\Sumpetit{\sum\limits}
\def\Prod{\displaystyle\prod}
\def\Prodpetit{\prod\limits}
\def\Int{\displaystyle\int}
\def\Lim{\displaystyle\lim}

\def\Im{\mathop{\mathrm{Im}}\nolimits}
\def\sg{\mathop{\mathrm{sg}}\nolimits}
\def\pr{\mathop{\mathrm{pr}}\nolimits}
\def\id{\mathop{\mathrm{id}}\nolimits}
\def\abs#1{{\left\vert#1\right\vert}}

\def\rg{\mathop{\mathrm{rg \,}}\nolimits}
\def\tr{\mathop{\mathrm{tr}}\nolimits}
\def\Ker{\mathop{\mathrm{Ker}}\nolimits}
\def\Sp{\mathop{\mathrm{Sp}}\nolimits}

\def\moins0{\mathchoice{\setminus\!\{0\}}{\setminus\!\{0\}}
                       {\setminus\{0\}}{\setminus\{0\}}}

\def\egdef{:=} 

\def\diagd#1#2{
\mathchoice
{\begin{array}{cc}   {#1} &   \\   &  {#2} \end{array}}
{\begin{smallmatrix} \textstyle {#1} &   \\
   & \textstyle {#2} \end{smallmatrix}}
{\begin{smallmatrix} \scriptstyle {#1} &   \\
   & \scriptstyle {#2} \end{smallmatrix}}
{\begin{smallmatrix} \scriptscriptstyle {#1} &   \\
   & \scriptscriptstyle {#2} \end{smallmatrix}}
}
\def\diagt#1#2#3{
\mathchoice
{\begin{array}{ccc}   {#1}& &  \\  &{#2}& \\  & &{#3} \end{array}}
{\begin{smallmatrix} \textstyle{#1}& &  \\  &\textstyle{#2}& \\
   & &\textstyle{#3} \end{smallmatrix}}
{\begin{smallmatrix} {#1}& &  \\  &{#2}& \\
   & &{#3} \end{smallmatrix}}
{\begin{smallmatrix} \scriptscriptstyle{#1}& &  \\  &\scriptscriptstyle{#2}& \\
   & &\scriptscriptstyle{#3} \end{smallmatrix}}
}
\def\diagq#1#2#3#4{
\mathchoice
{\begin{array}{cccc} {#1}& & & \\ &{#2}& & \\ & &{#3}& \\ & & &{#4} \end{array}}
{\begin{smallmatrix} \textstyle{#1}& & & \\ &\textstyle{#2}& & \\
   & &\textstyle{#3}& \\  & & &\textstyle{#4} \end{smallmatrix}}
{\begin{smallmatrix} {#1}& & & \\ &{#2}& & \\
   & &{#3}& \\ & & &{#4} \end{smallmatrix}}
{\begin{smallmatrix} \scriptscriptstyle{#1}& & &\\ &\scriptscriptstyle{#2}& &\\
   & &\scriptscriptstyle{#3}& \\ & & &\scriptscriptstyle{#4} \end{smallmatrix}}
}

\def\Rot#1{
\mathchoice
{\begin{array}{cc}   0 & -{#1} \\ {#1} &  0 \end{array}}
{\begin{smallmatrix} \textstyle 0 & \textstyle -{#1} \\
   \textstyle {#1} & \textstyle 0 \end{smallmatrix}}
{\begin{smallmatrix}  0 & -{#1} \\
   {#1} &  0 \end{smallmatrix}}
{\begin{smallmatrix} \scriptscriptstyle 0 & \scriptscriptstyle -{#1} \\
   \scriptscriptstyle {#1} & \scriptscriptstyle 0 \end{smallmatrix}}
}
\def\Rotneg#1{
\mathchoice
{\begin{array}{cc}   0 & {#1} \\ -{#1} &  0 \end{array}}
{\begin{smallmatrix} \textstyle 0 & \textstyle {#1} \\
   \textstyle -{#1} & \textstyle 0 \end{smallmatrix}}
{\begin{smallmatrix}  0 & {#1} \\
   -{#1} &  0 \end{smallmatrix}}
{\begin{smallmatrix} \scriptscriptstyle 0 & \scriptscriptstyle {#1} \\
   \scriptscriptstyle -{#1} & \scriptscriptstyle 0 \end{smallmatrix}}
}

\def\sym#1{\mathfrak{S}_{#1}}
\def\interieur{\mathop{\mathrm{int}}\nolimits}
\def\centregr#1{\mathrm{Z}(#1)}
\def\centrealg#1{\mathrm{Z}(#1)}
\def\derivegr#1{\mathrm{D}(#1)}
\def\derivealg#1{\mathrm{D}(#1)}
\def\centragr#1#2{\mathrm{C}_{#1}(#2)}
\def\centraalg#1#2{\mathrm{C}_{#1}(#2)}
\def\normagr#1#2{\mathrm{N}_{#1}(#2)}
\def\normaalg#1#2{\mathrm{N}_{#1}(#2)}

\newcommand{\me}{\mathrm{e}} 
\newcommand{\mi}{\mathrm{i}} 

\def\Log{\mathop{\mathrm{Log}}\nolimits}
\def\Supp{\mathop{\mathrm{Supp}}\nolimits}
\def\diff{\mathrm{d}}
\def\norm#1{{\left\Vert#1\right\Vert}}

\def\restriction#1#2{\mathchoice
              {\setbox1\hbox{${\displaystyle #1}_{\scriptstyle #2}$}
              \restrictionaux{#1}{#2}}
              {\setbox1\hbox{${\textstyle #1}_{\scriptstyle #2}$}
              \restrictionaux{#1}{#2}}
              {\setbox1\hbox{${\scriptstyle #1}_{\scriptscriptstyle #2}$}
              \restrictionaux{#1}{#2}}
              {\setbox1\hbox{${\scriptscriptstyle #1}_{\scriptscriptstyle #2}$}
              \restrictionaux{#1}{#2}}}
\def\restrictionaux#1#2{{#1\,\smash{\vrule height .8\ht1 depth .85\dp1}}_{\,#2}}

\def\fonctioncar#1{\mbox{\rm 1\hspace{-.25em}l}_{#1}}

\def\ad{\mathop{\mathrm{ad}}\nolimits}
\def\Ad{\mathop{\mathrm{Ad}}\nolimits}
\def\Ind{\mathop{\mathrm{Ind}}\nolimits}
\def\Car{\mathop{\mathrm{Car}}\nolimits}

\def\ssreg#1{#1_{\scriptscriptstyle\mathit{ss\,reg}}}
\def\ssreggr#1{#1_{\!\scriptscriptstyle\mathit{ss\,reg}}}

\def\chiinfty#1#2{\chi_{#1}^{U \smash{{#2}_\C}}}

\def\F+{{\calf^+}}
\def\Fh+{{\calf^+_{\!\gothh}}}
\def\Fp+{{{{\calf}'}^+}}
\def\Fep+{{{{\calf}_{\!e}'}^+}}
\def\Fphp+{{{\calf_{\!\gothh'}'}^+}}
\def\XInd{X^{\textit{Ind}}}
\def\Xirr{X^{\textit{irr}}}
\def\Xirrp{X^{\textit{irr,+}}}
\def\Xfin{X^{\textit{final,+}}}

\let\wt\widetilde
\def\lambdatilde{{\mathchoice
             {\wt{\lambda}}
             {\wt{\lambda}}
             {\tilde{\lambda}}
             {\tilde{\lambda}}}}
\def\mutilde{{\mathchoice
             {\wt{\mu}}
             {\wt{\mu}}
             {\tilde{\mu}}
             {\tilde{\mu}}}}
\def\wtlie#1#2{\wt{{#1}^{\vrule width 0ex height .4ex \smash{*}}}_{\!\!\!#2}}

\def\gstregtilde{\wtlie{\gothg}{\textit{reg}}}
\def\gstregGtilde{\wtlie{\gothg}{\textit{reg},G}}
\def\gestregtilde{\wtlie{\gothg(e)}{\textit{reg}}}
\def\gejstregtilde{\wtlie{\gothg(e_j)}{\textit{reg}}}
\def\mpestregtilde{\wtlie{\gothm'(e)}{\textit{reg}}}

\def\gstssfondG{\gothg^*_{\textit{ssfond},G}}
\def\gstfondtilde{\wtlie{\gothg}{\textit{fond}}}
\def\mpstfondMptilde{\wtlie{\gothm'}{\textit{fond},M'}}

\def\gstssI{\gothg^*_{\textit{ssI}}}
\def\gstItilde{\wtlie{\gothg}{\textit{I}}}
\def\gestItilde{\wtlie{\gothg(e)}{\textit{I}}}
\def\gezerostItilde{\wtlie{\gothg(e_0)}{\textit{I}}}
\def\gstssIG{\gothg^*_{\textit{ssI},G}}
\def\gstIGtilde{\wtlie{\gothg}{\textit{I},G}}
\def\mstIMtilde{\wtlie{\gothm}{\textit{I},M}}

\def\gstssInc{\gothg^*_{\textit{ssInc}}}
\def\gstInctilde{\wtlie{\gothg}{\textit{Inc}}}
\def\gestInctilde{\wtlie{\gothg(e)}{\textit{Inc}}}
\def\gstssIncG{\gothg^*_{\textit{ssInc},G}}
\def\gstIncGtilde{\wtlie{\gothg}{\textit{Inc},G}}

\def\gstreg{\gothg^*_{\textit{reg}}}
\def\glzerostreg{\gothg(l_0)^*_{\textit{reg}}}
\def\glambdastreg{\gothg(\lambda)^*_{\textit{reg}}}

\def\gstss{\gothg^*_{\textit{ss}}}
\def\gestss{\gothg(e)^*_{\textit{ss}}}
\def\gstssreg{\ssreg{\gothg^*}}
\def\gestssreg{\ssreg{\gothg(e)^*}}
\def\mpstssreg{\ssreg{{\gothm'}^*}}

\def\gssreg{\ssreg{\gothg}}
\def\gessreg{\ssreg{\gothg(e)}}
\def\mpessreg{\ssreg{\gothm'(e)}}
\def\Gssreg{\ssreggr{G}}
\def\MpGzerossreg{\ssreggr{(M'G_0)}}
\def\Mpssreg{\ssreggr{M'}}

\makeatletter
\def\cases#1{\left\{\vcenter{\normalbaselines\m@th
    \ialign{$##\hfil$&\quad##\hfil\crcr#1\crcr}}\,\right.}
\makeatother

\begin{document}
%

\vbox{\vskip2.5truecm}
\bfone
\centerline{M\'ethode des orbites et formules du caract\`ere}
\centerline{pour les repr\'esentations temp\'er\'ees d'un groupe}
\centerline{alg\'ebrique r\'eel r\'eductif non connexe}
\vskip.8truecm
\centerline{\bf Jean-Yves Ducloux}\vskip.8truecm\rm

\begin{abstract}
Let $G$ be a non-connected reductive real Lie group.
In this paper, I parametrize the set of irreductible tempered characters of
$G\!$.
After\-wards, I describe these characters by means of some ``Kirillov's
formulas'', using the descent method near each elliptic element in $G\!$.

If $G$ is linear and connected, the parameters that I use are ``final basic''
parameters in the sense of Knapp and Zuckerman (cf. \cite [p.$\!$~453] {KZ82}).
\end{abstract}
\bigskip

%
\centerline{\bf Table des mati\`eres}
%

\def\ldotfill{\leaders\hbox spread6pt{\hss .\hss}\hfill}

\bigskip\noindent
{{\bf Introduction et notations g\'en\'erales}
\ldotfill{\bf \pageref{RepIntro}}}

\medskip\noindent
{{\bf I. Les param\`etres << forme lin\'eaire >>}
\ldotfill{\bf \pageref{RepI}}}

\hspace{-2em}
{1\quad Les param\`etres $\,\lambdatilde \in \gstregGtilde$
\ldotfill\pageref{Rep1}}

\hspace{-2em}
{2\quad Les mesures $\beta_{G \cdot \lambdatilde}$
\ldotfill\pageref{Rep2}}

\hspace{-2em}
{3\quad Points fix\'es par un \'el\'ement elliptique
\ldotfill\pageref{Rep3}}

\medskip\noindent
{{\bf II. Les param\`etres << repr\'esentation projective >>}
\ldotfill{\bf \pageref{RepII}}}

\hspace{-2em}
{4\quad Rappels sur les groupes sp\'ecial-m\'etalin\'eaire et
m\'etaplectique
\ldotfill\pageref{Rep4}}

\hspace{-2em}
{5\quad Les param\`etres $\,\tau \in \XInd_G(\lambdatilde)$
\ldotfill\pageref{Rep5}}

\hspace{-2em}
{6\quad Condition d'int\'egrabilit\'e
\ldotfill\pageref{Rep6}}

\medskip\noindent
{{\bf III. Construction de repr\'esentations}
\ldotfill{\bf \pageref{RepIII}}}

\hspace{-2em}
{7\quad Cas $G$ connexe
\ldotfill\pageref{Rep7}}

\hspace{-2em}
{8\quad Les repr\'esentations $T_{\lambdatilde,{\gotha^*}^+\!,\tau_+}^G\!$
\ldotfill\pageref{Rep8}}

\hspace{-2em}
{9\quad L'injection
$\,G \cdot (\lambdatilde,\tau) \mapsto T_{\lambdatilde,\tau}^G$
de $G\,\backslash\,\XInd_G$ dans $\widehat{G}$
\ldotfill\pageref{Rep9}}

\medskip\noindent
{{\bf IV. Caract\`eres des repr\'esentations}
\ldotfill{\bf \pageref{RepIV}}}

\hspace{-2em}
{10\quad Formule de restrictions des caract\`eres
\ldotfill\pageref{Rep10}}

\hspace{-2em}
{11\quad Passage de $M'$ \`a $G$
\ldotfill\pageref{Rep11}}

\hspace{-2em}
{12\quad Translation au sens de G.~Zuckerman
\ldotfill\pageref{Rep12}}

\hspace{-2em}
{13\quad Passage pour $M'$ au cas r\'egulier
\ldotfill\pageref{Rep13}}

\medskip\noindent
{{\bf Index des notations}
\ldotfill{\bf \pageref{RepIndex}}}

\medskip\noindent
{{\bf R\'ef\'erences}
\ldotfill{\bf \pageref{RepRef}}}
\vfill\eject

%
\part*{
Introduction et notations g\'en\'erales
}\label{RepIntro}
\addcontentsline{toc}{section}{Introduction et notations g\'en\'erales}
%

Dans cet article, je vais d\'ecrire par la m\'ethode des orbites le dual
temp\'er\'e d'un groupe de Lie r\'eel r\'eductif, en m'appuyant sur les quatre
travaux suivants :

\noindent\quad
-- la description du dual temp\'er\'e d'un groupe de Lie r\'eel r\'eductif
connexe par J.~Adams, D.~Barbasch et D.~Vogan dans \cite[ch. 11]{ABV92} ;

\noindent\quad
-- la param\'etrisation par M.~Duflo des classes d'\'equivalence des
repr\'esentations d'un groupe de Lie r\'eel r\'eductif qui sont
temp\'er\'ees irr\'eductibles avec un caract\`ere infinit\'esimal r\'egulier,
obtenue via la << th\'eorie de Mackey >> dans \cite[III]{Df82a} ;

\noindent\quad
-- la formule du caract\`ere << \`a la Kirillov >> de A.~Bouaziz dans \cite
{Bo87}, qui exprime le  caract\`ere de ces repr\'esentations \`a l'aide des
transform\'ees de Fourier d'orbites coadjointes semi-simples r\'eguli\`eres ;

\noindent\quad
-- les formules de W.~Rossmann dans \cite {Ro82} (o\`u  << principal >> se
traduit par << r\'egulier >>, et << regular >> par << semi-simple r\'egulier >>)
qui relient les transform\'ees de Fourier des orbites r\'eguli\`eres \`a celles
des orbites semi-simples r\'eguli\`eres.

\medskip
Un th\'eor\`eme de M.~Duflo dans \cite [p.$\!$~189] {Df82b} ram\`ene la
classification des duaux unitaires  des groupes lin\'eaires alg\'ebriques
r\'eels \`a la classification des duaux unitaires des groupes r\'eductifs
<< presque alg\'ebriques \`a noyau fini >>.
Ce sont donc ces groupes r\'eductifs qui m'int\'eressent.

\hypothese{

On se donne un groupe de Lie r\'eel \`a base d\'enombrable $G$\labind{G aa}
d'alg\`ebre de Lie not\'ee $\gothg$\labind{g}, un groupe lin\'eaire alg\'ebrique
$\G $\labind{G  aabarre} d\'efini sur $\R $, et un morphisme de groupes de Lie de
$G$ dans $\G (\R )$ de noyau fini central dont l'image est ouverte pour la
topologie usuelle.
Ainsi, tout \'el\'ement de $G$ (respectivement $\gothg$) a une << d\'ecomposition
de Jordan r\'eelle >> en composantes elliptique, positivement hyperbolique et
unipotente (respectivement composantes infinit\'esimalement elliptique,
hyperbolique et nilpotente) d\'ecrite dans
\cite [haut p.$\!$~36 et lem. 31 p.$\!$~38] {DV93}.
On suppose que $\gothg$ est r\'eductive avec un centre form\'e d'\'el\'ements
semi-simples.

}

\medskip
M.~Duflo a param\'etr\'e dans \cite [lem. 8 p.$\!$~173] {Df82a} une partie
du dual unitaire $\widehat{G}$\labind{G  aachapeau} de $G$ (pr\'ecis\'ee
ci-dessus), en termes d'orbites coadjointes semi-simples r\'eguli\`eres, en se
ramenant par induction aux s\'eries discr\`etes de certains sous-groupes de $G$.
Quand $G$ est connexe, cette partie de $\widehat{G}$ avait \'et\'e d\'ecrite par
Harish-Chandra (cf. \cite [th. 1 p.$\!$~198] {Ha76}).
Je param\`etre ici dans mon \emph{th\'eor\`eme \ref{injection dans le dual} (b)}
p.~\!\pageref{injection dans le dual}, en termes d'orbites coadjointes
r\'eguli\`eres << g\'en\'eralis\'ees >>, la partie plus grosse de $\widehat{G}$,
\'egale au dual temp\'er\'e de $G$, obtenue en rempla\c{c}ant au d\'ebut de
l'induction les repr\'esentations des s\'eries discr\`etes par les
repr\'esentations limites de s\'eries discr\`etes.
Quand $G$ est connexe, mon \'enonc\'e reproduit une partie du th\'eor\`eme 11.14
de \cite [p.$\!$~131] {ABV92} (voir aussi \cite [th. 14.76 p.$\!$~598] {Kn86})
avec l'apport suivant : l'\'egalit\'e de deux repr\'esentations se traduit
exactement par la conjugaison sous $G$ de couples de << bons param\`etres >>
qui leur sont associ\'es.

\smallskip
L'expression << bons param\`etres >> renvoie aux contraintes impos\'ees \`a ces
param\`etres.
En voici le principe heuristique (cf. \cite[th. 1 p.$\!$~193]{Df82a} et
\cite[th. 19 p.$\!$~211]{Df82b}).
Une repr\'esentation unitaire irr\'eductible $T$ de $G$ sera param\'etr\'ee dans
les cas favorables par l'orbite sous $G$ d'un couple~$(\lambdatilde,\tau)$.
Le terme $\lambdatilde$ doit \^etre une sorte de forme
lin\'eaire dont la composante semi-simple $l$ correspond au caract\`ere
$\chiinfty{\mi\,l}{\gothg}$ par lequel le centralisateur $(U \gothg_\C )^G$ de
$G$ dans $U\gothg_\C $ agit sur l'espace des vecteurs $C^\infty$ de l'espace de
$T$ (cf. \ref{caracteres canoniques} (a)).
Dans cet article, ce premier param\`etre est reli\'e \`a une limite de s\'erie
discr\`ete dont $T$ va en gros \^etre une induite.
Pour cette raison, il s'\'ecrira sous la forme $\lambdatilde = (\lambda,\F+)$,
o\`u $\F+\!$ est une chambre de Weyl pour des racines imaginaires (cf.
\cite [th. 5.7 p.$\!$~305] {Zu77}).
Sa << composante semi-simple >> est la forme lin\'eaire $\lambda$, en un sens
compatible avec les d\'efinitions concernant le cas des orbites coadjointes (cf.
\ref{bijection entre orbites}~(a)).
On note $G(\lambdatilde)$ le stabilisateur de $\lambdatilde$ dans $G$.
Le terme $\tau$ doit \^etre une repr\'esentation projective de
$G(\lambdatilde)/G(\lambdatilde)_0$, remont\'ee en une repr\'esentation unitaire
d'un rev\^etement de degr\'e $2$ de $G(\lambdatilde)$, telle que
$\,\diff_1\tau=\mi\lambda\id\,$ et le caract\`ere central de $\tau$ prolonge
celui de $T$ (cf. \ref{caracteres canoniques} (b)).
Dans cet article, ce second param\`etre est une repr\'esentation unitaire
(peut-\^etre non irr\'eductible) d'un analogue ad\'equat pour $\lambdatilde$ du
rev\^etement << de Duflo >> du stabilisateur $G(f)$ dans $G$ d'un \'el\'ement $f$
de $\gothg^*\!$.
En fait, mes constructions utiliseront un couple de param\`etres
$((\lambdatilde,{\gotha^*}^+\!),\tau_+)$ au lieu de $(\lambdatilde,\tau)$, o\`u
${\gotha^*}^+$ est une certaine chambre de Weyl pour des racines restreintes, et
$\tau_+$ est une repr\'esentation projective irr\'eductible de
$\normagr{G(\lambdatilde)}{{\gotha^*}^+}/G(\lambdatilde)_0$.
Ce nouveau couple de param\`etres satisfait encore les conditions
d\'ecrites ci-dessus.
Le couple $(\lambdatilde,\tau)$ se r\'ev\`elera \^etre << plus canonique >> quand
le caract\`ere de $T$ aura \'et\'e calcul\'e, mais il n'aura qu'un r\^ole
secondaire.
Les $G$-orbites de $(\lambdatilde,{\gotha^*}^+)$ et 
$((\lambdatilde,{\gotha^*}^+),\tau_+)$ s'identifient respectivement \`a celles
de $\lambdatilde$ et $(\lambdatilde,\tau)$.
Dans l'\'enonc\'e de mes r\'esultats ci-dessous, je ferai intervenir
l'ensemble $\XInd_G$ des couples $(\lambdatilde,\tau)$.
Il sera d\'efini plus pr\'ecis\'em\'ent dans la partie \ref{RepII}, en relation
avec la notion de << caract\`ere final >> (cf. \ref{parametres adaptes} (c)).

\smallskip
Le r\'esultat principal de cet article est mon
\emph{th\'eor\`eme \ref{formule du caractere}}
p.~\!\pageref{formule du caractere} qui d\'ecrit le caract\`ere des
repr\'esentations temp\'er\'ees $T$ de $G$ en terme de transform\'ees de Fourier
des mesures canoniques sur certaines orbites coadjointes r\'eguli\`eres reli\'ees
\`a un param\`etre $\lambdatilde$ attach\'e \`a $T$.
Ces formules du caract\`ere montrent la n\'ecessit\'e d'adapter le point de vue
<< orbites nilpotentes >> de la m\'ethode des orbites en oubliant l'aspect
<< orbite coadjointe >> de $G \cdot \lambdatilde$.
Par exemple les orbites coadjointes associ\'ees aux deux limites de la s\'erie
discr\`ete de $SL(2,\R )$, qui sont l'un ou l'autre des deux demi-c\^ones
nilpotents ouverts, n'ont pas de point fixe sous l'action de la << rotation >>
$\Ad^* \!\left( {\scriptstyle \Rot{1}} \right)$
(cf. remarque \ref{insuffisance de la methode des orbites} (3)) ;
elles doivent \^etre remplac\'ees par les ensembles des demi-droites incluses
dans l'int\'erieur de leurs enveloppes convexes respectives, au produit par
$\mi $ pr\`es.
On verra aussi dans la remarque \ref{autre insuffisance de la methode des
orbites} que ce nouveau point de vue permet de comprendre pourquoi certaines
orbites nilpotentes << admissibles >> au sens de M.~Duflo ne correspondent \`a
aucune repr\'esentation.

\medskip
Mes deux th\'eor\`emes se r\'esument comme suit (voir l'index des notations).

\resultats{

Il existe une bijection << canonique >>
$\; G \cdot (\lambdatilde,\tau) \mapsto T_{\lambdatilde,\tau}^G \;$ de
$\,G\,\backslash\,\XInd_G\,$ sur le dual temp\'er\'e de $G$.
Elle est d\'etermin\'ee par des formules du type :

\smallskip\centerline{
$k_e(X) \; \tr T_{\lambdatilde,\tau}^G (e \exp X)
\;\;=\;\;
\Sum_{\substack
{\dot{\lambdatilde'} \,\in\,
G(e) \backslash \,G \cdot \lambdatilde \, \cap \, \gstregtilde(e)\, \\
\textrm{tel que } \,\lambdatilde'\![e] \,\in\, \gestregtilde}}\!\!\!
c_{\widehat{e'},\lambdatilde}\;
\tr\tau(\widehat{e'})\;
{\widehat{\beta}}_{G(e) \cdot \lambdatilde'\![e]} (X)$
}

\smallskip\noindent
o\`u
$\; (\lambdatilde,\tau) \in \XInd_G\!$,
$e$ d\'ecrit l'ensemble des \'el\'ements elliptiques de $G$, \\
$\,X \,$ est un \'el\'ement d'un certain voisinage ${\cal V}_{\!e}$ de $0$ dans
$\mathfrak{g}(e)$ tel que $\; e \exp X$ est semi-simple r\'egulier dans $G$
(auquel cas $X$ est semi-simple r\'egulier dans $\mathfrak{g}(e)\!$), \\
$\; \widehat{e'}$ est un \'el\'ement du rev\^etement double
$\,G(\lambdatilde)^{\mathfrak{g}/\!\mathfrak{g}(\lambda)(\mi \rho_\F+)}\,$
(en posant $\lambdatilde = (\lambda,\F+)$)
qui se projette sur $\, e' = g^{-1}eg \,$
pour un $\, g \in G \,$ tel que $\; \lambdatilde' = g\lambdatilde$, \\
et, $\,k_e(X)\,$ et $\,c_{\widehat{e'},\lambdatilde}$
sont certains nombres complexes non nuls.

}

\medskip
Les id\'ees qui m'ont permis d'obtenir ces r\'esultats sont les suivantes.
Les constructions de M.~Duflo reprises \'etape par \'etape vont fournir une
bijection d'un ensemble d'orbites sur le dual temp\'er\'e de $G$, gr\^ace \`a
mon \emph{lemme \ref{resultat inattendu} (b)} p.~\!\pageref{resultat inattendu}.
Pour passer du cas connexe au cas non connexe, j'aurai besoin d'une
g\'en\'eralisation d'un r\'esultat de D.~Vogan utilis\'e par M.~Duflo pour une
construction homologique.
Celle-ci se trouve \`a la page 555 du livre \cite
{KV95} de A.~Knapp et D.~Vogan.
Au vu du cas connexe (cf. \cite [p.$\!$~64] {Ro80}), il est ensuite naturel de
chercher \`a r\'ecup\'erer les caract\`eres des repr\'esentations temp\'er\'ees
irr\'eductibles par passage \`a la limite \`a partir de ceux associ\'es aux
orbites semi-simples r\'eguli\`eres.
Dans le cas non connexe, le calcul du caract\`ere des repr\'esentations se
d\'ecomposera avant tout en deux \'etapes dont la premi\`ere est peu commode.

\smallskip\noindent\quad
\emph{\'Etape 1.}
On induit une repr\'esentation limite de la s\'erie discr\`ete d'un
sous-groupe de $G$ \`a la composante neutre de la composante de Levi $M'$ (en
g\'en\'eral hors de la classe d'Harish-Chandra) d'un certain << sous-groupe
parabolique >> $M'U$ de $G$.
L'analogue de cette repr\'esentation induite dans le cas connexe \'etait
seulement une repr\'esentation limite de la s\'erie discr\`ete de la composante
de Levi $M$ d'un sous-groupe parabolique cuspidal $MAN$ de $G$.
La m\'ethode ici comme dans le cas connexe, est d'appliquer le foncteur de
translation de Zuckerman (adapt\'e au cas du groupe non connexe $M'$) pour se
ramener aux repr\'esentations de $M'$ dont on conna\^{\i}t le caract\`ere.
Une difficult\'e est qu'on passe de $M'_0$ \`a $M'$ avec une \'etape homologique.
Sur le conseil d'A.~Bouaziz, je me suis inspir\'e des \'etapes (i) et (ii) de la
page 550 de son article \cite {Bo84} pour calculer l'action d'un groupe
$M'(\lambdatilde)^{\gothm'\!/\gothh}$ dans un espace poids de l'homologie d'un
$M'_0$-module translat\'e au sens de G.~Zuckerman.

\smallskip\noindent\quad
\emph{\'Etape 2.}
On induit ensuite de $M'U$ \`a $G$.
Le terrain aura \'et\'e pr\'epar\'e dans la partie \ref{RepI} par des rappels et
compl\'em\'ents concernant des r\'esultats de W.~Rossmann pour le calcul des
limites de transform\'ees de Fourier des mesures canoniques sur des orbites
semi-simples r\'eguli\`eres.
En particulier, dans la section \ref{Rep3}, j'aurai mis en \'evidence ce que
donne le passage de $G$ \`a $G(e)$ au niveau des param\`etres
$\lambdatilde$.
Le  reste de cette \'etape consiste \`a utiliser les m\'ethodes et r\'esultats
d'A.~Bouaziz.

\medskip
Les d\'efinitions qui suivent vont me permettre, d'abord de pr\'eciser mes
notations g\'en\'erales, et ensuite d'introduire les notations relatives aux
groupes r\'eductifs que j'utiliserai le plus fr\'equemment.

\Introdefinition{
\label{notations generales}

{\bf (a)}
On note $\abs{A}$ le cardinal d'un ensemble $A$ et $\dot{a}$ la classe d'un
\'el\'ement $a$ de $A$ modulo une relation d'\'equivalence $\sim$ sur $A$,
$\,\restriction{f}{B'}: B' \to C\,$ la restriction d'une application
$\,f: B \to C\,$ \`a une partie $B'$ de $B$,
$\,\overline{\!F}$ l'adh\'erence d'un sous-ensemble $F$ d'un espace topologique
$E$ (dans les sections \ref{Rep2} et \ref{Rep12}),
$l^\C $~la complexifi\'ee d'une application $\R $-lin\'eaire~$l$,
$\overline{w}$~et $\overline{W}$ les conjugu\'es d'un vecteur $w$ et d'un
sous-espace vectoriel $W$ dans le complexifi\'e $V_\C $ d'un espace vectoriel
r\'eel~$V$ (dans les sections \ref{Rep1}, \ref{Rep3}, \ref{Rep4}, \ref{Rep5},
\ref{Rep8} et \ref{Rep9}),
$u_Y$ et $u_{X/Y}$ les endomorphismes induits sur $Y$ et $X/Y$ par un
endomorphisme $u$ d'un espace vectoriel $X$ laissant invariant un sous-espace
vectoriel $Y$ de $X$.

\smallskip
{\bf (b)}
Soient $A$ un groupe de Lie r\'eel \`a base d\'enombrable et $\gotha$ son
alg\`ebre de Lie.
On note $1$ l'\'el\'ement neutre de $A$,
$A_0$ la composante neutre de $A$,
$\exp_A$ (ou $\exp$) l'application exponentielle de $\gotha$ dans $A$,
$\Ad^*$ et $\ad^*$ les repr\'esentations coadjointes de $A$ et de $\gotha$,
$\interieur \gotha$ le sous-groupe du groupe lin\'eaire de $\gotha$ engendr\'e
par les \'el\'ements $\exp(\ad X)$ quand $X$ d\'ecrit $\gotha$,
$\centregr{A}$ et $\derivegr{A}$ le centre et le groupe d\'eriv\'e de $A$,
$\centrealg{\gotha}$ et $\derivealg{\gotha}$ le centre et l'alg\`ebre
d\'eriv\'ee de $\gotha$,
$\centrealg{U \gotha_\C }$ le centre de l'alg\`ebre enveloppante $U \gotha_\C $
de $\gotha_\C $,
$\widehat{A}$ le dual unitaire de $A$.
\`A toute sous-alg\`ebre de Lie $\gothb$ de $\gotha$, on associe les
sous-groupes de Lie

\centerline{
$\centragr{A}{\gothb}
= \{ x \in A \mid \restriction{(\Ad x)}{\gothb} = \id \} \;$
et
$\; \normagr{A}{\gothb}
= \{ x \in A \mid \Ad x . \gothb \subseteq \gothb \}$}

\noindent
de $A$, dont les alg\`ebres de Lie sont

\centerline{
$\centraalg{\gotha}{\gothb}
= \{ X \in \gotha \mid \restriction{(\ad X)}{\gothb} = 0 \} \;$
et
$\; \normaalg{\gotha}{\gothb}
= \{ X \in \gotha \mid \ad X . \gothb \subseteq \gothb \}$.
}

\smallskip
{\bf (c)}
Soit $M$ une vari\'et\'e $C^\infty\!$ s\'epar\'ee \`a base d\'enombrable de
dimension $m$.
Une densit\'e $C^\infty\!$ sur $M$ est une mesure de Radon complexe $\rho$ sur
$M$ qui se lit dans toute carte de $M$ centr\'ee en un point $x$ sous la forme
$\,a\,\diff x_1 \!\cdots \diff x_m$, o\`u $\,a\,$ est une fonction $C^\infty\!$
sur un voisinage ouvert de $0$ dans $\R ^m\!$.
Elle << s'identifie >> \`a la famille form\'ee des applications $\rho(x)$ avec
$x \in M$ qui envoient un \'el\'ement de $\bigwedge^m T_x M\,\moins0$
d'image $\,\omega_0\,$ par la carte pr\'ec\'edente, sur
$\,a(0)\,\abs{(\diff x_1 \wedge \cdots \wedge \diff x_m)(\omega_0)}$.
Une fonction g\'en\'eralis\'ee sur $M$ est une forme lin\'eaire continue sur
l'espace des densit\'es $C^\infty\!$ \`a support compact sur $M$, muni de la
topologie de Schwartz.

\smallskip
{\bf (d)}
Soient $A$ un groupe de Lie r\'eel \`a base d\'enombrable et $\,\diff_A$ une
mesure de Haar \`a gauche sur $A$.
Une repr\'esentation continue $T$ de $A$ dans un espace de Hilbert complexe est
dite tra\c{c}able si les op\'erateurs $\;T(\varphi \, \diff_A)$ avec
$\, \varphi \in C^\infty_c(A) \,$ sont tra\c{c}ables ;
dans ce cas
$\; \tr T : \varphi \, \diff_A \mapsto \tr T(\varphi \, \diff_A) \;$
est une fonction g\'en\'eralis\'ee sur $A$. Pour toute repr\'esentation unitaire
continue $\pi$ d'un sous-groupe ferm\'e $B$ de $A$ muni d'une mesure de Haar \`a
gauche
$\,\diff_B$, on note
$\Ind_B^A \pi$ la repr\'esentation unitaire de $A$ << induite >> \`a partir de
$\pi$ comme dans
\cite [p.$\!$~99] {B.72}.

\smallskip
{\bf (e)}
Soit $V$ un espace vectoriel r\'eel de dimension finie.
La transform\'ee de Fourier d'une distribution temp\'er\'ee $\mu$ sur $V^*$ est
la fonction g\'en\'eralis\'ee $\widehat{\mu}$ sur $V$ d\'efinie par
l'\'egalit\'e
$\;\; \widehat{\mu}(v) = \Int_{V^*} \me^{\mi \,l(v)} \diff \mu(l) \;\;$
de fonctions g\'en\'eralis\'ees en $v \in V$.

}

\Introdefinition{
\label{notations frequentes}

{\bf (a)}
On note $\Car \gothg$\labind{Carg} l'ensemble des sous-alg\`ebres de Cartan de
$\gothg$.
Soit $\gothh \in \Car \gothg$.
On fixe un syst\`eme de racines positives
$\, R^+(\gothg_\C ,\gothh_\C )$\labind{R+(gC,hC)},
arbitrairement (sauf indication contraire), dans l'ensemble
$\,R(\gothg_\C ,\gothh_\C )\,$\labind{R(gC,hC)}
des racines de $\gothh_\C $ dans $\gothg_\C $.

On note
$W(G,\gothh)$\labind{W(G,h)} le groupe fini
$\,\normagr{G}{\gothh}/\centragr{G}{\gothh}$,
$\gotht$\labind{t} et $\gotha$\labind{a} les composantes infinit\'esimalement
elliptique et hyperbolique de $\gothh$,
$\, \gothh_{(\R )} \!= \mi \, \gotht \oplus \gotha$\labind{hr},
$T_0=\exp \gotht$\labind{T0} et $A=\exp \gotha$\labind{A aa},
$\rho_{\gothg,\gothh}$\labind{zzr g,h} la demi-somme des \'el\'ements de
$R^+(\gothg_\C ,\gothh_\C )$, et $H_{\alpha}$\labind{Halpha} la racine duale
d'une $\alpha$ dans le syst\`eme de racines $R(\gothg_\C ,\gothh_\C )$.

Une racine $\, \alpha \in R(\gothg_\C ,\gothh_\C )\,$ est dite
complexe (respectivement r\'eelle, imaginaire, ou compacte) quand sa conjugu\'ee
$\;\overline{\alpha}: X \in \gothh_\C \mapsto \overline{\alpha(\overline{X})}\;$
v\'erifie $\,\overline{\alpha} \notin \{\alpha,-\alpha\}\,$
(respectivement $\,\overline{\alpha}=\alpha$, $\,\overline{\alpha}=-\alpha$, ou
$\; (\C  H_{\alpha} \oplus \gothg_\C ^\alpha
\oplus \gothg_\C ^{-\alpha}) \cap \gothg \simeq {\goth su}(2)$).

\smallskip
{\bf (b)}
On fixe une forme bilin\'eaire $G$-invariante non d\'eg\'en\'er\'ee
${\scriptstyle \langle} ~,\!~ {\scriptstyle \rangle}$\labind{<,>ev}
sur $\gothg$ dont la complexifi\'ee
(encore not\'ee ${\scriptstyle \langle} ~,\!~ {\scriptstyle \rangle}$)
se restreint en un produit scalaire sur chaque $\gothh_{(\R )}$,
$\gothh \in \Car \gothg$.
On note
$\, {\scriptstyle \perp}^
{\!{\scriptscriptstyle \langle} \!~,\!~ {\scriptscriptstyle \rangle}} \,$
la relation de
${\scriptstyle \langle} ~,\!~ {\scriptstyle \rangle}$-orthogonalit\'e.
Soient $x \in G$ et $X \in \gothg$ semi-simples.
On note
$G(x)$\labind{G (x)aa} (resp. $G(X)$\labind{G (Xaa)aa})
et
$\gothg(x)$\labind{g (x)} (resp. $\gothg(X)$\labind{g (Xaa)})
les commutants de $x$ (resp. $X$) dans $G$ et $\gothg$.
Le $G(x)$-module $\gothg(x)^*$ (resp. le $G(X)$-module $\gothg(X)^*$) est
canoniquement isomorphe par restriction \`a l'ensemble
$\gothg^*(x)$\labind{g*(x)} (resp. $\gothg^*(X)$\labind{g*(Xaa)}) des
\'el\'ements de $\gothg^*$ fixes sous $\Ad ^*x$ (resp. annul\'es par
$\,\ad^*X$), c'est-\`a-dire nuls sur
$\,\gothg(x)^{{\scriptstyle \perp}^
{\!{\scriptscriptstyle \langle} \!~,\!~ {\scriptscriptstyle \rangle}}}\!$
(resp.
$\,\gothg(X)^{{\scriptstyle \perp}^
{\!{\scriptscriptstyle \langle} \!~,\!~ {\scriptscriptstyle \rangle}}}\!$).

\smallskip
{\bf (c)}
Soit $f \in \gothg^* \!$.
On note $G(f)$\labind{G (f)aa} le stabilisateur de $f$ dans $G$ et
$\gothg(f)$\labind{g (f)} l'alg\`ebre de Lie de $G(f)$.
Dans la mesure du possible, je d\'esignerai par $\mu$\labind{zzm},
$\nu$\labind{zzn}, $\lambda = \mu\!+\!\nu$\labind{zzl}, et $\xi$\labind{zznx} les
<< composantes infinit\'esimalement elliptique, hyperbolique, semi-simple et
nilpotente de $f$ >> (notions issues de~$\gothg$ \`a l'aide de
${\scriptstyle \langle} ~,\!~ {\scriptstyle \rangle}$).
On identifie $\gothg(\lambda)^*$ \`a l'ensemble des \'el\'ements de $\gothg^*$
nuls sur
$\,\gothg(\lambda)^{{\scriptstyle \perp}^
{\!{\scriptscriptstyle \langle} \!~,\!~ {\scriptscriptstyle \rangle}}}\!$.

\smallskip
{\bf (d)}
Soit $\gothm$ le commutant dans $\gothg$ d'un \'el\'ement semi-simple $X$ de
$\gothg$.
\'Etant donn\'ees des mesures de Haar $\diff _{\gothg}$ et $\diff _{\gothm}$ sur
$\gothg$ et $\gothm$, on pose pour tout $f \in \gothm^* \!$ :

\centerline{
$\bigl|\Pi_{\gothg,\gothm}\bigr|(f)
= \frac{1}{k!}\;\bigl|
\bigl( \restriction{B_f}{[\gothg,X]^2} \bigr)^{\!k}
(\omega)\bigr|
\,{\scriptstyle \times}\,
\left( \diff _{[\gothg,X]} (0) (\omega) \right)^{\!-1}\;\;$%
\labind{zzP g,m}
(cf. (b) et d\'ef. \ref{notations generales} (c)),
}

\noindent avec
$k = \frac {1}{2} \dim [\gothg,X]$,
$B_f = f([\cdot,\cdot])$,
$\omega \in \bigwedge ^{2k}[\gothg,X]\moins0$,
et
$\, \displaystyle \diff _{[\gothg,X]}
=
\diff _{\gothg} / \diff _{\gothm}$.

}

\medskip
Les autres notations utilis\'ees dans les \'enonc\'es sont accessibles au moyen
de l'index des notations qui se trouve \`a la fin de cet article.

%
\part{
Les param\`etres << forme lin\'eaire >>
}\label{RepI}
%

La m\'ethode des orbites, propos\'ee initialement par A.~Kirillov dans le cas
des groupes de Lie nilpotents simplement connexes, consiste \`a param\'etrer les
repr\'esentations unitaires irr\'eductibles d'un groupe de Lie \`a l'aide des
orbites de sa repr\'esentation coadjointe.
Quand la repr\'esentation est tra\c{c}able, son caract\`ere au voisinage de
l'\'el\'ement neutre doit \^etre reli\'e \`a la transform\'ee de Fourier de la
mesure canonique sur l'orbite associ\'ee.

Une petite modification de ce point de vue va \^etre n\'ecessaire ici.
\`A chaque repr\'esentation temp\'er\'ee irr\'eductible de $G$ sera associ\'ee
une orbite $\wt{\Omega}$ de $G$ dans un certain ensemble $\gstregtilde$ autre
que $\gothg^* \!$, et \`a cette orbite $\wt{\Omega}$ sera attach\'ee une somme
finie de mesures canoniques sur des orbites coadjointes.
En vue d'un prochain article, j'introduis aussi ci-dessous des parties
$\gstItilde$ et $\gstInctilde$ de $\gstregtilde$.

\medskip
Je vais exploiter dans cette partie certains r\'esultats de W.~Rossmann.

\section{\boldmath
Les param\`etres $\,\lambdatilde \in \gstregGtilde$
}\label{Rep1}

On construit dans cette section un ensemble $\gstregtilde$ muni d'une action de
$G$ dans lequel l'ensemble des $\lambda \!\in \gothg^*\!$ semi-simples
r\'eguli\`eres s'injectera de mani\`ere $G$-\'equivariante.

\definition{

On note $\gstreg$\labind{g*reg} l'ensemble des $f \!\in \gothg^*\!$
r\'eguli\`eres (c'est-\`a-dire l'ensemble des $f \!\in \gothg^*\!$ pour
lesquelles la dimension de $\gothg(f)$ est \'egale au rang de~$\gothg$),
$\gstss$\labind{g*ss} l'ensemble des $f \!\in \gothg^*\!$ semi-simples et
$\, \gstssreg = \gstss \cap \gstreg$\labind{g*ssregzz}.

On note aussi $\gstssI$\labind{g*ssI}
(respectivement : $\gstssInc$\labind{g*ssInc})
l'ensemble des $\lambda \! \in \gstss$ telles que $\gothg(\lambda)$ a une
sous-alg\`ebre de Cartan $\gothh$ pour laquelle les racines de
$(\gothg(\lambda)_\C ,\gothh_\C )$ sont imaginaires
(respectivement : imaginaires non compactes).

Soient $\lambda \! \in \gstss$ et $\gothh \in \Car\gothg(\lambda)$.
On note $C(\gothg(\lambda),\gothh)$\labind{Cgh} l'ensemble des chambres
(ouvertes) dans $\, \gothh_{(\R )}^{~~*}$ pour les racines imaginaires de
$(\gothg(\lambda)_\C ,\gothh_\C )$ et
$C(\gothg(\lambda),\gothh)_{reg}$\labind{Cghreg}
l'ensemble des $\F+\! \!\in C(\gothg(\lambda),\gothh)$ telles que les racines
imaginaires de $(\gothg(\lambda)_\C ,\gothh_\C )$ qui sont simples relativement
\`a $\F+\!$ sont non compactes.

On pose
$\hfill \gstregtilde
\,=\, \Bigl\{
(\lambda,\F+) \,;
\lambda \!\in\! \gstss, \,
\gothh \!\in\! \Car\gothg(\lambda) \textrm{ et }
\F+\! \!\in C(\gothg(\lambda),\gothh)_{reg}
\Bigr\} \hfill$%
\labind{g*reg tilde}

\centerline{
et
$\quad\gstregGtilde
\,=\, \Bigl\{
(\lambda,\F+) \in \gstregtilde \mid
\forall Z \in \Ker \exp_{T_0} \;\,
\me^{(\mi \,\lambda+\rho_{\gothg,\gothh})(Z)}=1
\Bigr\}$%
\labind{g*reg,G tilde}
}

\noindent
ind\'ependamment du choix d'un syst\`eme de racines positives
$R^+(\gothg_\C ,\gothh_\C )$ associ\'e \`a l'\'el\'ement  $\gothh$ de
$\Car\gothg(\lambda)$ attach\'e \`a $\F+\!$.

On pose aussi

\centerline{
$\gstfondtilde
= \Bigl\{
(\lambda,\F+) \,;
\lambda \!\in\! \gstss, \,
\gothh \!\in\! \Car\gothg(\lambda) \textrm{ fondamentale et }
\F+\! \!\in C(\gothg(\lambda),\gothh)_{reg}
\Bigr\}$%
\labind{g*fond tilde}
}

\centerline{
et
$\quad\gstssfondG
\,=\, \Bigl\{ \lambda \!\in\! \gstss \mid
\forall Z \in \Ker \exp_{T_0} \;
\me^{(\mi \,\lambda+\rho_{\gothg,\gothh})(Z)}=1
\Bigr\}$%
\labind{g*ssfond,G}
}

\noindent
ind\'ependamment du choix d'une $\gothh \in \Car\gothg(\lambda)$ fondamentale et
d'un syst\`eme de racines positives $R^+(\gothg_\C ,\gothh_\C )$,
$\,\gstssIG = \gstssI \cap \gstssfondG$\labind{g*ssI,G} et
$\,\gstssIncG = \gstssInc \cap \gstssfondG$\labind{g*ssInc,G}.

\smallskip
On note ensuite $\gstItilde$\labind{g*I tilde},
$\gstInctilde$\labind{g*Inc tilde}, $\gstIGtilde\!$\labind{g*I,G tilde} et
$\gstIncGtilde\!$\labind{g*Inc,G tilde} les images r\'eciproques de $\gstssI$,
$\gstssInc$, $\gstssIG\!$ et $\gstssIncG\!$ par l'application canonique
(premi\`ere projection) de $\gstfondtilde\!$ dans $\gstss$.

}

\remarque{

{\bf (1)}
En appliquant \cite [th. p.$\!$~217] {Ro82} \`a $\gothg(\lambda)$ pour trois
sous-alg\`ebres de Cartan (une \'egale \`a $\gothh$, une sans racine imaginaire,
une fondamentale) et tenant compte de
\cite [lem. {\sc a} $(a) \!\Rightarrow\! (c)$ p.$\!$~220
et suppl. {\sc c} p.$\!$~218] {Ro82},
on constate que les images des applications canoniques de $\gstregtilde$ et
$\gstfondtilde$ dans $\gothg^*$ sont form\'ees des $\lambda \! \in \gstss$ tels
que $\gothg(\lambda)$ a une sous-alg\`ebre de Cartan sans racine imaginaire. 

\smallskip
{\bf (2)}
D'apr\`es \cite [th. 6.74 p.$\!$~341 et th. 6.88 p.$\!$~344] {Kn96}, chaque
alg\`ebre de Lie semi-simple complexe a, \`a isomorphisme pr\`es, une unique
forme r\'eelle $\gothg_0$ qui poss\`ede une sous-alg\`ebre de Cartan $\gothh_0$
pour laquelle le syst\`eme de racines $\, R({\gothg_0}_\C ,{\gothh_0}_\C ) \,$ a
une base constitu\'ee de racines imaginaires non compactes.
Compte tenu de (1) ci-dessus
(ou de \cite [pb. 18 p.$\!$~369 et p.$\!$~557] {Kn96}),
l'image de l'application canonique de $\gstItilde$ dans $\gothg^*$ est form\'ee
des $\lambda \! \in \gstss$ tels que chacun des id\'eaux simples de
$\gothg(\lambda)$ est isomorphe \`a l'une des alg\`ebres de Lie suivantes : \\
$\; \goth{su}(p,p)$ avec $p \geq 1$ et $\goth{su}(p,p-1)$ avec $p \geq 2$,
$\goth{so}(p,p-1)$ avec $p \geq 3$,
$\goth{sp}(2n,\R )$ (not\'ee $\goth{sp}(n,\R )$ dans \cite {Kn96})
avec $n \geq 3$,
$\goth{so}(p,p)$ avec $p$ pair $\geq 4$
et $\goth{so}(p,p-2)$ avec $p$ pair $\geq 6$,
$E \, I\!I$, $E \, V$, $E \, V\!I\!I\!I$, $F \, I$, $G$.

\smallskip
{\bf (3)}
La somme de deux racines imaginaires non compactes n'est jamais une racine non
compacte d'apr\`es \cite [(6.99) p.$\!$~352] {Kn96}.
Ainsi, $\gstssInc$ est form\'e des $\lambda \! \in \gothg^*$ semi-simples telles
que les id\'eaux simples de l'alg\`ebre de Lie r\'eductive $\gothg(\lambda)$
sont isomorphes \`a $\,\goth{sl}(2,\R )$.
\cqfr

}

\medskip
\`A partir de certains param\`etres enrichis $(\lambdatilde,{\gotha^*}^+)$, je
vais construire canoniquement deux syst\`emes de racines positives pour
$(\gothg_\C ,\gothh_\C )$, qui auront les m\^emes racines non complexes.
Le syst\`eme de racines positives $R^+(\gothg_\C ,\gothh_\C )$ est attach\'e \`a
une certaine forme lin\'eaire r\'eguli\`ere $\lambda_+$.
Il permettra, par des passages \`a la limite (section~\ref{Rep11}) et par une
translation (section \ref{Rep13}), de construire dans la partie \ref{RepIII}
des repr\'esentations en se ramenant au cas d'un caract\`ere infinit\'esimal
r\'egulier.
Le syst\`eme de racines positives $R^+_{\lambdatilde,{\gotha^*}^+}\!$ a pour
sous-ensemble de racines complexes \`a conjugu\'ee n\'egative, une partie
invariante sous l'action du groupe $G(\lambdatilde)$ d\'efini ci-dessous.
Il interviendra dans la d\'efinition \ref{parametres adaptes} (c) et dans le
lemme \ref{resultat inattendu} (b).
On verra dans le lemme \ref{nouveaux parametres} une autre param\'etrisation
des repr\'esentations, \`a \'equivalence pr\`es, dans laquelle le param\`etre
$(\lambdatilde,{\gotha^*}^+)$ sera remplac\'e par $\lambdatilde$, et le r\^ole
de $(R^+_{\lambdatilde,{\gotha^*}^+},\lambda_+)$ sera tenu par un certain couple
$(R^+_{\lambdatilde},\lambda_{\textit{can}})$ d\'eduit de $\lambdatilde$.

\definition{
\label{choix de racines positives}

Soient $\lambda \! \in \gstss$, $\gothh \in \Car\gothg(\lambda)$ et
$\,\F+\! \!\in C(\gothg(\lambda),\gothh)$.
On note $\mu$ et $\nu$ (resp. $\gotht$ et $\gotha$) les composantes
infinit\'esimalement elliptique et hyperbolique de $\lambda$ (resp. $\gothh$).
On pose $\lambdatilde = (\lambda,\F+)$\labind{zzl tilde}.

\smallskip
{\bf (a)}
On note $G(\lambdatilde)$\labind{G (zzl tilde)aa} le normalisateur de $\F+\!$
dans $\, G(\lambda)$ et $\gothg(\lambdatilde)$\labind{g (zzl tilde)} son
alg\`ebre de Lie.
On associe \`a $\lambdatilde$ la demi-somme
$\rho_\F+\! \in \mi\gotht^*$\labind{zzr F+}
des $\, \alpha \in R(\gothg(\lambda)_\C ,\gothh_\C ) \,$ imaginaires
tels que $\,\F+(H_{\alpha}) \subseteq \R ^+ \moins0$.
(Ainsi $\,\gothh$ est une sous-alg\`ebre de Cartan de
$\gothg(\lambda)(\mi \rho_\F+)$ sans racine imaginaire,
et donc $G \cdot (\lambda,\mi \rho_\F+)$ d\'etermine $G \cdot \lambdatilde$.)

On fixe pour la suite de cette d\'efinition une chambre ${\gotha^*}^+\!$ de
$(\gothg(\lambda)(\mi \rho_\F+),\gotha)$.

\smallskip
{\bf (b)}
On introduit le syst\`eme de racines positives suivant :

\smallskip\centerline{
$R^+(\gothg(\lambda)(\mi \rho_\F+)_\C ,\gothh_\C )
= \{ \, \alpha \! \in \! R(\gothg(\lambda)(\mi \rho_\F+)_\C ,\gothh_\C ) \mid
{\gotha^*}^+ (H_{\alpha}) \subseteq \R ^+ \!\moins0 \, \}$.
}

\smallskip
On fixe $\epsilon > 0$ assez petit pour que, en posant
$\; \nu_+
= \nu + \epsilon \, \rho_{\gothg(\lambda)(\mi \rho_\F+),\gothh} \in \gotha^*$%
\labind{zzn +}
on ait :
$\;\, a \cdot \nu_+ \not= \nu_+$ pour tout automorphisme $a$ du syst\`eme de
racines $R(\gothg_\C ,\gothh_\C )$ tel que $a \cdot \nu \not= \nu$.
On choisit m\^eme $\epsilon$ de fa\c{c}on que la propri\'et\'e pr\'ec\'edente
reste valable quand on remplace $\epsilon$ par $t\epsilon$,
$\,t \!\in\! \left]0,1\right]$.
(On a donc :
$\; \gothg(\nu_+) = \gothg(\nu)(\rho_{\gothg(\lambda)(\mi \rho_\F+),\gothh})$.)

On utilisera les syst\`emes de racines positives suivants :

\centerline{
$R^+(\gothg(\nu_+)_\C ,\gothh_\C )
= \left\{\, \alpha \in R(\gothg(\nu_+)_\C ,\gothh_\C ) \;\,\Big|\;\,
\mi \mu(H_{\alpha}) > 0
\textrm{ ou }
\Big\{\begin{smallmatrix}
\mi \mu(H_{\alpha}) = 0 \\ \textrm{et} \hfill \\ \rho_\F+\!(H_{\alpha}) > 0
\end{smallmatrix}
\,\right\}$
}

\smallskip\noindent
et
$\quad R^+(\gothg_\C ,\gothh_\C )
= \{ \alpha \in R(\gothg_\C ,\gothh_\C ) \mid
\nu_+ (H_{\alpha}) \! > \! 0 \}
\, \cup \, R^+(\gothg(\nu_+)_\C ,\gothh_\C )$%
\labind{R+(gC,hC)noncan}.

On pose
$\, \mu_+ = \mu - 2\mi \rho_{\gothg(\nu_+),\gothh} \in \gotht^*$\labind{zzm +}
et
$\;  \lambda_+ = \mu_+\!+\nu_+$\labind{zzl +}
\\
(donc
$\,R^+(\gothg(\nu_+)_\C ,\gothh_\C )
= \{\, \alpha \in R(\gothg(\nu_+)_\C ,\gothh_\C ) \mid
\mi \mu_+(H_{\alpha}) > 0 \}$).

\smallskip
Dans certains cas on notera
$\lambda_{\gothg,\lambdatilde,{\gotha^*}^+\!,\epsilon}$,
$\mu_{\gothg,\lambdatilde,{\gotha^*}^+}\!$ et
$\nu_{\gothg,\lambdatilde,{\gotha^*}^+\!,\epsilon}$
pour $\,\lambda_+\!$, $\mu_+\!$ et~$\nu_+\!$.

\smallskip
{\bf (c)}
On introduit l'ensemble

\centerline{
$R^+_\lambdatilde
= \left\{ \alpha \in R(\gothg_\C ,\gothh_\C ) \;\Big|\;
\nu (H_{\alpha}) \! > \! 0
\textrm{ ou }
\Big\{ \begin{smallmatrix}
\nu(H_{\alpha}) = 0 \\ \textrm{et} \hfill \\ \mi \mu(H_{\alpha}) > 0
\end{smallmatrix}
\textrm{ ou }
\Big\{\begin{smallmatrix}
\lambda(H_{\alpha}) = 0 \\ \textrm{et} \hfill \\ \rho_\F+\!(H_{\alpha}) > 0
\end{smallmatrix}
\right\}$%
\labind{R+lambdatilde},
}

\noindent
la demi-somme $\rho_{\textit{can}}$ des \'el\'ements de
$R^+_\lambdatilde \,\cap\, R(\gothg(\nu)_\C ,\gothh_\C )$,
et le syst\`eme de racines positives
$\; R^+_{\lambdatilde,{\gotha^*}^+}
= R^+_\lambdatilde \,\cup\,
R^+(\gothg(\lambda)(\mi \rho_\F+)_\C ,\gothh_\C ) \,$%
\labind{R+(gC,hC)can}
de $R(\gothg_\C ,\gothh_\C )$.
 
On pose
$\, \mu_{\textit{can}} = \mu - 2\mi \rho_{\textit{can}} \in \gotht^*$
et
$\; \lambda_{\textit{can}} = \mu_{\textit{can}}\!+\nu$\labind{zzl can}
\\
(donc en utilisant \cite [cor. 4.69 p.$\!$~271] {KV95}, on constate que
$\,\gothg(\lambda_{\textit{can}}) = \gothg(\lambda)(\mi \rho_\F+)\,$
et
$\; R^+_\lambdatilde \,\cap\, R(\gothg(\nu)_\C ,\gothh_\C )
= \{\, \alpha \in R(\gothg(\nu)_\C ,\gothh_\C ) \mid
\mi \mu_{\textit{can}}(H_{\alpha}) > 0 \}$).

}

\medskip
Voici maintenant le lemme clef qui permettra de se ramener de $\gstregtilde$
\`a $\gstssreg\!$.

\lemme{
\label{lemme clef}

On se place dans les conditions de la d\'efinition pr\'ec\'edente. \\
Soit $\,t \!\in\! \left]0,1\right]$.
On pose
$\; \nu_t = \nu + t\epsilon \, \rho_{\gothg(\lambda)(\mi \rho_\F+),\gothh}$,
$\mu_t = \mu - 2\mi t\rho_{\gothg(\nu_+),\gothh}$
et
$\lambda_t = \mu_t+\!\nu_t$.

On a : $\;\, \gothg(\lambdatilde) = \gothh$, $\; \lambda_t \in \gstssreg \,$
et $\; \normagr{G(\lambdatilde)}{{\gotha^*}^+} = G(\lambda_t)$.

Donc $\, G(\lambda_+) = G(\lambdatilde) \,$ quand $\, {\gotha^*}^+ = \gotha^*$
(par exemple quand $\, \lambdatilde \in \gstItilde$).

}

\dem{D\'emonstration du lemme}{

On a
$\;\, G(\lambdatilde) \subseteq \normagr{G}{\gothh}$,
donc
$\; \gothg(\lambdatilde) \subseteq \gothh \;$
puis
$\; \gothg(\lambdatilde) = \gothh$.

Par construction de $\epsilon$, on a
$\; \gothg(\nu_t)
= \gothg(\nu)(\rho_{\gothg(\lambda)(\mi \rho_\F+),\gothh})
=\gothg(\nu_+)$.
En outre, on a
$\; \mi \mu_t(H_{\alpha})
= \mi \mu(H_{\alpha}) + 2t\rho_{\gothg(\nu_+),\gothh}(H_{\alpha})
> 0 \;$
pour tout $\, \alpha \in R^+(\gothg(\nu_+)_\C ,\gothh_\C )$.
Donc $\; \lambda_t \in \gstssreg$.
Un argument de continuit\'e en $t$ permet d'en d\'eduire que :

\centerline{
$R^+(\gothg_\C ,\gothh_\C )
= \left\{\, \alpha \in R(\gothg_\C ,\gothh_\C ) \;\,\Big|\;\,
\nu_t(H_{\alpha}) > 0
\textrm{ ou }
\Big\{\begin{smallmatrix}
\nu_t(H_{\alpha}) = 0 \\ \textrm{et} \hfill \\ \mi \mu_t(H_{\alpha}) > 0
\end{smallmatrix}
\,\right\}$.
}

L'inclusion
$\, \normagr{G(\lambdatilde)}{{\gotha^*}^+} \subseteq G(\lambda_t) \,$
est imm\'ediate.
Soit $g \in G(\lambda_t)$. \\
On a
$\, g \cdot \gothh = g \cdot \gothg(\lambda_t) = \gothh \,$
et $\restriction{(\Ad g)}{\gothh}$ est un automorphisme du syst\`eme de
racines $R(\gothg_\C ,\gothh_\C )$ tel que $\, g \cdot \nu_t = \nu_t$,
donc
$\, g \cdot (\nu_t - \nu) = (\nu_t - \nu)$
puis
$\, g \cdot \nu_+ = \nu_+$.
Comme $g$ laisse stable
$R(\gothg(\nu_+)_\C ,\gothh_\C ) \cap R^+(\gothg_\C ,\gothh_\C )$,
on a :
$\, g \cdot \rho_{\gothg(\nu_+),\gothh} = \rho_{\gothg(\nu_+),\gothh}$
puis
$\, g \cdot \mu = \mu$. \\
Ensuite, comme $g$ laisse stable l'ensemble des racines imaginaires de
$R(\gothg(\lambda)_\C ,\gothh_\C ) \cap R^+(\gothg_\C ,\gothh_\C )$,
on a :
$\, g \cdot \F+ =\F+$.
Enfin, comme $g$ laisse stable l'ensemble des $\restriction{\alpha}{\gotha}$
avec
$\alpha \in
R(\gothg(\lambda)(\mi \rho_\F+)_\C ,\gothh_\C )
\cap R^+(\gothg_\C ,\gothh_\C )$,
on a :
$\, g \cdot {\gotha^*}^+ = {\gotha^*}^+$.
\cqfd

}

\medskip
Le lemme qui suit expliquera pourquoi, dans certaines formules que j'\'ecrirai
plus tard dans un autre article, et dans lesquelles apparaissent en facteur les
valeurs absolues de Pfaffiens
$\, \bigl|\Pi_{\gothg,\centraalg{\gothg(\mu)}{\gotha}}\bigr|$,
seuls interviennent les \'el\'ements de $\gstssI$.

Pour rendre ces Pfaffiens plus accessibles au calcul, je pr\'ecise d'abord une
formule de \cite[p.$\!$~301] {DV88}.
La forme bilin\'eaire ${\scriptstyle \langle} ~,\!~ {\scriptstyle \rangle}$ de
la d\'efinition \ref{notations frequentes} (b) d\'etermine une mesure de Haar
$\diff _V$ sur tout sous-espace vectoriel $V$ de $\gothg$ sur lequel elle a une
restriction non d\'eg\'en\'er\'ee, par la condition
$\; \diff _V (0) (v_1 \wedge \cdots \wedge v_n)
= |\! \det ( \langle v_i,v_j \rangle )_{_{1 \leq i,j \leq n}}|^{\frac{1}{2}}\;$
pour toute base $(v_1,\dots,v_n)$ de~$V$.
Soit $\, \gothh \in \Car \gothg$.
On fixe $\lambda_0 = \mu_0 + \nu_0\in \gothh^*$ r\'eguli\`ere dans $\gothg^*\!$.
On lui associe le syst\`eme de racines positives

\centerline{
$\; R_0^+(\gothg_\C ,\gothh_\C )
= \left\{\, \alpha \in R(\gothg_\C ,\gothh_\C ) \;\,\Big|\;\,
\nu_0(H_{\alpha}) > 0
\textrm{ ou }
\Big\{\begin{smallmatrix}
\nu_0(H_{\alpha}) = 0 \\ \textrm{et} \hfill \\ \mi \mu_0(H_{\alpha}) > 0
\end{smallmatrix}
\,\right\}$.
}

\noindent
On prend pour $\diff _{\gothg}$ et $\diff _{\gothh}$ les mesures de Haar sur
$\gothg$ et $\gothh$ d\'eduites de
${\scriptstyle \langle} ~,\!~ {\scriptstyle \rangle}$.
Pour tout $\lambda \in \gothh^*$ et tout
$\, \omega_0 \in \bigwedge ^{2k}[\gothg,\gothh]\backslash \{0\} \,$
\'el\'ement de l'orientation
$\scalo (B_{\lambda_0})_{\gothg/\gothh}$
(cf. \ref{geometrie metaplectique} (b) et \ref{parametres de Duflo} (a))
de $(\gothg/\gothh,B_{\lambda_0})$, on a

\smallskip\centerline{
{\bf (\boldmath$*$)}
$\; \frac{1}{k!}
\left( \restriction{B_{\lambda}}{[\gothg,\gothh]^2} \right)^{\!k}
\!(\omega_0)
\,{\scriptstyle\times}\,
\left( \diff _{[\gothg,\gothh]} (0) (\omega_0) \right)^{\!\!-1} \!\!
= \,
\mi^{n_0}
\!\!\!
\Prod_{\alpha \in R_0^+(\gothg_\C ,\gothh_\C )} \!\!
\langle \lambda,\alpha \rangle$
}

\noindent
avec
$\, k = \frac {1}{2} \dim [\gothg,\gothh]$,
$\, \displaystyle \diff _{[\gothg,\gothh]}
= \diff _{\gothg} / \diff _{\gothh} \,$
et
$\, n_0
=\abs{\{\alpha \in R_0^+(\gothg_\C ,\gothh_\C ) \mid
\overline{\alpha} \notin R_0^+(\gothg_\C ,\gothh_\C )\}}$,
en notant encore
${\scriptstyle \langle} ~,\!~ {\scriptstyle \rangle}$\labind{<,>dual}
la forme bilin\'eaire sur $\gothh_\C ^{~*}$ duale de la restriction \`a
$\gothh_\C $ de~${\scriptstyle \langle} ~,\!~ {\scriptstyle \rangle}$.

\lemme{
\label{description d'ensembles de formes lineaires}

Soit $\, f \!\in \gothg^* \,$ de composantes infinit\'esimalement elliptique,
hyperbolique, semi-simple et nilpotente $\mu$, $\nu$, $\lambda=\mu\!+\!\nu$ et
$\xi$.
On fixe $\, \gothh \in \Car\gothg \,$ telle que $\, \lambda \in \gothh^* \!$.
On note $\gotha$ la composante hyperbolique de $\gothh$.

\smallskip
{\bf (a)}
On a $\; f \!\in \gstreg$
si et seulement si
$\;\, \xi \!\in \glambdastreg$.

\smallskip
{\bf (b)}
On a $\; \lambda \in \gstssI\,$ et
$\, \gothh \in \Car\gothg(\lambda) \,$ est fondamentale, \\
si et seulement si
$\;\; \bigl|\Pi_{\gothg,\centraalg{\gothg(\mu)}{\gotha}}\bigr|(\lambda)
\not= 0$.

}

\dem{D\'emonstration du lemme}{

{\bf (a)}
Il est clair que l'\'egalit\'e
$\dim \gothg (f) \!=\! \rg \gothg$
\'equivaut \`a
$\dim \gothg (\lambda) (\xi) \!=\! \rg \gothg (\lambda)$.

\smallskip
{\bf (b)}
Pour tous $X \in \gothh$ et $\lambda' \in \gothh^*\!$, les sous-espaces
vectoriels suppl\'ementaires
$\, [\gothg,\gothg(X)] = [\gothg,X] \,$
et
$\, [\gothg(X),\gothh] \subseteq \gothg(X) \,$
de $[\gothg,\gothh]$ \'etant \`a la fois orthogonaux pour $B_{\lambda'}$ et pour
${\scriptstyle \langle} ~,\!~ {\scriptstyle \rangle}$, on a :
$\;\, \bigl|\Pi_{\gothg,\gothh}\bigr| (\lambda')
= \bigl|\Pi_{\gothg,\gothg(X)}\bigr| (\lambda')
\,{\scriptstyle \times}\,
\bigl|\Pi_{\gothg(X),\gothh}\bigr|(\lambda')$.
En choisissant $X \in \gothh$ tel que
$\gothg(X) = \centraalg{\gothg(\mu)}{\gotha}$, on trouve d'apr\`es ($*$) que :

\centerline{
$\bigl|\Pi_{\gothg,\centraalg{\gothg(\mu)}{\gotha}}\bigr| (\lambda) \;\,
= \!\!\Prod_{\substack
{\alpha \in R(\gothg_\C ,\gothh_\C ) \\
\alpha \notin R(\centraalg{\gothg(\mu)}{\gotha}_\C ,\gothh_\C )}} \!\!\!
\abs{\langle \lambda,\alpha \rangle}^{\frac{1}{2}}$.
}

\negmedskip\noindent
Cela donne le r\'esultat.
\cqfd

}

\section{\boldmath
Les mesures $\beta_{G \cdot \lambdatilde}$
}\label{Rep2}

On va maintenant attacher \`a chaque orbite
$\,\wt{\Omega} \!\in\! G \backslash \gstregtilde$ une mesure de Radon
$\beta_{\wt{\Omega}}$ sur $\gothg^*\!$ de fa\c{c}on que
$\beta_{G\cdot\lambdatilde}$ co\"{\i}ncide avec la mesure canonique sur
$G\!\cdot\!\lambda$ quand $\,\lambda \in \gstssreg$, $\gothh=\gothg(\lambda)$ et
$\lambdatilde = (\lambda,\gothh_{(\R )}^{~~*})$.

\definition{

Soient $\lambda \!\in \gstss$, $\gothh \in \Car\gothg(\lambda)$
et $\, \F+\! \!\in C(\gothg(\lambda),\gothh)$.

{\bf (a)}
On note $\, \Fh+ = \gotha^* + {\gotht^*}^+ \!$\labind{Fh+}, o\`u $\,\gotht$ et
$\gotha$ sont les composantes infinit\'esimalement elliptique et hyperbolique de
$\gothh$ et $\F+\!$ s'\'ecrit $\, \F+\! = \gotha^* + \mi {\gotht^*}^+\!$ avec
${\gotht^*}^+ \subseteq \gotht^*\!$.

\smallskip
{\bf (b)}
On note
$\, \Supp_{\gothg^*}(G \cdot (\lambda,\F+))\,$\labind{SuppgstOmegatilde}
la  r\'eunion dans $\gstreg$ des $\, G \cdot (\omega+\{\lambda\})$,
o\`u $\omega$ d\'ecrit l'ensemble des orbites nilpotentes r\'eguli\`eres de
$G(\lambda)$ dans  $\gothg(\lambda)^*$ incluses dans l'adh\'erence de
$G(\lambda) \cdot \Fh+$
(cf. le lemme \ref{description d'ensembles de formes lineaires} (a) et la
remarque \ref{explication de la notation support} (1)).

}

\medskip
La proposition suivante est donn\'ee dans le but de faciliter l'utilisation
de ma formule du caract\`ere \ref{formule du caractere}.
Elle permettra aussi, dans certains cas, de passer directement des caract\`eres
des repr\'esentations aux param\`etres associ\'es (cf. \ref{recuperation des
parametres}).

\proposition{
\label{bijection entre orbites}

{\bf (a)}
L'ensemble $\, G \,\backslash \gstfondtilde \,$ est en bijection avec
$\, G \,\backslash \gstreg \,$ par l'application  $R_G$\labind{RG} qui envoie
$\, G \cdot \lambdatilde \,$ sur $\, \Supp_{\gothg^*}(G \cdot \lambdatilde)$.

\smallskip
{\bf (b)}
Soient $\lambda \in \gstss$, $\gothh \in \Car\gothg(\lambda)$
et $\, \F+\! \in C(\gothg(\lambda),\gothh)$. \\
La partie $\, \Supp_{\gothg^*}(G \!\cdot\! (\lambda,\F+)) \,$ de $\gstreg$ est
la r\'eunion des $R_G(G \!\cdot\! (\lambda,\Fp+))$
avec $\gothh' \in \Car\gothg(\lambda)$ fondamentale et
$\Fp+\! \in C(\gothg(\lambda),\gothh')_{reg}$ v\'erifiant
$\, G(\lambda) \!\cdot\! \rho_\F+ \cap \overline{\Fp+} \not= \emptyset$.

}

\dem{D\'emonstration de la proposition}{

{\bf (a)}
Soit $\, (\lambda,\F+) \in \gstfondtilde$.
D'apr\`es \cite [th. p.$\!$~217 et suppl. {\sc a}(b) p.$\!$~218] {Ro82}
appliqu\'e \`a $G_0(\lambda)$, il existe une unique orbite nilpotente
r\'eguli\`ere $\omega_0$ de $G_0(\lambda)$ dans
$\; \overline{G_0(\lambda) \cdot \Fh+}$.
Donc $\; \omega \egdef G(\lambda) \!\cdot\! \omega_0 \;$
est une orbite nilpotente r\'eguli\`ere de $G(\lambda)$ dans
$\; \overline{G(\lambda) \cdot \Fh+}$.
Toute orbite nilpotente r\'eguli\`ere $\omega'$
de $G(\lambda)$ dans $\; \overline{G(\lambda) \cdot \Fh+}$ coupe
$\, u \cdot \overline{G_0(\lambda) \cdot \Fh+}$
pour un certain $u \in G(\lambda)$,
donc  v\'erifie $\, \omega' \supseteq u \cdot \omega_0 \,$
puis $\omega' = \omega$. \\
Cela permet de d\'efinir $R_G$.

\smallskip
On abandonne maintenant les notations du d\'ebut de cette d\'emonstration.

Soit $\, f \in \gstreg$ de composantes semi-simple et nilpotente $\lambda$ et
$\xi$.

On fixe $\gothh \in \Car\gothg(\lambda)$ fondamentale.
D'apr\`es le lemme \ref{description d'ensembles de formes lineaires} (a) et
\cite [lem. {\sc a} p.$\!$~220] {Ro82}, on peut passer par une suite finie de
transformations de Cayley inverses (cf. \cite [p.$\!$~218] {Ro82}) de $\gothh$
\`a une sous-alg\`ebre de Cartan de $\gothg(\lambda)$ sans racine imaginaire.
D'apr\`es
\cite [suppl. {\sc a}(a) et {\sc b} et {\sc c} p.$\!$~218, et th. p.$\!$~221]
{Ro82}, il existe $\F+\! \!\in C(\gothg(\lambda),\gothh)_{reg}$  unique \`a
l'action de $\normagr{G_0(\lambda)}{\gothh}$ pr\`es pour lequel
$\, \omega_0 \egdef G_0(\lambda) \!\cdot\! \xi \,$
est l'orbite nilpotente r\'eguli\`ere de $G_0(\lambda)$
dans $\; \overline{G_0(\lambda) \!\cdot\! \Fh+}$.
Ainsi, $G \!\cdot\! f$ est l'image de
$\, G \!\cdot\! (\lambda,\F+) \,$ par $R_G$.

Soit $\, G \cdot (\lambda',\Fp+) \,$ un ant\'ec\'edent de $G \cdot f$ par $R_G$.
On a $\, \Fp+\! \in C(\gothg(\lambda'),\gothh')_{reg} \,$ pour une certaine
$\gothh' \in \Car\gothg(\lambda')$ fondamentale.
On note $\omega'_0$ l'orbite nilpotente r\'eguli\`ere de $G_0(\lambda')$ dans
$\; \overline{G_0(\lambda') \!\cdot\! \Fphp+}$.
Vu le d\'ebut de cette d\'emonstration, il existe $g \in G\,$ tel que
$g\xi \in \omega'_0$ et $g\lambda = \lambda'$.
Les sous-alg\`ebres de Cartan $\gothh$ et $g^{-1}\gothh'$ de $\gothg(\lambda)$
\'etant fondamentales, il existe $x \in G(\lambda)_0$ tel que
$g^{-1}\gothh' = x\gothh$.
Par cons\'equent, $\omega_0$ qui est \'egal \`a $g^{-1} \omega'_0$ est aussi
l'orbite nilpotente r\'eguli\`ere de $G_0(\lambda)$ dans
$\; \overline{G_0(\lambda) \!\cdot\! x^{-1}g^{-1}\Fp+\!}$.
D'o\`u :
$\;\, x^{-1}g^{-1} \!\cdot\! (\lambda',\Fp+) \in
\normagr{G_0(\lambda)}{\gothh} \!\cdot\! (\lambda,\F+)$.
Il s'ensuit que $G \!\cdot\! (\lambda,\F+)$ est l'unique ant\'ec\'edant
de $G \!\cdot\! f$ par $R_G$.

{\bf (b)}
Avec un argument du d\'ebut de la d\'emonstration du (a), on constate que
$\, \Supp_{\gothg^*}(G \cdot (\lambda,\F+)) \,$ est le satur\'e sous $G$ de
$\, \Supp_{\gothg^*}(G_0 \cdot (\lambda,\F+))$.
On est donc ramen\'e \`a prouver le (b) pour $G_0$ avec
$\lambda = 0$.
On se place dans ce cas.
On fixe une sous-alg\`ebre de Cartan fondamentale $\gothh'$ de $\gothg$ qui
contient la forme lin\'eaire infinit\'esimalement elliptique $-\mi \,\rho_\F+$.
On pose $\, l_0 = -\mi \,\rho_\F+$ dans~$\Fh+\!$.
On va reprendre rapidement les id\'ees de
\cite [proof of suppl. {\sc c} p.$\!$~228] {Ro82}.

D'apr\`es \cite [th. 23 p.$\!$~50] {Va77} et
\cite [$(L1)\!\Leftrightarrow\!(L2)\!\Leftrightarrow\!(L3)$ p.$\!$~216]{Ro82},
on a :

\centerline{
$\Sum_{\Omega \,\in\, G_0 \backslash \Supp_{\gothg^*}(G_0 \cdot (0,\F+))}
\beta_{\Omega} (\varphi)
\,=\,
\lim_{t \to 0^+}\,
\Sum_{\Omega' \in L_t}
\beta_{\Omega'} (\varphi) \;\,$
pour tout $\, \varphi \in C^\infty_c(\gothg^*)$,
}

\noindent
o\`u on utilise \ref{mesures de Liouville} (c) et pose
$\; L_t
= \Lim_{\substack {l' \in \Fh+ \cap \, \gstreg \\ l' \to tl_0}}\!
G_0 \cdot l' \;$
dans $\, G_0 \backslash \gstreg$ quand $t>0$.

Par \cite [lignes 2 \`a 9 p.$\!$~217] {Ro82},
l'application canonique de $G_0(l_0) \backslash \glzerostreg$ dans
$G_0 \backslash \gstreg$ se restreint en une bijection de
$\; \Lim_{\substack {l'' \in \F+[l_0]_\gothh \cap \, \glzerostreg \\
l'' \to tl_0}}\!
G_0(l_0) \cdot l'' \;$
sur $L_t$ pour tout $t>0$,
o\`u la limite est prise dans $G_0(l_0) \backslash \glzerostreg$ et $\F+[l_0]$
est l'\'el\'ement de $C(\gothg(l_0),\gothh)$ contenant $\F+\!$.

Les racines de $(\gothg(l_0)_\C ,\gothh_\C )$ sont non imaginaires et celles de
$(\gothg(l_0)_\C ,\gothh'_\C )$ sont non r\'eelles.
D'apr\`es \cite [suppl. {\sc c} p.$\!$~218] {Ro82}, on en d\'eduit que :

\centerline{
$\Lim_{\substack {l'' \in \F+[l_0]_\gothh \cap \, \glzerostreg \\
l'' \to tl_0}}\!
G_0(l_0) \cdot l''
\;=
\!\bigcup_{\dot{\Fp+} \,\in\, W(G_0(l_0),\gothh') \backslash \cale}\;
\Lim_{\substack {l'' \in \Fp+[l_0]_{\gothh'} \cap \, \glzerostreg \\
l'' \to tl_0}}\!
G_0(l_0) \cdot l''$
}

\noindent
o\`u
$\, \cale \egdef
\big\{\, \Fp+ \! \in C(\gothg,\gothh') \mid l_0 \in \overline{\Fphp+} \,\Big\}
\,$ et  $\Fp+[l_0]$ est l'\'el\'ement de $C(\gothg(l_0),\gothh')$
contenant $\Fp+\!$.
En outre, la r\'eunion pr\'ec\'edente est une r\'eunion d'ensembles deux \`a
deux disjoints au vu de \cite [suppl. {\sc b} p.$\!$~218] {Ro82}.

Les mesures de Radon $\beta_{\Omega}$ avec
$\, \Omega \in G_0 \backslash \gstreg \,$
(cf. \ref{mesures de Liouville} (c))
sont lin\'eairement ind\'ependantes.
\`A l'aide de \cite [th. p.$\!$~217] {Ro82}, on trouve que
$\Supp_{\gothg^*}(G_0 \cdot (0,\F+))$ est r\'eunion disjointe des
$R_{G_0}(G_0 \cdot (0,\Fp+))$ o\`u
$\dot{\Fp+} \!\in\, W(G_0(l_0),\gothh') \backslash C(\gothg,\gothh')_{reg}$
et $\rho_\F+ \in \overline{\Fp+}$.
Cela conduit au r\'esultat.
\cqfd

}

\definition{
\label{mesures de Liouville}

{\bf (a)}
On note $\,\cald_{\gothg}$\labind{Dg} la fonction sur $\gothg$ dont la valeur en
$\,X \in \gothg$ est le coefficient de $T^r$ dans $\det(T\id - \ad X)$, o\`u
$\,r$ est le rang de $\gothg$.
Elle est polyn\^omiale et invariante sous $\Ad G$.

\smallskip
{\bf (b)}
On note $\,\gssreg\!$\labind{g ssreg} l'ensemble des $X \in \gothg$ semi-simples
tels que $\gothg(X)$ est commutative.
Donc $\,\gssreg\!$ est l'ouvert dense de compl\'ementaire n\'egligeable de
$\gothg$ form\'e des points o\`u $\,\cald_{\gothg}$ ne s'annule pas (cf. \cite
[(2) et lem. 1 p.$\!$~9] {Va77}).

\smallskip
{\bf (c)}
Soit $\Omega \in G\,\backslash \gothg^*\!$.
On note $\beta_{\Omega}$\labind{zzb Omega} la mesure positive image sur
$\gothg^*$ de la mesure de Liouville de la vari\'et\'e symplectique $\Omega$
(cf. \cite [2.6 p.$\!$~20] {B.72}).
C'est une mesure de Radon sur $\gothg^*\!$ (cf. \cite [th. 2 p.$\!$~509] {Ra72}).

\smallskip
{\bf (d)}
Soit $\; \wt{\Omega}_0 \in G_0\,\backslash \gstregtilde$.
On pose
$\;\; \beta_{\wt{\Omega}_0}
= \Sumpetit_{\Omega_0 \,\in\,
G_0 \backslash \Supp_{\gothg^*}(\wt{\Omega}_0)} \beta_{\Omega_0}$.

{\bf (e)}
Soit $\; \wt{\Omega} \in G\,\backslash \gstregtilde$.
On pose
$\;\; \beta_{\wt{\Omega}}
= \Sumpetit_{\wt{\Omega}_0 \,\in\, G_0
\backslash\wt{\Omega}}\beta_{\wt{\Omega}_0}$%
\labind{zzb Omegatilde}
(cf. fin \ref{limites de mesures de Liouville} (a)).
\negsmallskip

}

\remarque{
\label{explication de la notation support}

{\bf (1)}
Soit $\, \wt{\Omega} \in G\,\backslash \gstregtilde$.
D'apr\`es le d\'ebut de la d\'emonstration de la proposition
\ref{bijection entre orbites} (b), $\Supp_{\gothg^*}(\wt{\Omega})$ est la
r\'eunion des $\Omega \in G\,\backslash \gstreg$ tels que
$\beta_{\wt{\Omega}}(\Omega) \not= 0$.

Il d\'ecoule aussi de la d\'efinition ci-dessus que
$\;\; \beta_{\wt{\Omega}}
= \Sumpetit_{\Omega_0 \,\in\,G_0\,\backslash \Supp_{\gothg^*}(\wt{\Omega})}
m_{\Omega_0} \, \beta_{\Omega_0} \;$
o\`u $\, m_{\Omega_0}$ est le nombre (non nul) de
$\wt{\Omega}_0 \in G_0 \backslash \wt{\Omega}$
v\'erifiant
$\Omega_0 \subseteq \Supp_{\gothg^*}(\wt{\Omega}_0)$.

\noindent
En particulier :
$\;\; \beta_{\wt{\Omega}} = \beta_{R_G(\wt{\Omega})} \;$
quand $\; \wt{\Omega} \in G\,\backslash \gstfondtilde$.

\smallskip
{\bf (2)}
L'exemple de $G$ \'egal au produit semi-direct de $\Z / 2\Z $ par $SL(2,\R )^2$
o\`u l'\'el\'ement non trivial de $\Z / 2\Z $ op\`ere sur $SL(2,\R )^2$ par
permutation des coordonn\'ees, et de $\wt{\Omega}$ de la forme
$G \cdot (0,\F+)$ pour une chambre $\F+$ associ\'ee \`a $\gothh \in \Car \gothg$
ni fondamentale ni d\'eploy\'ee, montre que les coefficients $m_{\Omega_0}$ du
(1) ci-dessus peuvent prendre des valeurs diff\'erentes de $1$.

\smallskip
{\bf (3)}
D'apr\`es (1) et \ref{bijection entre orbites} (a), l'application
$\wt{\Omega} \mapsto \beta_{\wt{\Omega}}$ est injective sur $G\,\backslash
\gstfondtilde$.

Quand $G = SL(3,\R )$, les deux \'el\'ements $\wt{\Omega}$ de
$G\,\backslash \gstregtilde$ de la forme $G \cdot (0,\F+)$ (dont l'un est dans
$G\,\backslash \gstfondtilde$) correspondent \`a une m\^eme mesure
$\beta_{\wt{\Omega}}$, \'egale \`a la mesure de Liouville de l'unique orbite
nilpotente r\'eguli\`ere de $G$ dans~$\gothg^*$.
\cqfr

}

\medskip
Voici un r\'esultat bien connu dont on peut trouver la d\'emonstration dans
\cite[prop. 15 p.$\!$~66, prop. 13 p.$\!$~65, th. 17 p.$\!$~66, th. 28
p.$\!$~95, et bas p.$\!$~105] {Va77}.

\proposition{
\label{mesures de Liouville temperees}

{\bf (a)}
La fonction $\, \abs{\cald_{\gothg}}^{-1/2}$ est localement int\'egrable sur
$\gothg$.

\smallskip
{\bf (b)}
Soit $\,\Omega \!\in\! G\,\backslash\, \gothg^*$ tel que $\beta_{\Omega}$ est
une mesure de Radon temp\'er\'ee.
La fonction g\'en\'eralis\'ee $\,{\widehat{\beta}}_{\Omega}$ est (la classe
modulo l'\'egalit\'e presque partout d') une fonction localement int\'egrable
sur~$\gothg$ et analytique sur $\gssreg$.

}

\medskip
La proposition suivante sera utile quand on voudra se ramener, par passage \`a la
limite ou par translation, \`a des formes lin\'eaires semi-simples
r\'eguli\`eres.

\proposition{
\label{limites de mesures de Liouville}

Soient $\lambda \in \gstss$, $\,\gothh \in \Car\gothg(\lambda)\,$ et
$\, \F+\! \in C(\gothg(\lambda),\gothh)$.
On pose $\lambdatilde = (\lambda,\F+)$.
On note $\gotha$ la composante hyperbolique de $\gothh$, et se donne une chambre
${\gotha^*}^+\!$ de $(\gothg(\lambda)(\mi \rho_\F+),\gotha)$.

\smallskip
{\bf (a)}
Les mesures positives $\, \beta_{\wt{\Omega}}$ avec
$\, \wt{\Omega} \!\in\! G\,\backslash \gstregtilde$
sont des mesures de Radon temp\'er\'ees.
Soit $\; \varphi\in \cals(\gothg^*)$.
On a

\negmedskip\centerline{
$\Lim_{\substack
{\lambda' \in \Fh+ \cap \, \gstreg \\ \lambda' \to \lambda}}
\bigl(
{\scriptstyle \abs{W(G,\gothh) \, (\lambda')}} \;\;
\beta_{G\cdot \lambda'}(\varphi) \bigr)
= \cases{
{\scriptstyle \abs{W(G,\gothh) \, (\lambda_+)}} \;\;
\beta_{G \cdot \lambdatilde}(\varphi)
& si $\; \lambdatilde \in \gstregtilde$
\cr \hfil 0
& sinon.
}$
}

\noindent
En particulier
$\;\; \beta_{G \cdot \lambdatilde}(\varphi)
= \Lim_{t \to 0^+} \; \beta_{G \cdot \lambda_t}(\varphi) \;\,$
quand $\; \lambdatilde \in \gstregtilde \;$ (cf. \ref{lemme clef}).

\smallskip
{\bf (b)}
Pour tout $\, X \in \gssreg$, on a en tenant compte de
\ref{mesures de Liouville temperees} (b)

\centerline{
$\Lim_{\substack
{\lambda' \in \Fh+ \cap \, \gstreg \\ \lambda' \to \lambda}}
\bigl( {\scriptstyle \abs{W(G,\gothh) \, (\lambda')}} \;\;
{\widehat{\beta}}_{G\cdot \lambda'}(X) \bigr)
= \cases{
{\scriptstyle \abs{W(G,\gothh) \, (\lambda_+)}} \;\;
{\widehat{\beta}}_{G \cdot \lambdatilde}(X)
& si $\; \lambdatilde \in \gstregtilde$
\cr \hfil 0
& sinon.
}$
}

\noindent
En particulier
$\;\; \beta_{G \cdot \lambdatilde}(X)
= \Lim_{t \to 0^+} \; \beta_{G \cdot \lambda_t}(X) \;\,$
quand $\; \lambdatilde \in \gstregtilde \;$ (cf. \ref{lemme clef}).

Plus pr\'ecis\'ement, \'etant donn\'es $\gothj \in \Car\gothg$,
$y \in \interieur (\gothg_\C )$ pour lequel $\, \gothj_\C  = y \gothh_\C $,
un syst\`eme de racines positives $R^+(\gothg_\C ,\gothj_\C )$ de
$R(\gothg_\C ,\gothj_\C )$,
une composante connexe $\Gamma$ de $\gothj \cap \gssreg$
et
une composante connexe $\,{\gothh^*}^+$ de $\, \Fh+ \cap \gstreg \,$ \`a
laquelle $\lambda$ est adh\'erent,
il existe une famille $\bigl(c_w\bigr)_{w \in W(\gothg_\C ,\gothj_\C )}\!$
de nombres complexes telle que, pour tout $X \!\in\! \Gamma\!$, on ait

\smallskip
$\Sum_{w \in W(\gothg_\C ,\gothj_\C )}
c_w \;\; \me^{\mi \,wy\lambda'\,(X)}
= \Prod_{\alpha \in R^+(\gothg_\C ,\gothj_\C )}
\!\! \alpha(X) \;\;
{\scriptstyle \times} \;\;
{\widehat{\beta}}_{G_0\cdot \lambda'}(X) \;\,$
quand $\, \lambda' \in {\gothh^*}^+$
\\
et

\negmedskip
$\Sum_{w \in W(\gothg_\C ,\gothj_\C )}
c_w \;\, \me^{\mi \,wy\lambda\,(X)}
= \cases{
\Prod_{\alpha \in R^+(\gothg_\C ,\gothj_\C )}
\!\! \alpha(X) \;\;
{\scriptstyle \times} \;\,
{\widehat{\beta}}_{G_0 \cdot \lambdatilde}(X)
& si $\; \lambdatilde \in \gstregtilde$
\cr \; \hfil 0
& sinon.
}$

}

\dem{D\'emonstration de la proposition}{

{\bf (a)}
Soit $\, \varphi \in C^\infty_c(\gothg^*)$ positive.
D'apr\`es
\cite [$(L1)\!\Leftrightarrow\!(L2)\!\Leftrightarrow\!(L3)$ p.$\!$~216,
lignes 9 \`a 12 p.$\!$~217, et th. p.$\!$~217] {Ro82}, on a

\centerline{
$\Lim_{\substack
{\lambda' \in \Fh+ \cap \, \gstreg \\ \lambda' \to \lambda}}
\beta_{G_0\cdot\lambda'}(\varphi)
= \cases{
\beta_{G_0 \cdot \lambdatilde}(\varphi)
& si $\; \lambdatilde \in \gstregtilde$
\cr \hfil 0
& sinon.
}$
}

\noindent
On fixe une norme $\norm{\cdot}$ sur $\gothg^*\!$.
D'apr\`es \cite [(i) p.$\!$~40] {Va77}, il existe $N \!\in \N $ tel que la
famille form\'ee des
$ \;\Int_{\gothg^*} (1 + \norm{l})^{-N} \diff \beta_{\Omega_0}(l) \;$
avec
$\, \Omega_0 \!\in\! G_0 \backslash \gstssreg \,$
est born\'ee.
On en d\'eduit que~$\, \beta_{G_0 \cdot \lambdatilde} \,$ est une mesure de
Radon temp\'er\'ee quand $\, \lambdatilde \in \gstregtilde$,
et que le passage \`a la limite ci-dessus reste valable en rempla\c{c}ant la
condition << $\varphi \in C^\infty_c(\gothg^*)$ positive >> par
<< $\varphi \in \cals(\gothg^*)$ >>.

On obtient ensuite, pour tout $\varphi \in \cals(\gothg^*)$ :

\centerline{
$\Lim_{\substack
{\lambda' \in \Fh+ \cap \, \gstreg \\ \lambda' \to \lambda}} \;
{\textstyle \Sum_{\dot{x} \in G/G_0} \!\! \beta_{G_0\cdot x\lambda'}
(\varphi)} = \cases{
\Sumpetit_{\dot{x} \in G/G_0} \!\! \beta_{G_0 \cdot x\lambdatilde} (\varphi)
& si $\; \lambdatilde \in \gstregtilde$
\cr \hfil 0
& sinon.
}$
}

\noindent
On transforme d'abord la somme de gauche pour $\lambda' \in \Fh+ \cap \gstreg$ :

\centerline{
$\Sum_{\dot{x} \in G/G_0} \!\! \beta_{G_0\cdot x\lambda'}\textstyle
\,=\, {\scriptstyle \abs{G(\lambda') G_0 / G_0}} \;\;
\beta_{G\cdot\lambda'}
\,=\, \frac{\abs{W(G,\gothh)(\lambda')}}
{\abs{\centragr{G}{\gothh} / \centragr{G_0}{\gothh}}} \;\;
\beta_{G\cdot\lambda'}$.
}

\noindent
Par ailleurs, les chambres de $(\gothg(\lambda)(\mi \rho_\F+),\gotha)$ sont
conjugu\'ees sous $\normagr{G(\lambda)(\mi \rho_\F+)_0}{\gotha}$ et {\em a
fortiori} sous $G_0(\lambdatilde)$, ce qui prouve que
$\;G(\lambdatilde) =
G_0(\lambdatilde)\;\normagr{G(\lambdatilde)}{{\gotha^*}^+}$. \\
Quand $\; \lambdatilde \in \gstregtilde$, on a donc d'autre part \`a l'aide du
lemme \ref{lemme clef} :

\centerline{
$\Sum_{\dot{x} \in G/G_0} \!\! \beta_{G_0\cdot x\lambdatilde}\textstyle
\,=\, {\scriptstyle \abs{G(\lambdatilde) G_0 / G_0}} \;\;
\beta_{G\cdot\lambdatilde}
\,=\, \frac{\abs{W(G,\gothh)(\lambda_+)}}
{\abs{\centragr{G}{\gothh} / \centragr{G_0}{\gothh}}} \;\;
\beta_{G\cdot\lambdatilde}$.
}

\noindent
Cela fournit (a).

\smallskip
{\bf (b)}
L'avant derni\`ere \'egalit\'e de cet \'enonc\'e est extraite de
\cite [th. 4 p.$\!$~108]{Va77}.
Elle montre que $\, {\widehat{\beta}}_{G_0 \cdot \lambda'}$ converge
simplement sur $\gssreg$ vers une fonction $G_0$-invariante continue
quand $\lambda' \to \lambda$ avec $\lambda' \in {\gothh^*}^+\!$.
Soit $\psi \in C^\infty_c(\gssreg)$.
D'apr\`es la d\'emonstration de (a), on a

\smallskip\centerline{
$\Lim_{\substack
{\lambda' \in {\gothh^*}^+ \\ \lambda' \to \lambda}}
\Int_{\gssreg}
{\widehat{\beta}}_{G_0 \cdot \lambda'}(X) \;
\psi(X) \,\diff _{\gothg}(X)
= \cases{
{\widehat{\beta}}_{G_0 \cdot \lambdatilde}(\psi\,\diff _{\gothg})
& si $\; \lambdatilde \in \gstregtilde$
\cr \hfil 0
& sinon.
}$
}

Compte tenu des majorations
$\, \abs{c_w \; \me^{\mi \,wy\lambda'\,(X)}} \leq \abs{c_w} \,$
de \cite [th. 7 p.$\!$~111] {Va77}, on voit que les fonctions $G_0$-invariantes
$\, \abs{\cald_{\gothg}}^{1/2}
{\scriptstyle \times} \,
{\widehat{\beta}}_{G_0 \cdot \lambda'} \,$
sur $\gssreg$ avec $\, \lambda' \in {\gothh^*}^+$ sont uniform\'ement
born\'ees.
La proposition \ref{mesures de Liouville temperees} (a) permet donc
d'appliquer le th\'eor\`eme de convergence domin\'ee de Lebesgue dans le passage
\`a la limite ci-dessus.
D'o\`u la derni\`ere \'egalit\'e de l'\'enonc\'e.

Le passage de $G_0$ \`a $G$ se fait comme dans la d\'emonstration du (a).
\cqfd

}

\section{
Points fix\'es par un \'el\'ement elliptique
}\label{Rep3}

On se donne un \'el\'ement elliptique $e$ de $G$.

L'alg\`ebre de Lie $\gothg(e)$ est r\'eductive car $\gothg(e)_\C $ a pour
forme r\'eelle $\gothc(e)$, o\`u $\gothc$ est l'alg\`ebre de Lie d'un
sous-groupe compact maximal de $\G (\C )$ qui contient la projection dans
$\G (\R )$ de l'\'el\'ement $e$ de $G$.
De plus toute sous-alg\`ebre de Cartan $\gothh_e$ de $\gothg(e)$ coupe
$\,\gssreg$, ce qui signifie que le commutant
$\gothh$ de $\gothh_e$ dans $\gothg$ est une sous-alg\`ebre de Cartan de
$\gothg$ (cf. \cite [lem. 1.4.1 p.$\!$~6] {Bo87}).
Par cons\'equent $G(e)$ v\'erifie l'hypoth\`ese de l'introduction portant sur
$G$ (cf. \cite [bas p.$\!$~38] {DV93}).

\definition{

On note $\,\gstreg(e)\,$\labind{g*reg (e)}
(respectivement $\,\gstItilde(e)$\labind{g*I tilde(e)},
$\,\gstInctilde(e)$\labind{g*Inc tilde(e)},
$\,\gstfondtilde(e)\,$\labind{g*fond tilde(e)} et
$\,\gstregtilde(e)$\labind{g*reg tilde(e)})
l'ensemble des \'el\'ements de $\gstreg$
(respectivement $\gstItilde$, $\gstInctilde$, $\gstfondtilde$ et $\gstregtilde$)
qui sont fixes sous l'action de $\,e$.

}

\medskip
Parmi les r\'esultats du lemme qui suit, seul le point (a) est important.
Le point (b) sera utilis\'e dans la remarque \ref{recuperation des parametres}
et le point (c) dans un article ult\'erieur.%

\lemme{
\label{descente pour les formes lineaires}

Soient $\lambda \in \gstss$, $\,\gothh \in \Car\gothg(\lambda)\,$ de
composantes infinit\'esimalement elliptique et hyperbolique $\gotht$ et
$\gotha$, et  $\, \F+\! \!\in C(\gothg(\lambda),\gothh) \,$ tels que $e$ fixe
$\; \lambdatilde \egdef (\lambda,\F+)$.
On suppose qu'il existe une chambre ${\gotha^*}^+\!$ de
$(\gothg(\lambda)(\mi \rho_\F+),\gotha)$ stable par $e$.

\smallskip
{\bf (a)}
On a $\,\lambda \in \gestss$, $\gothh(e) \in \Car\gothg(e)(\lambda)$,
et l'\'el\'ement $\rho_\F+\!$ de $\gothh(e)_{(\R )}^{~~*}$ est r\'egulier
pour les racines imaginaires de $(\gothg(e)(\lambda)_\C ,\gothh(e)_\C )$.

On note par d\'efinition $\F+[e]$\labind{F+[e]} l'\'el\'ement de
$C(\gothg(e)(\lambda),\gothh(e))$ contenant $\rho_\F+\!$ et
$\;\, \lambdatilde[e] = (\lambda,\F+[e])$\labind{zzl tilde[e]}.
(On dira que << $\lambdatilde[e]$ existe >> quand $e$ stabilise une chambre~%
${\gotha^*}^+\!$.)

\smallskip
{\bf (b)}
Il existe $\, e_0 \!\in e\exp\gotht(e)$ tel que :
si $\, \Fp+\! \!\in C(\gothg(\lambda),\gothh)$ et
$\, \rho_\Fp+\! \in \F+[e_0]$, alors $\, \Fp+ \!=\! \F+\!$.
Quand $\lambdatilde \in \gstItilde$ on peut m\^eme choisir $e_0$ tel
que $\, \lambdatilde[e_0] \!\in \gezerostItilde$.

\smallskip
{\bf (c)}
On suppose que $\, \lambdatilde \in \gstItilde$.
Le cardinal de l'ensemble des $\lambdatilde' \!\in \gstItilde(e)$ v\'erifiant
$\lambdatilde'\![e] \!=\! \lambdatilde[e]$ est \'egal \`a
$\frac{|W(\gothg(\lambda)_\C ,\gothh_\C )(\Ad^* e)|}
{|W(\gothg(e)(\lambda)_\C ,\gothh(e)_\C )|}\!$,
o\`u $\, W(\gothg(\lambda)_\C ,\gothh_\C )(\Ad^* e) \,$ d\'esigne le commutant
de $\, (\Ad^* e^\C )_{\gothh_\C ^*} \,$ dans
$\, W(\gothg(\lambda)_\C ,\gothh_\C )$.

}

\dem{D\'emonstration du lemme}{

{\bf (a)}
D'apr\`es le lemme \ref{lemme clef}, on a :
$\; \lambda_+ \in \gestssreg$,
puis $\, \gothh(e) \!=\! \gothg(e)(\lambda_+) \in \Car\gothg(e)(\lambda)$.
Soit $\, \alpha' \in R(\gothg(e)(\lambda)_\C ,\gothh(e)_\C ) \,$ imaginaire.
Il existe $\, \alpha \in R(\gothg(\lambda)_\C ,\gothh_\C ) \,$ tel que
$\, \alpha' = \restriction{\alpha}{\gothh(e)_\C }$.
On a :
$\,\langle \nu_+,\alpha \rangle = \langle \nu_+,\alpha' \rangle = 0\,$
et
$\,\langle \rho_\F+,\alpha \rangle = \langle \rho_\F+,\alpha' \rangle$.
Donc $\; \langle \rho_\F+,\alpha' \rangle \not= 0$,
car $R^+(\gothg(\nu_+)_\C ,\gothh_\C )$ est un syst\`eme de racines.

\smallskip
{\bf (b)}
Soit $B$ la base de $R(\centraalg{\gothg}{\gotha}(\lambda)_\C ,\gothh_\C )$
associ\'ee \`a $\F+\!$.
Pour chaque $\alpha \!\in\! B$, on note $m_{\dot{\alpha}}$ le cardinal de
l'orbite ${\dot{\alpha}}$ de $\alpha$ sous l'action dans $B$ du sous-groupe
$\langle e \rangle$ de $G$ engendr\'e par $e$, et on choisit
$\, t_{\dot{\alpha}} \in \R \,$ tel que
$\; ((\Ad e^\C )^{\,m_{\dot{\alpha}}})_{\gothg_\C ^\alpha}
= \me^{\mi \,t_{\dot{\alpha}}} \id$.
On prend $\, e_0 = e \exp (X_0)$, o\`u $X_0$ est l'unique \'el\'ement de
$\, \gothh_\C  \cap \derivealg{\centraalg{\gothg}{\gotha}\,(\lambda)_\C } \,$ tel
que $\, \alpha(X_0) = -i t_{\dot{\alpha}} / m_{\dot{\alpha}} \,$ pour tout
$\alpha \!\in\! B$.
Il est clair que $\lambdatilde[e_0]$ existe.
\\
Chaque $\restriction{\alpha}{\gothh(e)_\C }$ avec
$\,\alpha \in B \,$ est un poids de
$\gothh(e)_\C $ dans $\gothg(e_0)_\C  $
dont l'espace propre contient
$(1 + e_0 + \cdots + {e_0}^{\! m_{\dot{\alpha}} - 1}) \cdot
\gothg_\C ^{\alpha}$.
Soit $\Fp+\! \in C(\gothg(\lambda),\gothh)$ tel que $\rho_\Fp+\! \in \F+[e_0]$,
(donc  $e_0\Fp+\! = \Fp+\!$).
Pour tout $\alpha \in B$, on a :
$\, \langle \rho_\F+,\restriction{\alpha}{\gothh(e)_\C } \rangle
= \langle \rho_\F+,\alpha \rangle
> 0 \,$
et
$\, \langle \rho_\Fp+,\alpha \rangle
= \langle \rho_\Fp+,\restriction{\alpha}{\gothh(e)_\C } \rangle$,
donc la racine $\restriction{\alpha}{\gothh(e)_\C }\!$ de
$R(\centraalg{\gothg(e_0)}{\gotha(e)}(\lambda)_\C ,\gothh(e)_\C )$
est positive relativement \`a  $\F+[e_0]$
et ensuite
$\, \langle \rho_\Fp+,\alpha \rangle > 0$.
D'o\`u $\, \Fp+\! = \F+\!$.

On suppose maintenant que  $\lambdatilde \in \gstItilde$.
Toute racine de $(\gothg(e_0)(\lambda)_\C ,\gothh(e)_\C )$ est restriction
d'une racine de $(\gothg(\lambda)_\C ,\gothh_\C )$.
Les racines $\, \alpha' = \restriction{\alpha}{\gothh(e)_\C }$ avec
$\,\dot{\alpha} \in \langle e \rangle \backslash B \,$ d\'ecrivent donc
lorsque $\dot{\alpha}$ varie
(de fa\c{c}on injective car $\dim \gothg(e_0)_\C ^{\alpha'} = 1$)
l'ensemble des racines simples relativement \`a $\F+[e_0]$ dans
$R(\gothg(e_0)(\lambda)_\C ,\gothh(e)_\C )$.
Ces racines $\alpha'$ sont non compactes car
$\alpha([\, \cdot \, ,\overline{\vrule width 0ex height 1ex \, \cdot \,}])$ est
positive sur le sous-espace vectoriel $\gothg(e_0)_\C ^{\alpha'}$ de la somme
directe orthogonale
$\, \gothg_\C ^{\alpha} \oplus e_0 \gothg_\C ^{\alpha} \oplus \cdots
\oplus {e_0}^{\! m_{\dot{\alpha}} - 1} \gothg_\C ^{\alpha}$.

\smallskip
{\bf (c)}
L'application
$\,\pi_e : W(\gothg(\lambda)_\C ,\gothh_\C ) \cdot
\rho_\F+ \cap \gothh(e)_\C ^*
\to C(\gothg(e)(\lambda),\gothh(e))\,$
qui \`a une forme lin\'eaire associe
la chambre qui la contient, commute \`a l'action du groupe
$\;W(\gothg(e)(\lambda)_\C ,\gothh(e)_\C )$.
Les groupes $\;W(\gothg(\lambda)_\C ,\gothh_\C )(\Ad^* e)\;$ et
$\;W(\gothg(e)(\lambda)_\C ,\gothh(e)_\C )\;$ op\`erent simplement
transitivement respectivement sur les ensembles de d\'epart et d'arriv\'ee de
$\pi_e$.
Les fibres de $\pi_e$ ont donc un cardinal \'egal \`a
$\,\frac{|W(\gothg(\lambda)_\C ,\gothh_\C )(\Ad^* e)|}
{|W(\gothg(\lambda)(e)_\C ,\gothh(e)_\C )|}$.
L'ensemble des
$\lambdatilde' \!=\! (\lambda',\Fp+) \in \gstItilde(e)$
tels que
$\lambdatilde'\![e] = \lambdatilde[e]$,
dont les \'el\'ements v\'erifient $\lambda' = \lambda$
et
$\gothg(\lambdatilde')
= \centraalg{\gothg}{\gothg(e)(\lambdatilde'\![e])}
= \gothh$,
est en bijection avec
$\pi_e^{-1}(\{\F+[e]\})$ par l'application
$(\lambda',\Fp+) \mapsto \pi_1^{-1}(\Fp+)$.
Cela permet de calculer son cardinal.
\cqfd

}

\remarque{
\label{insuffisance de la methode des orbites}

{\bf (1)}
Il peut arriver qu'un $\, \lambdatilde \!\in\! \gstItilde(e) \,$ v\'erifie
$\, \lambdatilde[e] \!\notin\! \gestItilde$.
(Cependant la condition $\, \lambdatilde \in \gstInctilde(e) \,$ implique
$\, \lambdatilde[e] \in \gestInctilde$.)
Par exemple :
$\, G = Sp(4,\R )$,
$\, \lambdatilde = (0,\F+) \,$
tel que $\gothh \egdef \gothg(\lambdatilde)$ est \'egale \`a sa composante 
infinit\'esimalement elliptique et le syst\`eme de racines positives de
$\, (\gothg_\C ,\gothh_\C ) \,$ associ\'e \`a $\F+$ s'\'ecrit
$\, \{ \, \alpha_1,\alpha_2,\alpha_1+\alpha_2,2\alpha_1+\alpha_2 \, \}$
avec pour seule racine compacte $\alpha_1+\alpha_2$.
En prenant $\,e = \exp E\,$ o\`u $\,E \in \gothh\,$ est d\'etermin\'e par
$\; \alpha_1(E) = -\alpha_2(E) = \mi \pi$, l'alg\`ebre de Lie $\gothg(e)$ est ici
isomorphe \`a $\gothu(2)$ et donc
$\, C(\gothg(e),\gothh(e))_{reg} = \emptyset$.

\smallskip
{\bf (2)}
On peut aussi trouver des
$\, \lambdatilde_e \in \gestregtilde \,$
qui ne sont pas de la forme $\lambdatilde[e]$ pour un
$\, \lambdatilde \in \gstregtilde(e)$.
Par exemple :
$\, G = SU(2)$,
$\, e = \left( {\scriptstyle \Rot{1}} \right) \,$
et
$\, \lambdatilde_e = (0,\mi \,\goth{so}(2)^*)$.

\smallskip
{\bf (3)}
Pour certains
$\, \lambdatilde \in \gstfondtilde(e) \,$
il n'existe pas de $f \!\in\! \gstreg(e)$
tel que $\, G \cdot f = R_G(G \cdot \lambdatilde)$.
C'est le cas, avec $\, G = SL(2,\R ) \,$ et
$\, e = \left( {\scriptstyle \Rot{1}} \right)$, des
\'el\'ements  de $\, \gstfondtilde(e) \,$ de la forme
$\, (0,\F+) \,$ tels que
$\, \F+ \!\subseteq \mi \,\goth{so}(2)^* \!$.
\cqfr

}

%
\part{
Les param\`etres << repr\'esentation projective >>
}\label{RepII}
%

Les param\`etres de cette partie sont calqu\'es sur ceux introduits par
M.~Duflo, qui \'etaient adapt\'es au cas des orbites coadjointes semi-simples
r\'eguli\`eres de $G$.

\section{
Rappels sur les groupes sp\'ecial-m\'etalin\'eaire et m\'etaplectique
}\label{Rep4}

Ma r\'ef\'erence concernant le groupe m\'etaplectique et les fonctions
orientation de M.~Duflo et M.~Vergne qui lui sont attach\'ees, est leur article
\cite {DV93}.
Je choisis comme eux d'\'eviter l'utilisation du caract\`ere de la
repr\'esentation m\'etaplectique en suivant le point de vue de
\cite [conjecture p.$\!$~291]{Ve94} (au signe de la forme symplectique sur les
orbites coadjointes pr\`es).

\definition{

On consid\`ere un $\R $-espace vectoriel $V$ de dimension finie.

\smallskip
{\bf (a)}
On note $\, DL(V) \to SL(V)$\labind{DL(V)} << le >> rev\^etement double de
$SL(V)$, unique \`a isomorphisme de rev\^etements pr\`es, qui est connexe quand
$\dim V \geq 2$.
(Le (b) ci-dessous fournira une description canonique d'un tel rev\^etement,
comme ensemble de couples form\'es d'un \'el\'ement de $SL(V)$ et d'une
orientation d'un certain sous-espace vectoriel de $V\!$.)

\smallskip
{\bf (b)}
Soit $\hat{a} \in DL(V)$ au-dessus d'un $a \in SL(V)$ elliptique.

Quand $\dim V \geq 2$, on note

$\scalo (\hat{a})_{(1-a) \cdot V}
= (-1)^{\frac{\alpha_1 +\cdots+ \alpha_p}{2\pi}} \,
\sg (\sin (\frac{\beta_1}{2})\ldots\sin (\frac{\beta_q}{2})) \;
{\scriptstyle \times} \;
\R ^+ \moins0 \,
(w_1 \wedge \cdots \wedge w_{2q})$%
\labind{Oa},
\\
o\`u << $\sg$\labind{sg} >> repr\'esente la fonction signe sur $\R \moins0$,
ind\'ependamment du choix d'un $A \in \goth{sl} (V)$ infinit\'esimalement
elliptique tel que $\,\hat{a} = \exp_{\!DL(V)} \! A$, et d'une base
$(v_1,\dots,v_{2p},w_1,\dots,w_{2q})$ de $A \cdot V$ dans laquelle la
matrice de la restriction de $A$ est de la forme

\negmedskip\centerline{
{\bf (\boldmath$**$)}
$\left(\!\! \begin{smallmatrix}
\begin{smallmatrix}
\vrule width 0ex height 2.2ex
\smash{\left({\scriptscriptstyle\Rot{\alpha_1}}\right)}
 & & \\
 & \!\!\smash{\ddots}\vrule width 0ex height 1ex \!\! & \\
 & &
\smash{\left({\scriptscriptstyle\Rot{\alpha_p}}\right)}
\end{smallmatrix} \!\!
 & \\
 &
\!\! \begin{smallmatrix}
\vrule width 0ex height 2.2ex
\smash{\left({\scriptscriptstyle\Rot{\beta_1}}\right)}
 & & \\
 & \!\!\smash{\ddots}\vrule width 0ex height 1ex \!\! & \\
 & &
\vrule width 0ex depth 1.2ex
\smash{\left({\scriptscriptstyle\Rot{\beta_q}}\right)}
\end{smallmatrix} \!\!
\end{smallmatrix} \!\!\right)$
}

\smallskip\noindent
avec
$\; \alpha_1,\dots,\alpha_p \in 2\pi\Z \moins0 \,$
et
$\,\beta_1,\dots,\beta_q \in \R \setminus\! 2\pi\Z $.

\smallskip
Quand $\dim V \leq 1$, $\; \scalo (\hat{a})_{\{0\}}$ vaut $\R ^+ \moins0$ ou
$\R ^- \moins0$ suivant que $\hat{a}$ est trivial ou non.

\smallskip
{\bf (c)}
Soit $A \in \goth{sl} (V)$ infinit\'esimalement elliptique.
On note

$\scalo (A)_{A \cdot V}
= \sg (\beta_1\ldots\beta_q) \;
{\scriptstyle \times} \;
\R ^+ \moins0 \,
(w_1 \wedge \cdots \wedge w_{2q})$%
\labind{OA},
\\
ind\'ependamment du choix d'une base $(w_1,\dots,w_{2q})$ de $A \cdot V$ dans
laquelle la matrice de la restriction de $A$ est de la forme

\smallskip\centerline{
$\left(\! \begin{smallmatrix}
\vrule width 0ex height 2.2ex
\smash{\left({\scriptscriptstyle\Rot{\beta_1}}\right)}
 & & \\
 & \!\!\smash{\ddots}\vrule width 0ex height 1ex \!\! & \\
 & &
\vrule width 0ex depth 1.2ex
\smash{\left({\scriptscriptstyle\Rot{\beta_q}}\right)}
\end{smallmatrix} \!\!\right)$
}

\noindent
avec $\; \beta_1,\dots,\beta_q \in \R \moins0$.

}

\definition{
\label{geometrie metaplectique}

On consid\`ere un $\R $-espace vectoriel $V$ de dimension finie muni d'une
forme bilin\'eaire $B$ altern\'ee non d\'eg\'en\'er\'ee.
On d\'esigne encore par $B$ le prolongement bilin\'eaire complexe de cette forme
bilin\'eaire \`a $V_\C $.

\smallskip
{\bf (a)}
On note $\,Mp(V)\,$\labind{Mp(V)} l'image r\'eciproque de $Sp(V)$ dans $DL(V)$.

\smallskip
{\bf (b)}
On note $\scalo (B)_V$\labind{OB} l'orientation de $V$ sur laquelle
$B^{\frac{1}{2}\dim V}\!$ est positive. \\
On appelle << base symplectique >> de $(V,B)$ toute base
$(P_1,\dots,P_n,Q_1,\dots,Q_n)$ de $V$ telle que
$B(P_i,P_j) = B(Q_i,Q_j) = 0$ et $B(P_i,Q_j) = \delta_{i,j}$
(symbole de Kronecker) pour $1 \leq i,j \leq n$.

\smallskip
{\bf (c)}
On appelle << lagrangien >> de $(V_\C ,B)$ tout sous-espace vectoriel de
$V_\C $ \'egal \`a son orthogonal pour $B$.
On appelle << lagrangien positif >> de $(V_\C ,B)$ tout lagrangien de
$(V_\C ,B)$ sur lequel la forme sesquilin\'eaire hermitienne
$\; (v,w) \mapsto \mi \,B(v,\overline{w}) \;$ est positive.

Soit $\,\call$ un lagrangien de $(V_\C ,B)$.
On note $Mp(V)_\call$\labind{Mp(V)L} le normalisateur de $\,\call$ dans $Mp(V)$.
\`A chaque $\, \hat{x} \!\in\! Mp(V)_\call \,$ au-dessus d'un \'el\'ement $x$ de
$Sp(V)$, on associe le nombre $\,n_\call(x)\,$\labind{nL} de valeurs propres
compt\'ees avec multiplicit\'e dans $\left] 1,+\infty \right[$ de la restriction
de $x^\C $ \`a $\,\call$, et le nombre $\,q_\call(x)\,$\labind{qL} de
composantes strictement n\'egatives dans la matrice de la forme sesquilin\'eaire
hermitienne
$\; (v,w) \mapsto \mi \,B(v,\overline{w}) \;$ sur $(1-x^\C ) \cdot \call$
relativement \`a une base orthogonale.

\smallskip
{\bf (d)}
Soit $\,\call$ un lagrangien de $(V_\C ,B)$.
Pour tout $\, \hat{x} \in Mp(V)_\call \,$ de composante elliptique $\hat{x}_e$
(au sens de \cite [lem. 31 p.$\!$~38]{DV93}) au-dessus d'un \'el\'ement $x$ de
$Sp(V)$ de composante elliptique $x_e$, on pose

\smallskip\centerline{
$\rho_\call(\hat{x})
\,=\, (-1)^{q_\call(x_e)} \;
\frac
{\scalo (\hat{x}_e)_{(1-x_e) \cdot V}}
{\scalo (B)_{(1-x_e) \cdot V}} \;
\Prod_{1 \leq k \leq n} (\sqrt{r_k} \, \me^{\mi \,\theta_k/2})$%
\labind{zzr L},
}

\noindent
o\`u $\; r_1 \me^{\mi \,\theta_1},\dots,r_n \me^{\mi \,\theta_n}$ sont les
valeurs propres  compt\'ees avec multiplicit\'e de la restriction de $x^\C $ \`a
$\,\call$, avec $\; r_1,\dots,r_n \in \R ^+ \moins0 \,$ et
$\; \theta_1,\dots,\theta_n \in \left] -2\pi,0 \right]$.

\smallskip
{\bf (e)}
Soit $\hat{x} \in Mp(V)$ de composante elliptique $\hat{x}_e$ au-dessus d'un
\'el\'ement semi-simple $x$ de $Sp(V)$ de composante elliptique $x_e$.
On pose

\smallskip\centerline{
$\delta(\hat{x})
\,=\; \frac
{\scalo (\hat{x}_e)_{(1-x_e) \cdot V}}
{\scalo (B)_{(1-x_e) \cdot V}} \;
\Prod_{1 \leq k \leq n} \me^{\mi \,\theta_k/2}$%
\labind{zzd},
}

\noindent
ind\'ependamment (compte tenu de la d\'emonstration du (c) de la proposition
ci-dessous) du choix d'un $\, E \in \goth{sp} (V)$ infinit\'esimalement
elliptique tel que $\,x_e = \exp_{\!Sp(V)} \! E$, et d'une base symplectique
$(P_1,\dots,P_n,Q_1,\dots,Q_n)$ de $(V,B)$ pour laquelle la matrice
de $E$ relativement \`a $(P_1,Q_1,\dots,P_n,Q_n)$ est de la forme

\centerline{
$\left(\! \begin{smallmatrix}
\vrule width 0ex height 2.2ex
\smash{\left({\scriptscriptstyle\Rotneg{\theta_1}}\right)}
 & & \\
 & \!\!\smash{\ddots}\vrule width 0ex height 1ex \!\! & \\
 & &
\vrule width 0ex depth 1.2ex
\smash{\left({\scriptscriptstyle\Rotneg{\theta_q}}\right)}
\end{smallmatrix} \!\!\right)$
}

\noindent
avec $\; \theta_1,\dots,\theta_n \in \left] -2\pi,0 \right]$.

\smallskip
{\bf (f)}
Soit $\hat{x} \in Mp(V)$ de composante elliptique $\hat{x}_e$ au-dessus d'un
\'el\'ement semi-simple $x$ de $Sp(V)$ de composante elliptique $x_e$.
On pose

\smallskip\centerline{
$\Phi(\hat{x})
\,=\; \frac
{\scalo (\hat{x}_e)_{(1-x_e) \cdot V}}
{\scalo (B)_{(1-x_e) \cdot V}} \;\,
\mi^{-\frac{1}{2}\dim(1-x_e) \cdot V} \;\,
|\det \;(1-x)_{_{\scriptstyle (1-x) \cdot V}}|^{-1/2}$%
\labind{zzV}.
}

}

\medskip
On va maintenant voir, comme le sous-entend M.~Vergne dans
\cite [prop. p.$\!$~289]{Ve94} (avec la repr\'esentation de Weil
contragr\'ediente de celle utilis\'ee dans \cite{DHV84}, d\'eduite de celle de
\cite{DHV84} en rempla\c{c}ant au choix $B$ par $-B$ \,ou\, le caract\`ere
central du groupe de Heisenberg par son conjugu\'e), que les fonctions que je
viens d'introduire dans les points (d), (e), (f) sont les fonctions $\rho_\call$
et $\delta$ de \cite{Df82a} et la fonction $\Phi$ de \cite{DHV84}.

\proposition{
\label{morphisme rho et fonction delta}

On consid\`ere un $\R $-espace vectoriel $V$ de dimension finie muni d'une
forme bilin\'eaire $B$ altern\'ee non d\'eg\'en\'er\'ee.

\smallskip
{\bf (a)}
Le rev\^etement $\, Mp(V) \to Sp(V)$ est connexe quand $V \not= \{0\}$.

\smallskip
{\bf (b)}
Soit $\,\call$ un lagrangien de $(V_\C ,B)$.
La fonction $\rho_\call$ est un morphisme de groupes de Lie de
$Mp(V)_\call$ dans $\C \moins0$ tel que
$\; \diff _1 \rho_\call
= \frac{1}{2}\,\tr (\,\cdot\,^\C )_\call$.
Pour tout $\hat{x} \in Mp(V)_\call$ au-dessus d'un \'el\'ement semi-simple $x$
de $Sp(V)$, on a

\smallskip\centerline{
$\Phi(\hat{x})
= (-1)^{n_\call(x) \,+\, q_\call(x)} \;
\rho_\call(\hat{x}) \;\,
(\det \;(1-x^\C )_{_{\scriptstyle (1-x^\C ) \cdot \call}})^{-1}$.
}

\smallskip
{\bf (c)}
Soit $\hat{x} \in Mp(V)$ au-dessus d'un $x \in Sp(V)$ de composante
elliptique~$x_e$.
Il existe un lagrangien positif de $(V_\C ,B)$ stable par $x^\C $.
Pour tout lagrangien $\,\call$ de $(V_\C ,B)$ stable par $x^\C $, on a

\centerline{
$\rho_\call(\hat{x}) \;
\abs{\rho_\call(\hat{x})}^{-1}
\,=\; \delta(\hat{x}) \;
\Prod_{z \in \Sp(x_e)}
z^{q_z}$,
}

\noindent
o\`u $\;\Sp(x_e)$ est le spectre de $x_e$ et, pour chaque
$z \in \Sp(x_e)$, $(p_z,q_z)$ est la signature de la forme
sesquilin\'eaire hermitienne $\, (v,w) \!\mapsto\! \mi \,B(v,\overline{w}) \,$
sur le sous-espace propre de $({x_e\!}^\C )_\call$ associ\'e \`a~$z$.

}

\dem{D\'emonstration de la proposition}{

{\bf (a)}
On suppose $V$ non nul.
Soit $(P_1,\dots,P_n,Q_1,\dots,Q_n)$ une base symplectique de $(V,B)$.
On note $T_{Sp}$ le tore maximal de $Sp(V)$ form\'e des endomorphismes de $V$
dont les matrices dans la base $(P_1,Q_1,\dots,P_n,Q_n)$ sont diagonales par
blocs avec des blocs dans $SO(2)$, et $\gotht$ son alg\`ebre de Lie.
Le tore $\, T_{Mp} \egdef \exp_{\!Mp(V)} \! \gotht \,$ est un tore maximal de
$DL(V)$.

On choisit un \'el\'ement $A$ de $\Ker \exp_{T_{Sp}}\!$ tel que
$\; \scalo (\exp_{\!DL(V)} \! A)_{\{0\}} = \R ^- \! \moins0$ (il en existe).
Le chemin $\, t \in [0,1] \mapsto \exp_{T_{Mp}} \! (t A) \,$ joint dans $Mp(V)$
les deux points du noyau du morphisme de groupes canonique de $Mp(V)$ dans
$Sp(V)$.

\smallskip
{\bf (b)}
Soit $\hat{x} \in Mp(V)_\call$ au-dessus d'un \'el\'ement semi-simple $x$ de
$Sp(V)$ de composante elliptique $x_e$.

On note $\wt{\rho}_\call$ le morphisme de groupes de Lie de
$Mp(V)_\call$ dans $\C \! \moins0$ d\'efini dans
\cite [p.$\!$~107]{Df84}.
D'apr\`es \cite [(10) p.$\!$~108]{Df84}, on a

\centerline{
$\abs{\wt{\rho}_\call(\hat{x})}
= \abs{\det(x^\C )_\call}^{1/2}
= \abs{\rho_\call(\hat{x})} \;$
et
$\; \diff _1 \wt{\rho}_\call
= \frac{1}{2}\,\tr (\,\cdot\,^\C )_\call$.
}

\noindent
Ces \'egalit\'es permettent de prouver, compte tenu de
\cite [fin du lemme 30 p.$\!$~ 37]{DV93}, que les applications $\wt{\rho}_\call$
et $\rho_\call$ commutent aux prises de composantes elliptique, positivement
hyperbolique et unipotente sur les groupes $Mp(V)_\call$ et $GL(\C )$, et ont
m\^emes restrictions aux parties de $Mp(V)_\call$ form\'ees de ses
\'el\'ements positivement hyperboliques ou de ses \'el\'ements unipotents.

On note $\wt{\Phi}$ la fonction (ind\'ependante de $\call$) de l'ensemble des
\'el\'ements semi-simples de $Mp(V)$ dans $\C \! \moins0$ d\'efinie dans
\cite [p.$\!$~102]{DHV84}.
Elle s'\'ecrit

\smallskip\centerline{
$\wt{\Phi}(\hat{x})
= (-1)^{n_\call(x) \,+\, q_\call(x)} \;
\wt{\rho}_\call(\hat{x}) \;\,
(\det \;(1-x^\C )_{_{\scriptstyle (1-x^\C ) \cdot \call}})^{-1}$.
}

\noindent
En outre, l'application $\, \dot{v} \mapsto B(v,.) \,$ de
$\, (1-x^\C ) \cdot V_\C  \,/\, (1-x^\C ) \cdot \call \,$
dans $\, ((1-x^\C ) \cdot \call)^* \,$
est une bijection lin\'eaire qui commute \`a l'action de $\, 1 \!-\! x$.
Il s'ensuit que

\smallskip\centerline{
$|\wt{\Phi}(\hat{x})|
= |\det \;(1-x)_{_{\scriptstyle (1-x) \cdot V}}|^{-1/2}
= \abs{\Phi(\hat{x})}$.
}

Pour montrer que $\,\wt{\Phi}=\Phi\,$
(respectivement $\,\wt{\rho}_\call=\rho_\call$),
vu \cite [(9) p.$\!$~108]{Df84} il reste \`a prouver que
$\,\wt{\Phi}/|\wt{\Phi}|\,$ et $\,\Phi/\abs{\Phi}\,$
(respectivement $\,\wt{\rho}_\call\,$ et $\,\rho_\call$)
co\"{\i}ncident en un des deux points de $Mp(V)$ situ\'es au-dessus de $x$
(respectivement de $x_e$).

\smallskip
On note $V^z$ (respectivement $V_\C ^z$) le sous-espace
propre de $x$ (respectivement $x^\C $) associ\'e \`a une valeur propre
$z \!\in\! \R \! \moins0$
(respectivement $z \!\in\! \C \! \moins0$), et
$\,{\scriptstyle \perp}^{\!\scriptscriptstyle B}\,$ la relation de
$B$-orthogonalit\'e.
On consid\`ere maintenant une valeur propre $z \in \C \! \moins0$ de $x^\C $
telle que $\, \Im z \geq 0 \,$ et $\, \abs{z} \geq 1$.
On va lui associer une certaine base symplectique $\calb_z$ de
$\; (( V_\C ^z + V_\C ^{z^{-1}} \!\!+ V_\C ^{\bar{z}} + V_\C ^{{\bar{z}}^{-1}} )
\cap V,B)$.

\quad
Si $\,z \in \R \,$ et $\,\abs{z} \geq 1$ :\
quand $\,\abs{z} = 1$, on fixe une base symplectique
$\; \calb_z = (P_1^0,\dots,P_{n^0}^0,Q_1^0,\dots,Q_{n^0}^0) \;$
de $V^z$ ;
quand $\,\abs{z} > 1$,
on fixe une base $(P_1^0,\dots,P_{n^0}^0)$ de $V^z$, et en d\'eduit une
unique base symplectique
$\, \calb_z  = (P_1^0,\dots,P_{n^0}^0,Q_1^0,\dots,Q_{n^0}^0) \,$
de $\, V^z \oplus V^{z^{-1}}$ telle que
$\,Q_1^0,\dots,Q_{n^0}^0 \in V^{z^{-1}}\!$,
car $V^z$ et $V^{z^{-1}}$ sont en dualit\'e avec
$\, V^z \,{\scriptstyle \perp}^{\!\scriptscriptstyle B}\, V^z \,$
et
$\, V^{z^{-1}} {\scriptstyle \perp}^{\!\scriptscriptstyle B}\, V^{z^{-1}}$.

\quad
Si $\,\Im z > 0\,$ et $\,\abs{z} > 1$ :\
on fixe une base
$\,(P_1 \!+\! \mi P_2,\dots,P_{2n-1} \!+\! \mi P_{2n})\,$
de $V_\C ^z$ telle que $P_1,\dots,P_{2n} \in V$ ;
on en d\'eduit une unique base symplectique
$\; (P_1 \!+\! \mi P_2,\dots,P_{2n-1} \!+\! \mi P_{2n},
Q_1 \!-\! \mi Q_2,\dots,Q_{2n-1} \!-\! \mi Q_{2n}) \;$
de $\,V_\C ^z \oplus V_\C ^{z^{-1}}$ muni de $B/2$
telle que
$\,Q_1 \!-\! \mi Q_2,\dots,Q_{2n-1} \!-\! \mi Q_{2n} \in V_\C ^{z^{-1}}\,$
et $\,Q_1,\dots,Q_{2n} \in V\!$,
car $V_\C ^z$ et $V_\C ^{z^{-1}}$ sont en dualit\'e avec
$\; V_\C ^z \,{\scriptstyle \perp}^{\!\scriptscriptstyle B}\,
V_\C ^z \;$ et
$\; V_\C ^{z^{-1}} {\scriptstyle \perp}^{\!\scriptscriptstyle B}\,
V_\C ^{z^{-1}}$ ;
puis on obtient la base symplectique
$\, \calb_z \egdef (P_1,\dots,P_{2n},Q_1,\dots,Q_{2n}) \,$
de
$\, ( V_\C ^z \oplus V_\C ^{z^{-1}}\! \oplus V_\C ^{\bar{z}}
\oplus V_\C ^{{\bar{z}}^{-1}} ) \cap V\!$,
car
$\; V_\C ^z \,{\scriptstyle \perp}^{\!\scriptscriptstyle B} \,V_\C ^{\bar{z}}$,
$\; V_\C ^{z^{-1}} \,{\scriptstyle \perp}^{\!\scriptscriptstyle B}\,
V_\C ^{\bar{z}^{-1}} \,$
et
$\; V_\C ^z \,{\scriptstyle \perp}^{\!\scriptscriptstyle B}\,
V_\C ^{\bar{z}^{-1}}\!$.

\quad
Si $\,\Im z > 0\,$ et $\,\abs{z} = 1$ :\
on fixe une base
$(P'_1 \!+\! \mi Q'_1,\dots,P'_p \!+\! \mi Q'_p,P''_1 \!-\! \mi Q''_1,
$ $\!\!
\dots,P''_q \!-\!\mi Q''_q)$
de $V_\C ^z$ telle que
$\, P'_1,Q'_1,\dots,P'_p,Q'_p,P''_1,Q''_1,\dots,P''_q,Q''_q \in V\!$,
dans laquelle la forme sesquilin\'eaire hermitienne non d\'eg\'en\'er\'ee
$\, (v,w) \mapsto iB/2 \, (v,\overline{w}) \,$
a pour matrice
$\left(\!
{\scriptstyle \begin{smallmatrix} I_p & 0 \\ 0 & -I_p \end{smallmatrix}}
\!\right)$ ;
elle fournit la base symplectique
$\; \calb_z
\egdef (P'_1,\dots,P'_p,P''_1,\dots,
$ $\!\!
P''_q,Q'_1,\dots,Q'_p,Q''_1,\dots,Q''_q) \;$
de
$\; (V_\C ^z \oplus V_\C ^{\bar{z}}) \cap V\!$,
car
$\; V_\C ^z \,{\scriptstyle \perp}^{\!\scriptscriptstyle B}\, V_\C ^z \;$
et
$\; V_\C ^{\bar{z}} \,{\scriptstyle \perp}^{\!\scriptscriptstyle B}\,
V_\C ^{\bar{z}}$.

\smallskip
On constate que le sous-espace vectoriel $\,\call_x^+$ (respectivement
$\,\call_{x,e}^+$) de
$V_\C $ engendr\'e par les diff\'erents vecteurs
$P_1^0,\dots,P_{n^0}^0$,
ou,
$P_1,\dots,P_{2n}$,
ou,
$P'_1 \!+\! \mi Q'_1,
$ $\!\!
\dots,P'_p \!+\! \mi Q'_p,P''_1 \!+\! \mi Q''_1,\dots,P''_q \!+\! \mi Q''_q$
(respectivement
$P_1^0 \!+\! \mi Q_1^0,\dots,P_{n^0}^0 \!+\! \mi Q_{n^0}^0$,
ou,
$P_1 \!+\! \mi Q_1,\dots,P_{2n} \!+\! \mi Q_{2n}$,
ou,
$P'_1 \!+\! \mi Q'_1,\dots,P'_p \!+\! \mi Q'_p,
P''_1 \!+\! \mi Q''_1,\dots,P''_q\!+\!\mi Q''_q \,$)
lorsque $z$ varie, est un lagrangien positif de $(V_\C ,B)$ stable
par~$x^\C $ (respectivement ${x_e\!}^\C $).

On note $E$ (respectivement $H$) l'\'el\'ement de ${\goth sp} (V)$ stabilisant
les espaces vectoriels
$\; ( V_\C ^z + V_\C ^{z^{-1}}\!\!
+ V_\C ^{\bar{z}} + V_\C ^{\bar{z}^{-1}} ) \cap V \;$
avec $z = r \me^{\mi \,\theta}$, $r \geq 1$, $-2\pi < \theta \leq 0$,
et dont la restriction \`a un tel espace a dans la base $\calb_z$ une matrice
diagonale par blocs dont les blocs sont les matrices
\\
$\left( {\scriptstyle \Rotneg{\theta}} \right)$
(resp.
$\left(\!
{\scriptstyle \begin{smallmatrix} \ln r & 0 \\ 0 & -\ln r \end{smallmatrix}}
\!\right)$)
relativement \`a $(P_k^0,Q_k^0)$,
\\
ou,
$\left(\!\!{\scriptscriptstyle\diagd{\Rotneg{\theta}}{\Rot{\theta}}}\!\!\right)$
(resp.
$\left( \scriptscriptstyle \diagq{\ln r}{\ln r}{-\ln r}{-\ln r} \right)$)
relativement \`a $(P_{2k-1},P_{2k},Q_{2k-1},-Q_{2k})$,
\\
ou,
$\left( {\scriptstyle \Rotneg{\theta}} \right)$
(resp. $\left( \begin{smallmatrix}  0 & 0 \\ 0 &  0 \end{smallmatrix} \right)$)
relativement \`a $(P'_k,Q'_k)$ et
$\left( {\scriptstyle \Rot{\theta}} \right)$
(resp. $\left( \begin{smallmatrix}  0 & 0 \\ 0 &  0 \end{smallmatrix} \right)$)
relativement \`a $(P''_l,Q''_l)$.

Les \'el\'ements
$\, \hat{x}_{E,H} = \exp_{\!Mp(V)} \! E \, \exp_{\!Mp(V)} \! H$
et
$\, \hat{x}_E = \exp_{\!Mp(V)} \! E$
de $Mp(V)$ sont respectivement au-dessus de $x$ et $x_e$.
D'apr\`es \cite [(15) p.$\!$~109]{Df84}, on a :

\centerline{
${\displaystyle\frac{\wt{\rho}_{\call_x^+}}{|\wt{\rho}_{\call_x^+}|}}
(\hat{x}_{E,H})
= \wt{\rho}_{\call_{x,e}^+}(\exp_{\!Mp(V)} \! E)
=\, \me^{\frac{1}{2} \tr (E^\C )_{\call_{x,e}^+}}$.
}

\noindent
D'o\`u les \'egalit\'es suivantes, qui permettent d'obtenir le r\'esultat :

\smallskip\centerline{
${\displaystyle\frac{\wt{\Phi}}{|\wt{\Phi}|}}(\hat{x}_{E,H})
= (-1)^{n_{\call_x^+}(x)} \; \me^{\frac{1}{2}\,\tr (E^\C )_{\call_{x,e}^+}} \;
\biggl(
\frac{\det \;(1-x^\C )_{(1-x^\C ) \cdot \call_x^+}}
{\abs{\det \;(1-x^\C )_{(1-x^\C ) \cdot \smash{\call_x^+}}}}
\biggr)^{\!\!-1} \!
= \cdots
= {\displaystyle\frac{\Phi}{\abs{\Phi}}}(\hat{x}_{E,H})$
}

\centerline{
et
$\quad \wt{\rho}_\call(\hat{x}_E)
\!= (-1)^{q_\call(x_e)} \;
{\displaystyle\frac{\Phi}{\abs{\Phi}}}(\hat{x}_E) \;
\frac{\det \;(1-\exp E^\C )_{(1-\exp E^\C ) \cdot \call}}
{\abs{\det \;(1-\exp E^\C )_{(1-\exp E^\C ) \cdot \call}}}
= \cdots
= \rho_\call(\hat{x}_E).$
}

\smallskip
{\bf (c)}
L'existence d'un lagrangien positif de $(V_\C ,B)$ stable par $x^\C $ est
d\'emontr\'ee dans \cite [cor. p.$\!$~82]{B.72}.
On suppose $V$ non nul et note $\hat{x}_e$ la composante elliptique de
$\hat{x}$.
Par d\'efinition $\delta(\hat{x})$ est \'egal \`a
$\rho_{\C (P_1 \!+ \mi Q_1) +\cdots+ \C (P_n \!+ \mi Q_n)} \! (\hat{x}_e)$,
avec les notations de \ref{geometrie metaplectique} (e).
L'\'egalit\'e \cite [(15) p.$\!$~109]{Df84} prouve que ce nombre
complexe est ind\'ependant du choix de $(P_1,\dots,P_n,Q_1,\dots,Q_n)$.

On se donne un lagrangien $\call$ de $(V_\C ,B)$ stable par $x^\C $ et un
suppl\'ementaire $\cals$ de $\, \call \cap \overline{\call} \,$ dans $\call$ qui
est stable par ${x_e\!}^\C $.
On note $W$ le sous-espace symplectique de $(V,B)$ \'egal \`a l'orthogonal de
$\, \cals + \overline{\cals} \,$ dans $V\!$.
On fixe une base
$\, (P_1^0,\dots,P_{n^0}^0,P_1^\R ,\dots,P_{2m}^\R ) \,$
du lagrangien $\, \call \cap \overline{\call} \,$ de $(W_\C ,B)$,
form\'ee de vecteurs de $V\!$, telle que $\, P_1^0,\dots,P_{n^0}^0 \,$ sont des
vecteurs propres de $x_e$ associ\'es \`a des valeurs propres r\'eelles et
$\, P_1^\R \!+\! \mi P_2^\R ,\dots,P_{2m-1}^\R \!+\! \mi P_{2m}^\R \,$
sont des vecteurs propres de ${x_e\!}^\C $ associ\'es \`a des valeurs propres
non r\'eelles.
Elle se compl\`ete en une base symplectique
$(P_1^0,\dots,P_{n^0}^0,P_1^\R ,\dots,P_{2m}^\R ,
Q_1^0,\dots,Q_{n^0}^0,2Q_1^\R ,\dots,2Q_{2m}^\R )$
de $(W,B)$ telle que $Q_1^0,\dots,Q_{n^0}^0$ et
$Q_1^\R \!+\! \mi Q_2^\R ,\dots,Q_{2m-1}^\R \!+\! \mi Q_{2m}^\R $
sont des vecteurs propres de ${x_e\!}^\C $.
En effet, on peut construire $(Q_1^\R ,\dots,Q_{2m}^\R )$
par r\'ecurrence sur $m \!=\! \frac{1}{2} \dim (x_e^2\!-\!1) \!\cdot\! W\!$,
en associant au vecteur $\, P \egdef P_{2m-1}^\R \!+\! \mi P_{2m}^\R \,$ de
$\, (x_e^2\!-\!1) \cdot
$ $\!\!
(\call \cap \overline{\call}) \,$
un vecteur propre $\, Q = Q_{2m-1}^\R \!-\! \mi Q_{2m}^\R \,$
de ${x_e\!}^\C $ avec $\, Q_{2m-1}^\R ,Q_{2m}^\R \in (x_e^2-1) \!\cdot\! W \,$
tel que $B(P,Q) = 1$ et $B(Q,\overline{Q}) = 0$
(cette derni\`ere \'egalit\'e se r\'ealise en rempla\c{c}ant $Q$ par
$\, Q - \frac{1}{2} B(Q,\overline{Q}) \, \overline{P}$).
On fixe aussi une base
$\, (P'_1 \!+\! \mi Q'_1,
$ $\!\!
\dots,P'_p \!+\! \mi Q'_p,
P''_1 \!-\! \mi Q''_1,\dots,P''_q \!-\!\mi Q''_q) \,$
de $\,\cals$ form\'ee de vecteurs propres de ${x_e\!}^\C $ avec
$\, P'_1,Q'_1,\dots,P'_p,Q'_p,P''_1,Q''_1,\dots,P''_q,Q''_q \in V\!$,
dans laquelle la forme sesquilin\'eaire hermitienne non d\'eg\'en\'er\'ee
${x_e\!}^\C $-invariante
$\, (v,w) \mapsto iB/2 \, (v,\overline{w}) \,$
a pour matrice~%
$\left(\!
{\scriptstyle \begin{smallmatrix}  I_p & 0 \\ 0 & -I_p \end{smallmatrix}}
\!\right)$.

La base de $\call$ constitu\'ee des vecteurs propres
$\, P_1^0,\dots,P_{n^0}^0,P_1^\R \!+ \mi P_2^\R ,
$ $\!\!
P_1^\R \!- \mi P_2^\R ,
$ $\!\!
\dots,P_{2m-1}^\R \!+ \mi P_{2m}^\R ,P_{2m-1}^\R \!- \mi P_{2m}^\R \,$
et
$\, P'_1 \!+ \mi Q'_1,\dots,P'_p \!+ \mi Q'_p,
P''_1 \!- \mi Q''_1,\dots,P''_q \!-\mi Q''_q \,$
de ${x_e\!}^\C $ fournit une expression pour $\, \rho_\call(\hat{x}_e)$.
On calcule $\delta(\hat{x})$ \`a l'aide de la base symplectique
$\; (P_1^0,\dots,P_{n^0}^0,P_1^+,\dots,P_{2m}^+,P'_1,\dots,P'_p,
P''_1,\dots,P''_q,Q_1^0,\dots,Q_{n^0}^0,Q_1^+,\dots,Q_{2m}^+,
$ $\!\!
Q'_1,\dots,Q'_p,Q''_1,\dots,Q''_q) \;$
de $(V,B)$,
\\
o\`u
$\; P_{2k-1}^+ \!+\! \mi Q_{2k-1}^+
\egdef
(P_{2k-1}^\R \!+\! \mi P_{2k}^\R ) + \mi (Q_{2k-1}^\R \!+\! \mi Q_{2k}^\R ) \;$
\\
et
$\; P_{2k}^+ \!-\! \mi Q_{2k}^+
\egdef
-\mi (P_{2k-1}^\R \!+\! \mi P_{2k}^\R ) - (Q_{2k-1}^\R \!+\! \mi Q_{2k}^\R ) \;$
pour $1 \leq k \leq m$.

On en d\'eduit facilement le r\'esultat.
\cqfd

}

\section{\boldmath
Les param\`etres $\,\tau \in \XInd_G(\lambdatilde)$
}\label{Rep5}

Dans les deux sections qui suivent, on se donne
$\,\lambdatilde = (\lambda,\F+) \in \gstregtilde$.
On pose $\,\gothh = \gothg(\lambdatilde)$.
On note $\gotha$ la composante hyperbolique de $\gothh$ et fixe une
chambre ${\gotha^*}^+\!$ de $(\gothg(\lambda)(\mi \rho_\F+),\gotha)$.
On note aussi $\mu$ et $\nu$ les composantes infinit\'esimalement elliptique et
hyperbolique de $\lambda$.

\smallskip
Afin de fournir une param\'etrisation canonique d'un sous-ensemble de
$\widehat{G}$ qui permette de donner une description de la formule de
Plancherel, M.~Duflo a introduit les objets du (a) de la d\'efinition ci-dessous.
Quand $\lambda\in\gstssreg$, ces objets seront compatibles avec les miens, car
on pourra prendre
$t \!\in\! \left[0,1\right]$ en posant $\lambda_0 \!=\! \lambda$ dans le lemme
\ref{fonctions canoniques sur le revetement double}, et on aura
$\,\XInd_G(\lambdatilde) = \Xirr_G(\lambda)$ dans la d\'efinition
\ref{parametres adaptes} (c).

\definition{
\label{parametres de Duflo}

{\bf (a)}
Soit $f \in \gothg^*\!$.
On note
\\
$B_f$\labind{Bf} la forme bilin\'eaire altern\'ee non d\'eg\'en\'er\'ee
$\, (\dot{X},\dot{Y}) \mapsto f([X,Y]) \,$ sur $\, \gothg / \gothg(f)$,
\\
$G(f)^{\gothg/\gothg(f)} \,$\labind{G (f)g/g(f)} le sous-groupe de Lie de
$\, G(f) \times DL(\gothg / \gothg(f)) \,$ form\'e des couples $(x,\hat{a})$ tels
que $(\Ad x)_{\gothg / \gothg(f)}$ est l'image de $\hat{a}$ dans
$SL(\gothg / \gothg(f))$
(cf. \ref{morphisme rho et fonction delta} (a)),
\\
$\{1,\iota\}$\labind{zzi} et $G(f)^{\gothg/\gothg(f)}_0\!$\labind{G (f)g/g(f)0}
les images r\'eciproques de $\{1\}$ et $G(f)_0$ par le morphisme de groupes de
Lie surjectif canonique de $G(f)^{\gothg/\gothg(f)}\!$ dans $G(f)$,
\\
$\Xirr_G(f)$\labind{Xirr} l'ensemble fini des classes d'isomorphisme de
repr\'esentations unitaires irr\'eductibles $\tau$ de $G(f)^{\gothg/\gothg(f)}$
telles que
$\tau(\iota) = -\id$ et $\, \tau(\exp X) = \me^{\mi \,f(X)}\id \,$
pour $X \!\in\! \gothg(f)$ (dans ce cas $\tau$ est de dimension finie).

\smallskip
{\bf (b)}
On note
$G(\lambdatilde)^{\gothg/\!\gothg(\lambda)(\mi \rho_\F+)}\!$%
\labind{G (zzl tilde)g/glambdairho}
l'image r\'eciproque de  $G(\lambdatilde)$ dans
$G(\lambda_{\textit{can}})^{\gothg/\gothg(\lambda_{\textit{can}})}\!\!$.

\noindent
On posera aussi
$\,G(\lambdatilde)^{\gothg/\gothh}_0 \!= G(\lambda_+)^{\gothg/\gothh}_0$%
\labind{G (zzl tilde)g/h0}
(cf. $\;G(\lambdatilde)_0 = G(\lambda_+)_0$ et $\gothg(\lambda_+) = \gothh$).

Lorsque  $\,{\gotha^*}^+ = \gotha^*$, on a
$\;\gothg(\lambda)(\mi \rho_\F+) = \gothh\,$ et
$\,G(\lambdatilde) = G(\lambda_+)$, et on notera dans ce cas
$\,G(\lambdatilde)^{\gothg/\gothh}\!$\labind{G (zzl tilde)g/h}
pour d\'esigner $G(\lambdatilde)^{\gothg/\!\gothg(\lambda)(\mi \rho_\F+)}\!$
et $\,G(\lambda_+)^{\gothg/\gothh}\!$.

}

\medskip
Dans la suite, on identifiera les fonctions complexes sur $\, G(\lambda_+) \,$
aux fonctions complexes sur le rev\^etement $\, G(\lambda_+)^{\gothg/\gothh} \,$
qui sont constantes sur les fibres.

\lemme{
\label{fonctions canoniques sur le revetement double}

On fixe un r\'eel $\,t \!\in\! \left]0,1\right]$.
On utilise les syst\`emes de racines positives de la d\'efinition
\ref{choix de racines positives} (b) et la notation $\lambda_t$ du lemme
\ref{lemme clef}. \\
On pose
$\;\, \call_{\lambda_t}
= \gothh_\C  \oplus
\Sumpetit_{\alpha \in R^+(\gothg_\C ,\gothh_\C )} \gothg_\C ^{\alpha}\;$
et, $\, \gothm' = \gothg(\nu_+)\,$ en anticipant sur la partie \ref{RepIII}.

{\bf (a)}
On a
$\; G(\lambda_t)^{\gothg/\gothh} = G(\lambda_+)^{\gothg/\gothh}$
et $\, \call_{\lambda_t}/\gothh_\C $ est un lagrangien de
$(\gothg_\C /\gothh_\C ,B_{\lambda_t})$ stable par $G({\lambda_t})$.
La fonction complexe $\,\delta_{\lambda_t}^{\gothg/\gothh}$ sur
$\, G(\lambda_+)^{\gothg/\gothh}$ et le morphisme de groupes de Lie
$\;\rho_{\lambda_t}^{\gothg/\gothh}$ de
$\, G(\lambda_+)^{\gothg/\gothh} \,$ dans $\, \C \! \moins0$ qu'on en d\'eduit
(cf. \ref{geometrie metaplectique} (e) et
\ref{morphisme rho et fonction delta} (b))
sont ind\'ependants de $t$, et donc respectivement \'egaux \`a
$\,\delta_{\lambda_+}^{\gothg/\gothh}\!$\labind{zzd zzl+}
et
$\, \rho_{\lambda_+}^{\gothg/\gothh}\!$\labind{zzr zzl+}. \\
On a $\;\, \diff _1 \rho_{\lambda_+}^{\gothg/\gothh} = \rho_{\gothg,\gothh}$.

\smallskip
{\bf (b)}
Soit $\hat{e} \in G(\lambda_+)^{\gothg/\gothh}$ au-dessus d'un \'el\'ement
elliptique $e$ de $G(\lambda_+)$.
On a 

\smallskip\centerline{
$\rho_{\lambda_+}^{\gothg/\gothh}(\hat{e})
\;=\;
\delta_{\lambda_+}^{\gothg/\gothh}(\hat{e}) \;{\scriptstyle \times}\;
\det\left(\Ad e^\C \right) \!\!\!
_{\Sumpetit_{\scriptstyle\alpha' \in R^+_K(\gothg_\C ,\gothh_\C )}
\!\!\gothg_\C ^{\alpha'}}
\;\,{\scriptstyle \times}\,
\Prod_{\calo \in \langle e \rangle \backslash
\wt{R}^+_\C (\gothm'_\C ,\gothh_\C )}
\!\!\! \! (-1)^{m_\calo-1} \, u_\calo$
}

\negsmallskip\noindent
o\`u $\;R^+_K(\gothg_\C ,\gothh_\C )$ est l'ensemble des
$\, \alpha' \!\in\! R^+(\gothg_\C ,\gothh_\C )$ compactes,
$\wt{R}^+_\C (\gothm'_\C ,\gothh_\C )$ est l'ensemble des classes des
$\beta \!\in\! R^+(\gothm'_\C ,\gothh_\C )$ complexes modulo
l'identification de $\beta$ avec $-\,\overline{\!\beta}$,
et, pour chaque orbite $\calo$ d'un \'el\'ement
$\{\beta,\!-\,\overline{\!\beta}\}$ de
$\wt{R}^+_\C (\gothm'_\C ,\gothh_\C )$ sous l'action du sous-groupe
$\langle e \rangle$ de $G$ engendr\'e par $e$, on note $m_\calo$ le cardinal de
$\calo$ et $u_\calo$ l'unique valeur propre de $(\Ad e^\C )^{\,m_\calo}$
associ\'ee \`a un vecteur propre qui s'\'ecrive
$\, X_\beta - \overline{\!X_{-\beta}\!} \,$ pour au moins un
$\, X_\beta \in \gothg_\C ^{\beta}$ et un
$\, X_{-\beta} \in \gothg_\C ^{-\beta}$ v\'erifiant
$\; [X_\beta,X_{-\beta}] = H_\beta$.

}

\dem{D\'emonstration du lemme}{

{\bf (a)}
La premi\`ere \'egalit\'e provient du lemme \ref{lemme clef}.
Il est imm\'ediat que $\, \call_{\lambda_t} / \gothh_\C $ est un lagrangien de
$(\gothg_\C /\gothh_\C ,B_{\lambda_t})$ stable par~$G({\lambda_t})$.

On note $\varphi_t$ la forme sesquilin\'eaire hermitienne
$\, (v,w) \mapsto \mi \,B_{\lambda_t}(v,\overline{w})$ sur
$(\gothg/\gothh)_\C $.
On va se servir des notations de l'\'enonc\'e du (b).
L'espace vectoriel $\, \call_{\lambda_t} / \gothh_\C $ est somme directe des
projections des espaces vectoriels suivants qui sont deux \`a deux orthogonales
pour $\varphi_t$ avec une signature de forme hermitienne restreinte pr\'ecis\'ee
entre parenth\`ese :
les $\, \gothg_\C ^{\alpha'}$ avec
$\alpha' \!\in\! R^+_K(\gothg_\C ,\gothh_\C )$
(signature $(0,1)$),
les $\, \gothg_\C ^{\alpha''}$ avec
$\alpha'' \!\in\! R^+(\gothg_\C ,\gothh_\C )$ imaginaire non compacte
(signature $(1,0)$),
les $\, \gothg_\C ^{\beta}\oplus\gothg_\C ^{-\,\overline{\!\beta}}$
avec
$\, \{\beta,\!-\,\overline{\!\beta}\} \in
\wt{R}^+_\C (\gothm'_\C ,\gothh_\C ) \,$
(signature $(1,1)$),
et les $\, \gothg_\C ^{\gamma}$ avec
$\gamma \!\in\! R^+(\gothg_\C ,\gothh_\C )$ hors de
$R^+(\gothm'_\C ,\gothh_\C )$
(signature $(0,0)$).

Soit $e$ un \'el\'ement elliptique de $G(\lambda_+)$.
Par ce qui pr\'ec\`ede, la forme sesquilin\'eaire hermitienne $\varphi_t$ est
non d\'eg\'en\'er\'ee sur
$(\call_{\lambda_t} \cap \gothm'_\C ) / \gothh_\C $
et nulle sur son orthogonal.
De plus, $(\call_{\lambda_t} \cap \gothm'_\C ) / \gothh_\C $
est somme directe orthogonale pour $\varphi_t$ de l'image et du noyau de
$\, 1 - \Ad e^\C $.
Donc la restriction de $\varphi_t$ \`a chacun d'entre eux est
non d\'eg\'en\'er\'ee.
Un argument de continuit\'e prouve ensuite que les signatures de ces
restrictions, et \emph{a fortiori} $q_{\call_{\lambda_t}}(e)$, sont
ind\'ependants de $t$.
De m\^eme, l'orientation
$\, \scalo (B_{\lambda_t})_{(1-\Ad e) \cdot (\gothg/\gothh)} \,$
est ind\'ependante de $t$.
Cela peut \^etre pr\'ecis\'e comme dans la preuve de
\ref{espace symplectique canonique}.
On utilise ensuite la d\'efinition \ref{geometrie metaplectique} (d) et la
proposition \ref{morphisme rho et fonction delta} (c).

La formule pour $\, \diff _1 \rho_{\lambda_+}^{\gothg/\gothh} \,$
se d\'eduit de la proposition \ref{morphisme rho et fonction delta} (b).

\smallskip
{\bf (b)}
Il s'agit d'appliquer la proposition \ref{morphisme rho et fonction delta} (c)
en tenant compte de la d\'emonstration du (a).
On utilise la notation $\varphi_t$ de cette d\'emonstration.

Soit $\calo$ l'orbite sous $\langle e \rangle$ d'un
$\{\beta,\!-\,\overline{\!\beta}\} \!\in\!
\wt{R}^+_\C (\gothm'_\C ,\gothh_\C )$.
On fixe
$(X_\beta,X_{-\beta}) \!\in\!
\gothg_\C ^{\beta} \,{\scriptstyle \times}\, \gothg_\C ^{-\beta}$
tel que $\; [X_\beta,X_{-\beta}] = H_\beta$.
La condition << $X_\beta - \overline{\!X_{-\beta}\!} \,$ est vecteur propre de
$\, (\Ad e^\C )^{\,m_\calo}$ >> est satisfaite (pour une unique valeur
propre de $\, (\Ad e^\C )^{\,m_\calo}$) quand
$\, \me^{\,m_\calo}\,\beta = \beta$, et se r\'ealise (\`a nouveau avec
unicit\'e de la valeur propre) en rempla\c{c}ant $X_\beta$ et $X_{-\beta}$ par
certains de leurs multiples quand
$\, \me^{\,m_\calo}\,\beta = -\,\overline{\!\beta}$.
Sous cette condition, les vecteurs propres
$\, \Sumpetit_{0 \leq k \leq m_\calo-1}\!\! \zeta^{-k} (\Ad e^\C )^k \,
(X_\beta \!- \overline{\!X_{-\beta}\!}) \,$
de $\Ad e^\C \!$, o\`u $\zeta$ d\'ecrit l'ensemble des racines
$m_\calo^{\textrm{i\`emes}}\!$ de $u_\calo$,
ont des projections dans
$\, \call_{\lambda_t} / \gothh_\C $ deux \`a deux
orthogonales pour $\varphi_t$ avec des << carr\'es >> strictement n\'egatifs.
Cela conduit \`a la formule annonc\'ee.
\cqfd

}

\remarque{
\label{sous-algebre de Borel canonique}

Soit $\, f \in \Supp_{\gothg(\lambda)^*}(G(\lambda)_0 \cdot \lambdatilde)$.

On note $\xi$ la composante nilpotente de $f$.
Il existe une unique sous-alg\`ebre de Borel $\gothb$ de
$\, \gothg(\lambda)_\C  \,$ sur laquelle $\xi$ est nulle.
En effet, d'apr\`es le lemme
\ref{description d'ensembles de formes lineaires} (a) et
\cite [cor. 5.3 p.$\!$~997 et cor. 5.6 p.$\!$~1001] {Ko59},
l'\'el\'ement de $\derivealg{\gothg(\lambda)}$ auquel $\xi\,$ s'identifie \`a
l'aide de la forme bilin\'eaire
${\scriptstyle \langle} ~,\!~ {\scriptstyle \rangle}$
appartient \`a une unique sous-alg\`ebre form\'ee d'\'el\'ements nilpotents
maximale de $\gothg(\lambda)_\C $, dont l'orthogonal pour
${\scriptstyle \langle} ~,\!~ {\scriptstyle \rangle}$
convient.
L'unicit\'e de $\gothb$ montre que $\; \overline{\gothb} = \gothb$.

On note $\, \gothn_{\lambda} \,$ la somme des espaces propres pour l'action de
$\, \centrealg{\gothg(\lambda)}_\C  \,$ dans $\gothg_\C $, dont le
poids $\, \alpha \in \centrealg{\gothg(\lambda)}_\C ^* \,$
(toujours restriction d'un poids de $\gothh_\C $) v\'erifie :
$\, \langle \nu,\alpha \rangle >0 \,$ ou
$\, (\langle \nu,\alpha \rangle =0 \;$ et
$\; i\langle \mu,\alpha \rangle >0)$.
On pose
$\, \call_f = \gothb \oplus \gothn_{\lambda}$.
On constate que $\call_f$ est une sous-alg\`ebre de Borel de $\gothg_\C $
en d\'ecomposant $\call_f$ en sous-espaces propres relativement \`a une
sous-alg\`ebre de Cartan de $\gothg_\C $ incluse dans $\gothb$.

Les inclusions
$\; [\gothb,\gothb] \subseteq \gothb \subseteq \Ker \xi \;$
et
$\; [\call_f,\gothn_{\lambda}] \subseteq \gothn_{\lambda}
\subseteq (\gothg(\lambda)_\C )^{{\scriptstyle \perp}^
{\!{\scriptscriptstyle \langle} \!~,\!~ {\scriptscriptstyle \rangle}}} \,$
fournissent les \'egalit\'es
$\; f([\gothb,\gothb]) = 0 \;$
et
$\; f([\call_f,\call_f]) = 0$.
Donc $\, \gothg(f)_\C  \subseteq \gothb \,$ et
$\, \call_f / \gothg(f)_\C $ est un lagrangien de
$(\gothg_\C  / \gothg(f)_\C ,B_f)$, clairement stable par $G(f)$.
D'apr\`es \ref{morphisme rho et fonction delta} (b) le morphisme de groupes de
Lie $\rho_f^{\gothg / \gothg(f)}$ de $G(f)^{\gothg/\gothg(f)}$ dans
$\, \C \moins0$ qu'on en d\'eduit v\'erifie
$\;\, \diff _1 \rho_f^{\gothg / \gothg(f)}
= \frac{1}{2} \tr ({\ad \cdot\,}^\C )_{\call_f / \gothg(f)_\C }$.
\cqfr

}

\medskip
J'adapte maintenant les param\`etres de Duflo \`a la situation qui
m'int\'eresse.

\lemme{
\label{espace symplectique canonique}

{\bf (a)}
L'espace vectoriel
$\; \call_{\lambdatilde,{\gotha^*}^+}
= \gothh_\C  \oplus
\Sumpetit_{\alpha \in R^+_{\lambdatilde,{\gotha^*}^+}}\!\gothg_\C ^{\alpha} \;$
(cf. \ref{choix de racines positives} (c)) fournit le lagrangien
$\, \call_{\lambdatilde,{\gotha^*}^+}/\gothh_\C $ de
$(\gothg_\C / \gothh_\C ,B_{\lambda_+})$.
On note
$\rho_{\lambdatilde,{\gotha^*}^+}^{\gothg/\gothh}$\labind{zzr zzl tilde,a*+ g,h}
le caract\`ere complexe de $G(\lambda_+)^{\gothg/\gothh}$ qui lui est associ\'e
(cf. \ref{morphisme rho et fonction delta} (b)).
La restriction de 
$\rho_{\lambdatilde,{\gotha^*}^+}^{\gothg/\gothh}
(\delta_{\lambda_+}^{\gothg/\gothh})^{-1}\!$
\`a $\centragr{G_0}{\gothh}$ est ind\'ependante de ${\gotha^*}^+\!$.

\smallskip
{\bf (b)}
La restriction de $B_{\lambda_+}\!$ au sous-espace vectoriel
$\gothg(\lambda)(\mi \rho_\F+)/\gothh$ de $\gothg/\gothh$ est non
d\'eg\'en\'er\'ee.
L'orthogonal de $\gothg(\lambda)(\mi \rho_\F+)/\gothh$ dans $\gothg/\gothh$ muni
de $B_{\lambda_+}\!$ est la projection dans $\gothg/\gothh$ de la trace sur
$\gothg$ de la somme des $\gothg_\C ^{\alpha}$ lorsque $\alpha$ d\'ecrit
$R(\gothg_\C ,\gothh_\C )$ priv\'e de
$R(\gothg(\lambda)(\mi \rho_\F+)_\C ,\gothh_\C )$.
On l'identifiera \`a $\gothg/\gothg(\lambda)(\mi \rho_\F+)$.

Soit $e \in G(\lambda_+)$ elliptique.
Les orientations
$\,\scalo (B_{\lambda_+})_{(1-\Ad e)(\gothg/\gothg(\lambda)(\mi\rho_\F+))}\,$
et
$\,\scalo (B_{\lambda_{\textit{can}}})
_{(1-\Ad e)(\gothg/\gothg(\lambda)(\mi\rho_\F+))}\,$
sont \'egales.

\smallskip
{\bf (c)}
L'espace vectoriel
$\; \call_{\lambdatilde}
= \gothg(\lambda)(\mi \rho_\F+)_\C  \oplus
\Sumpetit_{\alpha \in R^+_\lambdatilde}\!\gothg_\C ^{\alpha} \;$
(cf. \ref{choix de racines positives} (c)) fournit le lagrangien
$\call_{\lambdatilde}/\gothg(\lambda)(\mi \rho_\F+)_\C \!$ de
$\gothg_\C / \gothg(\lambda)(\mi \rho_\F+)_\C $, \`a la fois pour
$B_{\lambda_{\textit{can}}}$ et pour~$B_{\lambda_+}$.
On note
$\,\rho_\lambdatilde^{\gothg/\!\gothg(\lambda)(\mi \rho_\F+)}\!$%
\labind{zzr zzl tilde g,h}
le caract\`ere complexe de
$G(\lambdatilde)^{\gothg/\!\gothg(\lambda)(\mi \rho_\F+)}\!$
qui est associ\'e \`a
$(\gothg_\C / \gothg(\lambda)(\mi \rho_\F+)_\C ,B_{\lambda_{\textit{can}}})$
(cf. \ref{morphisme rho et fonction delta} (b)).
Il prolonge le caract\`ere complexe de l'image r\'eciproque de $G(\lambda_+)$
dans $G(\lambdatilde)^{\gothg/\!\gothg(\lambda)(\mi \rho_\F+)}\!$ qui est
associ\'e \`a $(\gothg_\C / \gothg(\lambda)(\mi \rho_\F+)_\C ,B_{\lambda_+})$.%

}

\dem{D\'emonstration du lemme}{

{\bf (a)}
Soit $e \in \centragr{G_0}{\gothh}$ elliptique.
L'ind\'ependance de
$\big( \rho_{\lambdatilde,{\gotha^*}^+}^{\gothg/\gothh}
(\delta_{\lambda_+}^{\gothg/\gothh})^{-1} \big) (e)$
par rapport \`a ${\gotha^*}^+\!$ s'obtient en reprenant la d\'emonstration du
lemme \ref{fonctions canoniques sur le revetement double}~(b).

\smallskip
{\bf (b)}
Au d\'ebut, on utilise le fait que $\gothg_\C \!/\gothh_\C \!$ est somme directe
$B_{\lambda_+}\!$-orthogonale des projections des sous-espaces vectoriels
$\gothg_\C ^{\alpha} \oplus \gothg_\C ^{-\alpha}$ de $\gothg_\C $ avec
$\alpha \in R^+(\gothg_\C ,\gothh_\C )$.

On pose $V = \gothg/\gothg(\lambda)(\mi\rho_\F+)$
et $V(e) = \gothg(e)/\gothg(e)(\lambda)(\mi\rho_\F+)$.
L'espace vectoriel $V$ est somme directe des sous-espaces vectoriels $V(e)$ et
$(1-\Ad e)(V)$, qui sont orthogonaux simultan\'ement pour
$B_{\lambda_{\textit{can}}}$ et pour $B_{\lambda_+}$.
Pour chaque $\alpha \in R^+(\gothg_\C ,\gothh_\C )$, on fixe
$X_{\alpha} \in \gothg_\C ^{\alpha}$ et
$X_{-\alpha} \in \gothg_\C ^{-\alpha}$ tels que
$[X_{\alpha},X_{-\alpha}] = H_{\alpha}$.
On constate que les orientations $\scalo (B_{\lambda_+})_V$ et
$\scalo (B_{\lambda_{\textit{can}}})_V$ sont toutes deux dirig\'ees par le
produit ext\'erieur des vecteurs suivants :
les $X_{\alpha} \wedge X_{-\alpha}$ avec $\alpha \in R^+_\lambdatilde$ r\'eelle,
les $\mi (X_{\alpha} \wedge X_{-\alpha})$ avec $\alpha \in R^+_\lambdatilde$
imaginaire,
et les
$X_{\alpha} \wedge X_{-\alpha}
\wedge \overline{X_{\alpha}\!}\, \wedge \overline{X_{-\alpha}\!}\,$
o\`u $\dot{\alpha}$ d\'ecrit les classes des
$\beta \in R(\gothg_\C ,\gothh_\C )
\!\setminus\!R(\gothg(\lambda)(\mi\rho_\F+)_\C ,\gothh_\C )$
complexes modulo l'action du groupe engendr\'e par la conjugaison et le passage
\`a l'oppos\'e.
On prouve de m\^eme que les orientations $\scalo (B_{\lambda_+})_{V(e)}$ et
$\scalo (B_{\lambda_{\textit{can}}})_{V(e)}$ sont \'egales, en utilisant le
syst\`eme de racines $R(\gothg(e)_\C ,\gothh(e)_\C )$ et le couple
$(\lambda,\rho_\F+)$.
Donc
$\;\scalo (B_{\lambda_+})_{(1-\Ad e)(V)}
= \scalo (B_{\lambda_{\textit{can}}})_{(1-\Ad e)(V)}$.%

\smallskip
{\bf (c)}
Soient $e \in G(\lambda_+)$ elliptique et
$\widehat{e} \in G(\lambdatilde)^{\gothg/\!\gothg(\lambda)(\mi \rho_\F+)}\!$
au-dessus de $e$.
On note $\varphi_{\textit{can}}$ et $\varphi_+$ les formes sesquilin\'eaires
hermitiennes $(v,w) \mapsto \mi \,B_{\lambda_{\textit{can}}}(v,\overline{w})$ et
$(v,w) \mapsto \mi \,B_{\lambda_+}(v,\overline{w})$ sur
$\gothg_\C / \gothg(\lambda)(\mi \rho_\F+)_\C $.
L'espace vectoriel
$\call \egdef \call_{\lambdatilde}/\gothg(\lambda)(\mi \rho_\F+)_\C \!$
est somme directe des sous-espaces vectoriels
$\call(e) \egdef (\call_{\lambdatilde} \cap \gothg(e)_\C )
/\gothg(e)(\lambda)(\mi \rho_\F+)_\C \!$
et $(1-\Ad e)(\call)$, qui sont orthogonaux simultan\'ement pour
$\varphi_{\textit{can}}$ et pour $\varphi_+$.
On s'inspire du d\'ebut de la d\'emonstration du lemme
\ref{fonctions canoniques sur le revetement double} (b).
Le nombre de composantes strictement n\'egatives dans la matrice de l'une ou
l'autre des formes sesquilin\'eaires $\varphi_{\textit{can}}$ et $\varphi_+$ sur
$\call$, relativement \`a une base orthogonale, est \'egal \`a :
la moiti\'e du nombre de racines compactes dans
$R(\gothg_\C ,\gothh_\C )
\!\setminus\!R(\gothg(\lambda)(\mi\rho_\F+)_\C ,\gothh_\C )$
plus le quart du nombre de racines complexes dans
$R(\gothg(\nu)_\C ,\gothh_\C )
\!\setminus\!R(\gothg(\lambda)(\mi\rho_\F+)_\C ,\gothh_\C )$.
On a une propri\'et\'e analogue pour $\call(e)$.
D'apr\`es (b) et la d\'efinition \ref{geometrie metaplectique} (d), les fonctions
$\rho$ associ\'ees \`a $B_{\lambda_{\textit{can}}}$ et $B_{\lambda_+}$ prennent
donc les m\^emes valeurs en $\widehat{e}$.
\cqfd

}

\definition{
\label{parametres adaptes}

{\bf (a)}
Soit $\alpha \!\in\! R(\gothg_\C ,\gothh_\C )$ r\'eelle.
On note
$\, n_\alpha = \frac{1}{2} \Sumpetit_{\beta} \beta(H_{\alpha}) \,$
o\`u la somme porte sur les
$\beta \in R(\gothg_\C ,\gothh_\C )$
tels que
$\beta + \overline{\beta} \in \R ^+ \alpha$.
\`A deux vecteurs
$\, X_{\alpha} \in \gothg_\C ^{\alpha}\cap\gothg \,$
et
$\; X_{-\alpha} \in \gothg_\C ^{-\alpha}\cap\gothg \,$
tels que
$\, [X_{\alpha},X_{-\alpha}] = H_{\alpha}$,
on associe l'\'el\'ement elliptique
$\, \gamma_\alpha = \exp(\pi(X_{\alpha} - X_{-\alpha})) \,$
de $\centragr{G_0}{\gothh}$.
(On sait que $\, n_\alpha \in \N \,$ et que l'ensemble
$\{\gamma_\alpha,\gamma_\alpha^{-1}\}$ est ind\'ependant du choix de
$(X_{\alpha},X_{-\alpha})$,  cf. \cite [p.$\!$~327] {DV88}.)

\smallskip
{\bf (b)}
On note $\, \Xirrp_G(\lambdatilde,{\gotha^*}^+) \,$\labind{Xirr+} l'ensemble fini
des classes d'isomorphisme des repr\'esentations unitaires irr\'eductibles
$\tau_+$\labind{zzt plus} de $G(\lambda_+)^{\gothg/\gothh}$
telles que
$\; \tau_+(\iota) = -\id \;$ et
$\; \tau_+(\exp X) = \me^{\mi \,\lambda(X)}\id \;$ pour $X \in \gothh$
(une telle $\tau_+$ est de dimension finie), \\
et
$\Xfin_G(\lambdatilde,{\gotha^*}^+)$\labind{Xfinal} l'ensemble des
$\tau_+ \in \Xirrp_G(\lambdatilde,{\gotha^*}^+)$ tels que pour toute racine
$\,\alpha \in R(\gothg(\lambda)_\C ,\gothh_\C )\,$ r\'eelle,
$(\delta_{\lambda_+}^{\gothg/\gothh} \tau_+) (\gamma_\alpha)$ n'admette pas la
valeur propre $(-1)^{n_\alpha}$.
(Quand $G$ est connexe, la d\'emonstration de \ref{resultat inattendu} (b)
montrera que cette condition peut \^etre remplac\'ee par
$(\delta_{\lambda_+}^{\gothg/\gothh} \tau_+) (\gamma_\alpha)
\not= (-1)^{n_\alpha} \id$.)

\noindent
On note aussi $\,\Xfin_G$\labind{Xf} l'ensemble des
$(\lambdatilde',{{\gotha'}^*}^+\!,\tau_+')$
avec $\;\lambdatilde' = (\lambda',\Fp+) \in \gstregtilde$,
$\gothh' = \gothg(\lambdatilde')$,
$\gotha'$ est la composante hyperbolique de $\gothh'$,
${{\gotha'}^*}^+\!$ est une chambre de
$(\gothg(\lambda')(\mi \rho_\Fp+),\gotha')$,
et $\;\tau_+' \in \Xfin_G(\lambdatilde',{{\gotha'}^*}^+)$.

\smallskip
{\bf (c)}
On note $X_G(\lambdatilde)$
(resp. $\Xirr_G(\lambdatilde)$\labind{Xirr tilde})
l'ensemble des classes d'isomorphisme des repr\'esentations unitaires
(resp. unitaires irr\'eductibles)
$\tau$\labind{zzt} de
$G(\lambdatilde)^{\gothg/\!\gothg(\lambda)(\mi \rho_\F+)}\!$ telles que
$\,\tau(\iota) \!=\! -\id\,$ et
$\;\tau(\exp X) = \me^{\mi \,\lambda(X)}\id\;$ pour $X \in \gothh$,
\\
et $\; \XInd_G(\lambdatilde)$\labind{XInd} l'ensemble (ind\'ependant du choix de
${\gotha^*}^+\!$) form\'e des classes d'isomorphisme des repr\'esentations
unitaires $\tau$ de $G(\lambdatilde)^{\gothg/\!\gothg(\lambda)(\mi \rho_\F+)}$
telles que

\centerline{
$\displaystyle
\frac{\rho_\lambdatilde^{\gothg/\!\gothg(\lambda)(\mi \rho_\F+)}}
{|\rho_\lambdatilde^{\gothg/\!\gothg(\lambda)(\mi \rho_\F+)}|}\;
\tau
\,=\, \Ind_{G(\lambda_+)} ^{G(\lambdatilde)}
\Big(\frac{\rho_{\lambdatilde,{\gotha^*}^+}^{\gothg/\gothh}}
{|\rho_{\lambdatilde,{\gotha^*}^+}^{\gothg/\gothh}|}\;
\tau_+\Big)$
}

\noindent
pour un $\,\tau_+ \in \Xfin_G(\lambdatilde,{\gotha^*}^+)\,$
(donc $\,\XInd_G(\lambdatilde) \subseteq X_G(\lambdatilde)$).

\noindent
On pose aussi :
$\; \XInd_G
= \Bigl\{
(\lambdatilde',\tau') \,;\,
\lambdatilde' \in \gstregtilde
\,\textrm{ et }\,
\tau' \in \XInd_G(\lambdatilde')
\Bigr\}$%
\labind{XI}.

}

\section{
Condition d'int\'egrabilit\'e
}\label{Rep6}

On conserve les notations de la section pr\'ec\'edente.

Le lemme suivant est une adaptation de \cite [remarque 2 p.$\!$~154] {Df82a}.
Il ne sera pas utilis\'e dans la suite de cet article.

\lemme{
\label{condition d'integrabilite}

Les propri\'et\'es suivantes sont \'equivalentes :

\smallskip
{\bf (i)}
$\;\; \Xirrp_G(\lambdatilde,{\gotha^*}^+) \not= \emptyset$ ;

\smallskip
{\bf (ii)}
il existe un caract\`ere unitaire
$\, \chi_{\lambdatilde}^G \,$\labind{zzw zzl+G}
(unique) de $\, G(\lambdatilde)^{\gothg/\gothh}_0$ tel que :

\centerline{
$\chi_{\lambdatilde}^G (\iota) \, = \, -1 \;$
et
$\; \diff _1 \chi_{\lambdatilde}^G \, = \, \mi \restriction{\lambda}{\gothh}$ ;
}

{\bf (iii)}
$\;\; \lambdatilde \in \gstregGtilde$.

}

\dem{D\'emonstration du lemme}{

{\bf (i) \boldmath$\!\!\Leftrightarrow\!$ (ii)} \quad
L'implication << (i) $\!\Rightarrow\!$ (ii) >> est claire.
On suppose (ii) v\'erifi\'e.
La restriction \`a
$\, G(\lambdatilde)^{\gothg/\gothh}_0$
de la repr\'esentation unitaire
$\; \Ind_{G(\lambda_+)^{\gothg/\gothh}_0}^{G(\lambda_+)^{\gothg/\gothh}}
\chi_{\lambdatilde}^G \;$
de $\, G(\lambda_+)^{\gothg/\gothh}$
est non nulle (car induite \`a partir d'un espace non nul) et multiple de
$\, \chi_{\lambdatilde}^G$.
Le th\'eor\`eme 8.5.2 de \cite [p.$\!$~153]{Di64} montre donc l'existence d'un
\'el\'ement de $\, \Xirrp_G(\lambdatilde,{\gotha^*}^+)$.

\smallskip
{\bf (ii) \boldmath$\!\!\Leftrightarrow\!$ (iii)} \quad
D'apr\`es le lemme \ref{fonctions canoniques sur le revetement double} (a), en
multipliant $\, \chi_{\lambdatilde}^G \,$ par
$\;\rho_{\lambda_+}^{\gothg/\gothh} \,$ la condition (ii)
\'equivaut \`a l'existence d'un caract\`ere complexe du groupe de Lie
$\, G(\lambda_+)_0 \,$ de diff\'erentielle
$\, \mi \lambda + \rho_{\gothg,\gothh}$.
Cette derni\`ere propri\'et\'e s'\'ecrit :
$\; \lambdatilde \in \gstregGtilde$.
\cqfd

}

\remarque{
\label{autre insuffisance de la methode des orbites}

Soit $\, f \in \Supp_{\gothg(\lambda)^*}(G(\lambda)_0 \cdot \lambdatilde)$.
En particulier, on a $\, f \in \gstreg$.

Le groupe de Lie r\'eel $G(f)_0$ est commutatif
(cf. \cite [th. p.$\!$~17]{B.72}).
Le stabilisateur de la composante nilpotente de $f$ dans le groupe adjoint de
$G(\lambda)_0$ est lin\'eaire alg\'ebrique unipotent d'apr\`es
\cite [th. 5.9 (b) p.$\!$~138] {Sp66}.
Le tore maximal de $G(f)_0$ est donc \'egal \`a celui de
$\centregr{G(\lambda)_0}$.
Compte tenu de la remarque \ref{sous-algebre de Borel canonique}, les
arguments de la d\'emonstration du lemme pr\'ec\'edant permettent d'en d\'eduire
que :

\smallskip\centerline{
$\Xirr_G(f) \not= \emptyset \iff
\forall Z \in \Ker \exp_{\centregr{G(\lambda)_0}} \;
\me^{(\mi \,\lambda+\rho_{\gothg,\gothh})(Z)}=1$.
}

\noindent
Donc la condition
<< $\Xfin_G(\lambdatilde,{\gotha^*}^+) \not= \emptyset$ >>
implique la condition
<< $\Xirr_G(f) \not= \emptyset$ >>.%

La r\'eciproque est fausse, comme le montre le cas o\`u $\, G = PSL(2,\R )$,
$\lambda = 0$ et $\gothh = \goth{so}(2)$, donc
$\Xfin_G(\lambdatilde,\{0\}) = \emptyset$,
$\Xirr_G(f) \not= \emptyset$, et aussi $\Xirr_G(-f) \not= \emptyset$.
Cela traduit le fait que $PSL(2,\R )$ n'a pas de repr\'esentation limite
de sa s\'erie discr\`ete hors de sa s\'erie discr\`ete.
Pour g\'en\'eraliser la m\'ethode des orbites propos\'ee par M.~Duflo, il \'etait
naturel d'imaginer que les orbites coadjointes r\'eguli\`eres de $G$ susceptibles
d'\^etre l'orbite $G \!\cdot\! l$ associ\'ee \`a une repr\'esentation
$T_{l,\tau}^G$ \`a caract\`ere infinit\'esimal nul seraient $G \!\cdot\! f$ et
$G \!\cdot\! (-f)$ (cf. \ref{caracteres canoniques} (a), en
rempla\c{c}ant $(\lambdatilde,{\gotha^*}^+)$ et $\lambda$ par $l$).
Cette g\'en\'eralisation aurait cr\'e\'e des difficult\'es, car
$\,G \!\cdot\! f \not= G \!\cdot\! (-f)$ et le dual unitaire de $G$ n'a qu'un
seul \'el\'ement \`a caract\`ere infinit\'esimal nul.
\cqfr

}

%
\part{
Construction de repr\'esentations
}\label{RepIII}
%

Cette partie suit de pr\`es l'article \cite {Df82a} de M.~Duflo, pages 160 \`a
180.

\medskip
On se donne
$\,\lambdatilde = (\lambda,\F+) \in \gstregtilde$.
On pose $\gothh = \gothg(\lambdatilde)$.
On note $\mu$ et $\nu$ (resp. $\gotht$ et $\gotha$) les composantes
infinit\'esimalement elliptique
et hyperbolique de $\lambda$ (resp. de $\gothh$).
On fixe aussi une chambre
${\gotha^*}^+\!$ de $(\gothg(\lambda)(\mi \rho_\F+),\gotha)$
et $\, \tau_+ \in \Xfin_G(\lambdatilde,{\gotha^*}^+)$.
D'o\`u $\; R^+(\gothg_\C ,\gothh_\C )$ (cf. \ref{choix de racines positives}
(b)).
On pose $\, \mutilde = (\mu,\F+ \!\cap \mi\gotht^*)$\labind{zzm tilde}.

\section{\boldmath
Cas $G$ connexe
}\label{Rep7}

Dans toute cette section, on suppose $G$ connexe.

\smallskip
L'inclusion $\, \Ad G \subseteq \interieur \gothg_\C  \,$ assure que
$\, G(\lambda_+) \,$ est \'egal \`a $\, \centragr{G}{\gothh}$, et que
$\, \delta_{\lambda_+}^{\gothg/\gothh} \,$ est un caract\`ere unitaire de
$\, G(\lambda_+)^{\gothg/\gothh} \,$ (cf. \ref{morphisme rho et fonction delta},
et aussi \cite [haut p.$\!$~118] {Df84}).
En effet pour toute composante connexe ${\gotha^*}^+_0$ de l'ensemble des
\'el\'ements de $\gotha^*$ qui ne s'annulent sur aucun $H_{\alpha}$ avec
$\, \alpha \in R(\gothg_\C ,\gothh_\C ) \,$
non imaginaire, l'espace vectoriel

\smallskip\centerline{
$\call^+
\;\egdef\; \gothh_\C  \oplus\!
\Sum_{\substack
{\alpha \in R^+(\gothg_\C ,\gothh_\C ) \\
\textrm{$\alpha$ imaginaire non compacte}}} \!\!\!
\gothg_\C ^{\alpha}
\oplus
\Sum_{\substack
{\alpha \in R^+(\gothg_\C ,\gothh_\C ) \\ \textrm{$\alpha$ compacte}}} \!
\gothg_\C ^{-\alpha}
\oplus
\Sum_{\substack
{\alpha \in R(\gothg_\C ,\gothh_\C ) \\
{\gotha^*}^+_0(H_{\alpha}) \,\subseteq\, \R ^+\moins0}}
\! \gothg_\C ^{\alpha}$
}

\noindent
fournit le lagrangien positif $\call^+ / \gothh_\C$ de 
$(\gothg_\C / \gothh_\C ,B_{\lambda_+})$
stable par $G(\lambda_+)$.

\smallskip
On note
$M$\labind{M} l'intersection des noyaux des caract\`eres r\'eels positifs de
$\centragr{G}{\gotha}$ et $\gothm$\labind{m} son alg\`ebre de Lie
(donc $\,\mutilde \in \mstIMtilde$ par \ref{construction de rep, cas connexe} (a)
et $\Xirr_M(\mutilde) = \Xfin_M(\mutilde,\{0\})$),
\\
$R^+(\gothm_\C ,\gotht_\C )
= R(\gothm_\C ,\gotht_\C ) \cap R^+(\gothg_\C ,\gothh_\C ) \;$
(notation de \ref{choix de racines positives} (b) relative \`a $\gothm$ et
$(\mutilde,\{0\})$),
\\
$\mu_{+,\gothm} = \mu_{\gothm,\mutilde,\{0\}}$\labind{zzm +m}
l'\'el\'ement $\, \mu - 2\mi \rho_{\gothm,\gotht} \,$ de $\gotht^*$
(donc $M(\mutilde) \!=\! \centragr{M}{\gotht}$ et
$\delta_{\mu_{+,\gothm}}^{\gothm / \gotht}$ est un caract\`ere unitaire
de $M(\mutilde)^{\gothm / \gotht}$),
\\
${\gothn}_M
= \!\!\! \Sum_{\alpha \in R^+(\gothm_\C ,\gotht_\C )}\! \gothg_\C ^{\alpha} \;$
et
$\; \gothb_M = \gotht_\C  \oplus {\gothn}_M$,
\\
$q 
= \abs{\{ \alpha \in R^+(\gothm_\C ,\gotht_\C ) \mid
\textrm{$\alpha$ compacte} \}}$%
\labind{q},
\\
$\gothk_M
= \gotht \oplus
\bigl( \Sum_{\substack
{\alpha \in R^+(\gothm_\C ,\gotht_\C ) \\ \textrm{$\alpha$ compacte}}}
\! \gothg_\C ^{\alpha}
\; \bigr)
\cap \gothg \;$
et
$\, K_{M_0} = \exp \gothk_M$.

\definition{
\label{limites de series discretes}

\negsmallskip
{\bf (a)}
On note $T_{\mutilde}^{M_0}$ la classe de repr\'esentation << limite de la
s\'erie discr\`ete de $M_0$ >> appartenant \`a $\widehat{M_0}$, dont l'espace
des vecteurs $K_{M_0}$-finis est isomorphe \`a l'image du caract\`ere unitaire
de $T_0$ de diff\'erentielle $\, \mi \mu \!-\! \rho_{\gothm,\gotht} \,$ par le
foncteur  << d'induction cohomologique >>
$\; \calr_{M_0}^q
\egdef
\left({}^{u\!}\calr_{\gothb_M,T_0}^{\gothm_\C ,K_{M_0}}\right)^{\!q}
(\, {\scriptscriptstyle\bullet} \otimes
\bigwedge^{\scriptscriptstyle\mathrm{max}}{\gothn}_M ) \;$
d\'ecrit dans \cite [(11.73) p.$\!$~677] {KV95}.
(Elle est not\'ee $\pi_{\mi \mu,\gothb_M}$ dans \cite [bas p.$\!$~734] {KV95},
en tenant compte de \cite [proposition 11.180 p.$\!$~733] {KV95}.)

\smallskip
{\bf (b)}
On note
$\; T_{\mutilde,\sigma}^M
= \, \Ind_{\centragr{M}{\gothm} . M_0}^M
(\delta_{\mu_{+,\gothm}}^{\gothm / \gotht} \sigma
\otimes T_{\mutilde}^{M_0} ) \;$
pour tout $\, \sigma \in \Xirr_M(\mutilde)$, \\
o\`u on a pu d\'efinir une repr\'esentation
$\delta_{\mu_{+,\gothm}}^{\gothm / \gotht} \sigma \otimes T_{\mutilde}^{M_0}$
de  $\centragr{M}{\gothm} M_0$ par l'\'egalit\'e

\centerline{
$\left( \delta_{\mu_{+,\gothm}}^{\gothm / \gotht} \sigma
\otimes T_{\mutilde}^{M_0} \right) (xy)
\,=\,
(\delta_{\mu_{+,\gothm}}^{\gothm / \gotht} \sigma)(x)
\otimes T_{\mutilde}^{M_0}(y) \;$
pour $\, x \in \centragr{M}{\gothm} \,$ et $\, y \in M_0$
}

\noindent
(dans laquelle
$\, (\delta_{\mu_{+,\gothm}}^{\gothm / \gotht} \sigma)(x) = \sigma(x,1)$),
car le groupe $\centregr{M_0}$ op\`ere dans
<< l'espace >> de $T_{\mutilde}^{M_0}$ par le caract\`ere unitaire
$\, \restriction{(\delta_{\mu_{+,\gothm}}^{\gothm / \gotht} \sigma)}
{\centregr{M_0}}$
d'apr\`es \cite[(11.184c) p.$\!$~734]{KV95}.

}

\remarque{
\label{cohomologie pour $M$}

On fixe un sous-groupe compact maximal $K$ de $G$ dont l'involution de Cartan
normalise~$\gothh$.
Donc $\, K_M \egdef K \cap M \,$ est un sous-groupe
compact maximal de $M$ contenant $\centragr{M}{\gothm} K_{M_0}$.
Tout $\, \pi \in (\centragr{M}{\gothm} M_0)\widehat{~} \,$ est reli\'e \`a
certaines repr\'esentations
$\, \pi_1 \in (\centragr{M}{\gothm})\widehat{~} \,$ et
$\, \pi_2 \in (M_0)\widehat{~} \,$ par les \'egalit\'es
$\; \pi(xy) = \pi_1(x) \otimes \pi_2(y) \;$ pour $x \in \centragr{M}{\gothm}$ et
$y \in M_0$ (cf. \cite [prop. 13.1.8 p.$\!$~251]{Di64}).
Soit $\, \sigma \in \Xirr_M(\mutilde)$.
D'apr\`es
\cite [(11.187) p.$\!$~735 et dem. de prop. 11.192 (a) p.$\!$~737] {KV95}
et \ref{fonctions canoniques sur le revetement double} (b) pour passer de
$M_0$ \`a $\centragr{M}{\gothm} M_0$,
et \cite [prop. 11.57 p.$\!$~672 et fin (5.8) p.$\!$~332] {KV95} pour passer du
<< sous-groupe parabolique >> $\centragr{M}{\gothm} M_0$ de $M$ \`a $M$, le
$(\gothm_\C ,K_M)$-module associ\'e \`a $T_{\mutilde,\sigma}^M$ est isomorphe \`a
$\,\calr_M^q((\rho_{\mu_{+,\gothm}}^{\gothm / \gotht}\!)^{-1} \sigma)$, o\`u
$\calr_M^q$ est le foncteur d'induction cohomologique relatif \`a $\gothb_M$.
\cqfr

}

\lemme{
\label{construction de rep, cas connexe}

{\bf (a)}
On a
$\;\, G(\lambda_+) = \centragr{M}{\gothm} \exp \gothh \;$
et
$\; M(\mutilde) = \centragr{M}{\gothm} \exp \gotht$. \\
Il existe donc un unique $\, \tau_M \in \Xirr_M(\mutilde) \,$\labind{zzt M}
tel que
$\;\, \delta_{\mu_{+,\gothm}}^{\gothm / \gotht} \tau_M
= \restriction{(\delta_{\lambda_+}^{\gothg/\gothh} \tau_+)}{M(\mutilde)}$. \\
(Cette notation ne tient pas compte du fait que $\tau_M$ d\'epend de
${\gotha^*}^+\!$.)

\smallskip
{\bf (b)}
La classe de repr\'esentation
$\; \Ind_{M . A . N}^G
( T_{\mutilde,\tau_M}^M \!\otimes \me^{\mi \, \nu \circ \Log} \otimes
\fonctioncar{N} ) \;$
construite \`a partir du choix d'un sous-groupe $\,N\,$ de $\,G\,$ qui est
radical unipotent d'un sous-groupe parabolique de $G$ de composante d\'eploy\'ee
$A$, appartient \`a $\widehat{G}$ et est ind\'ependante de $N\!$.

}

\dem{D\'emonstration du lemme}{

{\bf (a)}
On a
$\, G(\lambda_+) = \centragr{G}{\gothh}
\subseteq \centragr{G}{\gotha} = MA$
et
$M(\mutilde) = \centragr{M}{\gotht}$.
Vu \cite [th. 17 p.$\!$~199] {Va77}, on en d\'eduit que
$\, M(\mutilde) = \centragr{M}{\gothm} \exp \gotht$,
puis
$\, G(\lambda_+) = \centragr{M}{\gothm} \exp \gothh$.

\smallskip
{\bf (b)}
Soit $N$ le radical unipotent d'un sous-groupe parabolique de $G$ de composante
d\'eploy\'ee $A$.
D'apr\`es \cite [th. 18 p.$\!$~289] {Va77}, il existe une composante
connexe ${\gotha^*}^+_0$ de l'ensemble des \'el\'ements de $\gotha^*$ qui ne
s'annulent sur aucun $H_{\alpha}$ avec
$\, \alpha \in R(\gothg_\C ,\gothh_\C ) \,$ non imaginaire, telle que

\centerline{
$N = \exp \gothn \;\,$ avec
$\;\, \gothn
= \bigl(
\Sum_{\alpha \in R(\gothg_\C ,\gothh_\C ) \textrm{ et }
{\gotha^*}^+_0(H_{\alpha}) \,\subseteq\, \R ^+\moins0}
\! \gothg_\C ^{\alpha}
\; \bigr) \cap\, \gothg$.
}

\noindent
Soient $\nu_0 \in {\gotha^*}^+_0$ et
$\, \lambda_0 \in (\gotht^* + \nu_0) \cap \gstreg$.
On utilise la proposition \ref{caracteres d'induites} (a) avec $\lambda_0$ \`a
la place de $\lambda$
(auquel cas << $(M',U)$ >> devient << $(\centragr{M}{\gothm} \, M_0 \, A,N)$ >>)
et
$\;\, \pi =
(\delta_{\mu_{+,\gothm}}^{\gothm / \gotht} \tau_M \otimes T_{\mutilde}^{M_0})
\otimes \me^{\mi \, \nu \circ \Log}$.
Comme la classe de
$\; \Ind_{M . A . N}^G
( T_{\mutilde,\tau_M}^M \!\otimes \me^{\mi \, \nu \circ \Log} \otimes
\fonctioncar{N} ) \;$
est d\'etermin\'ee par son caract\`ere (cf. \cite [p.$\!$~64] {Co88}), elle ne
d\'epend pas de $N$.

\smallskip
On choisit $z=1$ et $\; \Lambda = (\Lambda^{can},R^+_{\mi \, \R },R^+_\R ) \;$
dans \cite [ligne 7 du bas p.$\!$~121] {ABV92}, o\`u
$\; \Lambda^{can} \otimes \rho(R^+(\gothg_\C ,\gothh_\C ))
= \rho_{\lambda_+}^{\gothg/\gothh} \, \tau_+ \;$
(cf. \cite [lignes 9 et 10 du bas p.$\!$~129] {ABV92})
et les syst\`emes de racines positives $R^+_{\mi \, \R }$ et $R^+_\R $ sont
inclus dans $R^+(\gothg_\C ,\gothh_\C )$. \\
D'apr\`es le lemme \ref{fonctions canoniques sur le revetement double} (b), la
repr\'esentation $\wt{\Lambda}$ de \cite [p.$\!$~123] {ABV92} est
$\, \wt{\Lambda} = \delta_{\lambda_+}^{\gothg/\gothh} \tau_+$.
Vu la proposition \ref{caracteres d'induites} (a), on a
$\; \pi(\Lambda)
= \Ind_{M . A . N}^G
( T_{\mutilde,\tau_M}^M \!\otimes \me^{\mi \, \nu \circ \Log} \otimes
\fonctioncar{N} ) \;$
dans \cite [(11.2)(e) p.$\!$~122] {ABV92}.
Les id\'ees de la d\'emonstration du lemme 5 de \cite [p.$\!$~335] {DV88}
permettent de montrer que $\Lambda$ est << final >> au sens de \cite [def. 11.13
p.$\!$~130] {ABV92}.
On applique enfin \cite [th. 11.14 (a) p.$\!$~131] {ABV92}.
\cqfd

}

\remarque{

On se place dans le cas o\`u $G = SL(3,\R )$ et $\lambda$ est l'\'el\'ement de
$\gstssIncG$ image de $\left( {\scriptscriptstyle \diagt{1}{1}{-2}} \right)$ par
l'isomorphisme de $G$-modules qui envoie $A \in \gothg$ sur
$\tr (A\,.) \in \gothg^*\!$.
Le groupe $MA$ qui est ici \'egal \`a $G(\lambda)$, est l'image de
$GL(2,\R )$ par le plongement de groupes de Lie qui envoie
$\, x \in GL(2,\R ) \,$ sur
$\; \left( {\scriptstyle \diagd{x}{(\det x)^{-1}}} \right) \in G$.
Mais l'induite
$\, \Ind_{M . A . N}^G
( T_{\mutilde,\tau_M}^M \!\otimes \me^{\mi \, \nu \circ \Log} \otimes
\fonctioncar{N} ) \,$
n'est pas donn\'ee par << nondegenerate data >> au sens de
\cite [p.$\!$~473] {KZ82}, car le groupe de Weyl d'une sous-alg\`ebre de Cartan
fondamentale de $\goth{gl}(2,\R )$ se repr\'esente dans $GL(2,\R )$.
\cqfr

}

\definition{

On note
$\;\, T_{\lambdatilde,{\gotha^*}^+\!,\tau_+}^G\!
= \Ind_{M . A . N}^G
( T_{\mutilde,\tau_M}^M \!\otimes \me^{\mi \, \nu \circ \Log} \otimes
\fonctioncar{N} )$%
\labind{Tzzl tilde,a*,zzt,connexe},
ind\'ependamment du choix de $N$ comme dans le lemme ci-dessus.
 
Ainsi lorsque $M$ est connexe, $\Xirr_M(\mutilde)$
a pour seul \'el\'ement $\chi_{\mutilde}^M\!$, et les classes des
repr\'esentations $\, T_{\mutilde,\chi_{\mutilde}^M}^M$ et
$\, T_{\mutilde,\{0\},\chi_{\mutilde}^M}^M$ sont toutes deux \'egales \`a
$\, T_{\mutilde}^{M_0} \!$.

}

\section{\boldmath
Les repr\'esentations $T_{\lambdatilde,{\gotha^*}^+\!,\tau_+}^G\!$
}\label{Rep8}

On va appliquer la construction du cas connexe \`a $\,M'_0 \egdef G(\nu_+)_0$.
Un passage par l'homologie va permettre (gr\^ace \`a un r\'esultat de D.~Vogan)
de normaliser des op\'erateurs d'entrelacement de fa\c{c}on \`a prolonger
$T_{\lambdatilde}^{M'_0}$ en une repr\'esentation d'un certain groupe $M'$.

\smallskip
On note
$\, M' = G(\lambda_+) G(\nu_+)_0 \,$\labind{M'} et
$\, \gothm' = \gothg(\nu_+) \,$\labind{m'} son alg\`ebre de Lie
\\
(donc $\,\lambdatilde \in \mpstfondMptilde$
par \ref{construction equivalente, cas connexe} (a),
$\gothm'(\lambda)(\mi \rho_\F+) = \gothh\,$ et
$\,\Xirr_{M'}(\lambdatilde) = \Xfin_{M'}(\lambdatilde,\gotha^*)$),
\\
$R^+(\gothm'_\C ,\gothh_\C )
= R(\gothm'_\C ,\gothh_\C )
\cap R^+(\gothg_\C ,\gothh_\C ) \;$
(notation de \ref{choix de racines positives} (b) relative \`a $\gothm'$ et
$(\lambdatilde,\gotha^*)$),
\\
$\lambda_{+,\gothm'} = \lambda_{\gothm'\!,\lambdatilde,\gotha^*\!,1} \;$%
\labind{zzl +m'}
l'\'el\'ement $\mu_+ + \nu$ de $\gothh^*$ o\`u $\mu_+$ repr\'esente encore
$\mu_{\gothg,\lambdatilde,{\gotha^*}^+}\!$
\\
(donc $\,B_{\lambda_{+,\gothm'}}\!$ est restriction de $B_{\lambda_+}$,
$M'(\lambdatilde) = G(\lambda_+)$
et
$\,\rho_\lambdatilde^{\gothm'\!/\gothh}
= \rho_{\lambda_{+,\gothm'}}^{\gothm'\!/\gothh}$),
\\
$\gothn_{M'}
= \! \Sum_{\alpha \in R^+(\gothm'_\C ,\gothh_\C )} \!
\gothg_\C ^{\alpha} \;\,$%
\labind{nM'}
et
$\; \gothb_{M'} = \gothh_\C  \oplus \gothn_{M'}$\labind{bM'},
\\
$q'
= \abs{\{ \alpha \in R^+(\gothm'_\C ,\gothh_\C ) \mid
\textrm{$\alpha$ compacte} \}}
\,+\,
\frac{1}{2} \, \abs{\{ \alpha \in R^+(\gothm'_\C ,\gothh_\C ) \mid
\textrm{$\alpha$ complexe} \}}$%
\labind{q'},
\\
$M_{\nu_+}$ l'intersection des noyaux des caract\`eres r\'eels positifs de
$\centragr{M'_0}{\gotha}$,
\\
$T_{\mutilde}^{M_{\nu_+}}
\!= T_{\mutilde,\chi_{\mutilde}^M}^{M_{\nu_+}}\!$
(cf. \ref{limites de series discretes} (b))
et
$T_{\lambdatilde}^{M'_0}
\!= T_{\lambdatilde,\gotha^*,\chi_{\lambdatilde}^{M'}}^{M'_0}$%
\labind{Tzzl ztilde,M'0}
o\`u $\;M_{\nu_+}(\mutilde) = \exp \gotht$
et
$M'_0(\lambdatilde) = \exp \gothh$,
\\
$\gothu
=\! \bigl( \!\Sum_{\alpha \in R(\gothg_\C ,\gothh_\C ) \textrm{ et }
\nu_+(H_{\alpha})>0} \! \gothg_\C ^{\alpha} \bigr)
\cap \gothg \,$%
\labind{u}
et $U = \exp \gothu$\labind{U}
(donc $\interieur M' . U \!\subseteq\! U$ et $M' \cap U = \{1\}$).%

\proposition{
\label{construction de rep, cas de $M'$}

On note $\calh$\labind{H} << l'espace >> de $\, T_{\lambdatilde}^{M'_0}$ et
$\,\calh^{\infty}$ l'ensemble de ses vecteurs~$C^{\infty}\!$.

\smallskip
{\bf (a)}
Le sous-espace propre de $H_{q'} (\gothn_{M'},\calh^{\infty})^*$ de poids
$-(\mi \lambda + \rho_{\gothm'\!,\gothh})$ sous l'action de $\gothh_\C $
issue de l'action de $\gothh_\C $ dans
$\; \bigwedge \gothn_{M'} \otimes_\C  \calh^{\infty}\!$,
est de dimension $1$.
On le note
$\, (H_{q'} ( \, \gothn_{M'},\calh^{\infty} )^*)
_{-(\mi \,\lambda + \rho_{\gothm'\!,\gothh})}$.

\smallskip
{\bf (b)}
Il existe une unique repr\'esentation unitaire continue $S$\labind{S} de
$M'(\lambdatilde)^{\gothm'\!/\gothh}$ dans $\calh$ satisfaisant les conditions
(i) et (ii) suivantes :

\hspace{1cm} (i)
$\; S(\hat{x}) \, T_{\lambdatilde}^{M'_0}\!(y) \,
S(\hat{x})^{-1}
=\, T_{\lambdatilde}^{M'_0}\!(xyx^{-1}) \;\;$
pour $\, \hat{x} \in M'(\lambdatilde)^{\gothm'\!/\gothh}$
au-dessus de $\, x \in M'(\lambdatilde)$ et $\, y \in M'_0$ ;

\hspace{1cm} (ii)
l'action de $M'(\lambdatilde)^{\gothm'\!/\gothh}$ dans
$\, (H_{q'} ( \, \gothn_{M'},\calh^{\infty} )^*)
_{-(\mi \,\lambda + \rho_{\gothm'\!,\gothh})}$
issue de l'action de $M'(\lambdatilde)^{\gothm'\!/\gothh}$
dans $\, \bigwedge \gothn_{M'} \otimes_\C  \calh^{\infty}$ d\'eduite
de $S$, est $\, (\rho_{\lambda_{+,\gothm'}}^{\gothm'\!/\gothh})^{-1} \id$.

}

\dem{D\'emonstration de la proposition}{

{\bf (a)}
On suppose $G$ connexe et utilise les notations de la section \ref{Rep7}

On fixe un sous-groupe compact maximal $K$ de $G$ dont l'involution de Cartan
normalise $\gothh$.
Donc $\, K_{M'_0} \egdef K \cap M'_0 \,$
et $\, K_{M_{\nu_+}} \egdef K \cap M_{\nu_+} \,$
sont des sous-groupes compacts maximaux de $M'_0$ et $M_{\nu_+}$.
D'apr\`es la remarque \ref{cohomologie pour $M$}, le
$(\gothm_\C ,K_{M_{\nu_+}})$-module associ\'e \`a
$T_{\mutilde}^{M_{\nu_+}}$ est isomorphe \`a
$\, \calr_{M_{\nu_+}}^q
( (\rho_{\mu_{+,\gothm}}^{\gothm / \gotht}\!)^{-1} \chi_{\mutilde}^M)$,
o\`u $\, \calr_{M_{\nu_+}}^q$ est le foncteur d'induction cohomologique
relatif \`a $\gothb_M$.

Compte tenu de la d\'efinition \ref{geometrie metaplectique} (d), et du
lemme \ref{fonctions canoniques sur le revetement double} (b) avec ses
notations, pour tout $\, e \in \centragr{M}{\gothm} T_0 \,$ on a

\smallskip\centerline{
$\displaystyle \frac
{(\rho_{\lambda_{+,\gothm'}}^{\gothm'\!/\gothh}
\delta_{\lambda_{+,\gothm'}}^{\gothm'\!/\gothh})(e)}
{(\rho_{\mu_{+,\gothm}}^{\gothm / \gotht}
\delta_{\mu_{+,\gothm}}^{\gothm / \gotht})(e)}
= \det (\Ad e^\C )_{\gothn_{M'}\!/\gothn_M}
\,{\scriptstyle\times}\,
\frac
{((\rho_{\lambda_{+,\gothm'}}^{\gothm'\!/\gothh})^{-1}
\delta_{\lambda_{+,\gothm'}}^{\gothm'\!/\gothh})(e)}
{((\rho_{\mu_{+,\gothm}}^{\gothm / \gotht}\!)^{-1}
\delta_{\mu_{+,\gothm}}^{\gothm / \gotht})(e)}
=\, u_1 \,\cdots\, u_b$
}

\noindent
o\`u
$\, \{\beta_1,\!-\overline{\beta_1}\},\dots,\{\beta_b,\!-\overline{\beta_b}\}$
sont les \'el\'ements de $\wt{R}^+_\C (\gothm'_\C ,\gothh_\C )$
et $u_1,\dots,u_b$ sont les rapports des homoth\'eties
$(\Ad e^\C )_{\gothg_\C ^{\beta_1}},\dots,(\Ad e^\C )_{\gothg_\C ^{\beta_b}}$.

On applique \cite [th. 11.225 p.$\!$~759 et cor. 8.28 p.$\!$~566] {KV95} au
groupe $M'_0$ et \`a la sous-alg\`ebre parabolique $\gothb_{M'}$ de
$\gothm'_\C $.
Le $(\gothm'_\C ,K_{M'_0})$-module $\,\calh^f$ associ\'e \`a
$T_{\lambdatilde}^{M'_0}$ est donc irr\'eductible et isomorphe \`a
$\, \calr_{M'_0}^{q'}
( (\rho_\lambdatilde^{\gothm'\!/\gothh})^{-1} \chi_{\lambdatilde}^{M'})$,
o\`u $\, \calr_{M'_0}^{q'}$ est le foncteur
d'induction cohomologique relatif \`a $\gothb_{M'}$.
On obtient

\centerline{
$\dim \, \mathrm{Hom}_{\gothh_\C ,T_0}
\bigl( H_{q'} ( \gothn_{M'},\calh^f),
\rho_\lambdatilde^{\gothm'\!/\gothh} \chi_{\lambdatilde}^{M'} \bigr)
= 1$
}

\noindent
en prenant $\, X = \calr_{M'_0}^{q'}(Z) \,$ avec
$\, Z = (\rho_\lambdatilde^{\gothm'\!/\gothh})^{-1}
\chi_{\lambdatilde}^{M'}$
dans \cite [prop. 8.11 p.$\!$~555] {KV95}, et en tenant compte du lemme de Schur
(cf. \cite [lem. 3.3.2 p.$\!$~80] {Wa88}).
Comme le morphisme canonique de $\,(\gothh_\C ,T_0)$-modules de
$\; H_{q'} ( \, \gothn_{M'},\calh^f ) \;$ dans
$\, H_{q'} ( \, \gothn_{M'},\calh^{\infty} ) \,$ est bijectif d'apr\`es
\cite [lem. 4 p.$\!$~165] {Df82a}, cela donne le r\'esultat.%

\smallskip
{\bf (b)}
D'apr\`es \ref{injection dans le dual} (a), on a
$\; \interieur x \cdot T_{\lambdatilde}^{M'_0} \!= T_{\lambdatilde}^{M'_0}$
pour tout $\, x \in M'(\lambdatilde)$.
On reprend ensuite la d\'emonstration de \cite [lem. 6 p.$\!$~169] {Df82a},
quasiment mot \`a mot.
(L'unicit\'e de $S$ provient bien s\^ur de
\ref{construction de rep, cas connexe} (b) et du lemme de Schur.)
\cqfd

}

\lemme{
\label{construction equivalente, cas connexe}

On conserve les notations de la proposition pr\'ec\'edente.

\smallskip
{\bf (a)}
Il existe un unique $\, \tau_{M'} \in \Xirr_{M'}(\lambdatilde) \,$\labind{zzt M'}
tel que
$\, \delta_{\lambda_{+,\gothm'}}^{\gothm'\!/\gothh} \tau_{M'}\!
= \delta_{\lambda_+}^{\gothg/\gothh} \tau_+$.
On note $\, \tau_{M'} \otimes S T_{\lambdatilde}^{M'_0}$
la repr\'esentation de $M'$ d\'efinie (clairement sans ambigu\"{\i}t\'e)~par

\centerline{
$\bigl( \tau_{M'} \otimes S T_{\lambdatilde}^{M'_0} \bigr)\!(xy)
= \tau_{M'}(\hat{x}) \otimes S(\hat{x}) \, T_{\lambdatilde}^{M'_0}\!(y)$
}

\noindent
pour $\, \hat{x} \in M'(\lambdatilde)^{\gothm'\!/\gothh}$
au-dessus de $\, x \in M'(\lambdatilde)$ et $\, y \in M'_0$.

\smallskip
{\bf (b)}
Lorsque $G$ est connexe, on a :
$\, T_{\lambdatilde,{\gotha^*}^+\!,\tau_+}^G\!\!
= \Ind_{M' \!. U}^G
( (\tau_{M'} \otimes S \, T_{\lambdatilde}^{M'_0}) \otimes \fonctioncar{U} )$.%

}

\dem{D\'emonstration du lemme}{

\negsmallskip
{\bf (a)}
D'apr\`es le lemme \ref{fonctions canoniques sur le revetement double} (b), on a
$\;\, \displaystyle
\frac{\rho_{\lambda_{+,\gothm'}}^{\gothm'\!/\gothh}}
{|\rho_{\lambda_{+,\gothm'}}^{\gothm'\!/\gothh}|} \,
(\delta_{\lambda_{+,\gothm'}}^{\gothm'\!/\gothh})^{-1}
=\, \frac{\rho_{\lambda_+}^{\gothg/\gothh}}
{|\rho_{\lambda_+}^{\gothg/\gothh}|} \,
(\delta_{\lambda_+}^{\gothg/\gothh})^{-1}$,
o\`u $\rho_{\lambda_{+,\gothm'}}^{\gothm'\!/\gothh}\!$ et
$\rho_{\lambda_+}^{\gothg/\gothh}\!$ sont des morphismes de groupes de Lie.
Cela permet de d\'efinir~$\tau_{M'}$.

\smallskip
{\bf (b)}
On suppose $G$ connexe et utilise les notations de la section \ref{Rep7} \\
On a donc
$\, M' = \centragr{M}{\gothm} M'_0$,
$\, M'(\lambdatilde) = \centragr{M}{\gothm} \exp \gothh \,$
et $\delta_{\lambda_{+,\gothm'}}^{\gothm'\!/\gothh}$ est un caract\`ere unitaire
de $M'(\lambdatilde)^{\gothm'\!/\gothh}\!$.

On choisit une composante connexe ${\gotha^*}^+_0$ de l'ensemble des
\'el\'ements de~$\gotha^*$ qui ne s'annulent sur aucun $H_{\alpha}$ avec
$\, \alpha \in R(\gothg_\C ,\gothh_\C ) \,$ non imaginaire, dont l'adh\'erence
contient $\nu_+$.
On pose

\centerline{
$\gothn
= \bigl(
\Sum_{\alpha \in R(\gothg_\C ,\gothh_\C ) \textrm{ \ et \ }
{\gotha^*}^+_0(H_{\alpha}) \,\subseteq\, \R ^+\moins0}
\! \gothg_\C ^{\alpha}
\; \bigr) \,\cap\, \gothg$,
$\; N = \exp \gothn \;$ et $\; N_{\nu_+} = N \cap M'_0$.
}

\noindent
Le groupe $N$ (resp. $N_{\nu_+}$) est radical unipotent d'un
sous-groupe parabolique de $G$ (resp. $M'_0$) de composante d\'eploy\'ee
$A$.  De plus, le produit de $G$ se restreint en un diff\'eomorphisme de
$\, N_{\nu_+} \times U \,$ sur $N$.

On note $\calv$ et $\calw$ << les espaces >> de $\tau_{M'}$ et
$\, T_{\mutilde}^{M_0} \!$.
L'espace de $\, T_{\lambdatilde}^{M'_0}$ se r\'ealise sous la forme
$\, \calh
\egdef
\Ind_{M_0 . A . N_{\nu_+}}^{M'_0}
( \calw \otimes \me^{\mi \, \nu \circ \Log} \otimes
\fonctioncar{N_{\nu_+}} )$.
On identifie $\calv \otimes_\C  \calh$ \`a un ensemble de fonctions
sur $M'_0$ \`a valeurs dans $\calv \otimes_\C  \calw \!$.
L'application $\, \Phi \mapsto \restriction{\Phi}{M'_0}$ de l'espace de
$\, \Ind_{\centragr{M}{\gothm} . M_0 . A . N_{\nu_+}}^
{\centragr{M}{\gothm}\,M'_0}
(\delta_{\lambda_{+,\gothm'}}^{\gothm'\!/\gothh} \! \tau_{M'}
\otimes T_{\mutilde}^{M_0}
\otimes \me^{\mi \, \nu \circ \Log} \otimes \fonctioncar{N_{\nu_+}} ) \,$
dans $\calv \otimes_\C  \calh$ muni de $1 \otimes T_{\lambdatilde}^{M'_0}\!$ est
un isomorphisme unitaire de $M'_0$-modules.
L'action d'un $\, x \in \centragr{M}{\gothm}$ sur
$\calv \otimes_\C  \calh$ qui s'en d\'eduit s'\'ecrit
$\tau_{M'}(\hat{x}) \otimes S_0(\hat{x})$, o\`u
$\hat{x} \in M'(\lambdatilde)^{\gothm'\!/\gothh}$ est au-dessus de $x$ et
$\; S_0(\hat{x})\cdot\varphi
= \delta_{\lambda_{+,\gothm'}}^{\gothm'\!/\gothh}\!(\hat{x})
\,{\scriptstyle \times}\,
\varphi (x^{-1} {\scriptscriptstyle\bullet\;} x) \;$
pour $\, \varphi \in \calh$.
On consid\`ere un sous-groupe ferm\'e $\Gamma$ de $\centragr{M}{\gothm}$.
On pose $\, M_\Gamma = \Gamma M_{\nu_+}$ et $\, M'_\Gamma = \Gamma M'_0$.
On fixe une composante irr\'eductible $\tau_{M'\!,\Gamma}$ (dans
$\Xirr_{M'_\Gamma}(\lambdatilde)$) de la restriction de
$\tau_{M'}$ \`a $M'_\Gamma(\lambdatilde)^{\gothm'\!/\gothh}$.
On obtient comme ci-dessus un isomorphisme unitaire de $M'_0$-modules avec
transport de l'action de $\Gamma$, en rempla\c{c}ant $\centragr{M}{\gothm}$ par
$\Gamma$, $\tau_{M'}$ par $\tau_{M'\!,\Gamma}$ et $\calv$ par l'espace
$\calv_{\Gamma}$ de $\tau_{M'\!,\Gamma}$.
On d\'etermine un \'el\'ement $\tau_{M,\Gamma}$ de
$\Xirr_{M_\Gamma}(\mutilde)$ par l'\'egalit\'e
$\; \delta_{\mu_{+,\gothm}}^{\gothm / \gotht} \tau_{M,\Gamma}
= \restriction{(\delta_{\lambda_{+,\gothm'}}^{\gothm'\!/\gothh} \!
\tau_{M'\!,\Gamma})}{\Gamma T_0}$.

On va voir que $S_0 = S$.
L'\'egalit\'e annonc\'ee deviendra une cons\'equence du th\'eor\`eme d'induction
par \'etages et du fait que
$\; \Ind_{H_1}^{G_1} (\pi \circ \restriction{p}{H_1})
\simeq (\Ind_{H_2}^{G_2} \pi) \circ p \;$
quand $p: G_1 \to G_2$ est un morphisme de groupes de Lie r\'eels surjectif,
$H_1$ est l'image r\'eciproque par $p$ d'un sous-groupe ferm\'e
$H_2$ de $G_2$, et $\pi$ est une repr\'esentation unitaire continue de $H_2$.

On reprend les arguments et notations de la d\'emonstration de
\ref{construction de rep, cas de $M'$} (a).
Les groupes
$\, K_{M'_\Gamma} \!\egdef K \cap M'_\Gamma \,$
et
$\, K_{M_\Gamma} \!\egdef K \cap M_\Gamma$
sont des sous-groupes compacts maximaux de $M'_\Gamma$ et $M_\Gamma$.
On constate que le $(\gothm_\C ,K_{M_\Gamma})$-module associ\'e \`a la
repr\'esentation
$\, \Ind_{\Gamma . M_0}^{M_\Gamma} (\delta_{\mu_{+,\gothm}}^{\gothm / \gotht}
\tau_{M,\Gamma} \otimes T_{\mutilde}^{M_0}) \,$
du groupe $M_\Gamma$ est isomorphe \`a
$\, \calr_{M_\Gamma}^q ((\rho_{\mu_{+,\gothm}}^{\gothm / \gotht}\!)^{-1}
\tau_{M,\Gamma})$,
o\`u $\,\calr_{M_\Gamma}^q$ est le foncteur d'induction cohomologique relatif
\`a $\gothb_M$.
Ensuite, \`a l'aide du th\'eor\`eme d'induction par \'etages on voit que le
$(\gothm'_\C ,K_{M'_\Gamma})$-module associ\'e \`a
$\, \Ind_{\Gamma . M_0 . A . N_{\nu_+}}^{M'_\Gamma}
( \delta_{\lambda_{+,\gothm'}}^{\gothm'\!/\gothh} \! \tau_{M'}
\otimes T_{\mutilde}^{M_0}
\otimes \me^{\mi \, \nu \circ \Log} \otimes \fonctioncar{N_{\nu_+}} ) \,$
est isomorphe \`a
$\, \calr_{M'_\Gamma}^{q'}
( (\rho_\lambdatilde^{\gothm'\!/\gothh})^{-1} \tau_{M'\!,\Gamma})$,
o\`u $\,\calr_{M'_\Gamma}^{q'}$ est le foncteur d'induction cohomologique
relatif \`a $\gothb_{M'}$.
Enfin, on trouve que

\centerline{
$\dim \, \mathrm{Hom}_{\gothh_\C ,\Gamma T_0} \,
\bigl( \calv_{\Gamma} \otimes_\C  H_{q'} ( \gothn_{M'},\calh^f),
\rho_\lambdatilde^{\gothm'\!/\gothh} \tau_{M'\!,\Gamma} \bigr)
= 1$,
}

\noindent
o\`u chaque $y \in \Gamma \,$ projection d'un
$\, \hat{y} \in M'(\lambdatilde)^{\gothm'\!/\gothh}$
op\`ere sur
$\calv_{\Gamma} \otimes_\C  H_{q'} ( \gothn_{M'},\calh^f)$
par le produit tensoriel de $\, \tau_{M'\!,\Gamma}(\hat{y}) \,$ et de
l'endomorphisme d\'eduit de $\, S_0(\hat{y})$. 

On se donne $\, \hat{x} \in M'(\lambdatilde)^{\gothm'\!/\gothh}$ au-dessus d'un
$\, x \in \centragr{M}{\gothm}$.
On choisit pour $\Gamma$ le sous-groupe ferm\'e de $\centragr{M}{\gothm}$
engendr\'e par $x$.
Donc $\dim \calv_{\Gamma} = 1$ et le calcul qui pr\'ec\`ede prouve que l'action
de $\hat{x}$ sur
$\, (H_{q'} ( \, \gothn_{M'},\calh^{\infty} )^*)
_{-(\mi \,\lambda + \rho_{\gothm'\!,\gothh})}$
issue de $S_0(\hat{x})$ est bien l'homoth\'etie de rapport
$\rho_\lambdatilde^{\gothm'\!/\gothh} (\hat{x})^{-1} \!$.
\cqfd

}

\definition{

On garde les notations de la proposition \ref{construction de rep, cas de $M'$}
et du lemme \ref{construction equivalente, cas connexe}.
On note
$\;\, T_{\lambdatilde,{\gotha^*}^+\!,\tau_+}^G\!
= \Ind_{M' .\, U}^G
( ( \tau_{M'} \otimes S \, T_{\lambdatilde}^{M'_0} ) \otimes \fonctioncar{U} )$%
\labind{Tzzl tilde,a*,zzt}.

Pour simplifier la notation
$T_{\lambdatilde,{\gotha^*}^+\!,\tau_+}^G\!$%
\labind{Tzzl zr,M'0}\labind{Tzzl zr,sigmar,M'},
j'enl\`everai ${\gotha^*}^+\!$ quand $\, {\gotha^*}^+ = \gotha^*$
(compatible \`a \ref{notation pour les representations})
et $\tau_+$ quand $\Xirrp_G(\lambdatilde,{\gotha^*}^+)$ est un singleton, et je
remplacerai $\lambdatilde$ par $\lambda$ quand
$\, \lambdatilde = (\lambda,\gothh_{(\R )}^{~~*})$
(conventions appliqu\'ees pour
$\,T_{\mutilde}^{M_0}\!$, $T_{\mutilde,\sigma}^M$,
$\,T_{\mutilde}^{M_{\nu_+}}\!$
et $\,T_{\lambdatilde}^{M'_0}$).

}

\proposition{
\label{caracteres canoniques}

{\bf (a)}
Pour chaque $\,l\in\gothg_\C ^*$ semi-simple, on note
$\chiinfty{l}{\gothg}$\labind{zzw lUgC} le caract\`ere $\chi_l$ de
$\centrealg{U \gothg_\C }$ canoniquement associ\'e \`a l'orbite de $\mi \lambda$
sous l'action de $\interieur \gothg_\C $ (cf. \cite [(4.114) p.$\!$~297] {KV95}).
Le centralisateur $(U \gothg_\C )^G$ de $G$ dans $U \gothg_\C $ op\`ere dans
l'espace des vecteurs $C^{\infty}\!$ de
$T_{\lambdatilde,{\gotha^*}^+\!,\tau_+}^G\!$ par restriction du caract\`ere
$\chiinfty{\mi \lambda}{\gothg}$ de $\centrealg{U \gothg_\C }$.

\smallskip
{\bf (b)}
Le morphisme de groupes $z \mapsto (z,1)$ de $\centregr{G}$ dans
$G(\lambda_+)^{\gothg/\gothh}\!$ est injectif.
Le caract\`ere central de $T_{\lambdatilde,{\gotha^*}^+\!,\tau_+}^G\!$ est
<< restriction >> de celui de $\tau_+$.

}

\dem{D\'emonstration de la proposition}{

{\bf (a)}
Le corollaire 5.25 (b) de \cite [p.$\!$~344] {KV95} permet de calculer le
caract\`ere infinit\'esimal de $\tau_{M'} \otimes S \, T_{\lambdatilde}^{M'_0}$
(\'egal \`a celui de $T_{\lambdatilde}^{M'_0}$).
Pour passer ensuite \`a $G$, la d\'emonstration de
\cite [prop. 11.43 p.$\!$~665] {KV95} s'adapte imm\'ediatement.

\smallskip
{\bf (b)}
Soit $z\in\centregr{G}$.
L'\'egalit\'e $S(z,1)=\id$ permet d'obtenir le r\'esultat.
\cqfd

}

\section{\boldmath
L'injection $\,G \cdot (\lambdatilde,\tau) \mapsto T_{\lambdatilde,\tau}^G$
de $G\,\backslash\,\XInd_G$ dans $\widehat{G}$
}\label{Rep9}

Le th\'eor\`eme \ref{injection dans le dual} va g\'en\'eraliser l'essentiel de
celui \'ecrit par M.~Duflo dans \cite [lem. 8 p.$\!$~173] {Df82a}.
Ma d\'emonstration est similaire \`a la sienne, mis \`a part un r\'esultat
inattendu : le lemme \ref{resultat inattendu} (b).
Je rappelle qu'une classe d'\'equivalence de repr\'esentation unitaire
irr\'eductible tra\c{c}able d'un groupe de Lie r\'eel est d\'etermin\'ee par
son caract\`ere (cf. \cite [p.$\!$~64] {Co88}).

\smallskip
Soit $\tau$ la repr\'esentation induite \`a partir de $\tau_+$ de la
d\'efinition \ref{parametres adaptes} (c).

\lemme{
\label{nouveaux parametres}

Les param\`etres  $\,\lambdatilde\,$ et $\,\tau\,$ d\'eterminent
$\, T_{\lambdatilde,{\gotha^*}^+\!,\tau_+}^G\!$.

}

\dem{D\'emonstration du lemme}{

Ce r\'esultat d\'ecoulera du th\'eor\`eme \ref{formule du caractere}, en y
enlevant dans un premier temps toute r\'ef\'erence \`a
$T_{\lambdatilde,\tau}^G\!$ (pour la coh\'erence du raisonnement).
\cqfd

}

\definition{
\label{notation pour les representations}

{\bf (a)}
On note
$\; T_{\lambdatilde,\tau}^G\! = T_{\lambdatilde,{\gotha^*}^+\!,\tau_+}^G\!$%
\labind{Tzzl tilde,zzt}.
Dans la notation $T_{\lambdatilde,\tau}^G\!$, j'enl\`everai $\tau$ quand
$\XInd_G(\lambdatilde)$ est un singleton, et je remplacerai $\lambdatilde$ par
$\lambda$ quand $\, \lambdatilde = (\lambda,\gothh_{(\R )}^{~~*})$.

\smallskip
{\bf (b)}
Soit $\,a$ un automorphisme du groupe de Lie $G$.
Il induit des isomorphismes canoniques encore not\'es $\,a$,
d'espaces vectoriels de $\gothg / \gothg(\lambda)(\mi \rho_\F+)$ sur
$\gothg / \gothg(a\lambda)(\mi \rho_{a\F+})$,
d'alg\`ebres de Lie de
$\, \goth{sl} (\gothg / \gothg(\lambda)(\mi \rho_\F+))$ sur
$\, \goth{sl} (\gothg / \gothg(a\lambda)(\mi \rho_{a\F+}))$ par conjugaison,
de groupes de Lie de
$\, DL (\gothg / \gothg(\lambda)(\mi \rho_\F+))$ sur
$\, DL (\gothg / \gothg(a\lambda)(\mi \rho_{a\F+}))$ par int\'egration,
et enfin de groupes de Lie de
$\,G(\lambdatilde)^{\gothg / \gothg(\lambda)(\mi \rho_\F+)}$ sur
$\,G(a\lambdatilde)^{\gothg / \gothg(a\lambda)(\mi \rho_{a\F+})}$.

}

\lemme{
\label{resultat inattendu}

Soit $\,\tau_0 \in (G_0(\lambda_+)^{\gothg/\gothh})\widehat{~~}$ qui intervient
dans $\restriction{\tau_+}{G_0(\lambda_+)^{\gothg/\gothh}}\!$.

\smallskip
{\bf (a)}
On a :
$\; \tau_0 \in \Xfin_{G_0}(\lambdatilde,{\gotha^*}^+) \,$
et
$\, G_0(\lambda_+) = \centragr{G_0}{\gothh}$.

{\bf (b)}
On pose :
$\; \dot{\tau_0}
= \frac{\rho_{\lambdatilde,{\gotha^*}^+}^{\gothg/\gothh}}
{|\rho_{\lambdatilde,{\gotha^*}^+}^{\gothg/\gothh}|}\,
\tau_0 \,$
(cf. \ref{parametres adaptes} (a)).
L'action du groupe $G_0(\lambdatilde)$ sur la classe de repr\'esentation
$\dot{\tau_0}$ est triviale.

}

\dem{D\'emonstration du lemme}{

{\bf (a)}
La d\'efinition \ref{parametres adaptes} (b) assure que $\tau_0$ est final.
L'\'egalit\'e $\, G_0(\lambda_+) = \centragr{G_0}{\gothh} \,$ d\'ecoule du lemme
\ref{lemme clef}.

\smallskip
{\bf (b)}
On prouve que $\,\tr\dot{\tau_0}(x\gamma x^{-1}) = \tr\dot{\tau_0}(\gamma)\,$
pour $x \in G_0(\lambdatilde)$ et $\gamma \in G_0(\lambda_+)$.%

\smallskip
Soit $G_1$ le sous-groupe de $G_0$ form\'e des \'el\'ements qui fixent $\mu$ et
$\mi \rho_\F+$.
On note $\gothg_1$ son alg\`ebre de Lie.
Comme les formes lin\'eaires $\mu$ et $\mi \rho_\F+$ sont infinit\'esimalement
elliptiques, $\gothg_1$ est r\'eductive, $\gothh \in \Car \gothg_1$, et $G_1$
est connexe.
Le syst\`eme de racines $\, R({\gothg_1}_\C ,\gothh_\C ) \,$ n'a pas de racine
imaginaire.
On note $\, W(\gothg_1(\nu),\gotha) \,$ le groupe de Weyl du syst\`eme de
racines form\'e des racines restreintes de $(\gothg_1(\nu),\gotha)$.
D'apr\`es
\cite [th. 5.17 p.$\!$~125, (5.5) p.$\!$~126 et lemma 5.16 p.$\!$~124] {Kn86}
appliqu\'e \`a $G_1$ d'une part, et \cite [prop. 4.12 p.$\!$~81] {Kn86} d'autre
part, l'application canonique de $\,G_0(\lambdatilde)/G_0(\lambda_+)\,$ dans
$W(\gothg_1(\nu),\gotha)$ est bijective.
On se donne $\beta \in R({\gothg_1(\nu)}_\C ,\gothh_\C )$,
$X \in \gothg_1(\nu)^{\restriction{\beta}{\gotha}}$ et
$Y \in \gothg_1(\nu)^{-\restriction{\beta}{\gotha}}$ tels que $[X,Y]$ est la
racine duale de $\restriction{\beta}{\gotha}$ dans $\gotha$.
On pose  $\, g = \exp(\frac{\pi}{2}(X-Y))$.
Donc $g \in G_0(\lambdatilde)$ se projette dans $W(\gothg_1(\nu),\gotha)$ sur la
r\'eflexion $s_{(\restriction{\beta}{\gotha})}$ au vu de la d\'emonstration de
\cite [prop. 5.15 (c) p.$\!$~123] {Kn86}.
On se contente maintenant de montrer l'\'egalit\'e
$\tr\dot{\tau_0}(g\gamma g^{-1}) = \tr\dot{\tau_0}(\gamma)$
pour $\gamma \in G_0(\lambda_+)$.%

\smallskip
Un calcul dans $SL(2,\C )$ montre que pour toute
$\alpha \in R(\gothg_\C ,\gothh_\C )$ r\'eelle, on a
$\, \Ad \gamma_\alpha = \exp(\mi\pi \ad H_\alpha)$ et
$\Ad \gamma_\alpha^2 = \id \,$ dans $GL(\gothg )\subseteq GL(\gothg_\C )$.
Donc $\gamma_\alpha^2 \in \centregr{G_0}$.
D'apr\`es \cite [lemma 12.30 (c) p.$\!$~469] {Kn86} et (a), tout \'el\'ement
$\gamma$ de $G_0(\lambda_+)$ s'\'ecrit
$\;\gamma = h\,z\,\gamma_{\alpha_1}\!\dots\gamma_{\alpha_k}\;$ avec
$h \in \exp\gothh$, $z \in \centregr{G_0}$ et
$\alpha_1,\dots,\alpha_k \in R(\gothg_\C ,\gothh_\C )$ r\'eelles.
Compte tenu de la proposition \ref{morphisme rho et fonction delta} (b), la
diff\'erentielle en $1$ de
$\frac{\rho_{\lambdatilde,{\gotha^*}^+}^{\gothg/\gothh}}
{|\rho_{\lambdatilde,{\gotha^*}^+}^{\gothg/\gothh}|}$
est la demi-somme des $\alpha \in R(\gothg(\nu)_\C ,\gothh_\C )$ tels que :
$\mi \mu(H_{\alpha}) > 0$ ou,
$\mi \mu(H_{\alpha}) = 0$ et $\rho_\F+\!(H_{\alpha}) > 0$.
Donc $G_0(\lambdatilde)$ laisse invariants les caract\`eres unitaires par
lesquels $\exp\gothh$ et $\centregr{G_0}$ agissent dans l'espace de
$\dot{\tau_0}$.
Dans la suite de cette d\'emonstration, on prendra pour cette raison $\gamma$ de
la forme $\,\gamma = \gamma_{\alpha_1}\cdots\gamma_{\alpha_k}\,$
avec $\alpha_1,\dots,\alpha_k \in R(\gothg_\C ,\gothh_\C )$ r\'eelles.
\smallskip

Pour toute $\alpha \in R(\gothg_\C ,\gothh_\C )$ r\'eelle,
$\, \gamma_\alpha g \gamma_\alpha^{-1} \,$
(resp. $\, \gamma_\alpha g^{-1} \gamma_\alpha^{-1}$)
est \'egal \`a $g$ (resp. $g^{-1}$)
si $\beta(H_\alpha)$ est pair et \`a $g^{-1}$ (resp. $g$) sinon.
Par cons\'equent, on a : \\
$\; g (\gamma g^{-1} \gamma^{-1})
= \cases{\; 1 & si $\,\beta(H_{\alpha_1}+\cdots+H_{\alpha_k})\,$ est pair
\cr\; g^2 & sinon.}$%
\smallskip

On se place tout d'abord dans le cas o\`u $\beta$ est r\'eelle et
$\; g \gamma g^{-1} \gamma^{-1} \not= 1$.
On rend les choix de $g$ et de $\gamma_\beta$ compatibles en imposant
$\,X=X_\beta\,$ et $\,Y=X_{-\beta}\,$ (cf. \ref{parametres adaptes} (a)).
On note $\zeta$ le rapport de l'homoth\'etie $\dot{\tau_0}(\gamma_\beta^2)$.
On a $\,\gamma_\beta = g^2\,$ puis
$\,\gamma_\beta \gamma \gamma_\beta^{-1}
= \gamma_\beta^2 \gamma\,$
et
$\,\gamma_\beta (g \gamma g^{-1}) \gamma_\beta^{-1}
= \gamma_\beta^2 (g \gamma g^{-1})$.
Donc $\; \tr\dot{\tau_0}(\gamma) = 0 = \tr\dot{\tau_0}(g\gamma g^{-1}) \;$ quand
$\zeta \not= 1$.
On va raisonner autrement quand $\zeta = 1$.
L'\'el\'ement $\gamma_\beta+\gamma_\beta^{-1}$ de l'alg\`ebre du groupe
$G_0(\lambda_+)$ est central.
Donc $\dot{\tau_0}(\gamma_\beta)$ est \'egal \`a l'homoth\'etie
$\; (\dot{\tau_0}(\gamma_\beta)+\dot{\tau_0}(\gamma_\beta^{-1}))
(\dot{\tau_0}(\gamma_\beta^2)^{-1}+\id)^{-1} \;$
lorsque  $\zeta\not= -1$.
En particulier, comme $\tau_0$ est final, on a
$\; \dot{\tau_0}(g \gamma g^{-1}) \dot{\tau_0}(\gamma)^{-1}
= \dot{\tau_0}(\gamma_\beta) \;$
avec
$\; (\delta_{\lambda_+}^{\gothg/\gothh} \tau_0) (\gamma_\beta)
= (-1)^{n_\beta+1} \id \;$
quand $\zeta = 1$.
On note $F$ le sous-groupe de $G$ engendr\'e par les $\gamma_\alpha$ avec
$\alpha \in R(\gothg_\C ,\gothh_\C )$ r\'eelle.
D'apr\`es \cite [(42) p.$\!$~335] {DV88}, \`a corriger, il existe un caract\`ere
unitaire de $F$ qui envoie $\gamma_\alpha$ sur $(-1)^{n_\alpha+1}$ pour toute
$\alpha \in R(\gothg_\C ,\gothh_\C )$ r\'eelle.
Par ailleurs, la restriction de
$\,\frac{\rho_{\lambdatilde,{\gotha^*}^+}^{\gothg/\gothh}}
{|\rho_{\lambdatilde,{\gotha^*}^+}^{\gothg/\gothh}|}
\big(\delta_{\lambda_+}^{\gothg/\gothh}\big)^{\!-1}$
\`a $G_0(\lambda_+)$ est invariante sous $G(\lambdatilde)$ d'apr\`es
\ref{espace symplectique canonique} (a).
Ainsi, il existe un caract\`ere unitaire $G(\lambdatilde)$-invariant $\chi$ de
$F$, tel que quand $\zeta = 1$ on ait :
$\; \dot{\tau_0}(\gamma_\beta)
= \chi(\gamma_\beta) \id \;$
puis
$\; \dot{\tau_0}(\gamma_\beta)
= \chi(g\gamma g^{-1})\,\chi(\gamma)^{-1} \id
= \id$.
\smallskip

On suppose pour finir que $\beta$ n'est pas r\'eelle et que
$g \gamma g^{-1} \gamma^{-1} \not= 1$.
On note $R_\beta$ le syst\`eme de racines
$(\R \beta \oplus \R \overline{\beta}) \cap R(\gothg_1(\nu)_\C ,\gothh_\C )$
sans racine imaginaire et avec des racines non r\'eelles.
En examinant le diagramme de Satake de la sous-alg\`ebre de Lie semi-simple de
$\gothg_1(\nu)$ de complexifi\'ee
$\,(\C H_\beta \oplus \C H_{\overline{\beta}}) \oplus
\Sumpetit_{\alpha \in R_\beta} \gothg_\C ^{\alpha}$,
on constate qu'il existe un morphisme d'alg\`ebres de Lie injectif
$\; \eta : \gothg_0 \to \gothg_1(\nu)$, avec
$\,\gothg_0 = \restriction{\goth{sl}(2,\C )}{\R }\,$
ou
$\,\gothg_0 = \goth{su}(1,2)$,
qui envoie pour tout $\,(x,y) \in \R ^2\,$ et selon la valeur de $\gothg_0$,
la matrice
$\left(\begin{smallmatrix}
x+\mi y & 0 \\ 0 & -(x+\mi y)
\end{smallmatrix}\right)$
ou la matrice
$\left(\begin{smallmatrix}
\mi y & x & 0 \\ x & \mi y & 0 \\ 0 & 0 & -2\mi y
\end{smallmatrix}\right)$
sur
$\,x\,(H_\beta+H_{\overline{\beta}})\,+\,\mi y\,(H_\beta-H_{\overline{\beta}})$.
On choisit $X$ et $Y$ dans
$(\gothg_\C ^\beta \oplus \gothg_\C ^{\overline{\beta}}) \cap \gothg$
et
$(\gothg_\C ^{-\beta} \oplus \gothg_\C ^{-\overline{\beta}}) \cap \gothg$
\'egaux, selon $\gothg_0$, \`a :
$\; X
= \eta\big(\!\!\left(\begin{smallmatrix}
0 & 1 \\ 0 & 0
\end{smallmatrix}\right)\!\!\big) \;$
et
$\; Y
= \eta\big(\!\!\left(\begin{smallmatrix}
0 & 0 \\ 1 & 0
\end{smallmatrix}\right)\!\!\big)$,
ou,
$\; X
= \eta\big(\!\!\left(\begin{smallmatrix}
0 & 0 & 1 \\ 0 & 0 & 1 \\ 1 & -1 & 0
\end{smallmatrix}\right)\!\!\big) \;$
et
$\; Y
= \eta\big(\!\!\left(\begin{smallmatrix}
0 & 0 & 1 \\ 0 & 0 & -1 \\ 1 & 1 & 0
\end{smallmatrix}\right)\!\!\big)$.
Un calcul dans l'un des groupes simplements connexes $SU(2)$ ou
$\{1\} \times SU(2)$ montre ensuite que $g^2$ appartient \`a
$\,\exp(\gothh\cap\derivealg{\gothg_1(\nu)})$.
La d\'efinition de $\gothg_1$ assure que $W(\gothg_1(\nu)_\C ,\gothh_\C )$ fixe
$\diff_1\dot{\tau_0}$.
Donc $\dot{\tau_0}$ est trivial sur
$\,\exp(\gothh\cap\derivealg{\gothg_1(\nu)})$.
Ainsi, on obtient :
$\, \dot{\tau_0}(g\gamma g^{-1}\gamma^{-1}) = \dot{\tau_0}(g^2) = \id$.
\cqfd

}

\medskip
Dans la proposition qui suit, j'\'enonce les r\'esultats de la th\'eorie du
petit groupe de Mackey qui me seront utiles dans la d\'emonstration du prochain
th\'eor\`eme.
Je ne rappellerai pas les notions de  << groupe de type $I$ >> et << int\'egrale
hilbertienne >> (voir \cite{Di69}).
On se donne pour cette proposition un groupe localement compact \`a base
d\'enombrable $A$, et un sous-groupe ferm\'e distingu\'e $B$ de $A$ qui est de
type $I$.
Le groupe $A$ agit canoniquement sur $\widehat{B}$.

Voici quelques pr\'ecisions concernant le vocabulaire.
On appellera cocycle mesurable unitaire de $A$ une application mesurable $c$ de
$A \times A$ dans l'ensemble des nombres complexes de module $1$ qui v\'erifie
$\,c(1,1) = 1\,$
et
$\,c(xy,z)\,c(x,y) = c(x,yz)\,c(y,z)\,$
pour $x,y,z \in A$.
Dans ce cas, \'etant donn\'e un espace de Hilbert complexe s\'eparable $V$, on
appellera $c$-repr\'esentation projective de $A$ dans $V$ une application
$\wt{\pi}$ de $A$ dans le groupe unitaire de $V$, pour laquelle chacune des
fonctions $\, x \in A \mapsto \langle \wt{\pi}(x) \cdot v,w \rangle \in \C \,$
avec $v,w \in V$ est mesurable, et telle que
$\,\wt{\pi}(1) = \id\,$
et
$\,\wt{\pi}(xy) = c(x,y)\,{\scriptstyle \times}\,\wt{\pi}(x)\,\wt{\pi}(y)\,$
pour $x,y \in A$.

\proposition{\hspace{-15pt}{\bf{(Mackey)}}\quad 
\label{theorie de Mackey}
%
{\bf (a)}
Soit $\pi \in \widehat{B}$.
On note $A_\pi$ le stabilisateur de $\pi$ dans $A$ et $\pr$ la projection
canonique de $A_\pi$ sur $A_\pi/B$.
Le groupe $A_\pi$ est ferm\'e dans $A$.
Il existe un cocycle mesurable unitaire $c$ de $A_\pi/B$ et une classe
d'isomorphisme $\wt{\pi}$ de
$c{\scriptstyle \,\circ\,}(\pr{\scriptstyle \times}\pr)$-repr\'esentation
projective de $A_\pi$ qui prolonge $\pi$ dans le m\^eme
espace de Hilbert.
Par ailleurs, il existe une mesure $\sigma$-finie non nulle $m$ sur
$\widehat{B}$, unique \`a \'equivalence pr\`es, telle que l'orbite de $\pi$ sous
$A$ dans $\widehat{B}$ a un compl\'ementaire $m$-n\'egligeable et $a \cdot m$
est \'equivalente \`a $m$ pour tout $a \in A$.%

\smallskip
{\bf (b)}
On conserve les notations du (a).
L'application
$\; \wt{\eta} \mapsto
\Ind_{A_\pi}^A (\wt{\eta} \circ \pr \otimes \wt{\pi}) \;$
de l'ensemble des classes d'isomorphisme des $c^{-1}$-repr\'esentations
projectives irr\'eductibles de $A_\pi/B$, dans l'ensemble des classes
d'isomorphisme des repr\'esentations unitaires de $A$, est injective.
Son image est form\'ee des \'el\'ements de $\widehat{A}$ dont la restriction \`a
$B$ est multiple de $\, \int_{\widehat{B\;}} \rho\,\diff m (\rho)$.

\smallskip
{\bf (c)}
Quand l'espace mesurable quotient $A \backslash \widehat{B}$ est
d\'enombrablement s\'epar\'e
(c'est-\`a-dire qu'on peut trouver une suite $(E_n)_{n \in \N }$ de parties
mesurables de $E = A \backslash \widehat{B}$ telle que pour $x \not= y$ dans $E$,
il existe $m \in \N $ v\'erifiant $x\in E_m$ et $y\notin E_m$),
tout \'el\'ement de $\widehat{A}$ s'obtient comme au (b) pour une unique orbite
de $A$ dans $\widehat{B}$.%

}

\dem{D\'emonstration de la proposition}{

{\bf (a)}
L'hypoth\`ese << metrically smooth and of type $I$ >> de Mackey va se traduire
par << de type $I$ >> d'apr\`es
\cite [prop. 4.6.1 p.$\!$~95 et th. 9.1 p.$\!$~168] {Di64}.
Le groupe $A_\pi$ est ferm\'e dans $A$ et $m$ existe d'apr\`es
\cite [th. 7.5 p.$\!$~295] {Ma58}.
L'existence de $c$ et de $\wt{\pi}$ est donn\'ee par
\cite [th. 8.2 p.$\!$~298] {Ma58}.

\smallskip
{\bf (b)}
Ce r\'esultat est cit\'e dans \cite [th. 8.1 p.$\!$~297 et 8.3 p.$\!$~300]
{Ma58}.
Il se r\'ef\`ere \`a la classe de repr\'esentation de $B$ associ\'ee \`a $m$,
construite \`a la fin de sa section 7 p.$\!$~296, qui renvoie aux lignes 4 \`a
6 p.$\!$~273.

\smallskip
{\bf (c)}
R\'esulte de \cite [th. 9.1 p.$\!$~302] {Ma58} et de la fin de sa section 7
p.$\!$~296.
\cqfd

}

\medskip
Soit $K_0$ un sous-groupe compact maximal de $G_0$.
On dira ci-dessous qu'une repr\'esentation unitaire irr\'eductible de $G$
(reli\'ee \`a $G_0$ par \ref{theorie de Mackey} (c)) est << temp\'er\'ee >>
quand les coefficients matriciels issus de ses vecteurs $K_0$-finis sont
major\'es comme dans \cite [3.6 p.$\!$~124] {BW80} pour $G_0$.
Cette notion est ind\'ependante du choix de $K_0$.

\theoreme{
\label{injection dans le dual}

{\bf (a)}
On a :
$\; a \cdot T_{\lambdatilde,\tau}^G
=\, T_{a\cdot\lambdatilde,a\cdot\tau}^G \;$
pour tout un automorphisme $\,a$ du groupe de Lie $G$.

\smallskip
{\bf (b)}
L'ensemble  $G\,\backslash\,\XInd_G$
(resp. $G\,\backslash\,\Xfin_G\!$)
s'injecte dans $\widehat{G}$ par l'application qui envoie l'orbite de
$(\lambdatilde',\tau')$
(resp. $(\lambdatilde',{{\gotha'}^*}^+\!,\tau_+')$)
sur $T_{\lambdatilde',\tau'}^G\!$
(resp. $T_{\lambdatilde',{{\gotha'}^*}^+\!,\tau_+'}^G\!$).
Son image est form\'ee des classes des repr\'esentations temp\'er\'ees
irr\'eductibles de $G$.
(D'apr\`es le lemme \ref{condition d'integrabilite}, les
$\lambdatilde' \in \gstregtilde$ qui interviennent ici appartiennent \`a
$\gstregGtilde$).

}

\dem{D\'emonstration du th\'eor\`eme}{

On utilise les notations du d\'ebut de la partie \ref{RepIII} et des sections
\ref{Rep7} et \ref{Rep8}

\smallskip
{\bf (a)}
Au vu d'une propri\'et\'e de l'induction cit\'ee dans la d\'emonstration
du lemme \ref{construction equivalente, cas connexe} (b), on est ramen\'e \`a
prouver que 
$a \cdot T_{\mutilde}^{M_0} \!= T_{a \cdot \mutilde}^{a \cdot M_0}$
quand $G$ est connexe.
Cette \'egalit\'e provient des formules pour les caract\`eres de
$a \cdot T_{\mutilde}^{M_0}$ et $T_{a \cdot \mutilde}^{a \cdot M_0}$
obtenues (sans le th. \ref{injection dans le dual}) dans la derni\`ere partie de
la d\'emonstration du th\'eor\`eme \ref{formule du caractere}.

\smallskip
{\bf (b)}
On suppose dans un premier temps que $G$ est connexe.
La classe de repr\'esentation $T_{\lambdatilde,{\gotha^*}^+\!,\tau_+}^G\!$
est temp\'er\'ee d'apr\`es la preuve de l'implication $(3) \!\Rightarrow\! (1)$
dans \cite [proposition 3.7 p.$\!$~124] {BW80}, jointe \`a
\cite[th. 6.8.1 p.$\!$~202, prop. p.$\!$~142, prop. p.$\!$~139] {Wa88}
(voir aussi \cite [corollary 12.27 p.$\!$~461] {Kn86}).
J'utilise le dictionnaire propos\'e dans la d\'emonstration du lemme
\ref{construction de rep, cas connexe} (b).
Tout \'el\'ement du dual temp\'er\'e de $G$ s'\'ecrit sous la forme
$\,T_{\lambdatilde',{{\gotha'}^*}^+\!,\tau_+'}^G\!$ pour un
$(\lambdatilde',{{\gotha'}^*}^+\!,\tau_+') \in \Xfin_G\!$ d'apr\`es
\cite [th. 11.14 (b) p.$\!$~131] {ABV92}, \cite [(11.196) p.$\!$~740] {KV95} et
\cite [lemma 4.9 p.$\!$~130] {BW80}.
On se place dans le cas o\`u $(\lambdatilde',{{\gotha'}^*}^+\!,\tau_+')$
v\'erifie
$\; T_{\lambdatilde',{{\gotha'}^*}^+\!,\tau_+'}^G\!\!
= T_{\lambdatilde,{\gotha^*}^+\!,\tau_+}^G\!$.
On introduit les triplets $(\Lambda^{can},R^+_{\mi \, \R },R^+_\R )$ et
$({\Lambda'}^{can},{R'}^+_{\mi \, \R },{R'}^+_\R )$ attach\'es \`a $\tau_+$ et
$\tau_+'$ comme dans la d\'emonstration de
\ref{construction de rep, cas connexe} (b),
et pose $\gothh' = \gothg(\lambdatilde')$.
Les diff\'erentielles des caract\`eres centraux de
$\Lambda^{can}$ et ${\Lambda'}^{can}$ sont \'egales \`a
$\restriction{\mi\lambda}{\gothh}$ et \`a $\restriction{\mi\lambda'}{\gothh'}$.
D'apr\`es \cite [th. 11.14 (c) p.$\!$~131 et def. 11.6  p.$\!$~124]
{ABV92}, il existe $g \in G$ tel que, en notant $\zeta$ l'application qui
envoie $x \in G(\lambda_+)$ sur le d\'eterminant de la restriction de
$\Ad x^\C $ \`a la somme des $\gothg_\C ^{\alpha}$ avec $\alpha \in R^+_\R $ et
$g \alpha \notin {R'}^+_\R $, on ait :
$\, \lambdatilde' = g \lambdatilde \,$
et
$\, {\Lambda'}^{can}
=\, g \cdot \big(\Lambda^{can}\otimes \frac{\zeta}{\abs{\zeta}}\big)$.
Gr\^ace au lemme \ref{espace symplectique canonique} (a) et \`a un calcul fait
dans la d\'emonstration de \ref{construction de rep, cas connexe} (b), on
obtient :
$\, \frac{\rho_{\lambdatilde,{{\gotha'}^*}^+}^{\gothg/\gothh'}}
{|\rho_{\lambdatilde,{{\gotha'}^*}^+}^{\gothg/\gothh'}|}\,\tau_+'
= g \cdot \big(\frac{\rho_{\lambdatilde,{\gotha^*}^+}^{\gothg/\gothh}}
{|\rho_{\lambdatilde,{\gotha^*}^+}^{\gothg/\gothh}|}\,\tau_+\big)$.
Comme les chambres de $\,(\gothg(\lambda)(\mi \rho_\F+),\gotha)\,$ sont
conjugu\'ees sous $\,\normagr{G(\lambda)(\mi \rho_\F+)_0}{\gotha}$,
il existe $u \in G(\lambdatilde)$ tel que
$\, (\lambdatilde',{{\gotha'}^*}^+) = gu \cdot (\lambdatilde,{\gotha^*}^+)$.
\emph{A fortiori}, on a $\; \tau_+' = gu \cdot \tau_+ \;$ d'apr\`es le lemme
\ref{resultat inattendu} (b).

\smallskip
On ne suppose plus $G$ connexe.
Il s'agit d'adapter la d\'emonstration de M.~Duflo dans \cite [p.$\!$~176 \`a
p.$\!$~179] {Df82a}.
On remarque que les sous-groupes ferm\'es distingu\'es
$B_1 = G_0(\lambda_+)^{\gothg/\gothh}$ et $B_2 = G_0$, respectivement de $A_1 =
G(\lambda_+)^{\gothg/\gothh}$ et $A_2 = G$, sont de type $I$ d'apr\`es
\cite [prop. 2.1 p.$\!$~425] {Di69}.
En particulier, les espaces mesurables $\widehat{B_1}$ et $\widehat{B_2}$ sont
standards d'apr\`es \cite [prop. 4.6.1 p.$\!$~95 et th. 9.1 p.$\!$~168] {Di64}.
Leurs quotients sous les actions \`a gauche des groupes finis respectivement
\'egaux \`a $A_1/B_1$ et  $A_2/B_2$ sont donc d\'enombrablement s\'epar\'es.

\smallskip
On va introduire plus bas un certain
$\tau_0 \in (G_0(\lambda_+)^{\gothg/\gothh})\widehat{~~}$.
On notera $G(\lambda_+)_{\tau_0}^{\gothg/\gothh}$ le stabilisateur de
$\tau_0$ dans $G(\lambda_+)^{\gothg/\gothh}\!$,
$G(\lambda_+)_{\tau_0}$ l'image de
$G(\lambda_+)_{\tau_0}^{\gothg/\gothh}$ dans $G(\lambda_+)$ et
$M'(\lambdatilde)_{\tau_0}^{\gothm'\!/\gothh}$ l'image r\'eciproque de
$G(\lambda_+)_{\tau_0}$ dans $M'(\lambdatilde)^{\gothm'\!/\gothh}\!$.
On d\'esignera par la m\^eme notation $\pr$ les projections canoniques de
$G(\lambda_+)_{\tau_0}^{\gothg/\gothh}$, $G(\lambda_+)_{\tau_0} G_0$
et $M'(\lambdatilde)_{\tau_0}^{\gothm'\!/\gothh}$ sur
$G(\lambda_+)_{\tau_0} / G_0(\lambda_+)$.
D'apr\`es la proposition \ref{theorie de Mackey} (a) et (c), il existe une
sous-repr\'esentation irr\'eductible $\tau_0$ de
$\restriction{\tau_+}{G_0(\lambda_+)^{\gothg/\gothh}}\!$,
un cocycle mesurable unitaire $c$ de $G(\lambda_+)_{\tau_0} / G_0(\lambda_+)$,
une classe d'isomorphisme $\wt{\tau_0}$ de
$c{\scriptstyle \,\circ\,}(\pr{\scriptstyle \times}\pr)$-repr\'esentation
projective de $G(\lambda_+)_{\tau_0}^{\gothg/\gothh}$ qui prolonge $\tau_0$
dans son espace, et une classe d'isomorphisme $\wt{\eta}$ de
$c^{-1}$-repr\'esentation projective irr\'eductible de
$G(\lambda_+)_{\tau_0} / G_0(\lambda_+)$, tels que :

\centerline{
$\tau_+
= \Ind_{G(\lambda_+)_{\tau_0}^{\gothg/\gothh}}
^{G(\lambda_+)^{\gothg/\gothh}}
(\wt{\eta} \circ \pr \otimes \wt{\tau_0})$.
}

\smallskip
On construit maintenant trois isomorphismes de repr\'esentations unitaires.%
\\
Le premier, relatif au groupe
$M'(\lambdatilde)^{\gothm'\!/\gothh}\!$, va de l'espace de
$\, \Ind_{G(\lambda_+)_{\tau_0}^{\gothg/\gothh}}
^{G(\lambda_+)^{\gothg/\gothh}}
(\wt{\eta} \circ \pr \otimes \wt{\tau_0}) \,$
muni de $\tau_{M'}$ sur celui de
$\, \Ind_{M'(\lambdatilde)_{\tau_0}^{\gothm'\!/\gothh}}
^{M'(\lambdatilde)^{\gothm'\!/\gothh}}
(\wt{\eta} \circ \pr \otimes (\wt{\tau_0})_{M'})$,
o\`u $(\wt{\tau_0})_{M'}$ est la
$c{\scriptstyle \,\circ\,}(\pr{\scriptstyle \times}\pr)$-repr\'esentation
projective de $M'(\lambdatilde)_{\tau_0}^{\gothm'\!/\gothh}$ telle que
$\, \delta_{\lambda_{+,\gothm'}}^{\gothm'\!/\gothh} (\wt{\tau_0})_{M'}
= \delta_{\lambda_+}^{\gothg/\gothh} \wt{\tau_0} \,$
sur $G(\lambda_+)_{\tau_0}$.
Il s'\'ecrit $\varphi \!\mapsto\! \varphi_{M'}$ avec
$\, (\delta_{\lambda_{+,\gothm'}}^{\gothm'\!/\gothh})^{-1} \varphi_{M'}
\!=\! (\delta_{\lambda_+}^{\gothg/\gothh})^{-1} \varphi \,$
(cf. la d\'emonstration de
\ref{construction equivalente, cas connexe} (a)).%
\\
Le second, relatif \`a $M'\!$, va de l'espace de
$\, \bigl( \Ind_{M'(\lambdatilde)_{\tau_0}^{\gothm'\!/\gothh}}
^{M'(\lambdatilde)^{\gothm'\!/\gothh}}
(\wt{\eta} \circ \pr \otimes (\wt{\tau_0})_{M'}) \bigr)
\otimes S T_{\lambdatilde}^{M'_0}$
muni par transport de $T_{\lambdatilde,\tau_{M'}}^{M'}\!\!$
sur celui de
$\Ind_{G(\lambda_+)_{\tau_0} M'_0}^{M'}\!
\bigl( (\wt{\eta} \circ \pr \otimes
(\wt{\tau_0})_{M'}) \otimes S T_{\lambdatilde}^{M'_0}
\bigr)$
o\`u la repr\'esentation de $G(\lambda_+)_{\tau_0} M'_0$ \`a laquelle il
est fait allusion est construite comme $\, T_{\lambdatilde,\tau_{M'}}^{M'}\!$
dans
\ref{construction equivalente, cas connexe} (a).
Il envoie un \'el\'ement de la forme $\varphi \otimes v$ sur l'\'el\'ement $\psi$
v\'erifiant

\negsmallskip\centerline{
$\psi(xy)
= \varphi(\hat{x}) \otimes
T_{\lambdatilde}^{M'_0}\!(y)^{-1} S(\hat{x})^{-1} \!\cdot v$}

\noindent
pour $\hat{x} \!\in\! M'(\lambdatilde)^{\gothm'\!/\gothh}$
au-dessus de $x \in G(\lambda_+)$ et $y \in M'_0$.
\\
On note $(M'_{G_0},(\tau_0)_{M'})$ le << couple $(M',\tau_{M'})$ >> obtenu \`a
partir de $(G_0,\lambdatilde,{\gotha^*}^+\!,\tau_0)$ au lieu de
$(G,\lambdatilde,{\gotha^*}^+\!,\tau_+)$.
Il existe une unique
$c{\scriptstyle \,\circ\,}(\pr{\scriptstyle \times}\pr)$-repr\'esentation
projective $\wt{T_0}$ de $G(\lambda_+)_{\tau_0} G_0$ qui prolonge
$T_{\lambdatilde,{\gotha^*}^+\!,\tau_0}^{G_0}\!$ de fa\c{c}on que son action
sur un \'el\'ement $\phi$ de l'espace de
$\, \Ind_{M'_{G_0} . U}^{G_0}
(((\tau_0)_{M'} \otimes S T_{\lambdatilde}^{M'_0}) \otimes \fonctioncar{U}) \,$
v\'erifie

\centerline{
$(\wt{T_0}(xg) \!\cdot\! \phi)(h)
= \abs{\det(\Ad x)_{\gothu}}^{1/2}
\,{\scriptstyle \times}\,
((\wt{\tau_0})_{M'}(\hat{x}) \otimes S(\hat{x})) \,\cdot\,
\phi(g^{-1}x^{-1}hx)$
}

\noindent
pour $\hat{x} \in M'(\lambdatilde)_{\tau_0}^{\gothm'\!/\gothh}$ au-dessus de
$x \in G(\lambda_+)_{\tau_0}$ et $g,h \in G_0$.
\\
Le troisi\`eme isomorphisme, relatif \`a $G(\lambda_+)_{\tau_0} G_0$,
est fourni par l'application de res\-triction \`a $G_0$, de l'espace de
$\Ind_{(G(\lambda_+)_{\tau_0} M'_0) . U}^{G(\lambda_+)_{\tau_0} G_0}
\bigl( ((\wt{\eta} \circ \pr \otimes
(\wt{\tau_0})_{M'}) \otimes S T_{\lambdatilde}^{M'_0})
\otimes \fonctioncar{U} \bigr)$
sur celui de $\wt{\eta} \circ \pr \otimes \wt{T_0}$
(identifi\'e \`a un espace de fonctions sur $G_0$).

\noindent
En r\'ecapitulant, on trouve que :

\centerline{
$\; T_{\lambdatilde,{\gotha^*}^+\!,\tau_+}^G\!
= \Ind_{G(\lambda_+)_{\tau_0} G_0}^G
(\wt{\eta} \circ \pr \otimes \wt{T_0})$.
}

\smallskip
Il reste \`a appliquer la th\'eorie de Mackey.
D'apr\`es l'\'etude du cas connexe faite ci-dessus, on a
$\, T_{\lambdatilde,{\gotha^*}^+\!,\tau_0}^{G_0}\! \in \widehat{G_0} \,$ et le 
stabilisateur de $T_{\lambdatilde,{\gotha^*}^+\!,\tau_0}^{G_0}\!$ dans $G$ est
\'egal \`a $G(\lambda_+)_{\tau_0} G_0$.
La proposition \ref{theorie de Mackey} (b) fournit la bijection
$\,\wt{\sigma} \mapsto
\Ind_{G(\lambda_+)_{\tau_0} G_0}^G
(\wt{\sigma} \circ \pr \otimes \wt{T_0})$
de l'ensemble des classes d'isomorphisme des $c^{-1}$-repr\'esentations
projectives irr\'eductibles de $G(\lambda_+)_{\tau_0} / G_0(\lambda_+)$ sur
l'ensemble des \'el\'ements de $\widehat{G}$ dont la restriction \`a $G_0$ est
\'egal \`a la somme hilbertienne des
$\, g \cdot T_{\lambdatilde,{\gotha^*}^+\!,\tau_0}^{G_0}$
avec $\,\dot{g} \in G/G(\lambda_+)_{\tau_0} G_0$.
Ainsi, au vu du cas connexe, on a : la classe de repr\'esentation
$T_{\lambdatilde,{\gotha^*}^+\!,\tau_+}^G\!$ est temp\'er\'ee irr\'eductible,
toute $\, \pi \in \widehat{G} \,$ temp\'er\'ee est atteinte d'apr\`es la
proposition \ref{theorie de Mackey} (c), et par un calcul facile tout
$(\lambdatilde',{{\gotha'}^*}^+\!,\tau_+') \in \Xfin_G\!$ tel que
$\; T_{\lambdatilde',{{\gotha'}^*}^+\!,\tau_+'}^G\!\!
=\, T_{\lambdatilde,{\gotha^*}^+\!,\tau_+}^G$
est conjugu\'e \`a $(\lambdatilde,{\gotha^*}^+\!,\tau_+)$ sous $G$.
\emph{A fortiori} pour un tel triplet
$(\lambdatilde',{{\gotha'}^*}^+\!,\tau_+')$ on a
$\,(\lambdatilde',\tau') \in G \cdot (\lambdatilde,\tau)$,
o\`u $\tau'$ se d\'eduit de $\tau_+'$ comme dans la d\'efinition
\ref{parametres adaptes} (c).
\cqfd

}

\remarque{
\label{recuperation des parametres}

On suppose dans cette remarque que $\, \lambdatilde \in \gstItilde \,$
et pose $\tau = \tau_+$.
On va retrouver explicitement l'orbite de $(\lambdatilde,\tau)$ sous $G$ \`a
partir du caract\`ere de $\, T_{\lambdatilde,\tau}^G$.

Dans un premier temps, je vais r\'ecup\'erer $G \cdot \lambdatilde$.
Le th\'eor\`eme \ref{formule du caractere} fournit l'\'egalit\'e
$\,\bigl( \tr T_{\lambdatilde,\tau}^G \bigr)_{\!1}
= \dim \tau \,{\scriptstyle \times}\,
\restriction{~{\widehat{\beta}}_{G \cdot \lambdatilde}}{\calv_1}$
o\`u $\calv_1$ est un voisinage ouvert de $0$ dans $\gothg$, qui coupe donc
chaque composante connexe de $\gssreg$.
D'apr\`es la proposition \ref{mesures de Liouville temperees}~(b) et la fin
de la remarque \ref{explication de la notation support} (1), on en d\'eduit que
$T_{\lambdatilde,\tau}^G$ d\'etermine $\beta_{R_G(G \cdot \lambdatilde)}$.
Les mesures de Radon temp\'er\'ees $\beta_{\Omega}$ sur $\gothg^*$ avec
$\Omega \in G \backslash \gstreg$ sont lin\'eairement ind\'ependantes,
car chaque $\beta_{\Omega}$ est non nulle concentr\'ee sur $\Omega$.
Compte tenu de la proposition \ref{bijection entre orbites} (a), la classe de
repr\'esentation $T_{\lambdatilde,\tau}^G$ d\'etermine
l'orbite $G \!\cdot\! \lambdatilde$.

Je vais maintenant reconstruire $\tr \tau$ \`a partir de
$\, T_{\lambdatilde,\tau}^G \,$ et de $\lambdatilde$.
On sait d\'ej\`a que $\; \tau(\exp X) = \me^{\mi \,\lambda(X)}\id \;$ pour
$X \in \gothh$.
Soit $\, \hat{e} \in G(\lambda_+)^{\gothg/\gothh}$ au-dessus d'un
$\, e \in G(\lambda_+)$ elliptique.
On lui associe un $\, e_0 \in e\exp\gotht(e)$ comme dans le lemme
\ref{descente pour les formes lineaires} (b).
Le th\'eor\`eme \ref{formule du caractere} exprime
$\, \bigl( \tr T_{\lambdatilde,\tau}^G \bigr)_{\!e_0}$
comme combinaison lin\'eaire des classes de fonction localement int\'egrable
$\restriction
{~{\widehat{\beta}}_{G(e_0) \cdot \lambdatilde'_{e_0}}}{\calv_{e_0}} \!$
avec
\smash{$\dot{\lambdatilde'_{e_0}} \!\!\in G(e_0) \backslash \gezerostItilde$}
qui sont lin\'eairement ind\'ependantes.
Le coefficient de $\bigl( \tr T_{\lambdatilde,\tau}^G \bigr)_{\!e_0} \!$ suivant
$\restriction{~{\widehat{\beta}}_{G(e_0) \cdot \lambdatilde[e_0]}}{\calv_{e_0}}$
fournit $\tr\tau(\hat{e})$.
\cqfr

}

%
\part{
Caract\`eres des repr\'esentations
}\label{RepIV}
%

\`A nouveau, on consid\`ere un \'el\'ement
$\, \lambdatilde \!=\! (\lambda,\F+)$ de $\gstregtilde$
et pose $\gothh = \gothg(\lambdatilde)$.
On note $\mu$ et $\nu$ (resp. $\gotht$ et $\gotha$) les composantes
infinit\'esimalement elliptique et hyperbolique de $\lambda$ (resp. $\gothh$).
Soit ${\gotha^*}^+$ une chambre de $(\gothg(\lambda)(\mi \rho_\F+),\gotha)$.
J'utiliserai les autres notations de la partie \ref{RepIII}
(dont $M$, $\gothm$, $\tau_M$ pour la section \ref{Rep7}, et, $M'$, $\gothm'$,
$U$, $\tau_{M'}$ pour la section \ref{Rep8}) sans rappeler ce qu'elles
repr\'esentent.

\smallskip
On se donne dans cette partie un \'el\'ement elliptique $e$ de $G$ (\`a ne pas
confondre avec la base $\me $ des logarithmes n\'ep\'eriens).

\section{
Formule de restrictions des caract\`eres
}\label{Rep10}

La formule que je vais maintenant \'ecrire g\'en\'eralise celle d'A.~Bouaziz
dans \cite [th. 5.5.3 p.$\!$~52]{Bo87}
(voir aussi \cite [th. (7) p.$\!$~106]{DHV84}),
tout en s'y ramenant, et celle de W.~Rossmann \cite [p.$\!$~64]{Ro80} relative
au cas $e=1$.

\definition{
\label{notations pour les restrictions}

{\bf (a)}
On note $D_G$\labind{DGaa} la fonction sur $G$ dont la restriction \`a chaque
composante connexe $G^+$ de $G$ envoie $x \in G^+\!$ sur le coefficient de
$T^{r^{\!+}}\!$ dans $\det(T \id + \id - \Ad x)$, o\`u $\,r^+\!$ est le rang
commun aux alg\`ebres de Lie r\'eductives $\gothg(x^+)$ quand $\,x^+\!$ d\'ecrit
l'ensemble des \'el\'ements semi-simples de $G^+$
(cf. \cite [lem. 1.4.1 p.$\!$~6] {Bo87}).
Elle est analytique et invariante sous $\interieur G$.

\smallskip
{\bf (b)}
On note $\Gssreg$\labind{G ssregaa} l'ensemble des $x \in G$ semi-simples tels
que $\gothg(x)$ est commutative.
Donc $\Gssreg$ est l'ouvert dense de compl\'ementaire n\'egligeable de $G$
form\'e des points o\`u $D_G$ ne s'annule pas (cf. \cite [1.3 p.$\!$~5] {Bo87}).

\smallskip
{\bf (c)}
On note
$\;d_e = \frac {1}{2} \dim (1-\Ad e)(\gothg/\gothj)$\labind{de}
et
$D_e = \det (1-\Ad e)_{_{\scriptstyle (1-\Ad e)(\gothg/\gothj)}} > 0$%
\labind{Deaa},

\noindent
ind\'ependamment du choix de $\; \gothj_e \!\in\! \Car \gothg(e) \,$
auquel on associe
$\; \gothj \egdef \centraalg{\gothg}{\gothj_e} \in \Car \gothg$.

\smallskip
{\bf (d)}
On note

\centerline{
$\calv_{\!e}\,
=\, \{\, X\in\gothg(e) \mid \abs{\Im z} < \varepsilon_e \;
\textrm{ pour toute valeur propre $z$ de $\ad_{\gothg} X$} \,\}$%
\labind{Ve}
}

\centerline{
et
$\;\, \displaystyle k_e(X)
= \left( \!
\det { \left(
\frac{\me^{\ad X/2} - \me^{- \ad X/2}}{\ad X}
\right) }_{\!\!\gothg(e)}
\frac{\det (1 - \Ad (e \exp X))_{\gothg/\gothg(e)}}
{\det (1 - \Ad e)_{\gothg/\gothg(e)}}
\right)^{\!\!1/2} \!\! > 0$%
\labind{ke}
}

\noindent
pour $\,X \in \calv_{\!e}$,
o\`u
$\;\, \varepsilon_e
=
\inf \{\, \theta \in \, \left] 0,2\pi \right] \mid
\textrm{$\me^{\mi \,\theta}$ est valeur propre de $\Ad_G e$} \,\}\;
\leq\;2\pi$.

\smallskip
{\bf (e)}
\`A toute fonction g\'en\'eralis\'ee $\interieur G$-invariante $\Theta$ sur $G$,
on associe (cf. \cite [(6) p.$\!$~98] {DHV84}) la fonction g\'en\'eralis\'ee
$\,\Ad G(e)$-invariante $\Theta_e$\labind{zzHt e} sur $\calv_{\!e}$,
d\'etermin\'ee par l'\'egalit\'e
$\; \Theta_e(X)
=\, k_e(X) \; {\scriptstyle \times} \;
\restriction{\Theta}{e \exp \calv_{\!e}} \, ( e \exp X) \;$
de fonctions g\'en\'eralis\'ees en $\,X \in \calv_{\!e}$.

}

\medskip
L'\'egalit\'e
$\;\displaystyle
G = \!
\bigcup_{\substack
{e_0 \in G \\ \textrm{$e_0$ elliptique}}} \!\!\!
\interieur G \cdot (e_0 \exp \calv_{\!e_0}) \;$
de \cite [lem. 8.1.1 p.$\!$~72] {Bo87} (cf. \cite [lem. 40 p.$\!$~41] {DV93})
montre que dans la situation du (d), les $\Theta_{e_0}$ avec $e_0 \!\in\! G$
elliptique d\'eterminent $\Theta$.

\smallskip
D'apr\`es \cite [3.1 haut p.$\!$~21] {Bo87}, toute repr\'esentation unitaire
topologiquement irr\'eductible $\,\pi\,$ de $G$ est tra\c{c}able, avec pour
caract\`ere une fonction $\,\tr \pi \,$ localement int\'egrable sur $G$
invariante sous $\interieur G$ dont la restriction \`a $\Gssreg$ est
analytique.

\theoreme{
\label{formule du caractere}

\'Etant donn\'es $\tau \in \XInd_G(\lambdatilde)$ et
$\tau_+ \in \Xfin_G(\lambdatilde,{\gotha^*}^+)$
associ\'es comme dans \ref{parametres adaptes} (c), on a

\noindent
$\bigl( \tr T_{\lambdatilde,\tau}^G \bigr)_{\!e}
=\; D_e^{-\frac{1}{2}}\!\!\!
\Sum_{\scriptstyle \dot{\lambdatilde'_e} \,\in\, G(e) \backslash \gestregtilde}
\!\!
{\scriptstyle \abs{\{ {{\gotha'_e}^*}^+ \}}^{-1}}\!
\Biggl(
\Sum_{\substack
{\lambdatilde' \in\, G \cdot \lambdatilde \, \cap \, \gstregtilde(e) \\
\textrm{tel que } \,\lambdatilde'\![e] \textrm{ existe} \\
\hfill \textrm{et } \,\lambdatilde'\![e] \,=\, \lambdatilde'_e}}\!\!\!
\pm_{\widehat{e'}}\;\, \mi^{-d_{e',\lambdatilde}}\;
\tr\tau(\widehat{e'})\!
\Biggr)\!
\restriction{~{\widehat{\beta}}_{G(e) \cdot \lambdatilde'_e}}{\calv_{\!e}}$

\noindent
$= \bigl( \tr T_{\lambdatilde,{\gotha^*}^+\!,\tau_+}^G \bigr)_{\!e}
\!=\, \mi^{-d_e} D_e^{-\frac{1}{2}}\!\!\!
\Sum_{\scriptstyle \dot{\lambdatilde'_e} \,\in\, G(e) \backslash \gestregtilde}
\!\!\!\!\!
{\scriptstyle \abs{\{ {{\gotha'_e}^*}^+ \}}^{-1}}\!
\Biggl(
\Sum_{\substack{\lambda_+' \in\, G \cdot \lambda_+ \cap\, \gothg(e)^* \\
\textrm{tel que } \,\lambdatilde'\![e] \,=\, \lambdatilde'_e}}\!\!\!\!\!
\pm_{\widehat{e'_+}}\!
\tr\tau_+(\widehat{e'_+})\!
\Biggr)\!\!
\restriction{~{\widehat{\beta}}_{G(e) \cdot \lambdatilde'_e}}{\calv_{\!e}}$

\smallskip\noindent
o\`u\
on donne un sens \`a la somme portant sur $\lambdatilde'$ en choisissant
$\,g \in G\,$ tel que $\,\lambdatilde' = g\lambdatilde\,$ puis fixant
$\,\widehat{e'} \in G(\lambdatilde)^{\gothg/\!\gothg(\lambda)(\mi \rho_\F+)}\!$
au-dessus de $\,e' = g^{-1}eg$
\\
(resp. :
on donne un sens \`a la somme portant sur $\lambda'_+$ en choisissant
$\,g \in G\,$ tel que $\, \lambda'_+ = g\lambda_+$ puis posant
$\,\lambdatilde' = g\lambdatilde \,$ et fixant
$\, \widehat{e'_+} \in G(\lambda_+)^{\gothg/\gothh}$ au-dessus de
$\,e' = g^{-1}eg$),
\\
$\pm_{\widehat{e'}}$ est le signe tel que
$\; \scalo (\widehat{e'})
_{(1-\Ad e')(\gothg/\gothg(\lambda)(\mi \rho_\F+))}
=\, \pm_{\widehat{e'}}\;
\scalo (B_{\lambda_{\textit{can}}})
_{(1-\Ad e')(\gothg/\gothg(\lambda)(\mi\rho_\F+))}$
\\
(resp. :
$\pm_{\widehat{e'_+}}$ est le signe tel que
$\; \scalo (\widehat{e'_+})_{(1-\Ad e')(\gothg/\gothh)}
=\, \pm_{\widehat{e'_+}}\;\,
\scalo (B_{\lambda_+})_{(1-\Ad e')(\gothg/\gothh)}$),
\\
$d_{e'}^{\gothg(\lambda)(\mi \rho_{\F+})}\!$ est le coefficient $d_{e'}$
relatif \`a $\gothg(\lambda)(\mi \rho_{\F+})$ et
$\;d_{e',\lambdatilde} = d_{e'} - d_{e'}^{\gothg(\lambda)(\mi \rho_{\F+})}\!$,
\\
$\gothh'_e = \gothg(e)(\lambdatilde'\![e])$ et $\gotha'_e$ est la composante
infinit\'esimalement hyperbolique de $\gothh'_e$,
\\
($\lambda',\Fep+)=\lambdatilde'\![e]$
et $\,\{ {{\gotha'_e}^*}^+ \}$ est l'ensemble des chambres de
$(\gothg(e)(\lambda')(\mi \rho_{\Fep+}),\gotha'_e)$.

\noindent
Les sommations portant sur $G(e) \backslash \gestregtilde$ font intervenir un
nombre fini de termes somm\'es non nuls.
Les autres sommations portent sur des ensembles finis.
La condition << $\lambdatilde'\![e]$ existe >> peut \^etre oubli\'ee, car elle
est r\'ealis\'ee quand $\;\tr\tau(\widehat{e'}) \not= 0$.%

\smallskip
En particulier, pour tout $\,X \in {\calv_{\!1}}\,$ on a

\noindent
$\det\left( \frac{\me^{\ad X/2} - \me^{- \ad X/2}}{\ad X} \right) _{\!\!\gothg}
^{\!1/2}
\,{\scriptstyle \times}\; \tr T_{\lambdatilde,{\gotha^*}^+\!,\tau_+}^G(\exp X)
\;=\; \dim\tau_+ \; {\scriptstyle \times} \;
{\widehat{\beta}}_{G \cdot \lambdatilde}(X)$.

}

\remarque{

Voici deux exemples et une pr\'ecision importante.

\smallskip
{\bf (1)}
Premier exemple :
$\, G = GL(2,\R )$ et $\, e = \left({\scriptstyle\Rot{1}}\right)$,
$\, \lambdatilde = (0,\F+) \,$ avec
$\F+ \subseteq \mi \,\goth{so}(2)^* \oplus \centrealg{\gothg}^*$
et $\tau = \chi_{\lambdatilde}^G$.
Donc $\lambdatilde[e]$ existe dans $\gestregtilde$.
La somme portant sur $G \cdot \lambdatilde \cap \gstregtilde(e)$ fait ici
intervenir deux nombres complexes non nuls et oppos\'es.

\smallskip
{\bf (2)}
Deuxi\`eme exemple :
$G$ est produit semi-direct de $\Z / 2\Z $ par $SL(2,\R )^2$ o\`u
l'\'el\'ement non trivial de $\Z / 2\Z $ op\`ere sur $SL(2,\R )^2$ par
permutation des coordonn\'ees et
$\,e = (\left({\scriptstyle\Rot{1}}\right),\left({\scriptstyle\Rot{1}}\right))$,
$\lambdatilde = (0,\F+) \,$ avec $\, \F+ \subseteq (\mi \,\goth{so}(2)^*)^2$
non stable sous $\Z / 2\Z $ et $\tau = \chi_{\lambdatilde}^G$.
Donc $\lambdatilde[e]$ existe dans $\gestregtilde$.
Les deux \'el\'ements $\lambdatilde'$ de
$G \cdot \lambdatilde \cap \gstregtilde(e)$ sont conjugu\'es sous $G(e)$, et
v\'erifient :
$\lambdatilde'\![e]$ existe et $\lambdatilde'\![e]=\lambdatilde[e]$.
Les nombres complexes somm\'es correspondants sont non nuls, et bien s\^ur
\'egaux.%

\smallskip
{\bf (3)}
Il peut arriver que des mesures $\,\widehat{\beta}_{G(e) \cdot \lambdatilde'_e}$
non lin\'eairement ind\'ependantes interviennent avec des coefficients non nuls
dans la formule du caract\`ere.
En effet, pla\c{c}ons nous dans le cas suivant :
$G$ est le produit semi-direct canonique du groupe sym\'etrique $\sym{4}$ par
$SL(3,\R )^4$, $\lambda=0$, $\gothh$ est produit de $2$ copies d'une
sous-alg\`ebre de Cartan d\'eploy\'ee de $\goth{sl}(3,\R )$ avec
$2$ copies d'une sous-alg\`ebre de Cartan fondamentale de $\goth{sl}(3,\R )$,
${\gotha^*}^+$ est stable par permutation des $2$ premi\`eres composantes.
On note $\,e = (1\,2)$ et $g = (1\,3)(2\,4)$ dans $\sym{4}$.
Donc
$\, G \cdot \lambdatilde \cap \gstregtilde(e)
= G(e) \cdot \lambdatilde \,\cup\, G(e) \cdot g\lambdatilde$,
$\lambdatilde[e]$ et $(g\lambdatilde)[e]$ existent dans $\gestregtilde$,
$\, G(e) \cdot \lambdatilde[e] \not= G(e) \cdot (g\lambdatilde)[e] \,$
et
$\, \widehat{\beta}_{G(e) \cdot \lambdatilde[e]}
= \widehat{\beta}_{G(e) \cdot (g\lambdatilde)[e]} \,$
(cf. \ref{explication de la notation support} (3)).
En outre, il existe $\, \tau_+ \in \Xfin_G(\lambdatilde,{\gotha^*}^+)$ pour
lequel les coefficients de
$\restriction{~{\widehat{\beta}}_{G(e) \cdot \lambdatilde[e]}}{\calv_{\!e}}$
et de
$\restriction{~{\widehat{\beta}}_{G(e) \cdot (g\lambdatilde)[e]}}{\calv_{\!e}}$
sous formes de sommes portant sur $\lambdatilde'$ sont \'egaux et non nuls.
\cqfr

}

\remarque{

On suppose que $G$ est connexe et utilise les notations de la section \ref{Rep7}
Diff\'erents auteurs ont utilis\'es dans certaines situations des formules du
caract\`ere \'equivalentes \`a celle du th\'eor\`eme pr\'ec\'edent, mais
relatives \`a d'autres param\'etrisations.
Pour relier ces formules, je fixe un syst\`eme de racines positives
$\, R_0^+(\gothg_\C ,\gothh_\C ) \,$ de
$\, R(\gothg_\C ,\gothh_\C ) \,$ qui contient les conjugu\'ees de ses
racines non imaginaires.
On note $\rho_{\gothg,\gothh,0}$ et $\rho_{\gothm,\gotht,0}$ les demi-sommes
de racines positives associ\'ees \`a
$\, R_0^+(\gothg_\C ,\gothh_\C ) \,$ et \`a
$\, R_0^+(\gothg_\C ,\gothh_\C ) \cap
R(\gothm_\C ,\gotht_\C )$.

\smallskip
{\bf (a)}
L'application de $C(\gothg(\lambda),\gothh)$ dans l'ensemble des chambres de
$(\gothm_\C ,\gotht_\C )$ dans $\, \mi \,\gotht^* \,$ dont l'adh\'erence
contient $\mi \mu$, qui \`a un \'el\'ement $\Fp+$ de $C(\gothg(\lambda),\gothh)$
associe la chambre de $(\gothm_\C ,\gotht_\C )$ contenant
$\mi \mu_{\gothm,(\lambda,\Fp+),\{0\}}$, est bijective.
On note $\calc^+$ l'image de $\F+$ par cette application.
On pose

\centerline{
$\sg_0 (\F+,\mu)
= \Prod_{\substack
{\alpha \, \in \, R_0^+(\gothg_\C ,\gothh_\C ) \\
\textrm{$\alpha$ compacte}}}\!
\sg(\calc^+(H_\alpha)) \;\;
{\scriptstyle \times}
\Prod_{\substack
{\alpha \, \in \, R_0^+(\gothg_\C ,\gothh_\C ) \\
\textrm{$\alpha$ imaginaire non compacte}}}\!\!\!
\sg(-\calc^+(H_\alpha))$.
}

\smallskip
{\bf (b)}
On note $\, T = \centragr{M}{\gotht} = \centragr{M}{\gothm} \exp \gotht$.
L'application de $\Xirr_M(\mutilde)$ dans l'ensemble des caract\`eres unitaires
de $T$ de diff\'erentielle $\, \mi \,\mu - \rho_{\gothm,\gotht,0} \,$ qui \`a
$\sigma \in \Xirr_M(\mutilde)$ associe le caract\`ere
unitaire $\eta$ de $T$ tel que
$\; \eta (x\exp X)
= \sigma(x,1) \, \me^{(\mi \,\mu - \rho_{\gothm,\gotht,0})(X)} \;$
pour $\, x \in \centragr{M}{\gothm} \,$ et $\, X \in \gotht$, est bijective.

\smallskip
{\bf (c)}
Soit $\sigma \in \Xirr_M(\mutilde)$.
On note $\eta$ son image par l'application du (b).
La fonction $\Theta_\eta^{(\calc^+)}$ de \cite [th. 23 p.$\!$~260] {Va77}
relative au groupe $M$ s'\'ecrit :

\centerline{
$\Theta_\eta^{(\calc^+)}\!
= \,\sg_0 (\F+,\mu) \,{\scriptstyle \times}\, \tr T_{\mutilde,\sigma}^M$.
}

\noindent
En particulier, quand $M_0$ est semi-simple et << acceptable >>, la fonction
$\Theta_{\mi \,\mu,\calc^+}$ de \cite [p.$\!$~305] {Ha65} relative au groupe
$M_0$ s'\'ecrit :

\centerline{
$\Theta_{\mi \,\mu,\calc^+}
= \sg_0 (\F+,\mu) \,{\scriptstyle \times}\, \tr T_{\mutilde}^{M_0}$.
}

Par ailleurs, quand $G$ est lin\'eaire, $\,\tr T_{\mutilde,\sigma}^M$
est not\'e
$\, \Theta^M \bigl( \mi \mu,\calc^+\!,
\restriction{\sigma({\scriptscriptstyle\bullet},1)}{\centregr{M}}\bigr) \,$
dans \cite [p.$\!$~397] {KZ82} en identifiant $\gotht$ et $\gotht^*\!$.

\smallskip
{\bf (d)}
On se place dans le cas o\`u $G$ est la composante neutre du groupe des points
r\'eels d'un groupe lin\'eaire alg\'ebrique complexe d\'efini sur $\R $
semi-simple et simplement connexe.
On dispose ainsi d'un caract\`ere complexe $\xi_\rho$ de $\exp \gothh_\C $ de
diff\'erentielle $\rho_{\gothg,\gothh,0}$.
On fixe un caract\`ere unitaire $\eta$ de $T$ comme au (b).
Lorsque $\nu$ est $(\gothg,\gotha)$-r\'eguli\`ere, la fonction
$\theta(TA,\xi_\rho\,\eta,\nu)$ de \cite [p.$\!$~244] {Hr83}~s'\'ecrit :

\smallskip\centerline{
$\theta(TA,\xi_\rho\,\eta,\nu)
\,=\, \abs{C(\gothg(\lambda),\gothh)}^{-1} \;{\scriptstyle \times}\!
\Sum_{\Fp+ \in \, C(\gothg(\lambda),\gothh)_{reg}}\!\!
\sg_0 (\Fp+\!,\mu) \,{\scriptstyle \times}\,
\tr T_{\lambdatilde',\tau'}^G$}

\negsmallskip\noindent
o\`u $\;\lambdatilde' = (\lambda,\Fp+)$,
$\lambda'_+$ est associ\'e \`a $(\lambdatilde',\gotha^*)$
comme dans \ref{choix de racines positives} (b),
et $\tau' \in \Xirr_G(\lambdatilde')$ est d\'etermin\'e par l'\'egalit\'e
$\,\restriction{(\delta_{\lambda'_+}^{\gothg/\gothh} \tau')}{\centregr{M}}\!\!
=\, \restriction{\eta}{\centregr{M}}$.
\cqfr

}

\medskip
La d\'emonstration du th\'eor\`eme va se faire en trois \'etapes.
Les deux premi\`eres \'etapes suivent l'article de Bouaziz :
on passe de $G$ \`a $G(\lambda_+) G_0$, puis de $G(\lambda_+) G_0$ \`a~$M'\!$.
Dans la derni\`ere \'etape, une variante du foncteur de translation de
Zuckerman va me permettre de relier les caract\`eres des repr\'esentations de
$M'$ qui m'int\'eressent aux caract\`eres \'etudi\'es par Bouaziz.

\smallskip
\dem{D\'ebut de la d\'emonstration du th\'eor\`eme}{

Soient $\,\widehat{e'_+} \in G(\lambda_+)^{\gothg/\gothh}$ et
$\,\widehat{e'} \in G(\lambdatilde)^{\gothg/\!\gothg(\lambda)(\mi \rho_\F+)}$
au-dessus d'un m\^eme \'el\'ement elliptique $e'$ de $G(\lambda_+)$.
Les valeurs propres autres que $-1$ et $1$ de la restriction de ${\Ad e'}^\C $
\`a
$\,\gothv_\C \egdef
\Sumpetit_{\alpha \in R^+(\gothg(\lambda)(\mi \rho_\F+)_\C ,\gothh_\C )}
\!\gothg_\C ^\alpha\,$
sont non r\'eelles deux \`a deux conjugu\'ees de module $1$.
Comme $\call_{\lambdatilde,{\gotha^*}^+}/\gothh_\C $
(cf. \ref{espace symplectique canonique} (a))
est somme directe des sous-espaces ${\Ad e'}^\C $-invariants
$\call_{\lambdatilde}/\gothg(\lambda)(\mi \rho_\F+)_\C \!$
(cf. \ref{espace symplectique canonique} (c)) et $\gothv_\C $, on a
$\; q_{\,\call_{\lambdatilde,{\gotha^*}^+}/\gothh_\C }
((\Ad e')_{\gothg/\gothh})
= q_{\,\call_{\lambdatilde}/\gothg(\lambda)(\mi \rho_\F+)_\C \!}
((\Ad e')_{\gothg / \gothg(\lambda)(\mi \rho_\F+)\!}) \,$
relativement \`a $B_{\lambda_+}$,
puis d'apr\`es \ref{espace symplectique canonique} et
\ref{geometrie metaplectique} (d) :

\centerline{
$\pm_{\widehat{e'_+}}\;
\rho_{\lambdatilde,{\gotha^*}^+}^{\gothg/\gothh}(\widehat{e'_+})
\;=\;
\pm_{\widehat{e'}}\;
\rho_\lambdatilde^{\gothg/\!\gothg(\lambda)(\mi \rho_\F+)}(\widehat{e'})\;
{\scriptstyle \times}\;
\mi^{-1/2 \dim (1-\Ad e')(\gothg(\lambda)(\mi\rho_\F+) / \gothh)}$
}

\noindent
o\`u les signes $\;\pm_{\widehat{e'_+}}$ et $\pm_{\widehat{e'}}$ sont d\'efinis
comme dans le th\'eor\`eme.

Soient $\lambdatilde'_0 \in G \cdot \lambdatilde \cap \gstregtilde(e)$,
$g_0 \in G\,$ v\'erifiant $\,\lambdatilde'_0 = g_0\lambdatilde$, et
$\widehat{e'_0} \in G(\lambdatilde)^{\gothg/\!\gothg(\lambda)(\mi \rho_\F+)}\!$
au-dessus de $e'_0 = g_0^{-1}eg_0$.
Par la d\'efinition \ref{parametres adaptes} (c) et ce qui pr\'ec\`ede, on
obtient

\centerline{
$\pm_{\widehat{e'_0}}\;\,\mi^{-d_{e'_0,\lambdatilde}}
\tr\tau(\widehat{e'_0})
\;= \Sum_{\substack{\dot{x} \,\in\, G(\lambdatilde)/G(\lambda_+) \\
\mathrm{tel\ que\ } x^{-1}e'_0\,x \,\in\, G(\lambda_+)}}\!\!
\pm_{\widehat{e'_{+,x}}}\;\,\mi^{-d_e}\tr\tau_+(\widehat{e'_{+,x}})$
}

\noindent
ind\'ependamment du choix de
$\; \widehat{e'_{+,x}} \in G(\lambda_+)^{\gothg/\gothh} \,$ au-dessus de
$\, e'_x = x^{-1}e'_0\,x$.

\noindent
Ensuite, pour tout $\lambdatilde'_e \in \gestregtilde$ on obtient l'\'egalit\'e
de sommes finies

\centerline{
$\Sum_{\substack
{\lambdatilde'_0 \in\, G \cdot \lambdatilde \, \cap \, \gstregtilde(e) \\
\textrm{tel que } \,\lambdatilde'_0[e] \,=\, \lambdatilde'_e}}\!\!
\pm_{\widehat{e'_0}}\,\mi^{-d_{e'_0,\lambdatilde}}
\tr\tau(\widehat{e'_0})
\;\,=\,
\Sum_{\substack{\lambda_+' \in\, G \cdot \lambda_+ \cap\, \gothg(e)^* \\
\textrm{tel que } \,\lambdatilde'\![e] \,=\, \lambdatilde'_e}}\!\!
\pm_{\widehat{e'_+}}\,\mi^{-d_e}
\tr\tau_+(\widehat{e'_+})$.
}

\noindent
Les deux expressions du caract\`ere propos\'ees dans le th\'eor\`eme sont donc
\'egales.

\smallskip
Soit $\lambdatilde'\!$ (resp. $\lambda_+'$) un terme associ\'e \`a un 
$\lambdatilde'_e\!$ comme dans le th\'eor\`eme.
On a :
\\
$\{\; g\lambdatilde \in G \cdot \lambdatilde \cap \gstregtilde(e) \mid
(g\lambdatilde)[e] = \lambdatilde'_e \textrm{ et }
g\lambdatilde \in G(e) \cdot \lambdatilde' \,\}
= G(e)(\lambdatilde'_e) \cdot \lambdatilde'$
\\
(resp.
$\; \{\; g\lambda_+ \in G \cdot \lambda_+ \cap \gothg(e)^* \mid
(g\lambdatilde)[e] = \lambdatilde'_e \textrm{ et }
g\lambda_+ \in G(e) \cdot \lambda_+' \,\}
= G(e)(\lambdatilde'_e) \cdot \lambda_+'\,$).
\\
Le th\'eor\`eme \'equivaut donc \`a chacune des formules $(F)$ et $(F_+)$
suivantes :

\noindent
$\bigl( \tr T_{\lambdatilde,\tau}^G \bigr)_{\!e}
\stackrel{\raisebox{.5ex}{\boldmath$
{\scriptscriptstyle (}\scriptstyle F{\scriptscriptstyle )}$}}{=}\,
D_e^{-\frac{1}{2}}\!\!
\Sum_{\substack
{\dot{\lambdatilde'} \,\in\, G(e) \backslash\, 
G \cdot \lambdatilde \, \cap \, \gstregtilde(e) \\
\textrm{tel que } \,\lambdatilde'\![e] \,\in\, \gestregtilde}}\!\!
\pm_{\widehat{e'}}\; \mi^{-d_{e',\lambdatilde}}\,
\textstyle\frac
   {\abs{G(e)(\lambdatilde'\![e]) \cdot \lambdatilde'}}
   {\abs{\{ {{\gotha'_e}^*}^+ \}}}\;
\tr\tau(\widehat{e'})\,
\restriction{~{\widehat{\beta}}_{G(e) \cdot \lambdatilde'\![e]}}{\calv_{\!e}}$%
\labind{ce} ;

\noindent
$\bigl( \tr T_{\lambdatilde,{\gotha^*}^+\!,\tau_+}^G \bigr)_{\!e}\!
\stackrel{\raisebox{.5ex}{\boldmath$
{\scriptscriptstyle (}\scriptstyle F_+{\scriptscriptstyle )}$}}{=}\,
\mi^{-d_e} D_e^{-\frac{1}{2}}\!\!
\Sum_{\substack
{\dot{\lambda'_+} \,\in\, G(e) \backslash\,
G \cdot \lambda_+ \cap\, \gothg(e)^* \\
\textrm{tel que }\,\lambdatilde'\![e] \,\in\, \gestregtilde}}\!\!
\pm_{\widehat{e'_+}}\,
\textstyle\frac
   {\abs{G(e)(\lambdatilde'\![e]) \cdot \lambda'_+}}
   {\abs{\{ {{\gotha'_e}^*}^+ \}}}\,
\tr\tau_+(\widehat{e'_+})\,
\restriction{~{\widehat{\beta}}_{G(e) \cdot \lambdatilde'\![e]}}{\calv_{\!e}}$.

\noindent
Cela permet en particulier de pr\'eciser les r\'esultats \'enonc\'es dans
l'introduction.

\smallskip
D'apr\`es la d\'efinition \ref{mesures de Liouville} (e), on a
$\;{\widehat{\beta}}_{\omega}
=\!\! \Sumpetit_{\omega_0 \in G(e)_0 \backslash \omega}
{\widehat{\beta}}_{\omega_0}$
pour $\,\omega \in G(e) \backslash \gestregtilde$.
Ainsi le th\'eor\`eme \'equivaut \`a l'\'enonc\'e qu'on obtient en y 
rempla\c{c}ant $G(e)$ par $G(e)_0$, et donc aussi \`a celui qu'on obtient en
rempla\c{c}ant $G(e)$ par $G(e)_0$ dans $(F_+)$.
Soient $\lambda'_+$ et $\lambdatilde'$ comme dans $(F_+)$.
Le groupe $G(e)_0(\lambdatilde'\![e])$ op\`ere transitivement sur
$\{ {{\gotha'_e}^*}^+ \}$ par l'action de certains \'el\'ements de
$W(\gothg(e)_\C ,{\gothh'_e}_\C )$ (notations du th\'eor\`eme).
Comme $\lambda'_+$ appartient \`a ${\gothh'_e}^* \cap \gestssreg$, le
stabilisateur d'un \'el\'ement de $\{ {{\gotha'_e}^*}^+ \}$ est \'egal \`a
$G(e)_0(\lambdatilde'\![e])(\lambda'_+)$ d'apr\`es le lemme \ref{lemme clef}.
Donc
$\;\abs{G(e)_0(\lambdatilde'\![e]) \cdot \lambda'_+}
= \abs{\{ {{\gotha'_e}^*}^+ \}}$.%

\medskip
Par cons\'equent, les formules \`a d\'emontrer \'equivalent \`a la formule
suivante :

\smallskip\noindent
$\!\bigl( \tr T_{\lambdatilde,{\gotha^*}^+\!,\tau_+}^G \bigr)_{\!e}
=\, \mi^{-d_e} \, D_e^{-\frac{1}{2}}\!\!
\Sum_{\substack
{\dot{\lambda'_+} \in\, G(e)_0 \backslash\,
G \cdot \lambda_+ \cap\, \gothg(e)^* \\
\textrm{tel que }
(g\lambdatilde)[e] \,\in\, \gestregtilde}}\!\!\!\!
{\textstyle \frac
{\scalo (\widehat{e'_+})_{(1-\Ad e')(\gothg/\gothh)}}
{\scalo (B_{\lambda_+})_{(1-\Ad e')(\gothg/\gothh)}}} \;
\tr\tau_+(\widehat{e'_+}) \,
\restriction{~{\widehat{\beta}}_{G(e)_0 \cdot (g\lambdatilde)[e]}}{\calv_{\!e}}$

\noindent
o\`u
$\;\, g \in G \,$ v\'erifie $\; \lambda'_+ \!= g\lambda_+ \;$
et
$\; \widehat{e'_+} \in G(\lambda_+)^{\gothg/\gothh} \,$ est au-dessus de
$\, e' = g^{-1}eg$.

\noindent
On remarque pour la suite (section \ref{Rep11}) que dans le raisonnement
pr\'ec\'edent, on peut remplacer $G(e)_0$ par n'importe quel sous-groupe $L$ de
$G(e)$ qui contient $G(e)_0$ et dont le groupe adjoint est inclus dans
$\interieur\gothg(e)_\C $.

\smallskip
Les orbites sur lesquelles on somme \'etant semi-simples r\'eguli\`eres, on va
pouvoir reproduire la plupart des arguments de \cite {Bo87}.

\smallskip
On a
$\;\, T_{\lambdatilde,{\gotha^*}^+\!,\tau_+}^G\!
=\, \Ind_{G(\lambda_+) \, G_0}^G
T_{\lambdatilde,{\gotha^*}^+\!,\tau_+}^{G(\lambda_+) \, G_0} \;$
par << induction par \'etages >>, puis

\centerline{
$\displaystyle
\bigl( \tr T_{\lambdatilde,{\gotha^*}^+\!,\tau_+}^G \bigr) _{\!e}
\;= \Sum_{\dot{x} \,\in\, G / G(\lambda_+) \, G_0\,
\textrm{ tel que } \,x^{-1}ex \,\in\, G(\lambda_+) \, G_0} \!\!\!
\Ad x \cdot
\bigl(\tr T_{\lambdatilde,{\gotha^*}^+\!,\tau_+}^{G(\lambda_+) \, G_0}\bigr)
_{\!x^{-1}ex}$.
}

Comme l'ensemble $\;\, G \cdot \lambda_+ \cap\, \gothg(e)^* \;$ est r\'eunion
des ensembles disjoints
$\, \Ad^* x \cdot \bigl( G(\lambda_+) G_0
\cdot \lambda_+ \cap\, \gothg(x^{-1}ex)^* \bigr) \,$
o\`u $\, \dot{x} \in G / G(\lambda_+) G_0$ v\'erifie
$\, x^{-1}ex \in G(\lambda_+) G_0$, on peut supposer que
$\, G = G(\lambda_+) \, G_0$.
\cqfdpartiel

}

\section{\boldmath
Passage de $M'$ \`a $G$
}\label{Rep11}

Dans cette section, on va se ramener \`a montrer le th\'eor\`eme pour
$\, T_{\lambdatilde,\tau_{M'}}^{M'}$ en reprenant mot \`a mot la
preuve de Bouaziz concernant le cas semi-simple r\'egulier.

\smallskip
Les r\'esultats \'enonc\'es dans la proposition suivante sont d\'emontr\'es dans
\cite [lem. 7.1.3 p.$\!$~65 et lem. 4.2.1 p.$\!$~38]{Bo87}.

\proposition{\hspace{-15pt}{\bf{(Bouaziz)}}\quad 
\label{caracteres d'induites}
%
{\bf (a)}
Soit $\pi$ une repr\'esentation unitaire topologiquement irr\'eductible de $M'$.
On lui associe
$\; \Pi = \Ind_{M'U}^{M'G_0} (\, \pi \, \otimes \, \fonctioncar{U} \,)$.
La repr\'esentation unitaire $\Pi$ de $M'G_0$ est somme hilbertienne d'un nombre
fini de sous-repr\'esentations topologiquement irr\'eductibles
(cf. \cite [prop. 5.4.13 (i) p.$\!$~110] {Di64} ou
\cite [prop. p.$\!$~25] {Wa88}).

On suppose que $e \in M'G_0$.
Pour tout $\, \gothj_e \in \Car \gothg(e)$ et tout $\, X \in \gothj_e$
tel que $\, e \exp X \in \MpGzerossreg$, on~a

\smallskip\noindent
\parbox{\textwidth}{
$\quad\; \abs{D_G(e\exp X)}^{1/2} \;\; (\tr\Pi) (e\exp X)$

\smallskip
$\displaystyle
=\! \Sum_{ (M' \cap G_0)\,g \;\in\; M' \cap G_0 \,\backslash\,
\{ g_0 \in G_0 \,\mid\, g_0 \, e\exp \gothj_e\,g_0^{-1} \subseteq M' \} }
\!\!\!
|D_{M'}(g e \exp X g^{-1})|^{1/2} \; (\tr\pi) (g e \exp X g^{-1})$
}

\smallskip\noindent
o\`u la somme de droite porte sur un ensemble fini.

\smallskip
{\bf (b)}
Soient $\nu' \in \gothg^*$ une forme lin\'eaire hyperbolique, $\gothh'$ une
sous-alg\`ebre de Cartan fondamentale de $\gothg(\nu')$, et $M'_1$ un
sous-groupe de $G_0(\nu')$ qui contient le groupe
$\centragr{G_0}{\gothh'} G(\nu')_0$.
Pour tous $\lambda' \in {\gothh'}^* \cap \gstssreg$,
$\gothj \in \Car \gothg$ et $X \in \gothj \cap \gssreg$, on a%

\smallskip\centerline{
$\abs{\cald_{\gothg}(X)}^{1/2} \;\;
{\widehat{\beta}}_{G_0 \cdot \lambda'} (X) \;
= \; \Sum_{ M'_1\,x \;\in\; M'_1 \,\backslash\,
\{ x_0 \in G_0 \,\mid\, x_0 \gothj \,\subseteq\, \gothg(\nu') \} }
|\cald_{\gothg(\nu')}(x\,X)|^{1/2} \;\;
{\widehat{\beta}}_{M'_1 \cdot \lambda'} (x\,X)$
}

\smallskip\noindent
o\`u la somme de droite porte sur un ensemble fini.

}

\smallskip
\dem{Suite de la d\'emonstration du th\'eor\`eme}{

On suppose r\'ealis\'ee l'\'egalit\'e $\; G = G(\lambda_+) \, G_0$,
\`a laquelle on s'est ramen\'e. \\
On admet (provisoirement) que le th\'eor\`eme est vrai pour
$\, T_{\lambdatilde,\tau_{M'}}^{M'}\!$.

\smallskip
Pour obtenir la formule du caract\`ere pour
$\, \tr T_{\lambdatilde,{\gotha^*}^+\!,\tau_+}^G\!$, il suffit de la v\'erifier
sur l'ouvert inclus dans $\,\gessreg$ de compl\'ementaire n\'egligeable (cf.
\cite [avant lem. 1.4.1 p.$\!$~6] {Bo87}) et sur lequel les deux membres de
l'\'egalit\'e sont des fonctions analytiques, form\'e des $\; X \in \calv_{\!e}
\,$ tels que $\; e \exp X \in \Gssreg$.

On fixe $\, \gothj_e \in \Car \gothg(e)$ et
$\, X \in \gothj_e \cap \calv_{\!e} \,$ tels que
$\, e \exp X \in \Gssreg$.
On a :

\smallskip\centerline{
$k_e(X)
\,=\, |\det (1-\Ad e)_{\gothg/\gothg(e)}|^{-1/2} \,{\scriptstyle \times}\,
|\cald_{\gothg(e)}(X)|^{-1/2} \,{\scriptstyle \times}\,
|D_G(e \exp X)|^{1/2}$.
}

\medskip
En paraphrasant les calculs de Bouaziz, on fixe des repr\'esentants
$g_1,\dots,g_n$ des doubles classes dans
$\; M' \cap G_0 \,\backslash\,
\{ g_0 \in G_0 \mid g_0 e g_0^{-1} \in M' \} \,/\, G(e)_0$.
On pose $e_1\!=\!g_1e g_1^{-1},\dots,e_n\!=\!g_ne g_n^{-1}$.
La r\'eunion disjointe des ensembles

\smallskip\centerline{
$M' \cap G(e_j)_0 \,\backslash\,
\{ x_0 \in G(e_j)_0 \mid x_0 g_j \, \gothj_e \subseteq \gothm'(e_j) \} \;$
o\`u $\, 1 \leq j \leq n$
}

\noindent
est en bijection avec l'ensemble

\centerline{
$M' \cap G_0 \,\backslash\,
\{ g_0 \in G_0 \mid g_0 \, e \exp \gothj_e \, g_0^{-1} \subseteq M' \}$
}

\noindent
par l'application qui envoie $(j,\dot{x})$ sur $\dot{x g_j}$.
On \'ecrit
$\, \bigl( \tr T_{\lambdatilde,{\gotha^*}^+\!,\tau_+}^G \bigr)_{\!e} (X) \,$
\`a l'aide de la proposition \ref{caracteres d'induites} (a) en prenant
$\, \pi = T_{\lambdatilde,\tau_{M'}}^{M'}$
et en tenant compte de cette bijection.
On utilise ensuite l'expression de
$\, \bigl( \tr T_{\lambdatilde,\tau_{M'}}^{M'} \bigr)_{\!e_j} \,$
propos\'ee au d\'ebut de la d\'emonstration du th\'eor\`eme,
avec $\, L = M' \cap G(e_j)_0$.
Cela fournit trois sommations.
On intervertit les deux derni\`eres et tombe sur :

\smallskip\centerline{
$\Sum_{1 \leq j \leq n} \;\;\;
\Sum_{\substack
{\dot{g\lambda_+} \;\in\; M' \cap G(e_j)_0 \,\backslash\,
M' \cdot \lambda_+ \cap\, \gothm'(e_j)^* \\
\textrm{tel que } (g\lambdatilde)[e_j] \,\in\, \gejstregtilde}} \;\;
\Sum_{\dot{x} \;\in\; M' \cap G(e_j)_0 \,\backslash\,
\{ x_0 \in G(e_j)_0 \,\mid\, x_0 \, g_j \gothj_e \,\subseteq\, \gothm' \}}
\;\;\cdots\;\;$.
}

\smallskip
La derni\`ere somme se calcule en appliquant la proposition \ref{caracteres
d'induites} (b) apr\`es avoir remplac\'e $G$ par $G(e_j)$, $\nu'$ par $\nu_+$,
$\gothh'$ par $(g\gothh)(e_j)$, $M'_1$ par $M' \cap G(e_j)_0$, $\lambda'$ par
$g\lambda_t$ avec $\,t \in \left]0,1\right]$ (cf. \ref{lemme clef}), $\gothj$ par
$g_j\gothj_e$ ,et $X$ par $g_jX$, puis en passant \`a la limite $t \to 0^+$
(cf. \ref{limites de mesures de Liouville} (b)).
En outre, la r\'eunion disjointe des ensembles
$\; M' \cap G(e_j)_0 \,\backslash\, M' \cdot \lambda_+ \!\cap \gothm'(e_j)^*\;$
o\`u $\, 1 \leq j \leq n \,$ est en bijection avec l'ensemble
$\; G(e)_0 \,\backslash\, G \cdot \lambda_+ \!\cap \gothg(e)^*$
par l'application qui envoie $(j,\dot{l})$ sur $\dot{g_j^{-1} l}$.
Cela permet de regrouper les deux premi\`eres sommations en une seule.

\smallskip
Soient
$\; \widehat{e'_+} \in G(\lambda_+)^{\gothg/\gothh} \;$
et
$\; \widehat{e'_{\gothm'}} \in M'(\lambdatilde)^{\gothm'/\gothh} \;$
au-dessus d'un \'el\'ement elliptique $e'$ de $\, G(\lambda_+)$.
Les valeurs propres autres que $-1$ et $1$ de la restriction de ${\Ad e'}^\C $
\`a $\gothu_\C $ sont non r\'eelles deux \`a deux conjugu\'ees de module $1$.
Comme $\, \call_{\lambda_t}/\gothh_\C $
(cf. \ref{fonctions canoniques sur le revetement double})
est somme directe des sous-espaces ${\Ad e'}^\C $-invariants
$\gothb_{M'} \!/ \gothh_\C $ et $\gothu_\C $, on a
$\, q_{\,\call_{\lambda_t} / \gothh_\C }((\Ad e')_{\gothg/\gothh})
= q_{\,\gothb_{M'} \!/ \gothh_\C }((\Ad e')_{\gothm'/\gothh}) \,$
relativement \`a $B_{\lambda_+}$,
puis d'apr\`es \ref{geometrie metaplectique} (d) :

\centerline{
$\frac
{\scalo (\widehat{e'_+})_{(1-\Ad e')(\gothg/\gothh)}}
{\scalo (B_{\lambda_+})_{(1-\Ad e')(\gothg/\gothh)}} \;
\rho_{\lambda_+}^{\gothg/\gothh} (\widehat{e'_+})
=
\frac
{\scalo (\widehat{e'_{\gothm'}})_{(1-\Ad e')(\gothm'/\gothh)}}
{\scalo (B_{\lambda_{+,\gothm'}})_{(1-\Ad e')(\gothm'/\gothh)}} \;
\rho_{\lambda_{+,\gothm'}}^{\gothm'\!/\gothh} (\widehat{e'_{\gothm'}})
\;{\scriptstyle \times}\;
\mi^{-1/2 \dim (1-\Ad e')(\gothg/\gothm')}$.
}

\`A l'aide de l'\'egalit\'e \'ecrite dans la d\'emonstration du lemme
\ref{construction equivalente, cas connexe} (a), on en d\'eduit la formule du
caract\`ere pour $\tr T_{\lambdatilde,{\gotha^*}^+\!,\tau_+}^G\!$
(en l'admettant pour $\tr T_{\lambdatilde,\tau_{M'}}^{M'}$).
\cqfdpartiel

}

\section{
Translation au sens de G.~Zuckerman
}\label{Rep12}

Dans cette section, il s'agit de cr\'eer un outil pour la seule \'etape
d\'elicate : le calcul du caract\`ere de la repr\'esentation d\'efinie apr\`es
un passage dans l'homologie.

\smallskip
Je commence par rappeler quelques r\'esultats classiques.

\lemme{
\label{composantes primaires}

Soient $\,\pi$ une repr\'esentation lin\'eaire continue de $G_0$ dans un espace
de Hilbert complexe $V\!$ et $K_0$ un sous-groupe compact maximal de $G_0$.
Dans ces conditions, on notera $V_{K_0}$ la r\'eunion des sous-espaces
vectoriels de dimension finie de $V\!$ stables par $K_0$.

Pour chaque $N \!\in \N $, les propri\'et\'es (i) et (ii) suivantes sont
\'equivalentes :

{\bf (i)}
il existe une suite finie
$\; V_0=V \supseteq V_1 \supseteq \cdots \supseteq V_N=\{0\} \;$
de sous-$G_0$-modules ferm\'es de $V\!$ tels que les quotients $V_{j-1}/V_j$
avec $1 \leq j \leq N \,$ sont topologiquement irr\'eductibles sous $G_0$ et
poss\`edent un caract\`ere infinit\'esimal ;

{\bf (ii)}
les multiplicit\'es des \'el\'ements de $\widehat{K_0}$ dans $V_{K_0}$ sont
finies, ce qui implique que $\gothg$ op\`ere sur $V_{K_0}$ par d\'erivation de
l'action de $G_0$ (cf. \cite [th. 14 p.$\!$~313] {Va77}), et il existe une
suite finie
$\; (V_{K_0})_0=V_{K_0} \supseteq (V_{K_0})_1 \supseteq
\cdots \supseteq (V_{K_0})_N = \{0\} \;$
de sous-$(\gothg_\C ,K_0)$-modules de $V_{K_0}$ tels
que les quotients
$\,(V_{K_0})_{j-1} \,/\, (V_{K_0})_j\,$ avec $1 \leq j \leq N$ sont des
$(\gothg_\C ,K_0)$-modules irr\'eductibles.

Quand ces propri\'et\'es sont v\'erifi\'ees, $\pi$ est tra\c{c}able, et les
<< composantes $\centrealg{U \gothg_\C }$-primaires >> $V^\chi$\labind{Vchi} de
$V$ associ\'ees aux caract\`eres $\chi$ de l'alg\`ebre unif\`ere
$\centrealg{U \gothg_\C }$, d\'efinies comme les adh\'erences dans $V\!$ des
sous-espaces vectoriels

\centerline{
$(V^{\infty})^\chi
\egdef
\{\, v\in V^{\infty} \mid \forall u \in \centrealg{U \gothg_\C } \;\;
(u-\chi(u))^n \! \cdot v = 0 \, \} \;$
pour << $n$ grand >>,
}

\noindent
sont en somme directe de somme dense dans $V\!$ ;
de plus l'ensemble $P$ des poids de $\centrealg{U \gothg_\C }$ dans
$V^{\infty}\!$, form\'e des caract\`eres $\chi$ de $\centrealg{U \gothg_\C }$
tels que $\, V^\chi \not= \{0\}$, est~fini.%

}

\dem{D\'emonstration du lemme}{

On constate que pour tout sous-$G_0$-module ferm\'e $W$ de $V\!$, l'injection
canonique de $V_{K_0}/W_{K_0}$ dans $(V/W)_{K_0}$ est surjective.
En effet, en notant $\pr_{V/W}$ la projection canonique de $V$ sur $V/W\!$, pour
tout sous-$K_0$-module $F$ de dimension finie de $V/W\!$ le sous-espace
vectoriel dense $\, \pr_{V/W}(V_{K_0} \cap \pr_{V/W}^{-1}(F)) \,$ de $F$ est
\'egal \`a~$F$.
Par ailleurs, on dispose pour chaque $\delta \in \widehat{K_0}$ d'un isomorphime
de $K_0$-modules canonique de $V_\delta/W_\delta$ sur
$(V_{K_0}/W_{K_0})_\delta\,$
(en notant $E_\delta$ la composante isotypique de type $\delta$ d'un
$K_0$-module $E$).
Donc pour chaque $\delta \in \widehat{K_0}$, les $K_0$-modules
$V_\delta/W_\delta$ et $(V/W)_\delta\,$ sont isomorphes.

Quand (i) est v\'erifi\'e, les espaces vectoriels $(V_{j-1}/V_j)_\delta$
avec $1 \leq j \leq N$ et $\delta \in \widehat{K_0}$ sont de dimension finie
d'apr\`es \cite [prop. 16 p.$\!$~314 et th. 19 p.$\!$~316] {Va77}, puis les
espaces vectoriels $V_\delta$ avec $\delta \in \widehat{K_0}$ sont de dimension
finie ; la suite
$\; (V_0)_{K_0}=V_{K_0} \supseteq (V_1)_{K_0} \supseteq
\cdots \supseteq (V_N)_{K_0}=\{0\} \;$
satisfait la condition de (ii) d'apr\`es \cite [fin de th. 14 p.$\!$~313] {Va77}.

Quand (ii) est v\'erifi\'e, la suite
$\; \overline{(V_{K_0})_0\!}\,=V \supseteq \overline{(V_{K_0})_1\!}\, \supseteq
\cdots \supseteq \overline{(V_{K_0})_N\!}\,=\{0\} \;$
satisfait la condition de (i), \`a nouveau d'apr\`es \cite [th. 14 p.$\!$~313]
{Va77}.

\smallskip
On suppose dor\'enavant que les propri\'et\'es (i) et (ii) sont v\'erifi\'ees.
Le r\'esultat << $\pi$ est tra\c{c}able >> est d\^u \`a Harish-Chandra
(cf. \cite [8.1.2 p.$\!$~292] {Wa88}).

D'apr\`es (ii) et \cite [cor. 7.207 p.$\!$~530] {KV95}, on peut appliquer
\cite [prop. 7.20 p.$\!$~446] {KV95} \`a la fois \`a $V^{\infty}$ et \`a ses
composantes $\centrealg{U \gothg_\C }$-primaires (ici \emph{a priori} au sens de
\cite {KV95}).
Par cons\'equent, le membre droite de l'\'egalit\'e qui d\'efinit
$(V^{\infty})^\chi$ stationne en $n\/$ et $V^{\infty}$ est somme directe d'un
nombre fini de ses composantes $\centrealg{U \gothg_\C }$-primaires.
On fixe $n \in \N $ qui convienne relativement aux d\'efinitions de tous les
$(V^{\infty})^\chi$, $\chi \in P$.

Soit $Q \!\subseteq\! P\!$.
Les caract\`eres de $\centrealg{U \gothg_\C }$ \'etant lin\'eairement
ind\'ependants, on~a

\smallskip\centerline{
$\Sum_{\chi \in Q} (V^{\infty})^{\chi}
=
\Big\{\, v\in V^{\infty} \mid \forall u \in \centrealg{U \gothg_\C } \;\;
\Prod_{\chi \in Q}(u-\chi(u))^n \cdot v = 0 \, \Big\}$.
}

\noindent
On se donne une mesure de Haar $\diff _{G_0}$ sur $G_0$ et une suite
$\,(\varphi_k)_{k \in \N }\,$ d'\'el\'ements positifs de $C^\infty_c(G_0)$
d'int\'egrale $1$ pour $\diff _{G_0}$, dont les supports << rentrent dans
tout voisinage de $1$ dans~$G_0$ >>.
On a
$\;\, (\overline{\Sumpetit_{\chi \in Q} V^{\chi}\!}\,)^{\infty}\!
= \Sumpetit_{\chi \in Q} (V^{\infty})^{\chi} \,$
car pour tout $u \in \centrealg{U \gothg_\C }$ les endomorphismes continus
$\, \pi\bigl( \varphi_k \diff _{G_0} \ast
\Prodpetit_{\chi \in Q}(u-\chi(u))^n \bigr) \,$
de $V$ avec $\,k \in \N \,$ sont nuls sur
$\, \Sumpetit_{\chi \in Q} (V^{\infty})^{\chi} \,$
puis sur
$\, \overline{\Sumpetit_{\chi \in Q} V^{\chi}}$,
et d'autre part convergent simplement sur
$\, (\overline{\Sumpetit_{\chi \in Q} V^{\chi}\!}\,)^{\infty} \,$
quand $k \to +\infty$ vers l'op\'erateur associ\'e \`a
$\, \Prodpetit_{\chi \in Q}(u-\chi(u))^n$.

Pour chaque $\, \chi_0 \in P\!$, on a donc
$\, \bigl(
V^{\chi_0} \cap\, \overline{\!\Sumpetit_{\chi \in P \setminus \{\chi_0\}}
V^{\chi}\!} \,\bigr)^{\!\infty} = \{0\}$,
ce qui montre que les espaces vectoriels $V^{\chi}$ avec $\chi \in P$ sont en
somme directe.
\cqfd

}

\definition{

On appelle << $G_0$-module de longueur finie admissible >> (ou << $G_0$-module
de Harish-Chandra >>) tout espace de Hilbert complexe $V\!$ muni d'une
repr\'esentation lin\'eaire continue $\pi$ de $G_0$ en son sein qui v\'erifie
les propri\'et\'es \'equivalentes du lemme pr\'ec\'edent.

}

\medskip
L'objet de la proposition suivante est d'adapter le foncteur de translation de
Zuckerman (cf. \cite {Zu77}) au cas d'une repr\'esentation d'un groupe non
connexe.

\proposition{
\label{foncteur de translation}

On fixe un groupe de Lie r\'eel $\underline{G}$ dont la composante neutre est
\'egale \`a $G_0$ et un \'el\'ement $\underline{a}$ de $\underline{G}$ qui
op\`ere par automorphisme int\'erieur sur $G_0$ comme un certain \'el\'ement $a$
de $G$.
On suppose pour simplifier que le groupe $\underline{G}$ est engendr\'e
par $\underline{a} G_0$.

Soit $\pi$ une repr\'esentation lin\'eaire continue de $\underline{G}$ dans un
espace de Hilbert complexe $V$ dont la restriction \`a $G_0$ est un $G_0$-module
de longueur finie admissible.

\smallskip
{\bf (a)}
On d\'efinit une fonction g\'en\'eralis\'ee $(\tr\pi)_{a G_0}$ sur $a G_0$,
invariante sous $\interieur G_0$, localement int\'egrable sur $a G_0$ et
analytique sur $a G_0 \cap \Gssreg$, en posant

\centerline{
$(\tr\pi)_{a G_0} (\varphi \, \diff _G)
= \tr\pi (\underline{\varphi} \, \diff _{\underline{G}}) \;\;$
pour tout $\varphi \in C^\infty_c(a G_0)$,
}

\noindent
o\`u $\; \underline{\varphi} \in C^\infty_c(\underline{a} G_0) \;$ est d\'efinie
par
$\; \underline{\varphi}(\underline{a}x_0) = \varphi(ax_0) \;$ pour $x_0 \in G_0$
et, $\diff _G$ et $\diff _{\underline{G}}$ sont des mesures de Haar sur $G$ et
$\underline{G}$ dont les restrictions \`a $G_0$ co\"{\i}ncident.

Elle s'\'ecrit :
$\qquad (\tr\pi)_{a G_0}
= \!\Sumpetit_{\chi \in P \textrm{ tel que } \underline{a} \, \chi = \chi}\!
(\tr\pi_{\dot{\chi}})_{a G_0}$,

\noindent
o\`u $\,P$ est l'ensemble des poids de $\centrealg{U \gothg_\C }$ dans
$V^{\infty}\!$ et $\pi_{\dot{\chi}}$ d\'esigne pour chaque $\chi \in P$ tel que
$\underline{a} \, \chi = \chi$ la repr\'esentation de $\underline{G}$ dans
$V^{\chi}$ issue de $\pi$.

\smallskip
{\bf (b)}
On suppose qu'il existe un caract\`ere $\chi$ de $\centrealg{U \gothg_\C }$ tel
que $V\!=V^{\chi}$.

Soit $x \in a G_0 \cap \Gssreg$.
Donc $\, \gothj_x \egdef \gothg(x) \,$ est commutative
(cf. \ref{notations pour les restrictions} (b)) et
$\, \gothj \egdef \centraalg{\gothg}{\gothj_x}$ est une sous-alg\`ebre
de Cartan de $\gothg$  (cf. \cite [lem. 1.4.1 p.$\!$~6] {Bo87}).
On note $\, \gothj(x)_1$\labind{j (x)1} la composante connexe de $0$ dans
$\{ \, Y \in \gothj(x) \mid x\exp Y \in \Gssreg \}$.

Il existe des \'el\'ements $p_l$ de $S(\gothj(x)_\C ^*)$ associ\'es aux
$l \!\in \gothj(x)_\C ^*$ v\'erifiant ${\chi \!=\! \chiinfty{l}{\gothg}}\!$
(notation de \ref{caracteres canoniques} (a)), tels que pour tout
$\, Y \in \gothj(x)_1$ on ait

\smallskip\centerline{
$\abs{D_G(x\exp Y)}^{1/2} \;\, (\tr\pi)_{a G_0}(x\exp Y) \;
= \Sum_{l \in \gothj(x)_\C ^* \textrm{ tel que }
\chi = \chiinfty{l}{\gothg}}\!\!
p_l(Y) \;\, \me^{l(Y)}$.
}

De plus, quand $\chi$ est << r\'egulier >> (c'est-\`a-dire associ\'e \`a une
orbite semi-simple r\'eguli\`ere de $\interieur \gothg_\C $ dans $\gothg_\C ^*$)
chaque $p_l$ est un scalaire.

\smallskip
{\bf (c)}
On suppose que $\pi$ a un caract\`ere infinit\'esimal r\'egulier $\chi$.
On choisit $\Lambda \in \gothg_\C ^*$ semi-simple r\'eguli\`ere telle que
$\, \chi = \chiinfty{\Lambda}{\gothg}\!$.
On note $C$ l'unique chambre de $R(\gothg_\C ,\gothg_\C (\Lambda))$
dans $\,\gothg_\C (\Lambda)^* \cap \derivealg{\gothg_\C }^*$ telle que
$\,\Lambda \in C + \centrealg{\gothg_\C }^*\!$,
et $\,\overline{C}$ son adh\'erence.

Soit $\pi_F$ une repr\'esentation lin\'eaire de $\underline{G}$ dans un espace
de Hilbert complexe de dimension finie $F\!$, dont la restriction \`a $G_0$ est
irr\'eductible, et telle que le plus bas poids de $\gothg_\C (\Lambda)$ dans
$F$ relativement \`a $C\!$, not\'e $\Lambda_F\!$, v\'erifie
$\Lambda + \Lambda_F \in \overline{C} + \centrealg{\gothg_\C }^*\!$.

La restriction \`a $G_0$ de la repr\'esentation $\pi \otimes \pi_F$ de
$\underline{G}$ donne un $G_0$-module de longueur finie admissible et
$\underline{G}$ op\`ere dans sa composante $\centrealg{U \gothg_\C }$-primaire
$(V \otimes_\C  F)^{\chiinfty{_{\Lambda + \Lambda_F}}{\gothg}}\!$
au moyen d'une repr\'esentation $\pi_{Zuc}$ ayant encore cette propri\'et\'e.

Soit $x \in a G_0 \cap \Gssreg$.
On lui associe une sous-alg\`ebre de Cartan $\gothj$ de $\gothg$ et des nombres
complexes $p_l$ comme au (b).
Pour tout $Y \in \gothj(x)_1$, on a

\smallskip\centerline{
$\abs{D_G(x\exp Y)}^{1/2} \;\; (\tr\pi_{Zuc})_{a G_0}(x\exp Y) \;
= \Sum_{l \,\in\, \gothj(x)_\C ^* \textrm{ tel que }
\chi = \chiinfty{l}{\gothg}}\!\!
\underline{x}_{~\!l_F} \;\, p_l \;\, \me^{(l \,+\, l_F)(Y)}$,
}

\noindent
o\`u $\; l_F \in \gothj_\C ^*$ est d\'efini \`a partir de $l$ comme
$\,\Lambda_F$ \`a partir de $\Lambda$ et $\,\underline{x}_{~\!l_F}$ est le
scalaire par lequel l'\'el\'ement
$\;\underline{x} \egdef \underline{a}\,(a^{-1}\,x) \;$
de $\underline{G}$ (qui fixe $\,l_F$) agit sur l'espace propre (de dimension
$1$) de poids $\,l_F$ pour l'action de $\, \gothj(x)_\C $ dans $F\!$.

}

\dem{D\'emonstration de la proposition}{

On utilisera l'action canonique par op\'erateurs diff\'erentiels invariants
\`a gau\-che de l'alg\`ebre enveloppante d'un groupe de Lie r\'eel sur
l'espace vectoriel des fonctions g\'en\'eralis\'ees sur ce groupe de Lie.

\smallskip
{\bf (a)}
La repr\'esentation $\pi$ est tra\c{c}able d'apr\`es le lemme
\ref{composantes primaires}.
La fonction g\'en\'eralis\'ee $(\tr\pi)_{a G_0}$, \'egale \`a
$\delta_a \ast (\delta_{\underline{a}^{-1}} \ast
\restriction{\tr\pi}{\underline{a} G_0})$,
est invariante sous $\,\interieur G_0$.

\smallskip
Soient
$\, u_0 \in \centrealg{U \gothg_\C } \,$
et
$\, \underline{\varphi} \in C^\infty_c(\underline{a} G_0)$.
On pose
$\; u = \Prodpetit_{\chi \in P} (u_0 - \chi(u_0))$.
La d\'emonstration de \cite [cor. 7.207 p.$\!$~530] {KV95} fournit
l'existence d'une suite finie
$\; W_0=V \supseteq W_1 \supseteq \cdots \supseteq W_N=\{0\} \;$
de sous-$\underline{G}$-modules ferm\'es de $V$ tels que les
quotients $W_{k-1}/W_k$ avec $1 \leq k \leq N' \,$ sont topologiquement
irr\'eductibles sous $\underline{G}$.

Soit $k \!\in\! \{1,\dots, N'\}$.
On note $\pi_k$ la repr\'esentation de $\underline{G}$ dans $W_{k-1}/W_k$.
Tout sous-quotient d'un $G_0$-module de longueur finie admissible est un
$G_0$-module de longueur finie admissible.
Donc $W_{k-1}/W_k$ a un sous-$G_0$-modules ferm\'e topologiquement
irr\'eductible $E$ qui poss\`ede un caract\`ere infinit\'esimal.
L'action de $u$ sur $\; \Sumpetit_{n \in \Z }\underline{a}^n \, E^\infty \,$ est
nulle.
Donc $\; \pi_k (\underline{\varphi} \diff _{\underline{G}} \ast u) = 0$.

Ainsi :
$\;\, \partial_u (\tr\pi) \,
(\underline{\varphi} \diff _{\underline{G}})
= \tr\pi_1 (\underline{\varphi} \diff _{\underline{G}} \ast u) + \cdots
+ \tr\pi_{N'} (\underline{\varphi} \diff _{\underline{G}} \ast u)
= 0 \;\,$
(cf. \cite [lem. 8.1.3 p.$\!$~293] {Wa88}).
De plus, on a $\;\, \underline{a} \, u = u$.
Il en r\'esulte que :

\centerline{
$\partial_u ((\tr\pi)_{a G_0})
= \delta_a \ast ((\delta_{\underline{a}^{-1}} \ast
\restriction{\tr\pi}{\underline{a} G_0}) \ast
(\delta_{\underline{a}} \ast \check{u} \ast \delta_{\underline{a}^{-1}}))
= 0$.
}

\noindent
D'apr\`es \cite [th. 2.1.1 p.$\!$~10] {Bo87} et
\cite [th. 4.95 p.$\!$~286 et th. 7.30 (a) p.$\!$~450] {KV95},
$(\tr\pi)_{a G_0}$ est localement int\'egrable sur $a G_0$ et analytique sur
$a G_0 \cap \Gssreg$.

\smallskip
Soient $\, \chi \in P$ tel que $\underline{a} \, \chi \not= \chi$ et (\`a
nouveau) $\, \underline{\varphi} \in C^\infty_c(\underline{a} G_0)$.
On note $m$ le cardinal de l'orbite $\dot{\chi}$ de $\chi$ sous l'action du
sous-groupe $\langle \underline{a} \rangle$ de $\underline{G}$
engendr\'e par $\underline{a}$, et $\pi_{\dot{\chi}}$ la repr\'esentation de
$\underline{G}$ dans
$\, W \egdef \overline{\Sumpetit_{n \in \Z } \underline{a}^n \, V^{\chi}\!}$.
La restriction de
$\, \pi_{\dot{\chi}} (\underline{\varphi}\diff _{\underline{G}}) \,$
\`a la somme directe
$\, V^{\chi} \oplus \underline{a}V^{\chi} \oplus \cdots \oplus
\underline{a}^{m-1}V^{\chi} \,$
se d\'ecompose en blocs avec des blocs diagonaux nuls.
J'exploite cette propri\'et\'e en utilisant une suggestion de G.~Skandalis pour
me ramener au cas de la dimension finie.

La sous-alg\`ebre de $\call (W)$ form\'ee des endomorphismes qui stabilisent
simultan\'ement $V^{\chi}$,
$\underline{a}V^{\chi},\dots,\underline{a}^{m-1}V^{\chi}$ est ferm\'ee.
D'apr\`es \cite [remarque p.$\!$~55 et prop. 7 p.$\!$~47] {B67}, les
sous-espaces primaires $W^z$ de la restriction de
$\, \pi_{\dot{\chi}} (\underline{\varphi} \diff _{\underline{G}})^m \,$
\`a $W$ associ\'es aux $\, z \in \C \moins0$, qui
sont de dimension finie, sont tous sommes de leurs intersections avec
$V^{\chi}$ et $\cdots$ et $\underline{a}^{m-1}V^{\chi}$.
En outre, chaque $W^z$ avec $\, z \in \C \moins0$ est \'egal \`a la somme
des sous-espaces primaires de la restriction de
$\, \pi_{\dot{\chi}} (\underline{\varphi} \diff _{\underline{G}}) \,$ \`a
$W$ associ\'es aux racines $m^{\textrm{i\`emes}}\!$ de $z$.
La formule de Lidskij (<< la trace est la somme des valeurs propres >>) donne
ensuite l'\'egalit\'e
$\;\, \tr\pi_{\dot{\chi}} (\underline{\varphi} \, \diff _{\underline{G}})
= \Sumpetit_{z \in \C \moins0}
\tr \big(
\restriction{\pi_{\dot{\chi}} (\underline{\varphi} \,
\diff _{\underline{G}})} {W^z} \big)
= 0$.

D'o\`u, en utilisant cette fois-ci \cite {B67} relativement \`a $\call (V)$ (cf.
dem. de \ref{composantes primaires}) :

\centerline{
$(\tr\pi)_{a G_0}
=
\Sum_{\dot{\chi} \in \langle \underline{a} \rangle \backslash P}
(\tr\pi_{\dot{\chi}})_{a G_0}
= \Sum_{\chi \in P \textrm{ tel que } \underline{a} \, \chi = \chi}
(\tr\pi_{\dot{\chi}})_{a G_0}$.
}

\smallskip
{\bf (b)}
Je vais pr\'eciser ici les calculs de \cite [p.$\!$~33 et p.$\!$~34] {Bo87}.

\smallskip
On pose
$\;\, \theta(Y) = \abs{D_G(x\exp Y)}^{1/2} \; (\tr\pi)_{a G_0}(x\exp Y) \;\,$
pour tout $\, Y \!\in \gothj(x)_1$. \\
On note $W$ le groupe $W(\gothg_\C ,\gothj_\C )$ et $\pr$ la projection canonique
de $\, S\, \gothj_\C  \,$ sur $\, S\, \gothj(x)_\C  \,$ obtenue comme dans la
d\'efinition \ref{notations frequentes} (b).
La d\'emonstration du (a) montre que $(\tr\pi)_{a G_0}$ est vecteur propre de
$\centrealg{U \gothg_\C }$ associ\'e \`a $\chi$.
On fixe un \'el\'ement $l_0$ de $\gothj_\C ^*$ tel que
$\chi = \chiinfty{l_0}{\gothg}\!$.
D'apr\`es \cite [ligne 4 de 2.5 p.$\!$~20] {Bo87}, on a :

\smallskip\centerline{
$q(l_0) \; \theta = \partial_{\pr(q)} \theta \;\,$
pour tout $\, q \in (S\, \gothj_\C )^W\!$.
}

\medskip
On r\'ep\`ete dans ce qui suit les arguments de
\cite [p.$\!$~369 \`a p.$\!$~371] {Kn86}.

Soit $X \in \gothj_\C $.
On note $q_0,\dots,q_{\abs{W}-1}$ les \'el\'ements de $(S\, \gothj_\C )^W$ pour
lesquels on a l'\'egalit\'e suivante de polyn\^omes en l'ind\'etermin\'ee $T$ :

\centerline{
$\Prod_{w \in W} (T-wX)
= \; T^{\abs{W}} \,+\, q_{\abs{W}-1} T^{\abs{W}-1} \,+\, \cdots \,+\, q_0\;\;$
dans $\, (S\, \gothj_\C ) [T]$.
}

\noindent
\`A partir de l\`a, d'une part en rempla\c{c}ant $T$ par $X$, et d'autre
part en prenant la valeur de chaque membre en $l_0$ puis rempla\c{c}ant $T$ par
$X$, on obtient dans $\, S\, \gothj_\C $ :

\smallskip\centerline{
$X^{\abs{W}} \,+\, X^{\abs{W}-1} q_{\abs{W}-1} \,+\, \cdots \,+\, q_0 = 0$
}

\centerline{
et
$\;\; X^{\abs{W}} \,+\, X^{\abs{W}-1} q_{\abs{W}-1}(l_0) \,+\, \cdots \,+\,
q_0(l_0)
= \Prod_{w \in W} (X-l_0(wX))$.
}

\noindent
On utilise l'\'equation aux d\'eriv\'ees partielles v\'erifi\'ee ci-dessus par
$\theta$ avec $q$ successivement \'egal \`a $\, q_{\abs{W}-1},\dots,q_0$.
Elle donne l'\'egalit\'e suivante, qu'il reste \`a exploiter :

\smallskip\centerline{
$\Prodpetit_{w \in W} (\partial_{\pr(X)} - l_0(wX)) \cdot \theta = 0$.
}

On fait d\'ecrire \`a $X$ une base de $\gothj(x)_\C $.
D'apr\`s \cite [prop. 3 p.$\!$~58] {Va77}, il existe une unique famille
$(p_l)_{l \in \gothj(x)_\C ^*}$ d'\'el\'ements de $\, S(\gothj(x)_\C ^*)$ \`a
support fini telle que

\centerline{
$\theta(Y) \;= \Sum_{l \in \gothj(x)_\C ^*} p_l(Y) \; \me^{l(Y)} \;\,$
pour tout $\, Y \!\in \gothj(x)_1$.
}

\noindent
On va obtenir des renseignements plus pr\'ecis en substituant cette
expression de $\theta$ dans l'\'egalit\'e qui pr\'ec\'edait (et prolongeant
analytiquement \`a $\gothj(x)$).

\smallskip
Soit $l \in \gothj(x)_\C ^*$ tel que $p_l \not= 0$.
On note $d$ le degr\'e de $p_l$ dans $S(\gothj(x)_\C ^*)$ et $p_l^{[d]}$
la composante homog\`ene de $p_l$ de degr\'e $d$.
Soit $X \in \gothj_\C $.
Pour tout $w \in W$, on a

\centerline{
$(\partial_{\pr(X)} - l_0(wX)) \cdot p_l \; \me^l \,
= \; (l(X) - l_0(wX)) \, p_l \; \me^l + (\partial_{\pr(X)} p_l) \; \me^l$.
}

\noindent
Le polyn\^ome nul
$\; \me^{-l} {\scriptstyle \times}\!
\Prodpetit_{w \in W} (\partial_{\pr(X)} - l_0(wX)) \cdot p_l \; \me^l \;$
admet
$\; \Prodpetit_{w \in W} (l - w^{-1} l_0) (X) \,{\scriptstyle \times}\,
p_l^{[d]}$
comme composante homog\`ene de degr\'e $d$.
Vu que l'anneau $\, S\, \gothj_\C ^*$ est int\`egre, cela montre que
$\; l \in W \!\cdot l_0$, puis
$\; \chi = \chiinfty{l}{\gothg}\!$.

On suppose maintenant $p_l$ non scalaire.
On choisit $X \in \gothj_\C $ tel que, pour tout $\, w \in W \,$ hors de $W(l)$
on ait $\, (w\,l)(X) \not= l(X)$, et
$\; p_l^{[d]} (\pr(X)) \not= 0$.
On note $\,(\pr(X)^*)\,$ la base de $\,(\C \pr(X))^*\,$ duale de $(\pr(X))$.
La restriction du polyn\^ome~nul

\centerline{
$\me^{-l} \,{\scriptstyle \times} \!\!
\Prodpetit_{w \in W(l)} \!\!
(\partial_{\pr(X)} - l(wX)) \cdot
\bigl( \Prodpetit_{w \in W \textrm{ et } w \notin W(l)} \!\!
(\partial_{\pr(X)} - l(wX))
\cdot p_l \; \me^l \bigr)$
}

\negsmallskip\noindent
\`a $\,\C \pr(X)\,$ est somme d'un polyn\^ome de la forme
$\; \zeta \,{\scriptstyle\times}\,
(\partial_{\pr(X)})^{\abs{W(l)}} \cdot (\pr(X)^*)^d \;$
avec $\zeta \in \C \moins0$ et de mon\^omes de degr\'es strictement plus petits.
D'o\`u $\, \abs{W(l)} > 1$.

\smallskip
{\bf (c)}
Je m'inspire de \cite [p.$\!$~301] {Zu77}.

\smallskip
D'apr\`es \cite [cor. 7.207 p.$\!$~530] {KV95}, $\pi \otimes \pi_F$ et les
repr\'esentations de $G_0$ sur les composantes
$\centrealg{U \gothg_\C }$-primaires de $V \otimes_\C  F$ donnent des
$G_0$-modules de longueur finie admissibles.

Par hypoth\`ese, on a
$\, \underline{a} \, \chi = \chi$.
Il existe donc
$\, \sigma \in (\Ad\underline{a})^\C  \, \interieur \gothg_\C  \,$
tel que $\sigma(\Lambda) = \Lambda$.
On a : $\, \sigma(\gothg_\C (\Lambda)) = \gothg_\C (\Lambda)$ et
$\sigma(C) = C$, puis $\sigma(\Lambda_F) = \Lambda_F$.
On en d\'eduit que
$\; \underline{a} \, \chiinfty{_{\Lambda + \Lambda_F}}{\gothg}
= \sigma \, \chiinfty{_{\Lambda + \Lambda_F}}{\gothg}
= \chiinfty{_{\Lambda + \Lambda_F}}{\gothg}\!$.
Le groupe $\underline{G}$ laisse donc stable
$(V \otimes_\C  F)^{\chiinfty{_{\Lambda + \Lambda_F}}{\gothg}}\!$.

On note $P'$ l'ensemble des poids de $\centrealg{U \gothg_\C }$ dans
$(V \otimes_\C  F)^{\infty}$, $P(F,\gothj(x)_\C )$ l'ensemble des poids de
$\gothj(x)_\C $ dans $F$ et $F^\gamma$ l'espace propre de $\gothj(x)_\C $
associ\'e \`a un \'el\'ement $\gamma$ de $P(F,\gothj(x)_\C )$.
Pour tout $Y \in \gothj(x)_1$, on a gr\^ace au (a) :

\smallskip\centerline{
$(\tr\pi)_{a G_0}(x\exp Y) \,{\scriptstyle \times}\,
\tr\pi_F(\underline{x}\exp Y) \,
= \Sum_{\chi' \in P' \textrm{ tel que } \underline{a} \, \chi' = \chi'}
(\tr(\pi \otimes \pi_F)_{\dot{\chi}'})_{a G_0}(x\exp Y)$
}

\centerline{
et
$\;\; \tr\pi_F(\underline{x}\exp Y) \;
= \!\Sum_{\gamma \in P(F,\gothj(x)_\C )}\!\!
\tr(\restriction{\pi_F(\underline{x})}{F^\gamma}) \; \me^{\gamma(Y)}$.
}

\noindent
On applique le (b) \`a $\pi$, et aux $(\pi \otimes \pi_F)_{\dot{\chi}'}$ avec
$\chi' \!\in\! P'$ tel que $\underline{a} \, \chi'\! = \!\chi'$, parmi lesquels
se trouve $\pi_{Zuc}$.
En identifiant les facteurs polyn\^omiaux \'ecrits devant les fonctions
$\,\me^l\,$ avec $l \in \gothj(x)_\C ^*$, on trouve que pour tout
$Y \in \gothj(x)_1$ on a

\smallskip\centerline{
$\abs{D_G(x\exp Y)}^{1/2} \;\; (\tr\pi_{Zuc})_{a G_0}(x\exp Y)
= \Sum_{(l,\gamma) \,\in\, \cala}
\tr(\restriction{\pi_F(\underline{x})}{F^\gamma}) \;\,
p_l \;\, \me^{(l \,+\, \gamma)(Y)}$,
}

\noindent
o\`u on a pos\'e
$\;\, \cala
= \{ (l,\gamma) \in \gothj(x)_\C ^* \times P(F,\gothj(x)_\C ) \mid
\chi = \chiinfty{l}{\gothg} \textrm{ et }
\chiinfty{_{\Lambda + \Lambda_F}}{\gothg} = \chiinfty{l + \gamma}{\gothg} \}$. 

\smallskip
Pour finir, on v\'erifie que
$\; \cala
= \{ (l,\gamma) \in \gothj(x)_\C ^* \times \gothj_\C ^* \mid
\chi = \chiinfty{l}{\gothg} \textrm{ et } \gamma = l_F \}$,
et $\dim F^\gamma = 1$ pour $(l,\gamma) \in \cala$.
L'inclusion $\,\supseteq\,$ est claire.
On consid\`ere un $(l,\gamma) \in \cala$.
On modifie le choix de $\Lambda$ en prenant $\, \Lambda = l$.
On va montrer que $\, \gamma = \Lambda_F$, et calculer $\dim F^\gamma\!$.
On fixe un $\, w \in W(\gothg_\C ,\gothj_\C ) \,$ tel que
$\, l + \gamma = w (\Lambda + \Lambda_F)$.
\\
On a
$\; \gamma - \Lambda_F = w (\Lambda + \Lambda_F) - (\Lambda + \Lambda_F) \,$ et
$\, w (\Lambda + \Lambda_F) - (\Lambda + \Lambda_F) \leq 0 \,$
pour l'ordre associ\'e \`a $C$, car
$\, \Lambda + \Lambda_F \in \overline{C} + \centrealg{\gothg_\C }^* \!$.
D'autre part, il existe un poids $\zeta$ de $\gothj_\C $ dans $F$ tel que
$\; \gamma = \restriction{\zeta}{\gothj(x)_\C }$.
L'inclusion $\, \gothj(x)_\C ^* \subseteq \gothj_\C ^* \;$ identifie
$\gamma$ \`a l'isobarycentre de la partie finie
$\, \{ \Ad^*x^n \cdot \zeta \}_{n \in \Z } \,$
de l'ensemble des poids de $\gothj_\C $ dans $F\!$.
Donc $\, \gamma - \Lambda_F \geq 0 \,$ pour l'ordre associ\'e \`a $C$.
Par cons\'equent $\, \gamma = \Lambda_F$.
Comme $\gamma$ est un poids extr\^emal de $\gothj_\C $ dans $F$ qui est
isobarycentre d'un ensemble fini de poids de $\gothj_\C $ dans $F\!$, il
est \'egal en particulier \`a $\zeta$.
Ainsi le sous-espace propre de $F$ de poids $\gamma$ sous l'action
de $\gothj(x)_\C$ est de dimension $1$.
\cqfd

}

\section{\boldmath
Passage pour $M'$ au cas r\'egulier
}\label{Rep13}

Pour simplifier les notations, on pose dans cette section :
$\M '=\G (\nu_+)$,
$\rho = \rho_{\gothm',\gothh}$, $\sigma = \tau_{M'}$,
et $\lambda_r = \lambda_{+,\gothm'}$.
On a  donc $\lambda_r  = (\mu-2\mi \rho) + \nu \,\in\, \mpstssreg$.
On va voir qu'en gros, le caract\`ere
$\tr T_{\lambdatilde,\sigma}^{M'}$ se d\'eduit d'un caract\`ere de la forme
$\tr T_{\lambda_r,\sigma_r}^{M'}$ par translation au sens de Zuckerman.%

\smallskip
La proposition suivante va s'obtenir en prouvant une variante des r\'esultats de
Bouaziz dans \cite [(i) et (ii) p.$\!$~550]{Bo84}
(cf. aussi \cite [p.$\!$~547 et p.$\!$~548]{KV95}).
Elle a pour origine le lemme 3.1 p.$\!$~406 de \cite {KZ82}, qui fournit
l'\'egalit\'e des restrictions \`a $M'_0$ des caract\`eres consid\'er\'es
au (b) ci-dessous.

\proposition{
\label{passage au cas regulier}

On fixe $\,\hat{a} \in M'(\lambdatilde)^{\gothm'/\gothh}$ au-dessus d'un
\'el\'ement elliptique $\,a$ de $M'(\lambdatilde)$.
On note ${\underline{M'\!}\,}$ le produit semi-direct de $\Z $ par $\,M'_0$ tel
que $\, \underline{a} \egdef 1 \in \Z $ agisse sur $\, M'_0$ par $\interieur a$.

\smallskip
{\bf (a)}
Il existe une repr\'esentation lin\'eaire continue $\pi_F$ de
${\underline{M'\!}\,}$ dans un espace de Hilbert complexe de dimension finie
$F$, unique \`a isomorphisme pr\`es, dont la restriction \`a $M'_0$ est
irr\'eductible avec pour plus bas poids $-2\rho$ pour
l'action de $\gothh_\C $  et relativement \`a l'ordre d\'eduit de
$R^+(\gothm'_\C ,\gothh_\C )$, telle que $\underline{a}$ agisse trivialement sur
l'espace propre de poids $-2\rho$.

\smallskip
{\bf (b)}
On prolonge les repr\'esentations unitaires $T_{\lambdatilde}^{M'_0}\!$ et
$T_{\lambda_r}^{M'_0}\!$ de $M'_0$ << d'espaces >> not\'es $\calh$ et $V\!$, en
des repr\'esentations lin\'eaires continues $\pi_{\lambdatilde}$ et
$\pi_{\lambda_r}\!$ de ${\underline{M'\!}\,}$ par les conditions
$\, \pi_{\lambdatilde}(\underline{a}) = S_{\lambdatilde}(\hat{a}) \,$
et
$\, \pi_{\lambda_r}(\underline{a}) = S_{\lambda_r}(\hat{a})$,
o\`u $\, S_{\lambdatilde} \,$ et $\, S_{\lambda_r}$ sont les repr\'esentations
de $M'(\lambdatilde)^{\gothm'/\gothh}$ (qui est \'egal \`a
$M'(\lambda_r)^{\gothm'/\gothh}$) attach\'ees \`a $\,\lambdatilde$ et
$\lambda_r$ dans la proposition \ref{construction de rep, cas de $M'$} (b).

Les repr\'esentations de ${\underline{M'\!}\,}$ dans $\calh$ et
$(V \otimes_\C F)^{\chiinfty{\mi \,\lambda}{\gothm'}}\!$
sont tra\c{c}ables et ont m\^eme caract\`ere.

}

\dem{D\'emonstration de la proposition}{

{\bf (a)}
Il existe --- et on se donne --- une repr\'esentation lin\'eaire irr\'eductible
de $\gothm'_\C $ dans un espace de Hilbert complexe $F$ de dimension finie, de
plus bas poids $-2\rho$ pour l'action de $\gothh_\C $
relativement \`a l'ordre d\'eduit de $R^+(\gothm'_\C ,\gothh_\C )$.
Elle s'int\`egre en une repr\'esentation << du >> rev\^etement universel de
$M'_0$ avec un caract\`ere central trivial, puis en une repr\'esentation de
$M'_0$ not\'ee $\pi_{F,0}$.
Comme $\underline{a}$ normalise $\gothh$ et fixe $-2\rho$, il
existe un unique op\'erateur d'entrelacement $\Phi$ de
$(F,\interieur\underline{a} \cdot \pi_{F,0})$ sur $(F,\pi_{F,0})$ qui agit
trivialement sur le sous-espace propre de $F$ de poids
$-2\rho$ sous l'action de $\gothh_\C $, commun aux deux
repr\'esentations et de dimension~$1$.
On construit $\pi_F$ en prolongeant $\pi_{F,0}$ et choisissant
$\pi_F(\underline{a})$ \'egal \`a $\Phi$.

\smallskip
{\bf (b)}
On remarque d'abord que ${\underline{M'\!}\,}$ stabilise
$\, \calh_0 \egdef (V \otimes_\C  F)^{\chiinfty{\mi \,\lambda}{\gothm'}}\!$.
D'apr\`es le lemme \ref{composantes primaires} et la proposition
\ref{foncteur de translation} (c), les repr\'esentations de
${\underline{M'\!}\,}$ dans $\calh$ et $\calh_0$ sont tra\c{c}ables.

On fixe un sous-groupe compact maximal $K_{M'_0\!}$ de $M'_0$ stable par
$\interieur a$ et dont l'involution de Cartan normalise $\gothh$.
Par exemple $\exp (\gothc_{M'} \cap \gothm')$, o\`u $\gothc_{M'}$ est
l'alg\`ebre de Lie d'un sous-groupe compact maximal du produit semi-direct de
$\{1,c\}$ ($c$ : conjugaison de $\M '(\C )$) par $\M '(\C )$ qui contient \`a la
fois $c$ et le sous-groupe engendr\'e par le produit de la projection dans
$\M '(\R )$ de $a$ avec $\exp_{\M '(\C )} (\mi \,\gothh_{(\R )})$.
On note $\calh^f\!$, $\calh_0^f$ et $V^f$ les $(\gothm'_\C ,K_{M'_0})$-modules
stables par $\underline{a}$ associ\'es \`a $\calh$, $\calh_0$ et~$V\!$.
D'apr\`es \cite [prop. 10.5 p.$\!$~336]{Kn86} avec sa d\'emonstration (cf.
\cite [cor. 8.8 p.$\!$~211]{Kn86}), les repr\'esentations de
${\underline{M'\!}\,}$ dans $\calh$ et $\calh_0$ ont m\^eme caract\`ere s'il
existe un isomorphisme de $(\gothm'_\C ,K_{M'_0})$-modules de $\calh^f$ sur
$\calh_0^f$ compatible \`a l'action de~$\underline{a}$.

\smallskip
On a vu au cours de la d\'emonstration de \ref{construction de rep, cas de $M'$}
(a) que $\calh^f$ et $V^f$ sont irr\'eductibles, isomorphes \`a
$\, \calr_{M'_0}^{q'} (\C_{\,\mi \,\lambda-\rho})\,$
et
$\, \calr_{M'_0}^{q'} (\C_{\,\mi \,\lambda_r-\rho}) \,$
o\`u $\, \calr_{M'_0}^{q'}$ est le foncteur d'induction cohomologique
relatif \`a $\gothb_{M'}$ et $\C_{\Lambda}$ est pour chaque
$\,\Lambda \in \mi \,\gothh^*\,$ le $(\gothh_\C ,T_0)$-module $\C$ sur
lequel $\gothh$ agit par $\Lambda$.
D'apr\`es \cite [th. 7.237 p.$\!$~544] {KV95} et en utilisant la notation
$\,\psi_{l}^{l'}\,$ de \cite [(7.141) p.$\!$~493] {KV95},
on a

\centerline{
$\calr_{M'_0}^{q'} (\psi_{\,\mi \,\lambda_r-\rho}^{\,\mi \,\lambda-\rho}
(\C_{\,\mi \,\lambda_r-\rho}))
\simeq
\psi_{\,\mi \,\lambda_r}^{\,\mi \,\lambda} (\calr_{M'_0}^{q'}
(\C_{\,\mi \,\lambda_r-\rho}))$.
}

\noindent
Ainsi, les $(\gothm'_\C ,K_{M'_0})$-modules $\calh^f$ et $\calh_0^f$
sont isomorphes.

Compte tenu de la caract\'erisation de $\, S_{\lambdatilde}(\hat{a}) \,$
dans \ref{construction de rep, cas de $M'$} (b), du lemme de Schur de
\cite [lem. 3.3.2 p.$\!$~80] {Wa88} et de l'identification
$\, H_{q'} (\gothn_{M'},\calh^{\infty})^* = H_{q'} (\gothn_{M'},\calh^f)^* \,$
de \cite [lem. 4 p.$\!$~165] {Df82a}, pour obtenir la compatibilit\'e \`a
l'action de $\,\underline{a}$ d'un isomorphisme de
$(\gothm'_\C ,K_{M'_0})$-modules fix\'e de $\calh^f$ sur $\calh_0^f$, il suffit
de montrer que $\,\underline{a}\,$ agit \`a partir de
$\, \pi_{\lambda_r}(\underline{a}) \otimes \pi_F(\underline{a}) \,$ dans
$\, (H_{q'} (\gothn_{M'},\calh_0^f)^*)_{-(\mi \,\lambda + \rho)}$
par $\rho_\lambdatilde^{\gothm'\!/\gothh} (\hat{a})^{-1} \id$,
c'est-\`a-dire par $\rho_{\lambda_r}^{\gothm'\!/\gothh} (\hat{a})^{-1} \id$.

L'id\'ee de Bouaziz consiste en gros \`a v\'erifier la compatibilit\'e \`a
l'action de $\, \underline{a}$ de l'isomorphisme qu'on obtiendrait au niveau de
l'homologie en prenant $\gothq$ \'egal au conjugu\'e de $\gothb_{M'}$ et $X$
\'egal \`a $V^f$ dans \cite [th. 7.242 p.$\!$~546] {KV95}.
Malheureusement, l'hypoth\`ese << $\mi \,\lambda_r$ est au moins autant
singulier que $\mi \,\lambda$ >> n'est pas satisfaite.
Au lieu de passer de $V^f$ \`a $\calh_0^f$, on va passer de $\calh_0^f$ \`a
$V^f$ en utilisant une propri\'et\'e de la compos\'ee de deux foncteurs de
translation de Zuckerman.

\smallskip
D'apr\`es \cite [th. 7.220 p.$\!$~536] {KV95}, le morphisme de
$(\gothm'_\C ,K_{M'_0})$-modules compatible \`a l'action de $\,\underline{a}$
de $(\calh_0^f \otimes_\C  F^*)^{\chiinfty{\mi \,\lambda_r}{\gothm'}}$
dans $V^f\!$ obtenu par restriction \`a partir de l'application canonique de
$V^f\! \otimes_\C  F \otimes_\C  F^*$ dans $V^f\!$, est non nul donc surjectif.
Les composantes $\centrealg{U \gothh_\C }$-primaires associ\'ees \`a
$\chiinfty{\mi \,\lambda_r + \rho}{\gothh}\!$
dans la suite exacte de $\gothh_\C $-modules d'homologie issue de cette
surjection fournissent un morphisme de $\gothh_\C $-modules de
$\; H_{q'}
(\gothn_{M'},(\calh_0^f \otimes_\C  F^*)^{\chiinfty{\mi \,\lambda_r}{\gothm'}})
^{\chiinfty{\mi \,\lambda_r + \rho}{\gothh}}$
dans
$\, H_{q'} (\gothn_{M'},V^f)^{\chiinfty{\mi \,\lambda_r + \rho}{\gothh}}$
compatible \`a l'action de $\,\underline{a}$, qui est surjectif d'apr\`es
\cite [lem. 8.9 p.$\!$~553] {KV95}.

Soit $v^*_0$ un vecteur non nul de $F^*$ de poids $2\rho$
pour l'action de $\gothh_\C $.
On a $\; \gothn_{M'} \cdot v^*_0 = \{0\}$ et
$\, \underline{a} \cdot v^*_0 = v^*_0 \,$ d'apr\`es (a).
Il reste \`a prouver que chacun des deux morphismes de $\gothh_\C $-modules
compatibles \`a l'action de $\,\underline{a}\,$ d\'eduits des injections
canoniques de
$(\calh_0^f \otimes_\C  F^*)^{\chiinfty{\mi \,\lambda_r}{\gothm'}}$
et
$\calh_0^f \otimes_\C  \C v^*_0$ dans $\calh_0^f \otimes_\C  F^*\!$,
d\'efinis respectivement sur
$\, H_{q'} (\gothn_{M'},(\calh_0^f \otimes_\C  F^*)
^{\chiinfty{\mi \,\lambda_r}{\gothm'}})
^{\chiinfty{\mi \,\lambda_r + \rho}{\gothh}}$
et
$\, H_{q'} (\gothn_{M'},\calh_0^f)
^{\chiinfty{_{\mi \,\lambda + \rho}}{\gothh}}
\otimes_\C  \C v^*_0$,
et tous deux \`a valeurs dans
$\, H_{q'} (\gothn_{M'},\calh_0^f \otimes_\C  F^*)
^{\chiinfty{\mi \,\lambda_r + \rho}{\gothh}}$,
est bijectif. \\
En effet, dans ce cas et avec la notation $\C_{\Lambda}$ introduite plus haut
dans cette d\'emonstration, la droite form\'ee des morphismes de
$\gothh_\C $-modules de $\, H_{q'} (\gothn_{M'},V^f) \,$ dans
$\C_{\mi \,\lambda_r + \rho}$ s'injectera de mani\`ere
compatible \`a l'action de $\,\underline{a}$ dans celle form\'ee des morphismes
de $\gothh_\C $-modules de $\, H_{q'} (\gothn_{M'},\calh_0^f) \,$ dans
$\, \C_{\mi \,\lambda + \rho}$.

\smallskip
Le premier des deux morphismes pr\'ec\'edents est bijectif d'apr\`es
\cite [(7.243) p.$\!$~547, cf. prop. 7.166 p.$\!$~506] {KV95}.
Pour le second, on pose $\, W = F^*/\,\C v^*_0$.

La suite exacte de $\gothb_{M'}$-modules
$0 \!\to\! \calh_0^f \otimes_\C  \C v^*_0 \!\to\! \calh_0^f
\otimes_\C  F^* \!\to\! \calh_0^f \otimes_\C  W \!\to\! 0$
fournit le morceau suivant de la suite exacte de $\gothh_\C $-modules
d'homologie :
\\
\parbox{\textwidth}{
$H_{q'+1} (\gothn_{M'},\calh_0^f \otimes_\C  W)
\to H_{q'} (\gothn_{M'},\calh_0^f) \otimes_\C  \C v^*_0$

$\hfill
\to H_{q'} (\gothn_{M'},\calh_0^f \otimes_\C  F^*)
\to H_{q'} (\gothn_{M'},\calh_0^f \otimes_\C  W)$.
}

\noindent
On termine la d\'emonstration en montrant que
$\; H_{\scriptscriptstyle\bullet}(\gothn_{M'},\calh_0^f \otimes_\C  W)
^{\chiinfty{\mi \,\lambda_r + \rho}{\gothh}}
= \{0\}$.

Par le th\'eor\`eme de Lie, il existe --- et on fixe --- une suite
$\; \{0\}
= W^{(-1)} \subseteq W^{(0)} \subseteq \cdots \subseteq W^{(N)}
= W \;$
de sous-$\gothb_{M'}$-modules de $W$ avec $\dim W^{(j)}/\,W^{(j-1)} = 1$
pour $0 \leq j \leq N$.
D'apr\`es \cite [prop. 3.12 p.$\!$~188] {KV95}, pour
tout
$(\gothb_{M'},T_0)$-module $E$ et tout $n \in \N $, le
$\gothh_\C $-module $H_n (\gothn_{M'},E)$ est isomorphe \`a $P_n(E)$, o\`u
$P_n$ est le $n^{\textrm{i\`eme}}\!$ foncteur d\'eriv\'e du foncteur exact
\`a droite de la cat\'egorie des $(\gothb_{M'},T_0)$-modules dans celle des
$(\gothh_\C ,T_0)$-modules not\'e
$\, P_{\gothb_{M'},T_0}^{\gothh_\C ,T_0}$ dans
\cite [(2.8) p.$\!$~104] {KV95}.
D'apr\`es \cite [prop. D.57 (b) p.$\!$~887, cf. (2.123) p.$\!$~162] {KV95}, il
existe --- et on fixe --- une filtration
$\, \{0\}
\!=\! C_{\scriptscriptstyle\bullet}^{(-1)}
\subseteq C_{\scriptscriptstyle\bullet}^{(0)}
\subseteq \cdots
\subseteq C_{\scriptscriptstyle\bullet}^{(N)}
\!=\! C_{\scriptscriptstyle\bullet} \,$
d'un complexe de chaines $C_{\scriptscriptstyle\bullet}$ en
$(\gothh_\C ,T_0)$-modules nul en degr\'e strictement n\'egatif, dont chaque
$\gothh_\C $-modules d'homologie $H_n (C_{\scriptscriptstyle\bullet})$ est
isomorphe \`a $H_n (\gothn_{M'},\calh_0^f\! \otimes_\C  W)$
pour $n \in \N $, et dont la suite spectrale $(E^r)_{r \geq 0}$ fournit des
$\gothh_\C $-modules $E^1_{p,q}$ isomorphes \`a
$H_{p+q} (\gothn_{M'},\calh_0^f) \otimes_\C  (W^{(p)}/\,W^{(p-1)}) \,$
pour $0 \leq p \leq N$ et $p+q \geq 0$.
D'apr\`es \cite [prop. 7.56 p.$\!$~460] {KV95}, les composantes
$\centrealg{U \gothh_\C }$-primaires de chacun de ces $\gothh_\C $-modules
$E^1_{p,q}$ sont associ\'ees \`a des caract\`eres de la forme
$\chiinfty{_{w\,\mi \,\lambda + \rho + \,\gamma}}
{\gothh}$
avec $w \in W(\gothm'_\C ,\gothh_\C )$, o\`u on a not\'e $\gamma$ le
poids de $\gothh_\C $ dans $\, W^{(p)}/\,W^{(p-1)}$.
Ainsi, au vu de \cite [prop. 7.166 p.$\!$~506] {KV95}, on a
$\; (E^1_{p,q})
^{\chiinfty{\mi \,\lambda_r + \rho}{\gothh}}
= \{0\}$.
Par ailleurs, pour $n \in \N $ et $q \in \Z $, le $\gothh_\C $-module
$E^{\infty}_{n-q,q}$ est un sous-quotient de $E^1_{n-q,q}$ isomorphe au quotient
de
$\, \pr _{H_n (C_{\scriptscriptstyle\bullet})}
(H_n (C_{\scriptscriptstyle\bullet}^{(n-q)})) \,$
par
$\, \pr _{H_n (C_{\scriptscriptstyle\bullet})}
(H_n (C_{\scriptscriptstyle\bullet}^{(n-q-1)}))$,
o\`u on a utilis\'e la m\^eme notation
$\pr_{H_n(C_{\scriptscriptstyle\bullet})}$
pour les applications canoniques \`a valeurs dans
$H_n (C_{\scriptscriptstyle\bullet})$.
D'o\`u pour chaque $n \in \N $ les relations suivantes entre
$\gothh_\C $-modules :
\\
\parbox{\textwidth}{
$H_n(\gothn_{M'},\calh_0^f \otimes_\C  W)
^{\chiinfty{\mi \,\lambda_r + \rho}{\gothh}}
\simeq
\pr _{H_n (C_{\scriptscriptstyle\bullet})}
(H_n (C_{\scriptscriptstyle\bullet}^{(N)}))
^{\chiinfty{\mi \,\lambda_r + \rho}{\gothh}}$

$\hfill
= \cdots
= \pr _{H_n (C_{\scriptscriptstyle\bullet})}
(H_n (C_{\scriptscriptstyle\bullet}^{(-1)}))
^{\chiinfty{\mi \,\lambda_r + \rho}{\gothh}}
= \{0\}$.
}
\\
Cela ach\`eve la d\'emonstration.
\cqfd

}

\remarque{

On conserve les notations de la d\'emonstration pr\'ec\'edente.
En particulier, $K_{M'_0\!}$ est stable par $\interieur a$.
On fixe un op\'erateur d'entrelacement $\Phi$ de
$(\calh,\interieur a \cdot T_{\lambdatilde}^{M'_0})$ sur
$(\calh,T_{\lambdatilde}^{M'_0})$ et
un isomorphisme $\Psi$ de $(\gothm'_\C ,K_{M'_0})$-modules de
$\calh^f$ sur $\calh_0^f$.
On note $\Phi^f$ la restriction de $\Phi$ \`a $\calh^f$ et
$\pi_0(\underline{a})^f$ l'op\'erateur par lequel $\underline{a}$ agit sur%
~$\calh_0^f$.
D'apr\`es le lemme de Schur, il existe $z\in\C \moins0$ tel que 
$\,\Phi^f = z\;\Psi^{-1}\!\circ\pi_0(\underline{a})^f\!\circ\Psi$.
La d\'emonstration de la proposition \ref{passage au cas regulier} (b) montre que
la repr\'esentation $S$ de la proposition \ref{construction de rep, cas de $M'$}
(b) v\'erifie $S(\hat{a}) = z^{-1}\,\Phi$.
Par ailleurs, $K_{M'_0\!}$ est inclus dans un sous-groupe compact maximal de $M'$
qui contient $a$, et les sous-groupes compacts maximaux de $M'$ sont conjugu\'es
sous $M'_0$ (cf. \cite [th. 3.1 p.$\!$~180 et lignes avant le th. 3.7 p.$\!$~186]
{Ho65}).
Il est facile d'en d\'eduire que $z$ est ind\'ependant du choix de $K_{M'_0}$.
On aurait donc pu d\'efinir $S$ en se ramenant au cas semi-simple r\'egulier
trait\'e par M.~Duflo (pour les valeurs sur les \'el\'ements elliptiques) et en
s'arrangeant pour que
$\,S(\exp X) = \me^{-\mi \lambda(X)} \, T_{\lambdatilde}^{M'_0}(\exp X)\,$ quand
$X\in\gothh$.
Il resterait \`a prouver \`a partir de cette nouvelle d\'efinition que $S$ est
une repr\'esentation (unitaire) de $M'(\lambdatilde)^{\gothm'/\gothh}$.
\cqfr

}

\smallskip
\dem{Fin de la d\'emonstration du th\'eor\`eme}{

Il reste \`a prouver le th\'eor\`eme pour $\, T_{\lambdatilde,\sigma}^{M'}\!$.
La situation de r\'ef\'erence est celle de \cite [lem. 5.4 p.$\!$~304] {Zu77},
dans laquelle $M'$ est remplac\'e par $M_0$ et $e$ est remplac\'e par un
\'el\'ement semi-simple r\'egulier de $M_0$.

On suppose que $\, e \in M'$.
On note $\calv_{\!e}^{M'}\!$, $d_e^{M'}$ et $D_e^{M'}$ les objets analogues \`a
$\calv_{\!e}$, $d_e$ et $D_e$ attach\'es \`a $M'$.
On choisit $\, \hat{a} \in M'(\lambdatilde)^{\gothm'/\gothh}$ au-dessus d'un
$\, a \in M'(\lambdatilde) \,$ elliptique tel que $\, e \in a\,M'_0$.
On lui associe ${\underline{M'\!}\,}$, $\pi_{\lambdatilde}$ et
$\pi_{\lambda_r}\!$ comme dans la derni\`ere proposition.
On va utiliser les notations de la proposition \ref{foncteur de translation} (a)
relatives au groupe ${\underline{M'\!}\,}$.
On introduit l'\'el\'ement $\sigma_r$ de $\Xirr_{M'}(\lambda_r)$ tel que :

\centerline{
$\sigma_r(\hat{u})
= \det(\Ad u^\C )_{\gothn_{M'}} \; \sigma(\hat{u}) \;\,$
pour $\, \hat{u} \in M'(\lambdatilde)^{\gothm'/\gothh}$.
}

\noindent
\begin{tabular}{@{}rr@{$\,=\;$}l@{}}
Il v\'erifie : &
$\restriction{\tr T_{\lambdatilde,\sigma}^{M'}}{e \, M'_0}$ &
$\tr \sigma (\hat{a}) \,{\scriptstyle \times}\,
(\tr\pi_{\lambdatilde})_{a M'_0}$ \\
et &
$\restriction{\tr T_{\lambda_r,\sigma_r}^{M'}}{e \, M'_0}$ &
$\tr \sigma_r (\hat{a}) \,{\scriptstyle \times}\,
(\tr\pi_{\lambda_r})_{a M'_0}$.
\end{tabular}

\noindent
Pour terminer la d\'emonstration, on remarque que compte tenu de la proposition
\ref{passage au cas regulier}, la proposition \ref{foncteur de translation} va
ensuite permettre de relier
$\bigl( \tr T_{\lambdatilde,\sigma}^{M'} \bigr)_{\!e}$ \`a
$\bigl( \tr T_{\lambda_r,\sigma_r}^{M'} \bigr)_{\!e}$.

\smallskip
On v\'erifie la formule du caract\`ere pour $\,\tr T_{\lambdatilde,\sigma}^{M'}$
(sous la forme propos\'ee dans le d\'ebut de la d\'emonstration du th\'eor\`eme)
sur l'ouvert inclus dans $\,\mpessreg$ de compl\'ementaire n\'egligeable et sur
lequel les deux membres de l'\'egalit\'e sont des fonctions analytiques, form\'e
des $\; X \in \calv_{\!e}^{M'} \,$ tels que $\; e \exp X \in \Mpssreg$.

Soit $\, X \in \calv_{\!e}^{M'}\!$ tel que $\, x \egdef e \exp X \in \Mpssreg$.
On pose $\, \gothj_x \egdef \gothm'(x)$.
Donc $\, \gothj_x = \gothm'(e)(X) \in \Car \gothm'(e) \,$ (vu les dimensions) et
$\, \gothj \egdef \centraalg{\gothm'}{\gothj_x} \in \Car \gothm'$.
On note $\Gamma$ la composante connexe de $X$ dans
$\, \gothj(x) \cap \mpessreg$.
Pour poursuivre le calcul, on consid\`ere un $Y \in \gothj(x)_1$
(cf. \ref{foncteur de translation} (b) pour $G=M'$)
tel que $X\!+\!Y \in \Gamma\cap \calv_{\!e}^{M'}\!$.

\smallskip
D'apr\`es \cite [th. 5.5.3 p.$\!$~52]{Bo87} et gr\^ace \`a
\ref{morphisme rho et fonction delta} (b), la formule du caract\`ere est acquise
pour $\, \tr T_{\lambda_r,\sigma_r}^{M'}\!$.
Par cons\'equent, on a :

\noindent\parbox{\textwidth}{
$\quad\; \tr \sigma_r(\hat{a})
\;{\scriptstyle \times}\;
\abs{D_{M'}(x\exp Y)}^{1/2} \;\,
(\tr\pi_{\lambda_r})_{a M'_0} (x\exp Y)$
\\
$\displaystyle
=\; |\det \;(1-\Ad e)_{_{\scriptstyle \gothm'/\gothm'(e)}}|^{1/2} \;\,
|\cald_{\gothm'(e)}(X\!+\!Y)|^{1/2} \;\,
\bigl( \tr T_{\lambda_r,\sigma_r}^{M'} \bigr)_{\!e}(X\!+\!Y)$
}
\\
\parbox{\textwidth}{
$\stackrel
{\raisebox{.5ex}{\boldmath$
{\scriptscriptstyle (}\scriptstyle\star{\scriptscriptstyle )}$}}
{=}\;
\mi^{-d_e^{M'}}\;(D_e^{M'})^{-\frac{1}{2}} \;\,
|\det \;(1-\Ad e)_{_{\scriptstyle \gothm'/\gothm'(e)}}|^{1/2}
\Sum_{\dot{\lambda'_r} \, \in \, M'(e)_0 \backslash
M'_0 \cdot \lambda_r \cap\, \gothm'(e)^*}$

\negsmallskip
$\hfill
{\textstyle \frac
{\scalo (\widehat{e'})_{(1-\Ad e')(\gothm'/\gothh)}}
{\scalo (B_{\lambda_r})_{(1-\Ad e')(\gothm'/\gothh)}}} \;
\tr\sigma_r(\widehat{e'}) \;\,
|\cald_{\gothm'(e)}(X\!+\!Y)|^{1/2} \;\,
{\widehat{\beta}}_{M'(e)_0 \cdot \lambda'_r}(X\!+\!Y)$
}

\noindent
o\`u
$\;\, g \in M'_0 \,$ v\'erifie $\; \lambda'_r = g\lambda_r \;$ et
$\; \widehat{e'} \in M'(\lambdatilde)^{\gothm'/\gothh} \;$ est au-dessus de
$\, e' = g^{-1}eg$.

\smallskip
On se donne $\,g \in M'_0\,$ comme ci-dessus.
On note $\pr$ l'application canonique de $M'$ dans $\M '(\C )$.
On associe \`a $g$ un \'el\'ement $y^g$ de $\M '(\C )(\pr(e))_0$
v\'erifiant $\gothj(x)_\C  = y^g \cdot (g \gothh)(e)_\C $.
On fixe un syst\`eme de racines positives $R^+(\gothm'(e)_\C ,\gothj(x)_\C )$ de
$R(\gothm'(e)_\C ,\gothj(x)_\C )$.
D'apr\`es la proposition \ref{limites de mesures de Liouville} (b), il existe
des nombres complexes $c^g_w$ ind\'ependants de $Y$ et index\'es par les
$\; w \in W(\gothm'(e)_\C ,\gothj(x)_\C ) \,$ tels qu'on ait

\centerline{
$\Prod_{\alpha \in R^+(\gothm'(e)_\C ,\gothj(x)_\C )}\!\! \alpha(X\!+\!Y)
\;{\scriptstyle \times}\;
{\widehat{\beta}}_{M'(e)_0 \cdot g\lambda_r}(X\!+\!Y) \,
\stackrel {\raisebox{.5ex}{\boldmath$
{\scriptscriptstyle (}\scriptstyle\star\star{\scriptscriptstyle )}$}}
{=}
\Sum_{w \in W(\gothm'(e)_\C ,\gothj(x)_\C )}\!
c^g_w \; \me^{\mi \,wy^gg\lambda_r(X+Y)}$
}

\noindent
et
$\hfill\Prod_{\alpha \in R^+(\gothm'(e)_\C ,\gothj(x)_\C )}\!\!\alpha(X)
\;{\scriptstyle \times}\;
{\widehat{\beta}}_{M'(e)_0 \cdot (g\lambdatilde)[e]}(X) \,
= \Sum_{w \in W(\gothm'(e)_\C ,\gothj(x)_\C )}\!
c^g_w \;\; \me^{\mi \,wy^gg\lambda\,(X)}\hfill$

\noindent
quand $(g\lambdatilde)[e] \in \mpestregtilde$, et
$\;\Sum_{w \in W(\gothm'(e)_\C ,\gothj(x)_\C )} \!\!
c^g_w \; \me^{\mi \,wy^gg\lambda\,(X)}
= 0\;$
si $(g\lambdatilde)[e] \notin \mpestregtilde$.%

En outre, la fonction
$\; |\cald_{\gothm'(e)}|^{1/2} \,{\scriptstyle \times}\,
\bigl( \Prodpetit_{\alpha \in R^+(\gothm'(e)_\C ,\gothj(x)_\C )}
\! \alpha \bigr) ^{-1} \;$
sur $\Gamma$ est
\`a valeurs dans l'ensemble des racines quatri\`emes de l'unit\'e, donc
constante.
On d\'eduit des \'egalit\'es $(\star)$ et $(\star\,\star)$ une formule pour
$\; \tr \sigma_r(\hat{a}) \,{\scriptstyle \times}\,
(\tr\pi_{\lambda_r})_{a M'_0} (x\exp Y) \;$
qui fait intervenir cette constante.

D'apr\`es la proposition \ref{caracteres canoniques} (a), la repr\'esentation
$\,\pi_{\lambda_r}\!$ de ${\underline{M'\!}\,}$ admet pour caract\`ere
infinit\'esimal $\, \chiinfty{\mi \, \lambda_r}{\gothm'}$.
On applique la proposition \ref{foncteur de translation} (c) \`a
${\underline{M'\!}\,}$ avec $\pi = \pi_{\lambda_r}$,
$\Lambda = \mi\,\lambda_r$ et $\Lambda_F = -2\rho$
($\pi_F$ ci-dessous).
On dispose d'une repr\'esentation $(\pi_{\lambda_r})_{Zuc}$ de
${\underline{M'\!}\,}$ et, vu ce qui pr\'ec\`ede, d'une formule pour
$\tr \sigma_r(\hat{a}) \,{\scriptstyle \times}\,
(\tr(\pi_{\lambda_r})_{Zuc})_{a M'_0} (x\exp Y)$.
Par la proposition \ref{passage au cas regulier} (b), on a aussi
$\,\tr(\pi_{\lambda_r})_{Zuc} = \tr\pi_{\lambdatilde}\,$ par choix de $\pi_F$.

On choisit maintenant $Y = 0$, et r\'ecapitule.
On constate que
$\bigl( \tr T_{\lambdatilde,\sigma}^{M'} \bigr)_{\!e}(X)$
se d\'eduit de l'expression de
$\bigl( \tr T_{\lambda_r,\sigma_r}^{M'} \bigr)_{\!e}(X)$
obtenue ci-dessus,
en rempla\c{c}ant le coefficient
$\, c^g_w \, \me^{\mi \,wy^gg\,\lambda_r\!(X)}$ par
$\, (\det(\Ad a^\C )_{\gothn_{M'}})^{-1}
\,{\scriptstyle \times}\; \underline{x}_{~\!l_F} \;
c^g_w \, \me^{\mi \,wy^gg\,\lambda_r\!(X)}$
pour chaque $ w\in W(\gothm'(e)_\C ,\gothj(x)_\C )$,
o\`u $\;l_F \egdef -2\,wy^gg\,\rho$.

\smallskip
Soit $\, w \in W(\gothm'(e)_\C ,\gothj(x)_\C )$.
On va calculer le terme $\underline{x}_{~\!l_F}$ qui lui est associ\'e.
On se donne un repr\'esentant $\wt{w}$ de $w$ dans $\,\M '(\C )(\pr(e))_0$.
On note ${\underline{M'\!}\,}_\C $ le produit semi-direct de $\Z $ par
$\M '(\C )_0$ tel que $\, 1 \in \Z $ agisse sur $\M '(\C )_0$ par
$\interieur\pr(a)$, et $\,\underline{\pr}\,$ l'application canonique de
${\underline{M'\!}\,}$ dans ${\underline{M'\!}\,}_\C $.
Il existe une unique repr\'esentation holomorphe $\pi_{F,\C }$ de
${\underline{M'\!}\,}_\C $ dans $F$ telle que
$\; \pi_F = \pi_{F,\C} \circ \underline{\pr}$.
On fixe un vecteur non nul $v_0$ de $F$ de poids
$\,-2\rho\,$ pour l'action de $\gothh_\C $ dans $F$.
On s'int\'eresse au vecteur non nul $\, \wt{w}y^gg \cdot v_0 \,$ de $F$ qui a
pour poids $\,l_F$ sous l'action de $\gothj(x)_\C $.
On pose encore $\, e' = g^{-1}eg$.
En se pla\c{c}ant dans le produit semi-direct de $\Z $ par $\M '(\C )$ au
moyen de $\interieur \pr(a)$ et utilisant la d\'efinition de $\underline{x}$
(cf. \ref{foncteur de translation} (c)), on obtient : \\
$\underline{\pr}(\underline{x}) \, \wt{w} y^g \pr(g)
\begin{array}[t]{cl}
\!=\!&
\pr(x) \, \wt{w} y^g \pr(g) \pr(a)^{-1} \, \underline{\pr}(\underline{a}) \\
\!=\!&
\pr(\exp X) \, \wt{w} y^g \pr(g) \pr(e'a^{-1}) \,
\underline{\pr}(\underline{a})
\end{array}$ \\
puis
$\;\, \underline{x} \cdot (\wt{w}y^gg \cdot v_0)
= \exp X \cdot (\wt{w}y^gg \cdot ((e'a^{-1}) \cdot v_0))$. \\
La propri\'et\'e
$\;\, e'a^{-1} \!\in M'_0(\lambdatilde) = \exp {\gothh} \;$
permet d'en d\'eduire que :

\centerline{
$\underline{x}_{~\!l_F}
= \det(\Ad a^\C )_{\gothn_{M'}} \; (\det({\Ad e'}^\C )_{\gothn_{M'}})^{-1}
\; \me^{-2\,wy^gg\,\rho(X)}$.
}

\smallskip
Cela permet de conclure.
\cqfd

}
\vfill\eject

%
\label{RepIndex}%
\begin{index des notations}
%

  \item[] ${\scriptstyle \langle} ~,\!~ {\scriptstyle \rangle}$,
          \pageref{<,>ev}, \pageref{<,>dual}

  \indexspace

  \item[] $\gotha$, \pageref{a}

  \indexspace

   \item[] $B_f$, \pageref{Bf}
   \item[] $\gothb_{M'}$, \pageref{bM'}

  \indexspace

  \item[] $\Car \gothg$, \pageref{Carg}
  \item[] $c_{\widehat{e'},\lambdatilde}$ cf. $(F)$ p.$\!$~\pageref{ce}
  \item[] $C(\gothg(\lambda ),\gothh)$, \pageref{Cgh}
  \item[] $C(\gothg(\lambda ),\gothh)_{reg}$, \pageref{Cghreg}

  \indexspace

  \item[] $d_e$, \pageref{de}
  \item[] $D_e$, \pageref{Deaa}
  \item[] $\cald_{\gothg}$, \pageref{Dg}
  \item[] $D_G$, \pageref{DGaa}
  \item[] $DL(V)$, \pageref{DL(V)}

  \indexspace

  \item[] $\F+[e]$, \pageref{F+[e]}
  \item[] $\Fh+$, \pageref{Fh+}

  \indexspace

  \item[] $\gothg$, \pageref{g}
  \item[] $G$, \pageref{G  aa}
  \item[] $\G $, \pageref{G  aabarre}
  \item[] $\widehat G$, \pageref{G  aachapeau}
  \item[] $\gothg(f)$, \pageref{g (f)}
  \item[] $G(f)$, \pageref{G (f)aa}
  \item[] $G(f)^{\gothg/\gothg(f)}$, \pageref{G (f)g/g(f)}
  \item[] $\gothg(x)$, \pageref{g (x)}
  \item[] $G(x)$, \pageref{G (x)aa}
  \item[] $\gothg(X)$, \pageref{g (Xaa)}
  \item[] $G(X)$, \pageref{G (Xaa)aa}
  \item[] $\gothg(\lambdatilde)$, \pageref{g (zzl tilde)}
  \item[] $G(\lambdatilde)$, \pageref{G (zzl tilde)aa}
  \item[] $G(\lambda_+)^{\gothg / \gothh}$, \pageref{G (f)g/g(f)}
  \item[] $G(\lambda_+)^{\gothg/\gothh}_0$, \pageref{G (f)g/g(f)0}
  \item[] $G(\lambdatilde)^{\gothg/\!\gothg(\lambda)(\mi \rho_\F+)}$,
          \pageref{G (zzl tilde)g/glambdairho}
  \item[] $G(\lambdatilde)^{\gothg/\gothh}$, \pageref{G (zzl tilde)g/h}
  \item[] $G(\lambdatilde)^{\gothg/\gothh}_0$, \pageref{G (zzl tilde)g/h0}
  \item[] $\gssreg$, \pageref{g ssreg}
  \item[] $\Gssreg$, \pageref{G ssregaa}
  \item[] $\gothg^*(x)$, \pageref{g*(x)}
  \item[] $\gothg^*(X)$, \pageref{g*(Xaa)}
  \item[] $\gstfondtilde$, \pageref{g*fond tilde}
  \item[] $\gstfondtilde(e)$, \pageref{g*fond tilde(e)}
  \item[] $\gstItilde$, \pageref{g*I tilde}
  \item[] $\gstItilde(e)$, \pageref{g*I tilde(e)}
  \item[] $\gstIGtilde$, \pageref{g*I,G tilde}
  \item[] $\gstInctilde$, \pageref{g*Inc tilde}
  \item[] $\gstInctilde(e)$, \pageref{g*Inc tilde(e)}
  \item[] $\gstIncGtilde$, \pageref{g*Inc,G tilde}
  \item[] $\gstreg$, \pageref{g*reg}
  \item[] $\gstregtilde$, \pageref{g*reg tilde}
  \item[] $\gstreg(e)$, \pageref{g*reg (e)}
  \item[] $\gstregtilde(e)$, \pageref{g*reg tilde(e)}
  \item[] $\gstregGtilde$, \pageref{g*reg,G tilde}
  \item[] $\gstss$, \pageref{g*ss}
  \item[] $\gstssfondG$, \pageref{g*ssfond,G}
  \item[] $\gstssI$, \pageref{g*ssI}
  \item[] $\gstssIG$, \pageref{g*ssI,G}
  \item[] $\gstssInc$, \pageref{g*ssInc}
  \item[] $\gstssIncG$, \pageref{g*ssInc,G}
  \item[] $\gstssreg$, \pageref{g*ssregzz}

  \indexspace

  \item[] $\gothh$ cf. \ref{lemme clef}, \pageref{lemme clef}
  \item[] $\calh$, \pageref{H}
  \item[] $\gothh_{(\R )}$, \pageref{hr}
  \item[] $H_{\alpha}$, \pageref{Halpha}

  \indexspace

  \item[] $\gothj(x)_1$, \pageref{j (x)1}

  \indexspace

  \item[] $k_e$, \pageref{ke}

  \indexspace

  \item[] $\gothm$, \pageref{m}
  \item[] $M$, \pageref{M}
  \item[] $\gothm'$, \pageref{m'}
  \item[] $M'$, \pageref{M'}
  \item[] $Mp(V)$, \pageref{Mp(V)}
  \item[] $Mp(V)_\call$, \pageref{Mp(V)L}

  \indexspace

  \item[] $n_\call$, \pageref{nL}
  \item[] $\gothn_{M'}$, \pageref{nM'}

  \indexspace

  \item[] $\scalo (A)_{A \cdot V}$, \pageref{OA}
  \item[] $\scalo (\hat{a})_{(1-a) \cdot V}$, \pageref{Oa}
  \item[] $\scalo (B)_V$, \pageref{OB}

  \indexspace

  \item[] $q$, \pageref{q}
  \item[] $q'$, \pageref{q'}
  \item[] $q_\call$, \pageref{qL}

  \indexspace

  \item[] $R(\gothg_\C ,\gothh_\C )$, \pageref{R(gC,hC)}
  \item[] $R^+(\gothg_\C ,\gothh_\C )$, \pageref{R+(gC,hC)},
                                        \pageref{R+(gC,hC)noncan}
  \item[] $R_G$, \pageref{RG}
  \item[] $R^+_\lambdatilde$, \pageref{R+lambdatilde}
  \item[] $R^+_{\lambdatilde,{\gotha^*}^+}$, \pageref{R+(gC,hC)can}

  \indexspace

  \item[] $S$, \pageref{S}
  \item[] $\sg$, \pageref{sg}
  \item[] $\Supp_{\gothg^*}\wt{\Omega}$, \pageref{SuppgstOmegatilde}

  \indexspace

  \item[] $\gotht$, \pageref{t}
  \item[] $T_0$, \pageref{T0}
  \item[] $T_{\lambdatilde,{\gotha^*}^+\!,\tau_+}^G\!$,
          \pageref{Tzzl tilde,a*,zzt,connexe}, \pageref{Tzzl tilde,a*,zzt}
  \item[] $T_{\lambdatilde,\tau}^G$, \pageref{Tzzl tilde,zzt}
  \item[] $T_{\lambda_r,\sigma_r}^{M'}$, \pageref{Tzzl zr,sigmar,M'}
  \item[] $T_{\lambda_r}^{M'_0}$, \pageref{Tzzl zr,M'0}
  \item[] $T_{\lambdatilde}^{M'_0}$, \pageref{Tzzl ztilde,M'0}

  \indexspace

  \item[] $\gothu$, \pageref{u}
  \item[] $U$, \pageref{U}

  \indexspace

  \item[] $\calv_e$, \pageref{Ve}
  \item[] $V^\chi $, \pageref{Vchi}

  \indexspace

  \item[] $W(G,\gothh)$, \pageref{W(G,h)}

  \indexspace

  \item[] $\Xfin_G$, \pageref{Xf}
  \item[] $\Xfin_G(\lambdatilde,{\gotha^*}^+)$, \pageref{Xfinal}
  \item[] $\XInd_G$, \pageref{XInd}
  \item[] $\XInd_G(\lambdatilde)$, \pageref{XI}
  \item[] $\Xirr_G(f)$, \pageref{Xirr}
  \item[] $\Xirr_G(\lambdatilde)$, \pageref{Xirr tilde}
  \item[] $\Xirrp_G(\lambdatilde,{\gotha^*}^+)$, \pageref{Xirr+}

  \indexspace

  \item[] $\beta _{\Omega}$, \pageref{zzb Omega}
  \item[] $\beta _{\wt{\Omega}}$, \pageref{zzb Omegatilde}

  \indexspace

  \item[] $\delta$, \pageref{zzd}
  \item[] $\delta_{\lambda_+}^{\gothg/\gothh}$, \pageref{zzd zzl+}

  \indexspace

  \item[] $\Theta _e$, \pageref{zzHt e}

  \indexspace

  \item[] $\iota $, \pageref{zzi}

  \indexspace

  \item[] $\lambda $, \pageref{zzl}
  \item[] $\lambda_+=\lambda_{\gothg,\lambdatilde,{\gotha^*}^+\!,\epsilon}$,
          \pageref{zzl +}
  \item[] $\lambda_{+,\gothm'}$, \pageref{zzl +m'}
  \item[] $\lambdatilde$, \pageref{zzl tilde}
  \item[] $\lambdatilde[e]$, \pageref{zzl tilde[e]}
  \item[] $\lambda_{\textit{can}}$, \pageref{zzl can}

  \indexspace

  \item[] $\mu $, \pageref{zzm}
  \item[] $\mu_+=\mu_{\gothg,\lambdatilde,{\gotha^*}^+}$, \pageref{zzm +}
  \item[] $\mu_{+,\gothm}$, \pageref{zzm +m}
  \item[] $\mutilde$, \pageref{zzm tilde}

  \indexspace

  \item[] $\nu $, \pageref{zzn}
  \item[] $\nu_+=\nu_{\gothg,\lambdatilde,{\gotha^*}^+\!,\epsilon}$,
          \pageref{zzn +}

  \indexspace

  \item[] $\xi $, \pageref{zznx}

  \indexspace

  \item[] $\displaystyle \left |\Pi _{\gothg,\gothm}\right |$, \pageref{zzP g,m}

  \indexspace

  \item[] $\rho_\F+$, \pageref{zzr F+}
  \item[] $\rho_{\gothg,\gothh}$, \pageref{zzr g,h}
  \item[] $\rho_\call$, \pageref{zzr L}
  \item[] $\rho_{\lambda_+}^{\gothg/\gothh}$, \pageref{zzr zzl+}
  \item[] $\rho_\lambdatilde^{\gothg/\!\gothg(\lambda)(\mi \rho_\F+)}$,
          \pageref{zzr zzl tilde g,h}
  \item[] $\rho_{\lambdatilde,{\gotha^*}^+}^{\gothg/\gothh}$,
          \pageref{zzr zzl tilde,a*+ g,h}

  \indexspace

  \item[] $\tau$, \pageref{zzt}
  \item[] $\tau_+$, \pageref{zzt plus}
  \item[] $\tau_M$, \pageref{zzt M}
  \item[] $\tau_{M'}$, \pageref{zzt M'}

  \indexspace

  \item[] $\Phi$, \pageref{zzV}

  \indexspace

  \item[] $\chiinfty{l}{\gothg}$, \pageref{zzw lUgC}
  \item[] $\chi_{\lambdatilde}^G$, \pageref{zzw zzl+G}

\end{index des notations}
\vfill

%
%

\begin{thebibl}{AAA 00}{R\'ef\'erences}{0}
\addcontentsline{toc}{section}{R\'ef\'erences}
\label{RepRef}

\bibitem[ABV 92]{ABV92}
{Adams,~J., D.~Barbasch, et D.~A.~Vogan},
{``The {L}anglands classification and irreducible characters for real reductive
groups''},
{Birkha\"user, 1992}.

\bibitem[B.. 72]{B.72}
{Bernat,~P., N.~Conze, M.~Duflo, M.~L\'evy-Nahas, M.~Ra\"{\i}s, P.~Renouard, et
M.~Vergne},
{``Repr\'esentations des groupes de Lie r\'esolubles''},
{Dunod, 1972}.

\bibitem[BW 80]{BW80}
{Borel,~A., et N.~Wallach},
{``Continuous cohomology, discrete subgroups, and representations of reductive
groups''},
{Annals of Math. Studies 94, Princeton University Press, 1980}.
{Seconde {\'e}dition : Math. Surveys 67, A. M. S., 2000}.

\bibitem[Bou 84]{Bo84}
{Bouaziz,~A.},
{Sur les repr\'esentations des groupes de  Lie r\'eductifs non connexes},
{{\em Math. Ann.}, {\bf268} (1984), 539--555}.

\bibitem[Bou 87]{Bo87}
{Bouaziz,~A.},
{Sur les caract\`eres des groupes de  Lie r\'eductifs non connexes},
{{\em J. Funct. Anal.}, {\bf70} (1987), 1--79}.

\bibitem[B~67]{B67}
{Bourbaki,~N.},
{``Th\'eories spectrales''},
{Hermann, 1967}.

\bibitem[Cow 88]{Co88}
{Cowling,~M.},
{On the characters of unitary representations},
{{\em J. Austral. Math. Soc. Ser. A}, {\bf45} (1988), 62--65}.

\bibitem[Dix 64]{Di64}
{Dixmier,~J.},
{``Les C* alg\`ebres et leurs repr\'esentations''},
{Gauthier-Villars, 1964}.

\bibitem[Dix 69]{Di69}
{Dixmier,~J.},
{Sur la repr\'esentation r\'eguli\`ere d'un groupe localement compact connexe},
{{\em Ann. Sci. \'Ecole Norm. Sup.}, {\bf2} (1969), 423--436}.

\bibitem[Duf 82a]{Df82a}
{Duflo,~M.},
{Construction de repr\'esentations unitaires d'un groupe de  Lie},
{In ``Harmonic analysis and group representations'', pages 129--221.
Liguori, Naples, 1982.
Cours d'\'et\'e du C.I.M.E., Cortona 1980}.

\bibitem[Duf 82b]{Df82b}
{Duflo,~M.},
{Th\'eorie de Mackey pour les groupes de  Lie alg\'ebriques},
{{\em Acta Math.}, {\bf149} (1982), 153--213}.

\bibitem[Duf 84]{Df84}
{Duflo,~M.},
{On the Plancherel formula for almost algebraic real  Lie groups},
{In ``Lie group representations III'', num\'ero 1077 in Lect. Notes in Math.,
pages 101--165. Springer-Verlag, 1984}.

\bibitem[DHV 84]{DHV84}
{Duflo,~M., G.~Heckman, et M.~Vergne},
{Projection d'orbites, formule de Kirillov et formule de Blattner},
{{\em M\'em. Soc. Math. Fr. (N.S.)}, {\bf15} (1984), 65--128}.

\bibitem[DV 88]{DV88}
{Duflo,~M., et M.~Vergne},
{La formule de Plancherel des groupes de  Lie semi-simples r\'eels},
{In ``Representations of Lie groups'', num\'ero~14 in Adv. Stud. Pure Math.,
pages 289--336. Academic Press, 1988}.

\bibitem[DV 93]{DV93}
{Duflo,~M., et M.~Vergne},
{Cohomologie \'equivariante et descente},
{{\em Ast\'erisque}, {\bf215} (1993), 5--108}.

\bibitem[Har 65]{Ha65}
{Harish-Chandra},
{Discrete series for semisimple  Lie groups I. Construction of invariant
eigendistributions},
{{\em Acta Math.}, {\bf113} (1965), 241--318}.

\bibitem[Har 76]{Ha76}
{Harish-Chandra},
{Harmonic analysis on real reductive groups III. The Maass-Selberg relations
and the Plancherel formula},
{{\em {Ann. of Math.}}, {\bf104} (1976), 117--201}.

\bibitem[Her 83]{Hr83}
{Herb,~R.~A.},
{Discrete series characters and Fourier inversion on semisimple real Lie
groups},
{{\em Trans. Amer. Math. Soc.}, {\bf277} (1983), 241--262}.

\bibitem[Hoc 65]{Ho65}
{Hochschild,~G.},
{``The structure of  Lie groups''},
{Holden-Day, 1965}.

\bibitem[Kna 86]{Kn86}
{Knapp,~A.~W.},
{``Representation Theory of Semisimple Groups''},
{Princeton University Press, 1986}.

\bibitem[Kna 96]{Kn96}
{Knapp,~A.~W.},
{``Lie groups beyond an introduction''},
{Birkha\"user, 1996}.

\bibitem[KV 95]{KV95}
{Knapp,~A.~W., et D.~A.~Vogan},
{``Cohomological Induction and Unitary Representations''},
{Princeton University Press, 1995}.

\bibitem[KZ 82]{KZ82}
{Knapp,~A.~W., et G.~J.~Zuckerman},
{Classification of irreductible tempered representations of semisimple groups},
{{\em Ann. of Math.}, {\bf116} (1982), 389--501.
(Typesetter's correction: {\em Ann. of Math.}, {\bf119} (1984), page 639)}.

\bibitem[Kos 59]{Ko59}
{Kostant,~B.},
{The principal three-dimensional subgroup and the Betti numbers of a complex
simple  Lie group},
{{\em Amer. J. Math.}, {\bf81} (1959), 973--1032}.

\bibitem[Mac 58]{Ma58}
{Mackey,~G.~W.},
{Unitary representations of group extension I},
{{\em Acta Math.}, {\bf99} (1958), 265--311}.

\bibitem[Ran 72]{Ra72}
{Ranga Rao,~R.},
{Orbital integrals in reductive groups},
{{\em Ann. of Math.}, {\bf96} (1972), 505--510}.


\bibitem[Ros 80]{Ro80}
{Rossmann,~W.},
{Limit characters of reductive  Lie groups},
{{\em Invent. Math.}, {\bf61} (1980), 53--66}.

\bibitem[Ros 82]{Ro82}
{Rossmann,~W.},
{Limit orbit in reductive  Lie algebras},
{{\em Duke Math. J.}, {\bf49} (1982), 215--229}.

\bibitem[Spr 66]{Sp66}
{Springer,~T.~A.},
{Some arithmetical results on semi-simple Lie algebras},
{{\em Inst. Hautes \'Etudes Sci.}, {\bf30} (1966), 115--141}.

\bibitem[Var 77]{Va77}
{Varadarajan,~V.~S.},
{``Harmonic analysis on real reductive groups''},
{num\'ero 576 in Lect. Notes in Math. Springer-Verlag, 1977}.

\bibitem[Ver 94]{Ve94}
{Vergne,~M.},
{Geometric quantization and equivariant cohomology},
{In ``First European Congress of Mathematics I, Paris 1992'',
  num\'ero 119 in Progr. Math., pages 249--295. Birkha\"user, 1994}.

\bibitem[Wal 88]{Wa88}
{Wallach,~N.},
{``Real Reductive Groups I''},
{Academic Press, 1988}.

\bibitem[Zuc 77]{Zu77}
{Zuckerman,~G.},
{Tensor products of finite and infinite dimensional representations
of semisimple  Lie groups},
{{\em Ann. of Math.}, {\bf106} (1977), 295--308}.

\end{thebibl}

\par\vspace*{1cm}
\noindent
\begin{minipage}[t]{.4\hsize}
\eightrm\parindent0pt
\noindent
Jean-Yves Ducloux\\
Universit\'e Paris 7,\\
UFR de Math\'ematiques, UMR 7586,\\
2 place Jussieu, 75\,251 Paris Cedex 05, France.\\
\font\eighttt=cmtt8 scaled\magstep1
m\'el : \makebox[1in][l]{\eighttt
ducloux{@}math.jussieu.fr}
\end{minipage}

%
\end{document}